\documentclass[11pt,reqno]{amsart}

\usepackage[no-math]{fontspec}
\usepackage{mathtools}
\usepackage{amssymb}
\usepackage{amsthm}
\usepackage[mathscr]{eucal}
\usepackage{enumitem}
\usepackage{microtype}
\usepackage{float}
\usepackage{tikz-cd}
\usetikzlibrary{arrows.meta,positioning}
\usepackage[
  backend=biber,
  style=numeric,
  doi=true,
  sorting=nyt,
  eprint=true,
  giveninits=true,
  maxnames=99
]{biblatex}
\DeclareFieldFormat[article]{title}{\mkbibemph{#1}}
\renewbibmacro*{in:}{%
  \ifentrytype{article}
    {}
    {\printtext{\bibstring{in}\intitlepunct}}}
\usepackage[
  colorlinks=false,
  linkbordercolor={1 0 0},
  citebordercolor={0 1 0},
  urlbordercolor={0 1 1}
]{hyperref}
\usepackage[nameinlink,noabbrev]{cleveref}

\newcommand{\RealNumbers}{\mathbb R}
\newcommand{\ComplexNumbers}{\mathbb C} 

\newcommand{\GHSpace}{\mathscr{M}}
\newcommand{\UnitDiameterGHSpace}{\GHSpace_1} 
\newcommand{\DiameterNormalization}{\operatorname{scale}_1} 
\newcommand{\GHDistance}{\mathcal{GH}}
\newcommand{\IsometryClass}[2]{[#1,#2]} 
\newcommand{\MetricSpace}[2]{(#1,#2)} 
\newcommand{\MeasuredIsometryClass}[3]{[#1,#2,#3]}
\newcommand{\MeasuredSpace}[3]{(#1,#2,#3)} 
\newcommand{\MetricQuotient}[2]{\operatorname{Quot}(#1,#2)} 
\newcommand{\RestrictedMetric}[2]{#1|_{#2\times#2}} 
\newcommand{\HausdorffDistance}[1]{\mathcal{H}[#1]}
\newcommand{\ProbabilityMeasures}{\operatorname{\mathscr{P}\!rob}}
\newcommand{\AmbientProbabilitySet}[2]{P_{#1}(#2)}

\newcommand{\InvarianceDefect}{\Diamond_{\mathrm{inv}}} 
\newcommand{\CompactHyperspace}{\mathscr{K}}
\newcommand{\AmbientClassMap}[2]{\chi_{#1,#2}} 
\newcommand{\MeasuredSpaces}{\mathscr E}
\newcommand{\FullSupportMeasuredClasses}{\mathscr S_{\mathrm{full}}}
\newcommand{\IsometryInvariantMeasuredClasses}{\mathscr S_{\mathrm{inv}}}
\newcommand{\InvariantMeasuredSpaces}{\mathscr I}
\newcommand{\IdentityMap}{\operatorname{id}}
\DeclareMathOperator{\diam}{diam}
\DeclareMathOperator{\dis}{dis}
\DeclareMathOperator{\Cor}{Cor}
\DeclareMathOperator{\Card}{Card}
\DeclareMathOperator{\dist}{dist}
\DeclareMathOperator{\supp}{supp}
\DeclareMathOperator{\Isom}{Isom}
\DeclareMathOperator{\spec}{spec}
\DeclareMathOperator{\range}{range}
\DeclarePairedDelimiter{\abs}{\lvert}{\rvert}
\DeclarePairedDelimiter{\norm}{\lVert}{\rVert}

\newcommand{\BoundedContinuousFunctions}{C_b} 
\newcommand{\ContinuousFunctions}{C} 
\newcommand{\OpenCover}{\mathcal U} 
\newcommand{\SeparableHilbertSpace}{\ell^2} 
\newcommand{\SubspaceFunctional}{\ell_0} 
\newcommand{\HahnBanachExtension}{\ell} 
\newcommand{\HahnBanachTestVector}{z} 
\newcommand{\SpanCoefficient}{t} 
\newcommand{\NormingVector}{u_0} 
\newcommand{\HomogeneityScalar}{c} 
\newcommand{\NormalizedNormingVector}{\widehat{\NormingVector}} 
\newcommand{\Indicator}{\mathbf 1} 
\newcommand{\LebesgueSpace}{L} 
\newcommand{\OrthogonalGroup}{O} 
\newcommand{\UnitSphere}{\mathbb S} 
\newcommand{\SummableNorm}{\ell^1} 
\newcommand{\ApproximationMap}{\Phi} 
\newcommand{\DualNorm}{^*} 
\newcommand{\LawSectionAt}[2]{\LawSection\IsometryClass{#1}{#2}}
\newcommand{\LawSection}{\Sigma} 
\newcommand{\MeasureProjection}{\mathsf p} 
\newcommand{\SingletonSpace}{*} 
\newcommand{\Pushforward}{_*} 
\newcommand{\RealizationMap}{\Psi} 
\newcommand{\UnitInterval}{[0,1]} 

\newenvironment{acknowledgements}{%
  \medskip\noindent\textit{Acknowledgements.}\ }{\par}
\newenvironment{useofai}{%
  \medskip\noindent\textbf{Use of AI.}\ }{\par}

\newcommand{\GHball}{\mathscr{B}}
\newcommand{\MetricBall}{B}
\newcommand{\SelectedMeasureSymbol}{\mathfrak m}
\newcommand{\SelectedMeasure}[2]{\SelectedMeasureSymbol_{#1,#2}}
\newcommand{\MeasureFiber}[1]{\mathscr L_{#1}}

\newcommand{\Barycenter}[1]{\operatorname{bar}_{#1}}
\newcommand{\LawBaseMap}{\varpi}
\newcommand{\FiberLawSpace}{\Lambda_{\MeasureProjection}}
\newcommand{\FiberBarycenterMap}{\beta_{\MeasureProjection}}
\newcommand{\MeasureIncidenceSpace}[2]{A_{#1,#2}}
\newcommand{\LawIncidenceSpace}[2]{D_{#1,#2}}
\newcommand{\IncidenceLaw}[2]{\ell_{#1,#2}}
\newcommand{\RealizedLaw}[2]{\sigma_{#1,#2}}
\newcommand{\RealizedBarycenter}[2]{b_{#1,#2}}
\newcommand{\MeasureLawAssignment}{\Xi}
\newcommand{\FiberRealization}[2]{I_{#1,#2}} 
\newcommand{\ModelNeighborhood}{\mathscr W}
\newcommand{\SpectralRankDomain}{\ModelNeighborhood_0}
\newcommand{\AuxiliaryAR}{\mathscr A}
\newcommand{\BorelSets}{\mathcal B}

\newcommand{\DistanceOperatorSymbol}{T}
\newcommand{\SpectralProjectionSymbol}{P}
\newcommand{\SpectralSubspaceSymbol}{V}
\newcommand{\DistanceOperator}[2]{\DistanceOperatorSymbol_{#1,#2}}
\newcommand{\SpectralProjection}[2]{\SpectralProjectionSymbol_{#1,#2}}
\newcommand{\SpectralSubspace}[2]{\SpectralSubspaceSymbol_{#1,#2}}

\newcommand{\ActingGroup}{G}
\newcommand{\AmbientCarrier}{Z}

\newcommand{\AmbientDisplacement}{\gamma}
\newcommand{\AmbientOperator}{A}
\newcommand{\ResolventErrorOperator}{\mathcal D} 
\newcommand{\CarrierAmbientOperator}{B} 
\newcommand{\ClosedBallSymbol}{\overline B}
\newcommand{\BasisExtension}{\widehat\Eigenfunction} 
\newcommand{\AmbientProjection}{\Pi}
\newcommand{\RieszProjection}[2]{\AmbientProjection_{#2}(#1)} 
\newcommand{\EigenvalueProjection}[2]{\AmbientProjection_{#1,#2}} 
\newcommand{\GeneralizedEigenspace}[2]{\mathcal{G}_{#2}(#1)} 
\newcommand{\ApproximatingMap}{g}

\newcommand{\ApproximationDomain}{K}
\newcommand{\ApproximationHomotopyAt}[3]{\ApproximationHomotopy(\IsometryClass{#1}{#2},#3)}
\newcommand{\ApproximationAt}[2]{\ApproximationMap(\IsometryClass{#1}{#2})}
\newcommand{\ModelWeightAt}[3]{\PartitionWeight_{#1}\IsometryClass{#2}{#3}}
\newcommand{\ApproximationHomotopy}{Q} 

\newcommand{\BumpFunction}[4]{\mathrm{Bump}_{#1,#2,#3,#4}}
\newcommand{\NormalizedBump}[5]{\mathrm{NB}_{#1,#2,#3,#4}^{#5}}
\newcommand{\BallMassBound}{c}

\newcommand{\BanachAmbientSpace}{E}
\newcommand{\BaseCarrier}{X}
\newcommand{\BasePoint}{x}
\newcommand{\BasisFunction}{\psi}
\newcommand{\BasisIndex}{k}

\newcommand{\BorelSubset}{A}

\newcommand{\CircleConstant}{\pi}
\newcommand{\ClosedDomain}{A}
\newcommand{\CoefficientBound}{M}
\newcommand{\NormComparisonBound}{L_{\NormSymbol}}
\newcommand{\CoefficientVector}{c}
\newcommand{\CommonMapDomain}{D}
\newcommand{\CompactFeatureImage}{S}
\newcommand{\CompactOperator}{T}
\newcommand{\CompactSetPoint}{s}
\newcommand{\CompactTestFamily}{\mathcal F}
\newcommand{\ComparisonCarrier}{Y}

\newcommand{\ComparisonFirstFactor}{C}
\newcommand{\ComparisonMetric}{e}
\newcommand{\EquilateralMetric}{\Upsilon} 
\newcommand{\ComparisonPoint}{y}
\newcommand{\ComparisonSecondFactor}{D}
\newcommand{\ConstantPoint}{p}
\newcommand{\ContractionHomotopy}{H}
\newcommand{\CoordinateIndex}{k}
\newcommand{\CoordinateProjection}{\operatorname{pr}}
\newcommand{\Correspondence}{R}
\newcommand{\CorrespondenceDisplacement}{\AmbientDisplacement}

\newcommand{\CoverIndex}{i}
\newcommand{\CoverScale}{L}
\newcommand{\CoverScaleAt}[2]{\CoverScale\IsometryClass{#1}{#2}}

\newcommand{\DiameterBound}{D}
\newcommand{\DiracProbability}{\delta}
\newcommand{\DomainIndexSet}{I}
\newcommand{\DomainPoint}{z}
\newcommand{\DomainSpace}{S}
\newcommand{\EigenIndexSet}{\operatorname{Ind}}
\newcommand{\Eigenfunction}{\psi}
\newcommand{\Eigenvalue}{\lambda}
\newcommand{\EquilateralCarrier}{E}
\newcommand{\EquilateralSize}{m}
\newcommand{\ErrorControl}{\varepsilon_{\mathrm{app}}}
\newcommand{\ErrorTolerance}{\varepsilon}

\newcommand{\ExtensionDomain}{S}
\newcommand{\ContinuousExtension}{F} 
\newcommand{\ExtensionTarget}{T}
\newcommand{\TopologicalSpaceClass}{\mathscr S}

\newcommand{\FamilyIndex}{i}
\newcommand{\FeatureCoordinates}{\mathbf v}
\newcommand{\FeatureDisplacement}{b}
\newcommand{\FiberPoint}{\xi}
\newcommand{\FiniteSize}{n}

\newcommand{\FirstCompactSet}{S}

\newcommand{\FirstEmbedding}{\iota}
\newcommand{\FirstFactor}{A}
\newcommand{\FirstFactorPoint}{a}

\newcommand{\FirstMeasure}{\mu}

\newcommand{\FirstProbabilityLaw}{\sigma}
\newcommand{\FirstScale}{a}
\newcommand{\FirstTime}{s}
\newcommand{\FirstVector}{v}

\newcommand{\GramMatrix}{\operatorname{Gram}}

\newcommand{\GraphProbability}{\nu}
\newcommand{\GroupElement}{g}
\newcommand{\FirstIsometry}{\alpha}
\newcommand{\HaarProbability}{\operatorname{Haar}}

\newcommand{\HilbertCube}{Q}
\newcommand{\HilbertDomain}{H}
\newcommand{\HilbertFunction}{f}
\newcommand{\HilbertTestFunction}{\psi}
\newcommand{\HomotopyTime}{t}
\newcommand{\IdentifyingIsometry}{u}
\newcommand{\IdentityMatrix}{\IdentityMap}
\newcommand{\ImaginaryPart}{\operatorname{Im}}
\newcommand{\ImaginaryUnit}{\sqrt{-1}}
\newcommand{\InputIsometry}{u}
\newcommand{\InputMap}{f}
\newcommand{\IntegrationPoint}{z}
\newcommand{\ResolventParameter}{\zeta} 
\newcommand{\IsometricEmbeddingMap}{f}
\newcommand{\Isometry}{g}
\newcommand{\IsometryGroup}{G}
\newcommand{\ActiveBlocksAt}[2]{\mathcal I(#1,#2)}
\newcommand{\FeatureBlock}[3]{B_{#1,#2,#3}}

\newcommand{\JointFeatureAt}[2]{\JointFeatureMap_{#1,#2}}
\newcommand{\FeatureImageAt}[2]{\CompactFeatureImage_{#1,#2}}
\newcommand{\JointFeatureMap}{F}
\newcommand{\KernelExtensionOperator}{\mathsf S}
\newcommand{\KroneckerSymbol}{\delta}
\newcommand{\LargerNeighborhood}{U}
\newcommand{\LawBasePoint}{t}
\newcommand{\LawBaseSpace}{T}

\newcommand{\LawTotalSpace}{W}

\newcommand{\LinearSpan}{\operatorname{span}}
\newcommand{\LocalError}{\mathcal{E}\mathrm{rr}}
\newcommand{\KuratowskiError}[2]{\LocalError^{\mathrm{Ku}}_{#1,#2}}
\newcommand{\KernelComparisonError}{\beta} 
\newcommand{\Mapping}{f}

\newcommand{\MetricSymbol}{d}

\newcommand{\ModelDimension}{n}
\newcommand{\ModelDomain}{\ModelNeighborhood}
\newcommand{\ModelError}{\LocalError_{\oplus}} 

\newcommand{\EquivariantRetraction}{r}
\newcommand{\RetractionMap}{r}

\newcommand{\NetPoint}{z}
\newcommand{\NetSize}{m}

\newcommand{\NormEmbedding}{J}
\newcommand{\BlockNormEmbedding}[2]{\NormEmbedding_{#1,#2}}
\newcommand{\NormSymbol}{N}
\newcommand{\OpenNeighborhood}{U}
\newcommand{\OpenSubset}{O}
\newcommand{\OpenSurjection}{p}
\newcommand{\OrbitMetric}{\mathcal{OD}}
\newcommand{\OrbitProjection}{\mathsf q}
\newcommand{\OrthogonalChange}{g}
\newcommand{\OtherFirstFactorPoint}{u}
\newcommand{\OtherSecondFactorPoint}{v}

\newcommand{\PartitionWeight}{w}

\newcommand{\ProbabilityCarrier}{S}

\newcommand{\ProductApproximation}{\operatorname{Prod}}
\newcommand{\PreviousApproximationImages}{\mathscr C}

\newcommand{\ProductScale}{r}
\newcommand{\ProjectedBasisFunction}{u}
\newcommand{\Pseudometric}{q}
\newcommand{\PseudometricError}{\Delta}
\newcommand{\Radius}{r}

\newcommand{\RealPart}{\operatorname{Re}}
\newcommand{\RecognitionTarget}{T}

\newcommand{\RestrictionOperator}{\mathsf R}
\newcommand{\RootOrder}{m}

\newcommand{\SecondBasisIndex}{l}

\newcommand{\SecondCompactSet}{T}

\newcommand{\SecondEmbedding}{\kappa}
\newcommand{\SecondFactor}{B}
\newcommand{\SecondFactorPoint}{b}

\newcommand{\SecondGroupElement}{h}
\newcommand{\SecondIsometry}{\beta}

\newcommand{\SecondMeasure}{\nu}
\newcommand{\SecondParameter}{t}
\newcommand{\SecondProbabilityCarrier}{T}

\newcommand{\SecondScale}{b}
\newcommand{\SecondVector}{w}
\newcommand{\SeparationScale}{c}
\newcommand{\SequenceIndex}{j}
\newcommand{\SmallError}{\eta}
\newcommand{\SmallerNeighborhood}{V}
\newcommand{\SpectralContour}{\Gamma}
\newcommand{\SpectralCutoff}{a}
\newcommand{\SpectralOuterBound}{R_{\mathrm{sp}}}
\newcommand{\SpectralRank}{n}
\newcommand{\SpectralRegion}{\Omega}
\newcommand{\SelectedEigenvalues}{\Lambda} 
\newcommand{\SphereFunction}{f}
\newcommand{\SpherePoint}{\mathbf s}

\newcommand{\TestFunction}{\varphi}

\newcommand{\TransposeSymbol}{\mathsf T}

\newcommand{\IntegrationDifferential}{\mathrm{d}}
\newcommand{\LocalTolerance}{\varepsilon_{\mathrm{loc}}}

\newcommand{\NonnegativeIntegers}{\mathbb Z_{\geq0}}
\newcommand{\AmbientMetric}{h}

\newcommand{\FirstFactorMetric}{d_0}
\newcommand{\SecondFactorMetric}{d_1}
\newcommand{\ComparisonFirstFactorMetric}{e_0}
\newcommand{\ComparisonSecondFactorMetric}{e_1}
\newcommand{\SubspaceMetric}{h}
\newcommand{\TargetMetric}{r}
\newcommand{\BanachMetric}{\mathsf d}
\newcommand{\CountableIntersectionType}{G_\delta}

\newcommand{\ClosedTestSet}{F}
\newcommand{\InnerNeighborhood}{V}
\newcommand{\JoiningContinuum}{C}
\newcommand{\ClosedUnitBall}{\ClosedBallSymbol}
\newcommand{\SpectralNeighborhood}{U}

\newcommand{\NetSet}{A}

\newcommand{\UniformMetric}{d_\infty} 
\newcommand{\SecondMapping}{g} 
\newcommand{\FunctionFamily}{\mathcal F} 
\newcommand{\CompactSubset}{K} 
\newcommand{\OperatorNorm}{\mathrm{op}} 
\newcommand{\SigmaAlgebraPrefix}{\sigma} 
\newcommand{\NetCoordinateMap}{\Theta} 

\newcommand{\RepresentingFunctional}{L}
\newcommand{\CompactParameterSpace}{K}
\newcommand{\ParameterPoint}{k}
\newcommand{\Integrand}{F}
\newcommand{\ComparisonPseudometric}{r}

\newcommand{\ContractionFactor}{C}
\newcommand{\PyramidSpace}{\operatorname{Pyr}}

\newcommand{\GHPDistance}{\mathcal{GHP}}
\newcommand{\ProkhorovDistance}[1]{\mathcal{LP}[#1]} 

\newcommand{\MMClassSpace}{\mathcal X}
\newcommand{\Pyramid}{\mathcal P}
\newcommand{\BoxDistance}{\square}

\newcommand{\BoxParameterInterval}{I_{\BoxDistance}}
\newcommand{\BoxLebesgueMeasure}{\mathcal L^1}
\newcommand{\BoxFirstParameter}{\theta_X}
\newcommand{\BoxSecondParameter}{\theta_Y}
\newcommand{\BoxLargeSubset}{S_{\BoxDistance}}
\newcommand{\ConcentrationDistance}{d_{\mathrm{conc}}}
\newcommand{\KyFanDistance}{d_{\mathrm{KF}}}
\newcommand{\LipschitzObservables}{\operatorname{Lip}_1}

\newcommand{\NormComparisonError}{\varepsilon_{\NormSymbol}}
\newcommand{\CoordinateComparisonError}{\varepsilon_{\FeatureCoordinates}}

\newcommand{\FunctionBound}{M} 
\newcommand{\TietzeExtension}{F} 

\newcommand{\ProjectedBasisRow}{\boldsymbol u} 

\DeclareMathOperator{\FunctionGraph}{graph}
\newcommand{\FiniteNormSpace}[1]{\mathcal N_{#1}}
\newcommand{\CoordinateClassFamily}{\mathfrak C}
\newcommand{\CoordinateClassSpace}[2]{\CoordinateClassFamily_{#1}(#2)}
\newcommand{\InvariantCoordinateClasses}[3]{\CoordinateClassFamily_{#1}^{\mathrm{inv}}\IsometryClass{#2}{#3}}
\newcommand{\CoordinateClassOf}[3]{[#2,#3]_{\OrthogonalGroup(#1)}}
\newcommand{\CoordinateClassElement}{\mathfrak c}
\newcommand{\LocalModelClass}[2]{\CoordinateClassElement_{#1,#2}} 
\newcommand{\LocalCoordinates}[2]{\FeatureCoordinates_{#1,#2}}
\newcommand{\LocalNorm}[2]{\NormSymbol_{#1,#2}}
\newcommand{\BlockCoordinates}[3]{\FeatureCoordinates_{#1,#2,#3}}
\newcommand{\BlockNorm}[3]{\NormSymbol_{#1,#2,#3}}
\newcommand{\LocalModel}{\mathfrak M}
\newcommand{\ModelAssignment}{\mathcal L}
\newcommand{\ModelValue}[2]{\ModelAssignment\IsometryClass{#1}{#2}}
\newcommand{\BlockModelValue}[3]{\ModelAssignment_{#1}\IsometryClass{#2}{#3}}
\newcommand{\CoordinateRadius}{\operatorname{rad}_0}
\newcommand{\ClassPseudometric}[1]{\Pseudometric_{#1}}
\newcommand{\LocalPseudometric}[2]{\ClassPseudometric{#1,#2}}
\newcommand{\CoordinateGraphSpace}[3]{\mathscr G_{#1}(#2,#3)}
\newcommand{\CoordinateGraphMetric}{D}
\newcommand{\CoordinateProductMetric}[1]{#1_{\times}}
\newcommand{\CoordinateOrbit}[3]{\mathcal O_{#1,#2}(#3)}

\newcommand{\GluingPseudometric}{\omega}
\newcommand{\BlockPseudometric}[3]{\Pseudometric_{#1,#2,#3}}
\newcommand{\GluedPseudometric}[2]{\GluingPseudometric_{#1,#2}}
\newcommand{\HomotopyPseudometric}[3]{\Pseudometric_{#1,#2,#3}} 
\newcommand{\IsometryRepresentation}[2]{\rho_{#1,#2}}

\newcommand{\AveragingOperator}{\operatorname{Av}}

\newcommand{\UltrametricGHSpace}{\mathscr U}
\newcommand{\NonArchimedeanGHDistance}{\mathcal{GH}_{\mathrm{NA}}}
\newcommand{\UrysohnFunctionModel}{\mathbf{G}}
\newcommand{\UrysohnFunctionDistance}{\Delta_{\mathrm U}}

\newcommand{\PackingSeparation}[1]{s_{#1}}
\newcommand{\GraphInterpolation}{\operatorname{Interp}}

\newcommand{\GramNormalization}{\operatorname{orth}}
\newcommand{\SphereMultiplier}[1]{\vartheta_{#1}}

\newcommand{\AmbientMeasureRealization}[2]{\operatorname{real}_{#1,#2}}
\newcommand{\AmbientLawMap}[2]{\operatorname{law}_{#1,#2}}
\newcommand{\AmbientMeanMap}[2]{\operatorname{mean}_{#1,#2}}
\newcommand{\LowerBallMass}{\operatorname{mass}}
\newcommand{\BallEvaluation}{\operatorname{ball}}
\newcommand{\SampleGramMatrix}{\mathbf M}

\newcommand{\SphereFunctionSpace}[1]{\ContinuousFunctions_{\UnitSphere}(#1)}
\newcommand{\SphereMaxTransform}{\operatorname{Max}}
\newcommand{\BoundedLinearMaps}{\operatorname{Lin}_{\mathrm b}}
\newcommand{\BlockModelError}[1]{\LocalError_{#1}}
\newcommand{\CommonModulus}{\omega_{\mathrm{cont}}}
\newcommand{\BasisCoordinateMap}{\operatorname{coord}_{\EigenbasisColumn}}
\newcommand{\EigenbasisColumn}{\boldsymbol{\BasisFunction}}

\newcommand{\ProjectedProfile}[3]{u_{#1,#2,#3}}
\newcommand{\ModelBump}[3]{b_{#1,#2,#3}}

\theoremstyle{plain}
\newtheorem{theorem}{Theorem}[section]
\newtheorem{lemma}[theorem]{Lemma}
\newtheorem{proposition}[theorem]{Proposition}
\newtheorem{corollary}[theorem]{Corollary}

\theoremstyle{definition}
\newtheorem{definition}[theorem]{Definition}

\newtheorem{question}[theorem]{Question}

\theoremstyle{remark}
\newtheorem{remark}[theorem]{Remark}

\AddToHook{env/theorem/begin}{\crefalias{section}{theorem}}
\AddToHook{env/lemma/begin}{\crefalias{section}{lemma}}
\AddToHook{env/proposition/begin}{\crefalias{section}{proposition}}
\AddToHook{env/corollary/begin}{\crefalias{section}{corollary}}
\AddToHook{env/claim/begin}{\crefalias{section}{claim}}
\AddToHook{env/definition/begin}{\crefalias{section}{definition}}
\AddToHook{env/example/begin}{\crefalias{section}{example}}
\AddToHook{env/question/begin}{\crefalias{section}{question}}
\AddToHook{env/remark/begin}{\crefalias{section}{remark}}

\crefname{theorem}{theorem}{theorems}
\Crefname{theorem}{Theorem}{Theorems}

\crefname{lemma}{lemma}{lemmas}
\Crefname{lemma}{Lemma}{Lemmas}

\crefname{proposition}{proposition}{propositions}
\Crefname{proposition}{Proposition}{Propositions}

\crefname{corollary}{corollary}{corollaries}
\Crefname{corollary}{Corollary}{Corollaries}

\crefname{claim}{claim}{claims}
\Crefname{claim}{Claim}{Claims}

\crefname{definition}{definition}{definitions}
\Crefname{definition}{Definition}{Definitions}

\crefname{example}{example}{examples}
\Crefname{example}{Example}{Examples}

\crefname{question}{question}{questions}
\Crefname{question}{Question}{Questions}

\crefname{remark}{remark}{remarks}
\Crefname{remark}{Remark}{Remarks}

\numberwithin{equation}{section}

\newcounter{proofstep}[theorem]
\newcounter{proofsubstep}[proofstep]

\newcommand{\ProofStep}[1]{%
  \leavevmode
  \refstepcounter{proofstep}%
  \textbf{Step~\theproofstep. #1}}

\makeatletter
\@ifundefined{Hy@theorem@makelinktarget}{}{%
\IfSocketExistsTF{refstepcounter/target}{%
  \NewSocketPlug{refstepcounter/target}{paper/theorem-heading}{%
    \AssignSocketPlug{refstepcounter/target}{hyperref}%
    \hyper@makecurrent{\@currentcounter}%
    \global\let\Hy@dth@currentHref\@currentHref
    \Hy@theorem@makelinktarget{\Hy@dth@currentHref}}%
  \AddToHook{cmd/@thm/before}{%
    \AssignSocketPlug{refstepcounter/target}{paper/theorem-heading}}%
}{}%
}
\makeatother

\title{The topology of Gromov--Hausdorff space}
\author{Yoshito Ishiki}
\address{Department of Mathematical Sciences\\
Tokyo Metropolitan University\\
Minami-osawa, Hachioji, Tokyo 192-0397, Japan}
\email{ishiki-yoshito@tmu.ac.jp}

\date{}

\begin{document}
\raggedbottom

\begin{abstract}
We prove that the Gromov--Hausdorff space is homeomorphic to the real Hilbert space $\SeparableHilbertSpace$. We construct a continuous assignment of full-support probability measures that is equivariant under isometries and finite-dimensional local approximations that control all pairwise distances. These approximations yield the absolute retract property for all metrizable spaces. We also prove that any countable family of continuous maps from compact metrizable spaces can be approximated, with respect to a prescribed open cover, by maps whose images form a discrete family.
\end{abstract}

\maketitle
\tableofcontents

\section{Introduction}\label{sec:introduction}
The Gromov--Hausdorff distance $\GHDistance$ measures how closely two compact metric spaces resemble each other.
Write $\GHSpace$ for the set of isometry classes of nonempty compact metric spaces with the Gromov--Hausdorff distance $\GHDistance$.
We call $(\GHSpace, \GHDistance)$ the \emph{Gromov--Hausdorff space}.
Antonyan asked whether $\GHSpace$ is homeomorphic to $\SeparableHilbertSpace$ (\cite[p.~2]{Antonyan2020Euclidean}).
The following classification theorem answers this question affirmatively.

\begin{theorem}[{\Cref{thm:part-iv-main}}]\label{thm:main}
The space $\GHSpace$ is homeomorphic to the Hilbert space $\SeparableHilbertSpace$.
\end{theorem}

Zava proved that the subspace of $\GHSpace$ consisting of finite metric spaces admits no coarse embedding into any Hilbert space (\cite[arXiv v3, Theorem~B]{Zava2025Coarse}).
In particular, the homeomorphism in \Cref{thm:main} cannot be bi-Lipschitz.

\medskip
\noindent \textbf{Related work.}
To place \Cref{thm:main} in context, we recall a classification of subspaces of $\GHSpace$.
For each positive integer $\FiniteSize $, Antonyan identified the subspace of $\GHSpace$ represented by compact subsets of $\RealNumbers^\FiniteSize $ with the Hilbert cube minus one point (\cite[Corollary~5.3]{Antonyan2021Euclidean}).

Geodesics and optimal correspondences describe the metric structure underlying the topological classification of $\GHSpace$.
Ivanov, Nikolaeva, and Tuzhilin showed that any two points of $\GHSpace$ can be joined by a geodesic (\cite[arXiv v1, Theorem~1]{IvanovNikolaevaTuzhilin2015}).
Ivanov, Iliadis, and Tuzhilin proved the existence of optimal correspondences between compact metric spaces and used them to realize their Gromov--Hausdorff distance as the Hausdorff distance between isometric copies in a common metric space (\cite[arXiv v1, Theorem~2.6 and Corollary~2.7]{IvanovIliadisTuzhilin2016}).
Rigidity under surjective isometries provides another comparison.
Ivanov and Tuzhilin proved that every surjective self-isometry of $\GHSpace$ is the identity (\cite[Main Theorem]{IvanovTuzhilin2019Isometries}).
Tuzhilin surveyed the geometry of the Hausdorff and Gromov--Hausdorff distances (\cite{Tuzhilin2020}).

Earlier work studied the topology of $\GHSpace$ through continuous families and embedded compact spaces.
Our result determines the homeomorphism type of $\GHSpace$ itself.
The author constructed families of branching Gromov--Hausdorff geodesics continuously parametrized by the Hilbert cube (\cite[Theorem~1.3]{Ishiki2022Branching}), and embedded arbitrary compact metrizable spaces into the subspaces of $\GHSpace$ consisting of continua (\cite[Theorem~1.1]{Ishiki2022Continua}), of spaces of prescribed dimensions (\cite[Theorem~1.3]{Ishiki2023FractalDimensions}), and of compact metric trees (\cite[Theorem~1.1]{Ishiki2023MetricTrees}).
In each case one can prescribe finitely many distinct values of the embedding.
Byakuno proved that a countable product of closed intervals of positive summable lengths, with the supremum metric, embeds isometrically into $\GHSpace$ (\cite[arXiv v1, Corollary~1.4]{Byakuno2026Embedding}).

For comparison with the classification of $\GHSpace$, we recall the local description and dimension of subspaces of finite metric spaces.
Iliadis, Ivanov, and Tuzhilin proved that within the subspace of spaces of the same cardinality, a sufficiently small neighborhood of a generic finite metric space is isometric to an open subset of a finite-dimensional space with the maximum norm (\cite[arXiv v1, Theorem~4.1]{IliadisIvanovTuzhilin2017Local}).
Here \emph{generic} means that the nonzero distances are pairwise distinct and every triangle inequality on three distinct points is strict.
Nakajima, Yamauchi, and Zava proved that the subspace of spaces with at most $\FiniteSize$ points has topological dimension $\FiniteSize(\FiniteSize-1)/2$ (\cite[arXiv v1, Theorem~A]{NakajimaYamauchiZava2025}).
Compactifications provide a complementary approach to the global topology of $\GHSpace$.
Nakajima and Shioya constructed a compact metrizable space containing a dense topological copy of $\GHSpace$ (\cite[arXiv v1, Main Theorem~1.2]{NakajimaShioyaCompactification}).

In Part~\ref{part:measures} we continuously select probability measures over $\GHSpace$.
When a measure is part of the data, Kazukawa, Nakajima, and Shioya proved that the space of metric measure spaces is contractible and locally path connected for both the box and concentration topologies (\cite[arXiv v1, Theorems~1.4--1.6]{KazukawaNakajimaShioya2024}).
Our classification concerns compact metric spaces.
One extension of the Gromov--Hausdorff framework treats incomplete spaces together with their boundaries.
Shibahara introduced a Gromov--Hausdorff distance with boundary on isometry classes of incomplete, precompact, locally compact metric spaces (\cite[arXiv v1, Definition~3.2 and Theorem~1.1]{Shibahara2021Boundary}).
The \emph{boundary} is the complement of the original space in its metric completion.
The local models in Part~\ref{part:spectral} retain isometries between representatives.
Isometries between individual spaces are also retained in the moduli stack used by Yuji to study families of compact metric spaces (\cite[arXiv v2, Definition~5.1]{Yuji2026ModuliStack}).

In Part~\ref{part:topology} we study $\GHSpace$ as a global analogue of a hyperspace in a fixed ambient metric space.
The relation between these settings also appears in the geometry of geodesics.
M\'emoli and Wan proved that every Gromov--Hausdorff geodesic can be realized as a Hausdorff geodesic of compact subsets of a single compact metric space (\cite[arXiv v2, Theorem~1]{MemoliWan2023}).
For comparison with our absolute retract theorem, Curtis proved, in particular, that the hyperspace of nonempty compact subsets of a Banach space with the Hausdorff metric of its norm is an absolute retract (\cite[Theorem~1.6]{Curtis1980}).
We use a quotient of such a hyperspace by a compact group to approximate the identity on $\GHSpace$.

Toru\'nczyk's characterization, used in Part~\ref{part:hilbert}, has a precedent in the study of the Urysohn universal metric space.
In the early twenty-first century Uspenskij proved that this space is homeomorphic to $\SeparableHilbertSpace$ (\cite[arXiv v1, Theorem~2.1]{Uspenskij2004Urysohn}).
The Urysohn space is the complete separable metric space containing an isometric copy of every separable metric space in which every isometry between finite subsets extends to a surjective self-isometry.

\medskip
\noindent \textbf{Proof strategy and the four parts.}
We first prove that $\GHSpace$ is an absolute retract for metrizable spaces.
We then approximate sequences of maps from compact metrizable spaces so that the approximating images form a discrete family.
Since $\GHSpace$ is complete and separable, the absolute retract and discrete approximation properties imply \Cref{thm:main} by Toru\'nczyk's characterization (\cite[p.~248, assertion~(i)]{Torunczyk1981}) and its correction (\cite[Section~C]{Torunczyk1985}).
Parts~\ref{part:measures} through~\ref{part:topology} construct the approximations needed for the absolute retract property.
The discrete approximation argument in Part~\ref{part:hilbert} uses only the common preliminaries and a separate metric construction.
Figure~\ref{fig:proof-strategy} shows how the two arguments yield the classification.

A difficulty in the absolute retract argument is that the isometry group depends on the compact metric space.
The following density statement follows from Rouyer's theorem (\Cref{thm:rouyer-generic-metric-spaces}).
\[
 \overline{
 \{\MetricSpace{\BaseCarrier}{\MetricSymbol}\in\GHSpace
   \mid\Isom\MetricSpace{\BaseCarrier}{\MetricSymbol}
      =\{\IdentityMap_\BaseCarrier\}\}
 }=\GHSpace.
\]
This also implies that the assignment of isometry groups with their uniform metrics is discontinuous as a map from $\GHSpace$ to itself.
Our construction must also handle spaces with nontrivial isometry groups.
On each model neighborhood, we construct homomorphisms from the varying isometry groups into one fixed orthogonal group so that the local point maps are equivariant.
We then use a partition of unity to combine the local point maps in a Banach space.
A countable product of orthogonal groups acts on this Banach space.
Passing to the orbit space of nonempty compact subsets of this Banach space makes the resulting approximation independent of all coordinate choices.

\medskip
\noindent \textbf{Part~\ref{part:measures}.
Continuous invariant probability measures.}
In the first part, we construct measures for integration on varying compact spaces.
The main result of this part, \Cref{thm:invariant-measures}, assigns a Borel probability measure $\SelectedMeasure{\BaseCarrier}{\MetricSymbol}$ to every nonempty compact metric space $\MetricSpace{\BaseCarrier}{\MetricSymbol}$.
Each assigned measure has full support,
\[
 \supp(\SelectedMeasure{\BaseCarrier}{\MetricSymbol})=\BaseCarrier.
\]
For every isometry $\IdentifyingIsometry\colon
\MetricSpace{\BaseCarrier}{\MetricSymbol}\to
\MetricSpace{\ComparisonCarrier}{\ComparisonMetric}$, the assignment satisfies
\[
 \IdentifyingIsometry\Pushforward\SelectedMeasure{\BaseCarrier}{\MetricSymbol}
 =\SelectedMeasure{\ComparisonCarrier}{\ComparisonMetric}.
\]
The assignment of measures is continuous in the following sense.
Let $\BaseCarrier_\SequenceIndex$ and $\BaseCarrier$ be compact subsets of a compact metric space $\MetricSpace{\AmbientCarrier}{\AmbientMetric}$ with restricted metrics $\MetricSymbol_\SequenceIndex$ and $\MetricSymbol$.
If $\HausdorffDistance{\AmbientMetric}(\BaseCarrier_\SequenceIndex,\BaseCarrier)\to0$, the measures $\SelectedMeasure{\BaseCarrier_\SequenceIndex}{\MetricSymbol_\SequenceIndex}$, transported to $\AmbientCarrier$ by inclusion, converge weakly in $\ProbabilityMeasures(\AmbientCarrier)$ to $\SelectedMeasure{\BaseCarrier}{\MetricSymbol}$.

For a fixed compact space, averaging over its isometry group produces an invariant probability measure.
Our construction must also preserve continuity when the space varies.
To construct such measures, we prove that the space of invariant full-support measured compact spaces is Polish and that the map that forgets the measure is an open surjection onto $\GHSpace$.
The theorem of Banakh and Radul for Polish spaces implies that the forgetful map is a Milyutin map.
We state Valov's generalization in \Cref{thm:valov}.
We choose an associated continuous family of probability measures on its fibers and take their barycenters to obtain the desired invariant measures.
We use full support and continuity in the spectral approximation of Part~\ref{part:spectral}.

\medskip
\noindent \textbf{Part~\ref{part:spectral}.
Continuous finite-dimensional approximation.}
The main result of Part~\ref{part:spectral} is \Cref{thm:local-models}, which locally approximates compact metric spaces by finite-dimensional coordinates.
Given $\IsometryClass{X}{d}\in\GHSpace$ and $\LocalTolerance>0$, on some open neighborhood $\ModelNeighborhood$ one can continuously assign normed coordinate classes $\ModelValue{Y}{e}$ of fixed dimension $\ModelDimension$.
Such a class consists of a map from $Y$ into $\RealNumbers^\ModelDimension$ and a norm on $\RealNumbers^\ModelDimension$, up to a simultaneous orthogonal change of coordinates.
The induced pseudometric $\LocalPseudometric{Y}{e}$ is invariant under every self-isometry and satisfies
\[
 \LocalError\IsometryClass{Y}{e}=\max_{y,z\in Y}|e(y,z)-\LocalPseudometric{Y}{e}(y,z)|<\LocalTolerance
 \qquad (\IsometryClass{Y}{e}\in\ModelNeighborhood).
\]
The error is also continuous in $\IsometryClass{Y}{e}$.
For each representative pair, the self-isometries of the input space act on the coordinate image by orthogonal transformations preserving its norm (\Cref{prop:local-isometry-representation}).

To construct $\ModelNeighborhood$, we use the measures of Part~\ref{part:measures} to define the distance operator $\DistanceOperator{\BaseCarrier}{\MetricSymbol}$ on $\LebesgueSpace^2(\BaseCarrier,\SelectedMeasure{\BaseCarrier}{\MetricSymbol})$.
For every $f\in\LebesgueSpace^2(\BaseCarrier,\SelectedMeasure{\BaseCarrier}{\MetricSymbol})$ and $\BasePoint\in\BaseCarrier$, set
\[
 (\DistanceOperator{\BaseCarrier}{\MetricSymbol}f)(\BasePoint)
 =\int_\BaseCarrier
   \MetricSymbol(\BasePoint,\IntegrationPoint)f(\IntegrationPoint)
   \,\IntegrationDifferential
   \SelectedMeasure{\BaseCarrier}{\MetricSymbol}(\IntegrationPoint).
\]
Maria, Oudot, and Solomon studied these operators and their spectral coordinates (\cite[Section~3]{MariaOudotSolomon2020}).
For $\BasePoint\in\BaseCarrier$ we call $\MetricSymbol_\BasePoint=\MetricSymbol(\BasePoint,\cdot)\in\ContinuousFunctions(\BaseCarrier)$ the \emph{Kuratowski map} of $\BasePoint$.
Projecting a normalized bump at a point onto a finite-dimensional spectral subspace and then applying the distance operator yields a uniform approximation to the Kuratowski map of that point.
When the center of the local model has at least two points, we take the coefficients in a real orthonormal basis as coordinates and transfer the uniform norm to the coefficient space.
At a one-point center we use the zero map into $\RealNumbers$ with the absolute value norm throughout the neighborhood.
The radius of the bump and the spectral parameter control the accuracy.
We retain whole eigenspaces so that orthogonal changes of basis accommodate repeated eigenvalues.
We extend the distance kernel to a common compact metric space.
The spectral projection theorem of Anselone and Palmer (\Cref{thm:riesz-convergence}) then implies that the selected spectral subspace has constant dimension on a neighborhood.
Part~\ref{part:topology} uses both the compact coordinate images and their norms.
We hope that this use of the distance operator encourages further applications of geometric analysis and spectral geometry to the Gromov--Hausdorff space.

The author's interpolation theorem (\cite[Theorem~1.1]{Ishiki2026Interpolation}) provides a uniform error estimate for extending metrics prescribed on a discrete family of closed subsets of a fixed metrizable space.
In the absolute retract proof, the continuous coordinate classes of \Cref{thm:local-models} yield pseudometrics approximating the metrics on varying compact spaces.

\medskip
\noindent \textbf{Part~\ref{part:topology}.
The absolute retract property.}
In \Cref{thm:part-iii-main}, the main result of the third part, we prove that $\GHSpace$ is an absolute retract for all metrizable spaces.
In \Cref{thm:absolute-extensor}, we prove that $\GHSpace$ is an absolute extensor for all metrizable spaces.

Fix a sequence $(\ModelDimension_\CoverIndex)_{\CoverIndex\in\NonnegativeIntegers}$ in which every positive integer occurs infinitely often, and fix a continuous function $\ErrorControl\colon\GHSpace\to(0,\infty)$.
At each point choose a neighborhood carrying a local model of error less than $\ErrorControl$ throughout.
Take a countable locally finite partition of unity subordinate to this cover.
For each partition function choose a model whose domain contains its support.
Assign distinct indices $\CoverIndex$ to these pairs so that the model dimension equals $\ModelDimension_\CoverIndex$, and write $\PartitionWeight_\CoverIndex$ for the corresponding partition function.
For each model, linearly and isometrically embed $\RealNumbers^{\ModelDimension_\CoverIndex}$ with its model norm into $\ContinuousFunctions(\UnitSphere^{\ModelDimension_\CoverIndex-1},\RealNumbers)$ with the uniform norm, where $\UnitSphere^{\ModelDimension_\CoverIndex-1}$ is the Euclidean unit sphere.
We use the explicit Banach--Mazur construction of \Cref{lem:dual-representation}.
Let $\BanachAmbientSpace$ be the $\SummableNorm$-sum of these function spaces over all $\CoverIndex\in\NonnegativeIntegers$.
In the components where $\PartitionWeight_\CoverIndex$ is positive, use the point map multiplied by its weight, and set the other components to zero.
The representation of the norms is compatible with orthogonal changes of coordinates (\Cref{lem:dual-representation,lem:dual-representation-equivariance}).
The compact group
\[
 \ActingGroup
 =\prod_{\CoverIndex\in\NonnegativeIntegers}
   \OrthogonalGroup(\ModelDimension_\CoverIndex)
\]
acts continuously on $\BanachAmbientSpace$ by linear isometries.
Passing to the orbit of each compact coordinate image removes the choices of bases and isometric representatives.
By Antonyan's equivariant hyperspace theorem, $\CompactHyperspace(\BanachAmbientSpace)$ is a $\ActingGroup$-AR.
Since the compact group $\ActingGroup$ has a countable base, Antonyan's orbit space theorem implies that the quotient $\AuxiliaryAR=\CompactHyperspace(\BanachAmbientSpace)/\ActingGroup$ is an absolute retract (\Cref{lem:hyperspace-orbit-ar}).
By \Cref{thm:global-domination} there are continuous maps
\[
 \GHSpace\xrightarrow{\ \ApproximationMap\ }
 \AuxiliaryAR\xrightarrow{\ \RealizationMap\ }\GHSpace,
\]
and a continuous homotopy $\ApproximationHomotopy\colon\GHSpace\times\UnitInterval\to\GHSpace$.
For every $\MetricSpace{\ComparisonCarrier}{\ComparisonMetric}\in\GHSpace$, this homotopy satisfies
\begin{equation}\label{eq:intro-domination-1}
 \ApproximationHomotopyAt{\ComparisonCarrier}{\ComparisonMetric}{0}
 =\IsometryClass{\ComparisonCarrier}{\ComparisonMetric},
\end{equation}
and
\begin{equation}\label{eq:intro-domination-2}
 \ApproximationHomotopyAt{\ComparisonCarrier}{\ComparisonMetric}{1}
 =\RealizationMap\ApproximationAt{\ComparisonCarrier}{\ComparisonMetric}.
\end{equation}
The homotopy tracks satisfy
\begin{equation}\label{eq:intro-domination-3}
 \diam_{\GHDistance}
 \bigl\{\ApproximationHomotopyAt{\ComparisonCarrier}{\ComparisonMetric}{t}
       \mid t\in\UnitInterval\bigr\}
 <\ErrorControl\IsometryClass{\ComparisonCarrier}{\ComparisonMetric}.
\end{equation}
A retract of an ANR is an ANR.
Hanner's domination theorem (\Cref{thm:hanner-domination}) extends this fact by allowing the retraction identity to hold up to homotopy, with arbitrarily small point tracks.
The maps in \eqref{eq:intro-domination-1} and \eqref{eq:intro-domination-2} and the estimate \eqref{eq:intro-domination-3} satisfy its hypotheses.
Thus $\GHSpace$ is an ANR.
Here we first obtain the ANR property for separable metrizable spaces, and then use Hanner's theorem (\Cref{thm:anr-category}) to obtain it for all metrizable spaces.
Scaling all distances to $0$ contracts $\GHSpace$ to the one-point space.
Using the ANR/ANE equivalence in Dugundji's theorem (\Cref{thm:retract-extensor}), we apply Hanner's theorem (\Cref{thm:contractible-ane}) to conclude that $\GHSpace$ is an AE for metrizable spaces.
By the AE/AR equivalence in Dugundji's theorem (\Cref{thm:retract-extensor}), we obtain the absolute retract property.

\begin{figure}[H]
\centering
\begin{tikzpicture}[
 >=Stealth,
 font=\small,
 partbox/.style={
  draw=black,
  fill=white,
  rounded corners=3pt,
  line width=.65pt,
  text width=4.9cm,
  minimum height=2.35cm,
  align=center,
  inner sep=7pt
 },
 route/.style={->,line width=.85pt,draw=black}
]
\node[partbox] (measures) at (-2.95,0) {
 \textbf{Part~\ref{part:measures}}\\
 Continuous choice of invariant\\
 probabilities of full support\\
 \Cref{thm:invariant-measures}
};
\node[partbox] (models) at (2.95,0) {
 \textbf{Part~\ref{part:spectral}}\\
 Finite-dimensional local models\\
 $\Isom(Y,e)\to\OrthogonalGroup(n)$\\
 \Cref{thm:local-models}
};
\node[partbox,minimum height=2.95cm] (ar) at (2.95,-3.55) {
 \textbf{Part~\ref{part:topology}}\\
 $\ActingGroup=\prod_{i\in\NonnegativeIntegers}\OrthogonalGroup(n_i)$\\
 $\CompactHyperspace(\BanachAmbientSpace)/\ActingGroup$ is an AR\\
 Small tracks $\Longrightarrow$ ANR\\
 Contraction by scaling $\Longrightarrow$ AR\\
 \Cref{thm:part-iii-main}
};
\node[
 partbox,
 minimum height=2.95cm
] (discrete) at (-2.95,-3.55) {
 \textbf{Part~\ref{part:hilbert}}\\
 Products with finite\\equilateral spaces\\
 Discrete approximation\\
 \Cref{thm:discrete-approximation}
};
\node[
 partbox,
 text width=10.8cm,
 minimum height=1.9cm
] (hilbert) at (0,-7) {
 \textbf{Toru\'nczyk's characterization}\\
 AR property and discrete approximation,\\completeness and separability\\
 $\GHSpace$ is homeomorphic to real $\SeparableHilbertSpace$\\
 \textbf{Part~\ref{part:hilbert}.
 \Cref{thm:part-iv-main}}
};
\draw[route] (measures.east)--(models.west);
\draw[route] (models.south)--(ar.north)
 node[midway,left=5pt,font=\footnotesize,align=right]
 {Partition of unity\\in the Banach space $\BanachAmbientSpace$};
\draw[route] (ar.south)--(2.95,-5.5)--(0,-5.5)--(hilbert.north);
\draw[line width=.85pt,draw=black]
 (discrete.south)--(-2.95,-5.5)--(0,-5.5);
\fill[black] (0,-5.5) circle (1.5pt);
\end{tikzpicture}
\caption{The two arguments that determine the topology of $\GHSpace$. The arrows indicate dependencies. Parts~\ref{part:measures}--\ref{part:topology} establish the absolute retract property. The discrete approximation theorem in Part~\ref{part:hilbert} uses the common preliminaries and a separate metric construction. Together with completeness and separability, these two conclusions give the Hilbert-space classification.}
\label{fig:proof-strategy}
\end{figure}

\medskip
\noindent \textbf{Part~\ref{part:hilbert}.
Discrete approximation and Hilbert space topology.}
Part~IV establishes the approximation condition in Toru\'nczyk's characterization of Hilbert space (\Cref{thm:torunczyk}).
Its discrete approximation theorem, \Cref{thm:discrete-approximation}, applies to every sequence $\{\ApproximationDomain_\FamilyIndex\}_{\FamilyIndex\in\NonnegativeIntegers}$ of compact metrizable spaces.
If $\InputMap_\FamilyIndex\colon\ApproximationDomain_\FamilyIndex\to\GHSpace$ are continuous and $\OpenCover$ is an open cover of $\GHSpace$, then there exist continuous maps $\ApproximatingMap_\FamilyIndex\colon\ApproximationDomain_\FamilyIndex\to\GHSpace$ with the following two properties.
\[
 \forall\FamilyIndex\in\NonnegativeIntegers\
 \forall x\in\ApproximationDomain_\FamilyIndex\
 \exists U\in\OpenCover,\qquad
 \{\InputMap_\FamilyIndex(x),\ApproximatingMap_\FamilyIndex(x)\}\subset U.
\]
Moreover, the images $\{\ApproximatingMap_\FamilyIndex(\ApproximationDomain_\FamilyIndex)\}_{\FamilyIndex\in\NonnegativeIntegers}$ form a discrete family in $\GHSpace$.
This means that each point of $\GHSpace$ has an open neighborhood $V$ such that $V\cap\ApproximatingMap_\FamilyIndex(\ApproximationDomain_\FamilyIndex)\ne\varnothing$ for at most one index $\FamilyIndex$.

Taking products with sufficiently small finite equilateral spaces changes each input by a controlled amount in Gromov--Hausdorff distance.
We choose the cardinalities recursively to exceed the bounds on separated subsets in the previously constructed compact images.
This separates the images, and letting the cardinalities tend to infinity prevents accumulation in $\GHSpace$.
This construction is independent of the results of Parts~\ref{part:measures}--\ref{part:topology}.
If all domains are the Hilbert cube $\HilbertCube=\UnitInterval^{\NonnegativeIntegers}$, then we obtain the discrete approximation condition in Toru\'nczyk's theorem.
We combine this condition with completeness, separability, and the absolute retract property from Part~\ref{part:topology} to prove \Cref{thm:main}.

The main theorem also answers the question of local contractibility.
Here \emph{local contractibility} means that every neighborhood $\OpenNeighborhood $ of any point contains a neighborhood $\SmallerNeighborhood $ of that point whose inclusion in $\OpenNeighborhood $ is null-homotopic.
By \Cref{thm:main}, every point has a neighborhood base of open sets that strongly deformation retract onto the point.
At the end of the paper we ask whether sufficiently small Gromov--Hausdorff balls are contractible (\Cref{question:gh-balls-contractible}).

\medskip
\noindent \textbf{Formal verification.}
The proof of \Cref{thm:main} was formalized in Lean~4 with five explicit assumptions from the literature.
The first asserts that the pushforward of probability measures induced by a continuous open surjection between Polish spaces has a continuous right inverse.
It follows from the softness theorem of Banakh and Radul (\Cref{thm:valov}), which Valov extended to completely metrizable spaces.
The other four are Antonyan's two retract theorems for equivariant hyperspaces and orbit spaces, Hanner's domination theorem, and Toru\'nczyk's characterization of Hilbert space.
At the cited commit of the companion repository (\href{https://github.com/yoshito-ishiki-math/lean-gh-space-hilbert/tree/8b25d3fc5b9f2a61dffa59e1daa8c9c483f70a46}{\texttt{lean-gh-space-hilbert}}, commit \texttt{8b25d3fc5b9f}), Comparator accepted the formal theorem \texttt{PalomarPaperN.main}, and NanoDa, con-ron, and the Lean kernel accepted its proof.
The five results from the literature remain assumptions in the formalization.
The formalized result is the main theorem.
The explicit local coordinates obtained from spectral projections, the auxiliary AR common to all error functions, and the two strong deformation retractions are proved in this paper.

\medskip
\noindent \textbf{Organization.}
The introduction and \Cref{sec:preliminaries} provide the common background for four parts.
Each part also has its own preliminaries, in Sections~\ref{sec:prelim-measures}, \ref{sec:prelim-spectral}, \ref{sec:prelim-topology}, and~\ref{sec:prelim-hilbert}, respectively.
Part~\ref{part:measures}, introduced in \Cref{sec:intro-measures} and developed in \Cref{sec:measures}, constructs continuous invariant probabilities.
Part~\ref{part:spectral}, introduced in \Cref{sec:intro-spectral} and developed in \Cref{sec:spectral}, uses these probabilities to construct finite-dimensional local models.
Part~\ref{part:topology}, introduced in \Cref{sec:intro-topology} and developed in \Cref{sec:ar}, combines the models into global approximations and proves the absolute retract property.
Part~\ref{part:hilbert}, introduced in \Cref{sec:intro-hilbert} and developed in \Cref{sec:approximation}, proves the discrete approximation property and completes the proof of \Cref{thm:main}.
\Cref{sec:questions} concludes with questions about related spaces and the geometry of GH balls.

\medskip
\noindent \textbf{Conventions and notation.}
Write $\Card(\BaseCarrier)$ for the cardinality of a set $\BaseCarrier$.
Every metric space $\MetricSpace{\BaseCarrier }{\MetricSymbol }$ representing a point of $\GHSpace$ is nonempty and compact.
In $\GHSpace$ we identify such a pair with its isometry class.
The symbol $\BaseCarrier $ denotes the underlying set with the topology induced by $\MetricSymbol $.
A construction indexed as in $\SelectedMeasure{\BaseCarrier}{\MetricSymbol}$ refers to the specified pair $\MetricSpace{\BaseCarrier }{\MetricSymbol }$.
When an isometric embedding into $\MetricSpace{\AmbientCarrier }{\AmbientMetric }$ is specified, we use the same notation for its image with the restricted metric.
Denote the identity on $\BaseCarrier$ by $\IdentityMap_{\BaseCarrier}$ and omit the subscript when the domain is specified.
Banach and Hilbert spaces are real unless the complex field is indicated.
Write $\CompactHyperspace(\BaseCarrier )$ for the nonempty compact subsets of $\BaseCarrier $, equipped with the Hausdorff metric $\HausdorffDistance{\MetricSymbol }$ when a metric $\MetricSymbol $ is specified.
For a metrizable space $\DomainSpace$, denote by $\ProbabilityMeasures(\DomainSpace )$ the space of Radon probability measures on $\DomainSpace$ with the topology of weak convergence against bounded continuous functions.
An \emph{isometric embedding} $\IsometricEmbeddingMap \colon\MetricSpace{\BaseCarrier }{\MetricSymbol }\to\MetricSpace{\ComparisonCarrier }{\ComparisonMetric }$ satisfies $\ComparisonMetric (\IsometricEmbeddingMap (\BasePoint ),\IsometricEmbeddingMap (\ComparisonPoint ))=\MetricSymbol (\BasePoint ,\ComparisonPoint )$ for all $\BasePoint ,\ComparisonPoint \in \BaseCarrier $.
An \emph{isometry} is a surjective isometric embedding.
Unless otherwise specified, our sequences are indexed by $\NonnegativeIntegers=\{0,1,2,\ldots\}$.

\begin{useofai}
OpenAI Codex, based on GPT-6, was used to explore and develop the mathematical constructions and proofs, search the literature, and prepare the exposition and \LaTeX{} source.
The research and writing workflow used customized versions of \texttt{math-research-harness} (\url{https://github.com/haruhisa-enomoto/math-research-harness}) and \texttt{math-paper-skills} (\url{https://github.com/haruhisa-enomoto/math-paper-skills}), both made publicly available by Haruhisa Enomoto.
The author takes full responsibility for the mathematical content and the final manuscript.
\end{useofai}

\begin{acknowledgements}
The author was supported by JSPS KAKENHI Grant Number JP24KJ0182.
\end{acknowledgements}

\section{Common preliminaries}\label{sec:preliminaries}
We recall metric convergence and compactness, comparison of metric quotients, and weak convergence of probability measures.
Section~\ref{subsec:preliminaries-topology} fixes terminology for contractibility and absolute extension.
Section~\ref{subsec:measures-preliminaries-rouyer} states Rouyer's theorem used to prove density of spaces with full-support isometry-invariant probability measures (\Cref{thm:dense-invariant-measured-spaces}).
\subsection{Gromov--Hausdorff convergence and compactness}\label{subsec:preliminaries-metric}
We begin with the correspondence formula for the distance and compactness criteria.

For a metric space $\MetricSpace{\BaseCarrier}{\MetricSymbol}$, a point $\BasePoint\in\BaseCarrier$, and a radius $\Radius>0$, we write
\[
 \MetricBall(\BasePoint,\Radius;\MetricSymbol)
 =\{\ComparisonPoint\in\BaseCarrier\mid\MetricSymbol(\BasePoint,\ComparisonPoint)<\Radius\}.
\]
For a subset $\BorelSubset\subset\BaseCarrier$, we put
\[
 \dist_{\MetricSymbol}(\BasePoint,\BorelSubset)
 =\inf_{\ComparisonPoint\in\BorelSubset}\MetricSymbol(\BasePoint,\ComparisonPoint),
 \qquad \dist_{\MetricSymbol}(\BasePoint,\emptyset)=\infty.
\]
We denote by $\CompactHyperspace(\BaseCarrier)$ the set of all nonempty compact subsets of $\BaseCarrier$.
The \emph{Hausdorff distance} on this set is defined by
\[
 \HausdorffDistance{\MetricSymbol}\colon
 \CompactHyperspace(\BaseCarrier)\times\CompactHyperspace(\BaseCarrier)\to[0,\infty).
\]
For every $\FirstCompactSet,\SecondCompactSet\in\CompactHyperspace(\BaseCarrier)$, put
\[
 \HausdorffDistance{\MetricSymbol}(\FirstCompactSet,\SecondCompactSet)
 =\max\left\{
   \sup_{\BasePoint\in\FirstCompactSet}\dist_{\MetricSymbol}(\BasePoint,\SecondCompactSet),
   \sup_{\ComparisonPoint\in\SecondCompactSet}\dist_{\MetricSymbol}(\ComparisonPoint,\FirstCompactSet)
 \right\}.
\]
The \emph{Vietoris topology} on $\CompactHyperspace(\BaseCarrier)$ has a basis consisting of the sets obtained as follows.
For an integer $m\geq1$ and open subsets $U_1,\ldots,U_m$ of $\BaseCarrier$, take
\[
 \left\{\FirstCompactSet\in\CompactHyperspace(\BaseCarrier)\;\middle|\;
 \FirstCompactSet\subset\bigcup_{i=1}^{m}U_i,
 \quad \FirstCompactSet\cap U_i\ne\emptyset\ (1\leq i\leq m)\right\}.
\]
\begin{theorem}[{Michael~\cite[Theorem~3.3 and Proposition~4.5.1]{Michael1951Hyperspaces}, see also \cite[Exercises~2.2.11(a),(b) and~3.2.9]{Beer1993}}]\label{thm:compact-hyperspace-topology}
Let $\MetricSpace{\BaseCarrier}{\MetricSymbol}$ be a metric space.
Then the Hausdorff distance induces the Vietoris topology on $\CompactHyperspace(\BaseCarrier)$.
If $\BaseCarrier$ is separable, then $\CompactHyperspace(\BaseCarrier)$ is separable.
\end{theorem}

\begin{lemma}[{Hausdorff Lipschitz estimate, Hutchinson~\cite[Section~2.4]{Hutchinson1981} for self-maps, see \cite[Chapter~2, Proposition~5.3]{Wicks1991} for the general case}]\label{lem:hausdorff-lipschitz-images}
Let $\MetricSpace{\BaseCarrier}{\MetricSymbol}$ and $\MetricSpace{\ComparisonCarrier}{\ComparisonMetric}$ be metric spaces, let $L\geq0$, and let $\Mapping\colon\BaseCarrier\to\ComparisonCarrier$ be $L$-Lipschitz.
Then for all $\FirstCompactSet,\SecondCompactSet\in\CompactHyperspace(\BaseCarrier)$, we have $\HausdorffDistance{\ComparisonMetric}
 (\Mapping(\FirstCompactSet),\Mapping(\SecondCompactSet))
 \leq L\cdot \HausdorffDistance{\MetricSymbol}(\FirstCompactSet,\SecondCompactSet)$.
\end{lemma}


For $\ErrorTolerance>0$, a subset $\NetSet\subset\BaseCarrier$ is an \emph{$\ErrorTolerance$-net} if
\[
 \forall\BasePoint\in\BaseCarrier\quad
 \exists\NetPoint\in\NetSet\quad
 \MetricSymbol(\BasePoint,\NetPoint)<\ErrorTolerance.
\]
Every compact metric space admits a finite $\ErrorTolerance$-net for every $\ErrorTolerance>0$.

A \emph{correspondence} between nonempty sets $\BaseCarrier $ and $\ComparisonCarrier $ is a subset of $\BaseCarrier \times \ComparisonCarrier $ with both coordinate projections surjective.
Denote the set of correspondences between $\BaseCarrier$ and $\ComparisonCarrier$ by $\Cor(\BaseCarrier,\ComparisonCarrier)$.
The \emph{distortion} of a correspondence $\Correspondence\in\Cor(\BaseCarrier,\ComparisonCarrier)$ between metric spaces $\MetricSpace{\BaseCarrier }{\MetricSymbol }$ and $\MetricSpace{\ComparisonCarrier }{\ComparisonMetric }$ is
\[
 \dis(\Correspondence )=
 \sup\left\{
   \abs{\MetricSymbol (x,z)-\ComparisonMetric (y,w)}
   \,\middle|\,
   (x,y),(z,w)\in \Correspondence
 \right\}.
\]
\begin{theorem}[{Kalton and Ostrovskii~\cite[Theorem~2.1]{KaltonOstrovskii1999}, see also \cite[Theorem~7.3.25]{BuragoBuragoIvanov2001} and \cite[arXiv v1, Theorem~6.12]{Tuzhilin2020}}]\label{thm:gh-correspondence}
For nonempty compact metric spaces $\MetricSpace{\BaseCarrier }{\MetricSymbol }$ and $\MetricSpace{\ComparisonCarrier }{\ComparisonMetric }$ the \emph{Gromov--Hausdorff distance} satisfies
\begin{equation}\label{eq:correspondence}
 \GHDistance(\MetricSpace{\BaseCarrier }{\MetricSymbol },\MetricSpace{\ComparisonCarrier }{\ComparisonMetric })=\frac12\inf_{\Correspondence\in\Cor(\BaseCarrier,\ComparisonCarrier)}\dis(\Correspondence).
\end{equation}
Equivalently, it is the infimum of the Hausdorff distances between isometric copies of $\MetricSpace{\BaseCarrier }{\MetricSymbol }$ and $\MetricSpace{\ComparisonCarrier }{\ComparisonMetric }$ in a common metric space.
\end{theorem}

Kalton and Ostrovskii proved this formula using the union of the graph of a map $f\colon X\to Y$ and the transpose of the graph of a map $g\colon Y\to X$.
Every correspondence contains one of this form, so both infima of distortion agree.
The next lemma realizes a prescribed correspondence in a common metric space with control of the distance between related points.
In the proof of \Cref{lem:common-gh-embedding} we use \Cref{lem:correspondence-embedding} to embed a convergent sequence of compact metric spaces into a single compact metric space.

\begin{lemma}\label{lem:correspondence-embedding}
Let $\MetricSpace{\BaseCarrier}{\MetricSymbol}$ and $\MetricSpace{\ComparisonCarrier}{\ComparisonMetric}$ be nonempty compact metric spaces, let $\Correspondence\subset\BaseCarrier\times\ComparisonCarrier$ be a correspondence, and let $\SmallError>0$ satisfy $\dis(\Correspondence)\leq2\SmallError$.
Then there is a metric $\AmbientMetric$ on the disjoint union $\BaseCarrier\sqcup\ComparisonCarrier$ extending both metrics such that for every $(\BasePoint,\ComparisonPoint)\in\Correspondence$ we have
\[
 \AmbientMetric(\BasePoint,\ComparisonPoint)\leq\SmallError.
\]
In particular, $\HausdorffDistance{\AmbientMetric}(\BaseCarrier,\ComparisonCarrier)\leq\SmallError$.
\end{lemma}

\begin{proof}
As in Kalton and Ostrovskii (\cite[proof of Theorem~2.1]{KaltonOstrovskii1999}), define the distance for $(x,y)\in\BaseCarrier\times\ComparisonCarrier$ by $\inf_{(z,w)\in\Correspondence}\bigl(\MetricSymbol(x,z)+\SmallError+\ComparisonMetric(w,y)\bigr)$.
The assumption $\dis(\Correspondence)\leq2\SmallError$ implies the triangle inequality, and every $(x,y)\in\Correspondence$ has cross-distance at most $\SmallError$ (see \cite[proof of Theorem~7.3.25, pp.~257--258]{BuragoBuragoIvanov2001}).
\end{proof}

\begin{theorem}[{Consequence of Gromov's compactness theorem~\cite[Section~6]{Gromov1981}, see also \cite[arXiv v1, Theorem~7.18]{Tuzhilin2020}}]\label{thm:uniform-nets}
If $\MetricSpace{\BaseCarrier_\SequenceIndex}{\MetricSymbol_\SequenceIndex}
 \to\MetricSpace{\BaseCarrier}{\MetricSymbol}$ in $\GHSpace$, then for every $\ErrorTolerance>0$ there is a positive integer $\NetSize$ such that each $\BaseCarrier_\SequenceIndex$ admits an $\ErrorTolerance$-net of cardinality at most $\NetSize$.
\end{theorem}

Let $\MetricSpace{\BaseCarrier}{\MetricSymbol}$ be a metric space, and let $\SeparationScale>0$.
A subset $S\subset\BaseCarrier$ is called \emph{$\SeparationScale$-separated} if
\[
 \forall x,y\in S,\qquad
 x\ne y\Longrightarrow\MetricSymbol(x,y)\ge\SeparationScale.
\]
For an integer $m\geq2$ and $\xi=\IsometryClass{X}{d}\in\GHSpace$, denote the largest possible minimum separation of $m$ points by
\[
 \PackingSeparation{m}(\xi)=
 \max_{(x_0,\ldots,x_{m-1})\in X^m}\ \min_{0\leq a<b<m}d(x_a,x_b).
\]
The maximum exists by compactness of $X^m$.
We allow repetitions, so $\Card(X)<m$ implies $\PackingSeparation{m}(\xi)=0$.
This quantity is $2$ times the $m$th packing radius of Grove and Markvorsen (\cite[Introduction, formula~(1)]{GroveMarkvorsen1995Packing}).

\begin{corollary}\label{cor:uniform-separated-sets}
For each $m\geq2$ the function $\PackingSeparation{m}:\GHSpace\to[0,\infty)$ is $2$-Lipschitz, and $\PackingSeparation{m}\downarrow0$ uniformly on compact subsets of $\GHSpace$.
In particular, if $\xi_j\to\xi$, then for every $\eta>0$ the cardinalities of $\eta$-separated subsets of representatives of $\xi_j$ have a common bound independent of $j$.
\end{corollary}

\begin{proof}
Transferring $m$ points along a correspondence $R$ changes each pairwise distance by at most $\dis(R)$.
Comparing the minima and interchanging the spaces shows that the values of $\PackingSeparation{m}$ differ by at most $\dis(R)$.
Taking the infimum over correspondences yields the $2$-Lipschitz estimate by the formula of Kalton and Ostrovskii (\Cref{thm:gh-correspondence}).
Increasing the number of points cannot increase $\PackingSeparation{m}$.
For any $\eta>0$, if more points are chosen than there are points in a finite $\eta/3$-net of $X$, two lie within distance $2\eta/3$.
Hence $\PackingSeparation{m}(\xi)\downarrow0$ for every $\xi$.
The increasing open sets $\{\PackingSeparation{m}<\eta\}$ cover $\GHSpace$, so every compact subset lies in one of them.
Applying this to $\{\xi\}\cup\{\xi_j:j\in\NonnegativeIntegers\}$ yields the uniform cardinality bound.
\end{proof}

We use completeness and separability of the Gromov--Hausdorff space $\GHSpace$.

\begin{theorem}[{Edwards~\cite{Edwards1975Superspace}, see also \cite[arXiv v1, Corollary~7.24]{Tuzhilin2020}}]\label{lem:gh-basics}
The metric space $\MetricSpace{\GHSpace}{\GHDistance}$ is complete and separable.
\end{theorem}

Gromov constructed a common compact metric space into which every member of a family of compact metric spaces embeds isometrically, provided their diameters and their covering numbers at each positive radius are uniformly bounded (\cite[Section~6, pp.~64--65]{Gromov1981}, see also \cite[arXiv v1, Theorem~7.19]{Tuzhilin2020}).
The following theorem also realizes Hausdorff convergence of the embedded spaces.

\begin{theorem}\label{lem:common-gh-embedding}
Every Cauchy sequence $\{\MetricSpace{\BaseCarrier_\SequenceIndex}{\MetricSymbol_\SequenceIndex}\}_{\SequenceIndex\in\NonnegativeIntegers}$ in $\GHSpace$ admits isometric embeddings into a compact metric space $\MetricSpace{\AmbientCarrier}{\AmbientMetric}$ whose images converge in Hausdorff distance to a nonempty compact subset.
If its Gromov--Hausdorff limit is $\MetricSpace{\BaseCarrier}{\MetricSymbol}$, then we may also embed $\MetricSpace{\BaseCarrier}{\MetricSymbol}$ into the same ambient space.
Identifying these spaces with their images, we obtain
\[
 \HausdorffDistance{\AmbientMetric}(\BaseCarrier_\SequenceIndex,\BaseCarrier)\longrightarrow0.
\]
\end{theorem}

\begin{proof}
Let $\MetricSpace{\BaseCarrier}{\MetricSymbol}$ be the limit provided by Edwards's theorem (\Cref{lem:gh-basics}), and put
\[
 \ErrorTolerance_\SequenceIndex
 =\GHDistance(\MetricSpace{\BaseCarrier_\SequenceIndex}{\MetricSymbol_\SequenceIndex},
              \MetricSpace{\BaseCarrier}{\MetricSymbol})+2^{-\SequenceIndex}.
\]
Use the correspondence formula of Kalton and Ostrovskii (\Cref{thm:gh-correspondence}) and \Cref{lem:correspondence-embedding}, and glue along the common copy of $\BaseCarrier$ (\cite[author's manuscript, Theorem~27.2]{Villani2009}).
The resulting metric $\AmbientMetric$ on the union of $\BaseCarrier$ and all $\BaseCarrier_\SequenceIndex$ extends each metric and satisfies
\[
 \HausdorffDistance{\AmbientMetric}
 (\BaseCarrier_\SequenceIndex,\BaseCarrier)
 \leq\ErrorTolerance_\SequenceIndex\longrightarrow0.
\]
For every $\SmallError>0$, the tail lies in the $\SmallError/2$-neighborhood of $\BaseCarrier$.
A finite $\SmallError/2$-net in $\BaseCarrier$, together with finite $\SmallError$-nets in the finitely many remaining spaces, forms a finite $\SmallError$-net for the union.
The completion of this union is therefore compact and retains the displayed Hausdorff bounds.
\end{proof}

We also recall the compactness criterion used for families of functions.

\begin{theorem}[{Arzel\`a--Ascoli, \cite[Theorem~6.26]{Muscat2024}}]\label{thm:arzela-ascoli}
Let $\MetricSpace{\BaseCarrier}{\MetricSymbol}$ be nonempty and compact, and let $\MetricSpace{\ComparisonCarrier}{\ComparisonMetric}$ be complete.
Equip the set $\ContinuousFunctions(\BaseCarrier,\ComparisonCarrier)$ of continuous maps with the uniform metric
\[
 \UniformMetric(\Mapping,\SecondMapping)
 =\sup_{\BasePoint\in\BaseCarrier}
   \ComparisonMetric(\Mapping(\BasePoint),\SecondMapping(\BasePoint)).
\]
A family $\FunctionFamily\subset\ContinuousFunctions(\BaseCarrier,\ComparisonCarrier)$ has compact closure if and only if the union of its images $\bigcup_{\Mapping\in\FunctionFamily}\Mapping(\BaseCarrier)$ is totally bounded and $\FunctionFamily$ is equicontinuous.
Here \emph{equicontinuous} means that for every $\ErrorTolerance>0$ there exists $\SmallError>0$ such that for every $\Mapping\in\FunctionFamily$ and all $\BasePoint,\ComparisonPoint\in\BaseCarrier$ with $\MetricSymbol(\BasePoint,\ComparisonPoint)<\SmallError$,
\[
 \ComparisonMetric(\Mapping(\BasePoint),\Mapping(\ComparisonPoint))<\ErrorTolerance.
\]
\end{theorem}

The cited theorem characterizes total boundedness in the uniform metric.
Completeness of the target $\MetricSpace{\ComparisonCarrier}{\ComparisonMetric}$ implies completeness of the function space $\ContinuousFunctions(\BaseCarrier,\ComparisonCarrier)$ with the uniform metric $\UniformMetric$, which yields the formulation in terms of compact closure above.
\subsection{Pseudometrics and metric quotients}\label{subsec:preliminaries-quotients}
We prove estimates that also apply when distinct points are identified.

The diameter is continuous on $\GHSpace$ (\cite[arXiv v1, Example~6.29]{Tuzhilin2020}).
More precisely, for any nonempty compact metric spaces $\MetricSpace{\BaseCarrier}{\MetricSymbol}$ and $\MetricSpace{\ComparisonCarrier}{\ComparisonMetric}$ the cited result states that
\[
 \abs{\diam\MetricSpace{\BaseCarrier }{\MetricSymbol }-\diam\MetricSpace{\ComparisonCarrier }{\ComparisonMetric }}\le 2\GHDistance(\MetricSpace{\BaseCarrier }{\MetricSymbol },\MetricSpace{\ComparisonCarrier }{\ComparisonMetric }).
\]
A continuous pseudometric $\Pseudometric $ on a compact space $\BaseCarrier $ defines an equivalence relation.
For $\BasePoint,\ComparisonPoint\in\BaseCarrier$ write $\BasePoint\sim_\Pseudometric\ComparisonPoint$ if $\Pseudometric(\BasePoint,\ComparisonPoint)=0$.
Its \emph{metric quotient} is
\[
 {\MetricQuotient{\BaseCarrier }{\Pseudometric }}
 =\MetricSpace{\BaseCarrier /{\sim_\Pseudometric }}{\overline{\Pseudometric }},
\]
where for representatives $\BasePoint,\ComparisonPoint\in\BaseCarrier$ the \emph{quotient metric} satisfies
\[
 \overline{\Pseudometric }([\BasePoint ],[\ComparisonPoint ])=\Pseudometric (\BasePoint ,\ComparisonPoint ).
\]
The quotient is compact.
If $\Pseudometric $ is a metric, the quotient is isometric to $\MetricSpace{\BaseCarrier }{\Pseudometric }$.
For comparison of a pseudometric quotient with a metric on the same set, see the author's result (\cite[Proposition~2.12]{Ishiki2023FractalDimensions}).
We use the following formulation in terms of correspondences to compare two metric quotients.

\begin{proposition}\label{lem:pseudometric-comparison}
Let $\BaseCarrier$ and $\ComparisonCarrier$ be nonempty compact spaces.
Let $\Correspondence\subset\BaseCarrier\times\ComparisonCarrier$ be a correspondence.
For continuous pseudometrics $\Pseudometric$ on $\BaseCarrier$ and $\ComparisonPseudometric$ on $\ComparisonCarrier$, put
\begin{equation}\label{eq:correspondence-pseudometric-error}
 \PseudometricError_\Correspondence(\Pseudometric,\ComparisonPseudometric)
 =\sup_{(x,y),(z,w)\in\Correspondence}
 |\Pseudometric(x,z)
       -\ComparisonPseudometric(y,w)|.
\end{equation}
\begin{enumerate}[label=\textup{(Q\arabic*)},ref=\textup{(Q\arabic*)}]
\item\label{item:quotient-1}
The quotient spaces satisfy
\begin{equation}\label{eq:pseudometric-correspondence-bound}
 \GHDistance(\MetricQuotient{\BaseCarrier}{\Pseudometric},
              \MetricQuotient{\ComparisonCarrier}{\ComparisonPseudometric})
 \leq\tfrac12\PseudometricError_\Correspondence(\Pseudometric,\ComparisonPseudometric).
\end{equation}
\item\label{item:quotient-2}
If $\Pseudometric_0,\Pseudometric_1$ are continuous pseudometrics on $\BaseCarrier$ and $\ComparisonPseudometric_0,\ComparisonPseudometric_1$ are continuous pseudometrics on $\ComparisonCarrier$, then
\begin{equation}\label{eq:pseudometric-uniform-comparison}
 \bigl|\norm{\Pseudometric_0-\Pseudometric_1}_\infty
       -\norm{\ComparisonPseudometric_0-\ComparisonPseudometric_1}_\infty\bigr|
 \leq\PseudometricError_\Correspondence(\Pseudometric_0,\ComparisonPseudometric_0)
       +\PseudometricError_\Correspondence(\Pseudometric_1,\ComparisonPseudometric_1).
\end{equation}
\item\label{item:quotient-3}
Let the pseudometrics satisfy the hypotheses in \ref{item:quotient-2} and let $\HomotopyTime\in[0,1]$.
Then
\begin{equation}\label{eq:pseudometric-interpolation}
 \PseudometricError_\Correspondence((1-\HomotopyTime)\Pseudometric_0+\HomotopyTime\Pseudometric_1,
                   (1-\HomotopyTime)\ComparisonPseudometric_0+\HomotopyTime\ComparisonPseudometric_1)
 \leq(1-\HomotopyTime)\PseudometricError_\Correspondence(\Pseudometric_0,\ComparisonPseudometric_0)
       +\HomotopyTime\PseudometricError_\Correspondence(\Pseudometric_1,\ComparisonPseudometric_1).
\end{equation}
\item\label{item:quotient-4}
For continuous pseudometrics $\Pseudometric_0,\Pseudometric_1$ on the same space $\BaseCarrier$, take the diagonal correspondence $\{(x,x)\mid x\in\BaseCarrier\}$ in \ref{item:quotient-1} to obtain
\begin{equation}\label{eq:pseudometric-bound}
 \GHDistance(\MetricQuotient{\BaseCarrier}{\Pseudometric_0},\MetricQuotient{\BaseCarrier}{\Pseudometric_1})
 \leq\tfrac12\norm{\Pseudometric_0-\Pseudometric_1}_\infty.
\end{equation}
\end{enumerate}
\end{proposition}

\begin{proof}
\textbf{Proof of \ref{item:quotient-1}.
Correspondence estimate.}
Let $\pi_{\BaseCarrier,\Pseudometric}\colon\BaseCarrier\to\BaseCarrier/{\sim_\Pseudometric}$ and $\pi_{\ComparisonCarrier,\ComparisonPseudometric}\colon\ComparisonCarrier\to\ComparisonCarrier/{\sim_\ComparisonPseudometric}$ be the quotient maps.
Since both quotient maps and both coordinate projections of $\Correspondence$ are surjective, the following set is a correspondence.
\begin{equation}\label{eq:quotient-correspondence}
 \{(\pi_{\BaseCarrier,\Pseudometric}(\BasePoint),\pi_{\ComparisonCarrier,\ComparisonPseudometric}(\ComparisonPoint))
       \mid(\BasePoint,\ComparisonPoint)\in\Correspondence\}.
\end{equation}
For $(x,y),(z,w)\in\Correspondence$, by the definition in \eqref{eq:correspondence-pseudometric-error},
\begin{equation}\label{eq:quotient-correspondence-error}
 |\Pseudometric(x,z)-\ComparisonPseudometric(y,w)|
 \leq\PseudometricError_\Correspondence(\Pseudometric,\ComparisonPseudometric).
\end{equation}
By the correspondence formula of Kalton and Ostrovskii (\Cref{thm:gh-correspondence}), estimate \eqref{eq:quotient-correspondence-error} implies \eqref{eq:pseudometric-correspondence-bound}.

\textbf{Proof of \ref{item:quotient-2}.
Uniform comparison.}
Fix $(x,y),(z,w)\in\Correspondence$.
The triangle inequality on the real line implies
\begin{equation}\label{eq:pseudometric-pointwise-comparison}
 |\Pseudometric_0(x,z)
    -\Pseudometric_1(x,z)|
 \leq|\ComparisonPseudometric_0(y,w)
         -\ComparisonPseudometric_1(y,w)|
 +\PseudometricError_\Correspondence(\Pseudometric_0,\ComparisonPseudometric_0)
       +\PseudometricError_\Correspondence(\Pseudometric_1,\ComparisonPseudometric_1).
\end{equation}
Every pair in $\BaseCarrier\times\BaseCarrier$ has a corresponding pair in $\ComparisonCarrier\times\ComparisonCarrier$.
Taking the supremum in \eqref{eq:pseudometric-pointwise-comparison}, we obtain
\begin{equation}\label{eq:pseudometric-sup-comparison}
 \norm{\Pseudometric_0-\Pseudometric_1}_\infty
 \leq\norm{\ComparisonPseudometric_0-\ComparisonPseudometric_1}_\infty
      +\PseudometricError_\Correspondence(\Pseudometric_0,\ComparisonPseudometric_0)
      +\PseudometricError_\Correspondence(\Pseudometric_1,\ComparisonPseudometric_1).
\end{equation}
Applying \eqref{eq:pseudometric-sup-comparison} after interchanging the carriers proves \eqref{eq:pseudometric-uniform-comparison}.

\textbf{Proof of \ref{item:quotient-3}.
Interpolation.}
Fix $(x,y),(z,w)\in\Correspondence$.
The difference of the interpolated pseudometrics at these pairs is
\begin{equation}\label{eq:pseudometric-pointwise-interpolation}
 (1-\HomotopyTime)
   \bigl(\Pseudometric_0(x,z)
       -\ComparisonPseudometric_0(y,w)\bigr)
 +\HomotopyTime
   \bigl(\Pseudometric_1(x,z)
       -\ComparisonPseudometric_1(y,w)\bigr).
\end{equation}
The absolute value of \eqref{eq:pseudometric-pointwise-interpolation} is bounded by the weighted sum of the two errors, since both weights are nonnegative.
Taking the supremum and using \eqref{eq:correspondence-pseudometric-error} proves \eqref{eq:pseudometric-interpolation}.

\textbf{Proof of \ref{item:quotient-4}.
A common carrier.}
Use $\Correspondence=\{(\BasePoint,\BasePoint)\mid\BasePoint\in\BaseCarrier\}$.
By \eqref{eq:correspondence-pseudometric-error},
\begin{equation}\label{eq:pseudometric-diagonal-error}
 \PseudometricError_\Correspondence(\Pseudometric_0,\Pseudometric_1)
 =\norm{\Pseudometric_0-\Pseudometric_1}_\infty.
\end{equation}
Substituting \eqref{eq:pseudometric-diagonal-error} into \eqref{eq:pseudometric-correspondence-bound} proves \eqref{eq:pseudometric-bound}.
\end{proof}

\begin{lemma}[{Orbit metric of a compact group}]\label{lem:compact-isometric-orbit-distance}
Let a compact group $G$ act continuously by isometries on a metric space $\MetricSpace{S}{d}$, and let $\OrbitProjection\colon S\to S/G$ be the orbit map.
Equip $S/G$ with the Hausdorff distance between orbits, the \emph{orbit metric} $\OrbitMetric(Gx,Gy)=\HausdorffDistance{d}(Gx,Gy)$.
Then for all $x,y\in S$ and $r>0$,
\[
 \OrbitMetric(Gx,Gy)=\min_{g\in G}d(x,gy),\qquad
 \MetricBall(\OrbitProjection(x),r;\OrbitMetric)=\OrbitProjection\bigl(\MetricBall(x,r;d)\bigr).
\]
In particular, $\OrbitMetric$ is a metric inducing the quotient topology.
\end{lemma}
\begin{proof}
Orbits are compact, so distinct orbits are at positive Hausdorff distance.
Palais (\cite[proof of Theorem~4.3.4, p.~319]{Palais1961}) and Antonyan (\cite[p.~4, (2.1)]{Antonyan2020Euclidean}) also use this orbit metric.
For each $g\in G$ we have $\dist_d(gx,Gy)=\dist_d(x,Gy)$.
Replacing a group element by its inverse shows that $\dist_d(y,Gx)$ has the same value.
Hence $\HausdorffDistance{d}(Gx,Gy)=\dist_d(x,Gy)$, and compactness of $Gy$ ensures attainment of the minimum.
Thus $Gy$ belongs to the open ball about $\OrbitProjection(x)$ of radius $r$ precisely when $Gy$ meets $\MetricBall(x,r;d)$.
This proves the open-ball identity.
Moreover, $\OrbitMetric(Gx,Gy)\leq d(x,y)$ implies continuity of $\OrbitProjection$, and the open-ball identity implies openness.
Hence $\OrbitMetric$ induces the quotient topology.
\end{proof}

\subsection{Probability measures and weak convergence}\label{subsec:preliminaries-probability}
We collect the measure-theoretic tools used for averaging probabilities and comparing distance operators on varying compact spaces.

For a metrizable space $\ProbabilityCarrier $, write $\BorelSets(\ProbabilityCarrier )$ for its \emph{Borel $\SigmaAlgebraPrefix$-algebra}, the $\SigmaAlgebraPrefix$-algebra generated by the open subsets of $\ProbabilityCarrier$.
A finite Borel measure $\FirstMeasure$ on $\ProbabilityCarrier$ is \emph{Radon} if it is \emph{inner regular} by compact sets, that is, if for every $\BorelSubset\in\BorelSets(\ProbabilityCarrier)$,
\[
 \FirstMeasure(\BorelSubset)
 =\sup\{\FirstMeasure(\CompactSubset)\mid
          \CompactSubset\subset\BorelSubset\text{ is compact}\}.
\]
Write $\ProbabilityMeasures(\ProbabilityCarrier)$ for the Radon measures of total mass one.
A \emph{probability law} on $\ProbabilityCarrier$ is an element of $\ProbabilityMeasures(\ProbabilityCarrier)$.
For Radon probabilities $\FirstMeasure\in\ProbabilityMeasures(\ProbabilityCarrier)$ and $\SecondMeasure\in\ProbabilityMeasures(\SecondProbabilityCarrier)$, write $\FirstMeasure\otimes\SecondMeasure$ for their \emph{product probability measure} on $\ProbabilityCarrier\times\SecondProbabilityCarrier$.
For all $\BorelSubset\in\BorelSets(\ProbabilityCarrier)$ and $\ClosedTestSet\in\BorelSets(\SecondProbabilityCarrier)$ it satisfies
\[
 (\FirstMeasure\otimes\SecondMeasure)(\BorelSubset\times\ClosedTestSet)
 =\FirstMeasure(\BorelSubset)\SecondMeasure(\ClosedTestSet).
\]
These identities characterize the product probability measure.
Write $\supp(\FirstMeasure)$ for the \emph{support} of $\FirstMeasure\in\ProbabilityMeasures(\ProbabilityCarrier)$.
A point $\BasePoint\in\ProbabilityCarrier$ belongs to $\supp(\FirstMeasure)$ if $\FirstMeasure(\OpenSubset)>0$ for every open neighborhood $\OpenSubset$ of $\BasePoint$.
We say $\FirstMeasure$ \emph{has full support} if $\supp(\FirstMeasure)=\ProbabilityCarrier$.
Equivalently,
\[
 \forall\OpenSubset\subset\ProbabilityCarrier\text{ open and nonempty}\quad
 \FirstMeasure(\OpenSubset)>0.
\]
If $\Mapping \colon \ProbabilityCarrier \to \SecondProbabilityCarrier $ is continuous and $\FirstMeasure \in\ProbabilityMeasures(\ProbabilityCarrier )$, its \emph{pushforward} is the probability measure $\Mapping \Pushforward\FirstMeasure $ defined for every $\BorelSubset\in\BorelSets(\SecondProbabilityCarrier)$ by
\[
 (\Mapping\Pushforward\FirstMeasure)(\BorelSubset)
 =\FirstMeasure(\Mapping^{-1}(\BorelSubset)).
\]
More generally, if $\Mapping\colon\ProbabilityCarrier\to\SecondProbabilityCarrier$ is Borel measurable, write $\Mapping\Pushforward\FirstMeasure$ for the Borel measure defined by the same preimage formula.
For a Borel measurable function $\TestFunction\colon\SecondProbabilityCarrier\to[0,\infty]$,
\begin{equation}\label{eq:pushforward-integration}
 \int_{\SecondProbabilityCarrier}\TestFunction(y)\,
       \IntegrationDifferential(\Mapping\Pushforward\FirstMeasure)(y)
 =\int_{\ProbabilityCarrier}\TestFunction(\Mapping(x))\,
       \IntegrationDifferential\FirstMeasure(x).
\end{equation}
For indicators this is the definition of the pushforward.
Linearity proves it for nonnegative simple functions, and monotone convergence proves the general nonnegative case.
For a real- or complex-valued Borel function, $\TestFunction$ is integrable with respect to $\Mapping\Pushforward\FirstMeasure$ if and only if $\TestFunction\circ\Mapping$ is integrable with respect to $\FirstMeasure$, and the same identity then holds.

For a nonempty metrizable space $\ProbabilityCarrier$, write
\[
 \BoundedContinuousFunctions(\ProbabilityCarrier)
 =\{\TestFunction\colon\ProbabilityCarrier\to\RealNumbers\mid
       \TestFunction\text{ is bounded and continuous}\},
\]
\[
 \norm{\TestFunction}_\infty=\sup_{\BasePoint\in\ProbabilityCarrier}|\TestFunction(\BasePoint)|.
\]
For a compact metric space $\BaseCarrier$ put
\[
 \ContinuousFunctions(\BaseCarrier)=\ContinuousFunctions(\BaseCarrier,\RealNumbers)
 =\{\TestFunction\colon\BaseCarrier\to\RealNumbers\mid\TestFunction\text{ is continuous}\},
\]
\[
 \norm{\TestFunction}_\infty=\max_{\BasePoint\in\BaseCarrier}|\TestFunction(\BasePoint)|.
\]
This is a real Banach space.
Continuous functions with values in $\ComplexNumbers$ form the complex Banach space $\ContinuousFunctions(\BaseCarrier,\ComplexNumbers)$.

Equip $\ProbabilityMeasures(\ProbabilityCarrier)$ with the \emph{weak topology}, the coarsest topology for which, for each $\TestFunction\in\BoundedContinuousFunctions(\ProbabilityCarrier)$, the map
\[\FirstMeasure\longmapsto\int_\ProbabilityCarrier\TestFunction\,
 \IntegrationDifferential\FirstMeasure
\]
defined by
\[
 \ProbabilityMeasures(\ProbabilityCarrier)\to\RealNumbers
 \]
is continuous.
We write $\FirstMeasure_\SequenceIndex\to\FirstMeasure$ for \emph{weak convergence}, meaning that for every $\TestFunction\in\BoundedContinuousFunctions(\ProbabilityCarrier)$,
\[
 \int_\ProbabilityCarrier\TestFunction\,\IntegrationDifferential\FirstMeasure_\SequenceIndex
 \longrightarrow
 \int_\ProbabilityCarrier\TestFunction\,\IntegrationDifferential\FirstMeasure.
\]
On a compact space all continuous functions are bounded.
For bounded complex-valued continuous test functions the same convergence follows by applying the definition to their real and imaginary parts.

We construct averages of measures by specifying their integrals of continuous functions.
The Riesz--Markov--Kakutani theorem states that every positive normalized linear functional on the real continuous functions on a nonempty compact metric space is integration against a unique Radon probability measure.

\begin{theorem}[{Riesz--Markov--Kakutani, \cite[Theorem~7.10.4 and Corollary~7.10.5]{Bogachev2007II}}]\label{thm:riesz-markov-kakutani}
Let $\MetricSpace{\BaseCarrier}{\MetricSymbol}$ be a nonempty compact metric space.
Let
\[
 \RepresentingFunctional\colon\ContinuousFunctions(\BaseCarrier,\RealNumbers)\to\RealNumbers
\]
be a linear functional.
Assume that for every $\TestFunction\in\ContinuousFunctions(\BaseCarrier,\RealNumbers)$ with $\TestFunction\geq0$ we have
\[
 \RepresentingFunctional(\TestFunction)\geq0.
\]
Let $1$ denote the constant function with value one and assume that
\[
 \RepresentingFunctional(1)=1.
\]
Then there is a unique Radon probability $\FirstMeasure\in\ProbabilityMeasures(\BaseCarrier)$ such that for every $\TestFunction\in\ContinuousFunctions(\BaseCarrier,\RealNumbers)$, we have
\[
 \RepresentingFunctional(\TestFunction)
 =\int_\BaseCarrier\TestFunction\,\IntegrationDifferential\FirstMeasure.
\]
\end{theorem}

Positivity and normalization imply $|\RepresentingFunctional(\TestFunction)|\leq\norm{\TestFunction}_\infty$, so the functional $\RepresentingFunctional$ is automatically continuous.
The cited representation theorem applies, and evaluating at the constant function $1$ proves that the representing measure $\FirstMeasure$ has total mass one.

\begin{theorem}[{A.~D.~Alexandroff's weak convergence theorem (Portmanteau theorem), \cite[Theorem~2.1]{Billingsley1999} and \cite[Theorem~8.2.3, Corollary~8.2.4(a), and p.~184]{Bogachev2007II}}]\label{thm:portmanteau}
Let $\ProbabilityCarrier$ be a metrizable space and let $\FirstMeasure,\FirstMeasure_\SequenceIndex\in\ProbabilityMeasures(\ProbabilityCarrier)$ for $\SequenceIndex\in\NonnegativeIntegers$.
Then the following conditions are equivalent.
\begin{enumerate}[label=\textup{(P\arabic*)},ref=\textup{(P\arabic*)}]
\item\label{item:portmanteau-weak}
For every $\TestFunction\in\BoundedContinuousFunctions(\ProbabilityCarrier)$,
\[
 \lim_{\SequenceIndex\to\infty}\int_{\ProbabilityCarrier}\TestFunction\,\IntegrationDifferential\FirstMeasure_\SequenceIndex
 =\int_{\ProbabilityCarrier}\TestFunction\,\IntegrationDifferential\FirstMeasure.
\]
\item
For every closed subset $\ClosedTestSet\subset\ProbabilityCarrier$,
\[
 \limsup_{\SequenceIndex\to\infty}\FirstMeasure_\SequenceIndex(\ClosedTestSet)
 \leq\FirstMeasure(\ClosedTestSet).
\]
\item
For every open subset $\OpenSubset\subset\ProbabilityCarrier$,
\[
 \FirstMeasure(\OpenSubset)
 \leq\liminf_{\SequenceIndex\to\infty}\FirstMeasure_\SequenceIndex(\OpenSubset).
\]
\end{enumerate}
\end{theorem}

We also use the following consequence for measures on varying compact sets.

\begin{lemma}\label{lem:limit-support}
Let $\FirstCompactSet_\SequenceIndex,\FirstCompactSet$ be nonempty compact subsets of a compact metric space $\MetricSpace{\AmbientCarrier}{\AmbientMetric}$ with $\HausdorffDistance{\AmbientMetric}(\FirstCompactSet_\SequenceIndex,\FirstCompactSet)\to0$.
If $\FirstMeasure_\SequenceIndex\to\FirstMeasure$ weakly in $\ProbabilityMeasures(\AmbientCarrier)$ and $\FirstMeasure_\SequenceIndex(\FirstCompactSet_\SequenceIndex)=1$, then $\FirstMeasure(\FirstCompactSet)=1$.
\end{lemma}

\begin{proof}
Define the function $\TestFunction\colon\AmbientCarrier\to\RealNumbers$ by
\[
 \TestFunction(\IntegrationPoint)
 =\dist_{\AmbientMetric}(\IntegrationPoint,\FirstCompactSet).
\]
Since $\AmbientCarrier$ is compact, $\TestFunction$ is bounded and continuous.
Weak convergence and $\FirstMeasure_\SequenceIndex(\FirstCompactSet_\SequenceIndex)=1$ imply
\[
 \begin{aligned}
 0\leq\int_\AmbientCarrier\TestFunction(\IntegrationPoint)
       \,\IntegrationDifferential\FirstMeasure(\IntegrationPoint)
 &=\lim_\SequenceIndex\int_\AmbientCarrier\TestFunction(\IntegrationPoint)
       \,\IntegrationDifferential\FirstMeasure_\SequenceIndex(\IntegrationPoint)\\
 &\leq\lim_\SequenceIndex
       \HausdorffDistance{\AmbientMetric}(\FirstCompactSet_\SequenceIndex,\FirstCompactSet)
 =0.
 \end{aligned}
\]
Since $\TestFunction^{-1}(\{0\})=\FirstCompactSet$, we obtain $\FirstMeasure(\FirstCompactSet)=1$.
\end{proof}

Let $\MetricSpace{\AmbientCarrier}{\AmbientMetric}$ be a compact metric space.
For $\BorelSubset\in\BorelSets(\AmbientCarrier)$ and $\SmallError>0$, put
\[
 \BorelSubset^{\SmallError}
 =\{\BasePoint\in\AmbientCarrier\mid
       \dist_{\AmbientMetric}(\BasePoint,\BorelSubset)<\SmallError\}.
\]
For probabilities $\FirstMeasure,\SecondMeasure\in\ProbabilityMeasures(\AmbientCarrier)$, we define the \emph{L\'evy--Prokhorov distance} by
\[
 \ProkhorovDistance{\AmbientMetric}(\FirstMeasure,\SecondMeasure)
 =\inf\left\{\SmallError>0\;\middle|\;
 \begin{gathered}
 \text{for every }\BorelSubset\in\BorelSets(\AmbientCarrier),\\
 \FirstMeasure(\BorelSubset)\leq\SecondMeasure(\BorelSubset^{\SmallError})+\SmallError,\\
 \SecondMeasure(\BorelSubset)\leq\FirstMeasure(\BorelSubset^{\SmallError})+\SmallError
 \end{gathered}\right\}.
\]
By Prokhorov's theorem (\cite[Section~1]{Prokhorov1956Weak}), this is a metric on $\ProbabilityMeasures(\AmbientCarrier)$ inducing weak convergence (\cite[Section~6, pp.~72--73, Theorem~6.8]{Billingsley1999}, \cite[author's manuscript, Definition~1.12, Proposition~1.13 and Lemma~1.15]{ShioyaMetricMeasure}).

\begin{lemma}[{Khezeli~\cite[arXiv v5, Example~2.1(iii)]{Khezeli2023Framework}}]\label{lem:prokhorov-isometric-embedding}
Let $\MetricSpace{\BaseCarrier}{\MetricSymbol}$ and $\MetricSpace{\AmbientCarrier}{\AmbientMetric}$ be compact metric spaces, and let $\FirstEmbedding\colon\MetricSpace{\BaseCarrier}{\MetricSymbol}\to\MetricSpace{\AmbientCarrier}{\AmbientMetric}$ be an isometric embedding.
For all $\FirstMeasure,\SecondMeasure\in\ProbabilityMeasures(\BaseCarrier)$, we have
\[
 \ProkhorovDistance{\AmbientMetric}
 (\FirstEmbedding\Pushforward\FirstMeasure,\FirstEmbedding\Pushforward\SecondMeasure)
 =\ProkhorovDistance{\MetricSymbol}(\FirstMeasure,\SecondMeasure).
\]
\end{lemma}


By Ulam's tightness theorem, every Borel probability measure on a Polish space is Radon (\cite[Theorem~7.1.7]{Bogachev2007II}).
Varadarajan proved metrizability of the weak topology on probability measures over a separable metric space.
Prokhorov proved complete metrizability over a complete separable metric space and characterized relative compactness by tightness.
In particular, when the underlying space is compact, all probability measures are uniformly tight and the space of probability measures is compact.


\begin{theorem}[{Varadarajan~\cite{Varadarajan1958Weak} and Prokhorov~\cite[Section~1]{Prokhorov1956Weak}, see also \cite[Theorems~8.9.3 and~8.9.4]{Bogachev2007II}}]\label{thm:probability-topology}
Let $\ProbabilityCarrier$ be a separable metrizable space.
Then $\ProbabilityMeasures(\ProbabilityCarrier)$ is separable and metrizable in the weak topology.
If $\ProbabilityCarrier$ is Polish, then $\ProbabilityMeasures(\ProbabilityCarrier)$ is Polish.
If $\ProbabilityCarrier$ is compact, then $\ProbabilityMeasures(\ProbabilityCarrier)$ is compact.
\end{theorem}


Weak convergence is uniform over compact families of continuous test functions.

\begin{lemma}\label{lem:uniform-weak-integrals}
Let $\AmbientCarrier $ be a nonempty compact metric space, let $\FirstMeasure _\SequenceIndex \to\FirstMeasure $ weakly in $\ProbabilityMeasures(\AmbientCarrier )$, and let $\CompactTestFamily \subset\ContinuousFunctions(\AmbientCarrier )$ be nonempty and compact in the uniform norm.
Then
\begin{equation}\label{eq:uniform-test-integrals}
 \sup_{\Mapping \in\CompactTestFamily }\left|\int_\AmbientCarrier  \Mapping \,\IntegrationDifferential \FirstMeasure _\SequenceIndex -\int_\AmbientCarrier  \Mapping \,\IntegrationDifferential \FirstMeasure \right|\to0.
\end{equation}
The same conclusion holds for complex continuous functions.
\end{lemma}

\begin{proof}
For each $\SequenceIndex$, choose $\Mapping_\SequenceIndex\in\CompactTestFamily$ attaining the supremum in \eqref{eq:uniform-test-integrals}.
By compactness of $\CompactTestFamily$, every subsequence has a further subsequence $\{\Mapping_{n_k}\}_{k\in\NonnegativeIntegers}$ converging uniformly to some $\Mapping\in\CompactTestFamily$.
For these indices,
\[
 \left|\int_\AmbientCarrier\Mapping_{n_k}\,
           \IntegrationDifferential(\FirstMeasure_{n_k}-\FirstMeasure)\right|
 \leq2\norm{\Mapping_{n_k}-\Mapping}_\infty
      +\left|\int_\AmbientCarrier\Mapping\,
           \IntegrationDifferential(\FirstMeasure_{n_k}-\FirstMeasure)\right|
 \longrightarrow0.
\]
Thus the whole sequence of suprema tends to $0$.
The same argument applies to complex functions, since weak convergence controls both their real and imaginary parts.
\end{proof}

We use \Cref{lem:uniform-weak-integrals} to prove uniform convergence of integrals whose integrands vary uniformly with compact parameters.

\begin{lemma}\label{lem:parameter-integrals}
Let $\CompactParameterSpace$ and $\AmbientCarrier$ be nonempty compact metric spaces.
Let $\FirstMeasure_\SequenceIndex\to\FirstMeasure$ weakly in $\ProbabilityMeasures(\AmbientCarrier)$.
Let $\Integrand_\SequenceIndex,\Integrand\colon\CompactParameterSpace\times\AmbientCarrier\to\ComplexNumbers$ be continuous real-valued or complex-valued functions.
Assume that $\norm{\Integrand_\SequenceIndex-\Integrand}_\infty\to0$.
Then
\[
 \sup_{\ParameterPoint\in\CompactParameterSpace}
 \left|\int_\AmbientCarrier\Integrand_\SequenceIndex(\ParameterPoint,\IntegrationPoint)\,\IntegrationDifferential\FirstMeasure_\SequenceIndex(\IntegrationPoint)
 -\int_\AmbientCarrier\Integrand(\ParameterPoint,\IntegrationPoint)\,\IntegrationDifferential\FirstMeasure(\IntegrationPoint)\right|\to0.
\]
\end{lemma}

\begin{proof}
The triangle inequality yields
\begin{equation}\label{eq:parameter-integral-error}
\begin{aligned}
&\sup_{\ParameterPoint\in\CompactParameterSpace}
 \left|
 \int_\AmbientCarrier
 \Integrand_\SequenceIndex(\ParameterPoint,\IntegrationPoint)
 \,\IntegrationDifferential\FirstMeasure_\SequenceIndex(\IntegrationPoint)
 -\int_\AmbientCarrier
 \Integrand(\ParameterPoint,\IntegrationPoint)
 \,\IntegrationDifferential\FirstMeasure(\IntegrationPoint)
 \right|\\
&\leq
 \norm{\Integrand_\SequenceIndex-\Integrand}_\infty
 +\sup_{\ParameterPoint\in\CompactParameterSpace}
 \left|
 \int_\AmbientCarrier
 \Integrand(\ParameterPoint,\IntegrationPoint)
 \,\IntegrationDifferential\FirstMeasure_\SequenceIndex(\IntegrationPoint)
 -\int_\AmbientCarrier
 \Integrand(\ParameterPoint,\IntegrationPoint)
 \,\IntegrationDifferential\FirstMeasure(\IntegrationPoint)
 \right|.
\end{aligned}
\end{equation}
The uniform error $\norm{\Integrand_\SequenceIndex-\Integrand}_\infty$ in \eqref{eq:parameter-integral-error} tends to $0$ by assumption.
The family $\{\Integrand(\ParameterPoint,\cdot)\mid\ParameterPoint\in\CompactParameterSpace\}$ is compact in the uniform norm, so \Cref{lem:uniform-weak-integrals} proves that the supremum of the integral differences on the right-hand side of \eqref{eq:parameter-integral-error} tends to $0$.
By \eqref{eq:parameter-integral-error}, the parameter integrals converge uniformly on $\CompactParameterSpace$.
\end{proof}

We use the following continuous functions to detect the support of a measure and approximate a distance function by an integral.

\begin{lemma}\label{lem:metric-bump}
Let $\MetricSpace{\BaseCarrier}{\MetricSymbol}$ be a nonempty compact metric space, $\FirstMeasure\in\ProbabilityMeasures(\BaseCarrier)$, and $\Radius>0$.
For $\BasePoint\in\BaseCarrier$, define the \emph{bump function}
\[
 \BumpFunction{\BaseCarrier}{\MetricSymbol}{\BasePoint}{\Radius}
 \colon\BaseCarrier\to[0,1]
\]
by
\begin{equation}\label{eq:metric-bump-definition}
 \BumpFunction{\BaseCarrier}{\MetricSymbol}{\BasePoint}{\Radius}(\ComparisonPoint)
 =\max\{1-\MetricSymbol(\BasePoint,\ComparisonPoint)/\Radius,0\}.
\end{equation}
Then for all $\BasePoint,\BasePoint_0,\ComparisonPoint,\ComparisonPoint_0\in\BaseCarrier$,
\begin{equation}\label{eq:metric-bump-lipschitz}
 \left|\BumpFunction{\BaseCarrier}{\MetricSymbol}{\BasePoint}{\Radius}(\ComparisonPoint)
 -\BumpFunction{\BaseCarrier}{\MetricSymbol}{\BasePoint_0}{\Radius}(\ComparisonPoint_0)\right|
 \leq\frac{\MetricSymbol(\BasePoint,\BasePoint_0)
              +\MetricSymbol(\ComparisonPoint,\ComparisonPoint_0)}{\Radius}.
\end{equation}
Moreover, for every $\BasePoint\in\BaseCarrier$,
\begin{equation}\label{eq:metric-bump-integral}
 \frac12\FirstMeasure\bigl(\MetricBall(\BasePoint,\Radius/2;\MetricSymbol)\bigr)
 \leq\int_\BaseCarrier
       \BumpFunction{\BaseCarrier}{\MetricSymbol}{\BasePoint}{\Radius}
       \,\IntegrationDifferential\FirstMeasure
 \leq\FirstMeasure\bigl(\MetricBall(\BasePoint,\Radius;\MetricSymbol)\bigr).
\end{equation}
\end{lemma}

\begin{proof}
The function $\SecondParameter\mapsto\max\{1-\SecondParameter/\Radius,0\}$ is $1/\Radius$-Lipschitz.
Thus the triangle inequality for $\MetricSymbol$ implies \eqref{eq:metric-bump-lipschitz}.
Integrating the following pointwise estimates yields \eqref{eq:metric-bump-integral}.
\[
 \tfrac12\Indicator_{\MetricBall(\BasePoint,\Radius/2;\MetricSymbol)}
 \leq\BumpFunction{\BaseCarrier}{\MetricSymbol}{\BasePoint}{\Radius}
 \leq\Indicator_{\MetricBall(\BasePoint,\Radius;\MetricSymbol)}\qedhere.
\]
\end{proof}

\begin{lemma}\label{lem:normalized-metric-bump}
Let $\MetricSpace{\BaseCarrier}{\MetricSymbol}$ be a nonempty compact metric space, let $\FirstMeasure$ be a full-support probability measure on $\BaseCarrier$, and let $\Radius>0$.
For every $\BasePoint\in\BaseCarrier$ the following \emph{normalized bump function} is well defined and continuous and has the three properties below.
\begin{equation}\label{eq:normalized-bump-definition}
 \NormalizedBump{\BaseCarrier}{\MetricSymbol}{\BasePoint}{\Radius}{\FirstMeasure}
 =\frac{\BumpFunction{\BaseCarrier}{\MetricSymbol}{\BasePoint}{\Radius/2}}
 {\displaystyle\int_\BaseCarrier
      \BumpFunction{\BaseCarrier}{\MetricSymbol}{\BasePoint}{\Radius/2}
      \,\IntegrationDifferential\FirstMeasure}
\end{equation}
\begin{equation}\label{eq:normalized-bump-properties-1}
 \NormalizedBump{\BaseCarrier}{\MetricSymbol}{\BasePoint}{\Radius}{\FirstMeasure}\geq0,
\end{equation}
\begin{equation}\label{eq:normalized-bump-properties-2}
 \int_\BaseCarrier
 \NormalizedBump{\BaseCarrier}{\MetricSymbol}{\BasePoint}{\Radius}{\FirstMeasure}
 \,\IntegrationDifferential\FirstMeasure=1,
\end{equation}
\begin{equation}\label{eq:normalized-bump-properties-3}
 \supp\left(\NormalizedBump{\BaseCarrier}{\MetricSymbol}{\BasePoint}{\Radius}{\FirstMeasure}\right)
 \subset\MetricBall(\BasePoint,\Radius;\MetricSymbol).
\end{equation}
\end{lemma}

\begin{proof}
Fix $\BasePoint\in\BaseCarrier$.
Since $\supp(\FirstMeasure)=\BaseCarrier$, \eqref{eq:metric-bump-integral} with radius $\Radius/2$ proves that the denominator in \eqref{eq:normalized-bump-definition} is positive.
Continuity and \eqref{eq:normalized-bump-properties-1}--\eqref{eq:normalized-bump-properties-2} follow from \eqref{eq:metric-bump-lipschitz} and \eqref{eq:normalized-bump-definition}.
By \eqref{eq:metric-bump-definition},
\[
 \supp\left(\NormalizedBump{\BaseCarrier}{\MetricSymbol}{\BasePoint}{\Radius}{\FirstMeasure}\right)
 \subset\{\ComparisonPoint\in\BaseCarrier\mid
              \MetricSymbol(\BasePoint,\ComparisonPoint)\leq\Radius/2\}
 \subset\MetricBall(\BasePoint,\Radius;\MetricSymbol).
\]
This proves \eqref{eq:normalized-bump-properties-3}.
\end{proof}

The next lemma compares product measures when one factor is a Dirac measure.
It will be used to average laws over a convergent family of compact spaces.

\begin{lemma}[{Consequence of Giry's continuity of product measures~\cite[pp.~78--79]{Giry1982Probability}, see also \cite[Chapter~III, Lemma~1.1]{Parthasarathy2005}}]\label{lem:dirac-product-convergence}
Let $\CompactParameterSpace$ be a compact metric space and $\ProbabilityCarrier$ a Polish space.
Suppose $\ParameterPoint_\SequenceIndex\to\ParameterPoint$ in $\CompactParameterSpace$ and $\FirstMeasure_\SequenceIndex\to\FirstMeasure$ weakly in $\ProbabilityMeasures(\ProbabilityCarrier)$.
For $\ParameterPoint\in\CompactParameterSpace$ write $\DiracProbability_\ParameterPoint$ for the \emph{Dirac measure} concentrated at $\ParameterPoint$.
Then the following weak convergence holds in $\ProbabilityMeasures(\CompactParameterSpace\times\ProbabilityCarrier)$.
\begin{equation}\label{eq:dirac-product-convergence}
 \DiracProbability_{\ParameterPoint_\SequenceIndex}\otimes\FirstMeasure_\SequenceIndex
 \longrightarrow\DiracProbability_\ParameterPoint\otimes\FirstMeasure.
\end{equation}
\end{lemma}


For comparisons with metric measure geometry, we also recall the box distance.

\begin{definition}\label{def:box-distance}
An \emph{mm-space} is a complete separable metric space $(X,d_X)$ equipped with a Borel probability measure $\mu_X$.
Write $\mathbf X=(X,d_X,\mu_X)$ for the triple.
Two mm-spaces are \emph{mm-isomorphic} if there is a measure-preserving isometry between their supports.
Denote the space of classes by $\MMClassSpace$.
Put $\BoxParameterInterval=[0,1)$ and denote Lebesgue measure on this interval by $\BoxLebesgueMeasure$.
A \emph{parameter} of $\mathbf X$ is a Borel map $\BoxFirstParameter\colon\BoxParameterInterval\to X$ satisfying $(\BoxFirstParameter)_*\BoxLebesgueMeasure=\mu_X$.
Every mm-space admits a parameter.
The \emph{box distance} $\BoxDistance(\mathbf X,\mathbf Y)$ between mm-spaces $\mathbf X$ and $\mathbf Y=(Y,d_Y,\mu_Y)$ is the infimum of the numbers $\ErrorTolerance>0$ for which there exist a parameter $\BoxFirstParameter$ of $\mathbf X$, a parameter $\BoxSecondParameter$ of $\mathbf Y$, and a Borel set $\BoxLargeSubset\subset\BoxParameterInterval$ such that
\[
 \BoxLebesgueMeasure(\BoxLargeSubset)\geq1-\ErrorTolerance,
\]
and such that for every $s,t\in\BoxLargeSubset$ we have
\[
 \left|d_X(\BoxFirstParameter(s),\BoxFirstParameter(t))
       -d_Y(\BoxSecondParameter(s),\BoxSecondParameter(t))\right|
 \leq\ErrorTolerance.
\]
It is a metric on $\MMClassSpace$.
The topology it induces is called the \emph{box topology} (see Kazukawa, Nakajima, and Shioya~\cite[arXiv v1, Definitions~2.1, 2.6, 2.7 and Theorem~2.8]{KazukawaNakajimaShioya2024}).
\end{definition}
\subsection{Partitions of unity, contractions and absolute retracts}\label{subsec:preliminaries-topology}
We recall partitions of unity, an extension theorem for continuous real functions, and the terminology for contractions and absolute retracts.

A family $\{A_i\}_{i\in I}$ of subsets of a topological space is \emph{locally finite} if every point has an open neighborhood $U$ such that $U\cap A_i\ne\varnothing$ for only finitely many $i\in I$.
One cover \emph{refines} another if every member of the first is contained in a member of the second.
A Hausdorff space is \emph{paracompact} if every open cover has a locally finite open refinement.

Let $X$ be a topological space and $I$ an index set.
A family of continuous functions $\PartitionWeight_i\colon X\to[0,1]$ indexed by $i\in I$ is a \emph{locally finite partition of unity} if the following supports form a locally finite family and the displayed sum equals one for every $x\in X$.
\[
 \supp\PartitionWeight_i=\overline{\{x\in X\mid\PartitionWeight_i(x)\ne0\}}
\]
\[
 \sum_{i\in I}\PartitionWeight_i(x)=1.
\]
It is \emph{subordinate} to an open cover if each support is contained in some member of the cover.

\begin{theorem}[{A.~H.~Stone~\cite[Corollary~1]{Stone1948Paracompactness} and Michael~\cite[Proposition~2]{Michael1953Paracompact}, see also \cite[Theorems~1.3.33 and~1.3.38]{CobzasMiculescuNicolae2019}}]\label{thm:metric-partition-unity}
By A.~H.~Stone's theorem every metrizable space is paracompact.
By Michael's characterization of paracompactness, every open cover admits a subordinate locally finite partition of unity.
\end{theorem}

\begin{theorem}[{Tietze--Urysohn, \cite[Theorem~2.1.8]{Engelking1989}}]\label{thm:tietze}
Let $\ExtensionDomain$ be a normal space, let $\ClosedDomain\subset\ExtensionDomain$ be closed, and let $\FunctionBound\geq0$.
Then every continuous function $\Mapping\colon\ClosedDomain\to[-\FunctionBound,\FunctionBound]$ has a continuous extension $\TietzeExtension\colon\ExtensionDomain\to[-\FunctionBound,\FunctionBound]$.
\end{theorem}

A continuous map $\Mapping \colon \ClosedDomain \to \ExtensionTarget $ is \emph{null-homotopic} if there exist a point $\ConstantPoint \in \ExtensionTarget $ and a continuous map $\ContractionHomotopy \colon \ClosedDomain \times\UnitInterval\to \ExtensionTarget $ such that for every $\FirstScale\in\ClosedDomain$,
\[
 \ContractionHomotopy (\FirstScale ,0)=\Mapping (\FirstScale ),
\]
\[
 \ContractionHomotopy (\FirstScale ,1)=\ConstantPoint.
\]
A nonempty space is \emph{contractible} if its identity map is null-homotopic.

The following proposition shows that scaling distances defines a contraction.
Write $\SingletonSpace=\MetricSpace{\{0\}}{0}$ for the one-point metric space.

\begin{proposition}\label{lem:scaling-contraction}
Let $\IsometryClass{\BaseCarrier}{\MetricSymbol},\IsometryClass{\ComparisonCarrier}{\ComparisonMetric}\in\GHSpace$.
For any real numbers $\FirstScale,\SecondScale\geq0$,
\begin{equation}\label{eq:scaling}
 \GHDistance\bigl(\MetricQuotient{\BaseCarrier}{\FirstScale\MetricSymbol},
                   \MetricQuotient{\ComparisonCarrier}{\SecondScale\ComparisonMetric}\bigr)
 \leq \FirstScale\GHDistance\bigl(\IsometryClass{\BaseCarrier}{\MetricSymbol},
                                   \IsometryClass{\ComparisonCarrier}{\ComparisonMetric}\bigr)
       +\tfrac12|\FirstScale-\SecondScale|\,\diam\MetricSpace{\ComparisonCarrier}{\ComparisonMetric}.
\end{equation}
Consequently the map
\[
 \ContractionHomotopy\colon\GHSpace\times\UnitInterval\longrightarrow\GHSpace,\qquad
 \ContractionHomotopy(\IsometryClass{\BaseCarrier}{\MetricSymbol},\HomotopyTime)
 =\MetricQuotient{\BaseCarrier}{(1-\HomotopyTime)\MetricSymbol}
\]
is a strong deformation retraction of $\GHSpace$ onto the singleton $\{\SingletonSpace\}$.
In particular, $\GHSpace$ is contractible.
\end{proposition}

\begin{proof}
The triangle inequality and the rescaling identities (\cite[arXiv v1, Examples~6.31 and~6.32]{Tuzhilin2020}) imply \eqref{eq:scaling}, including when either scale factor vanishes.
Since the diameter is continuous on $\GHSpace$, this bound proves continuity of the following map even at $\FirstScale=0$.
\[
 \GHSpace\times[0,\infty)\to\GHSpace,
 \qquad
 {(\IsometryClass{\BaseCarrier}{\MetricSymbol},\FirstScale)
       \longmapsto\MetricQuotient{\BaseCarrier}{\FirstScale\MetricSymbol}}.
\]
The map $\ContractionHomotopy(-,0)$ is the identity.
At time $1$ the metric quotient of the zero pseudometric is a singleton, so $\ContractionHomotopy(\IsometryClass{\BaseCarrier}{\MetricSymbol},1)=\SingletonSpace$.
For every $\HomotopyTime\in\UnitInterval$ we have \mbox{$\ContractionHomotopy(\SingletonSpace,\HomotopyTime)=\SingletonSpace$.}\qedhere
\end{proof}

For balls in Gromov--Hausdorff space we use a different symbol,
\[
 \GHball(\MetricSpace{\BaseCarrier}{\MetricSymbol},\Radius)
 =\{\MetricSpace{\ComparisonCarrier}{\ComparisonMetric}\in\GHSpace\mid
 \GHDistance(\MetricSpace{\BaseCarrier}{\MetricSymbol},\MetricSpace{\ComparisonCarrier}{\ComparisonMetric})<\Radius\}.
\]
The correspondence formula also implies $\GHDistance(\MetricSpace{\BaseCarrier }{\MetricSymbol },\SingletonSpace)=\diam\MetricSpace{\BaseCarrier }{\MetricSymbol }/2$.

We state extension and retraction properties relative to a class of topological spaces.

\begin{definition}\label{def:extension}
Let $\TopologicalSpaceClass$ be a class of topological spaces.
A topological space $\ExtensionTarget$ is an \emph{absolute extensor for $\TopologicalSpaceClass$}, or AE for $\TopologicalSpaceClass$, if for every space $\BaseCarrier\in\TopologicalSpaceClass$, every closed subset $\ClosedDomain\subset\BaseCarrier$, and every continuous map $\Mapping\colon\ClosedDomain\to\ExtensionTarget$, there exists a continuous map $\ContinuousExtension\colon\BaseCarrier\to\ExtensionTarget$ such that $\ContinuousExtension|_\ClosedDomain=\Mapping$.
It is an \emph{absolute neighborhood extensor for $\TopologicalSpaceClass$}, or ANE for $\TopologicalSpaceClass$, if for every space $\BaseCarrier\in\TopologicalSpaceClass$, every closed subset $\ClosedDomain\subset\BaseCarrier$, and every continuous map $\Mapping\colon\ClosedDomain\to\ExtensionTarget$, there exist an open neighborhood $\OpenNeighborhood$ of $\ClosedDomain$ in $\BaseCarrier$ and a continuous map $\ContinuousExtension\colon\OpenNeighborhood\to\ExtensionTarget$ such that $\ContinuousExtension|_\ClosedDomain=\Mapping$.
A space $\ExtensionTarget\in\TopologicalSpaceClass$ is an \emph{absolute retract for $\TopologicalSpaceClass$}, or AR for $\TopologicalSpaceClass$, if for every space $\BaseCarrier\in\TopologicalSpaceClass$ and every embedding $\FirstEmbedding\colon\ExtensionTarget\to\BaseCarrier$ with closed image, there exists a continuous map $\RetractionMap\colon\BaseCarrier\to\ExtensionTarget$ such that $\RetractionMap\circ\FirstEmbedding=\IdentityMap_\ExtensionTarget$.
A space $\ExtensionTarget\in\TopologicalSpaceClass$ is an \emph{absolute neighborhood retract for $\TopologicalSpaceClass$}, or ANR for $\TopologicalSpaceClass$, if for every space $\BaseCarrier\in\TopologicalSpaceClass$ and every embedding $\FirstEmbedding\colon\ExtensionTarget\to\BaseCarrier$ with closed image, there exist an open neighborhood $\OpenNeighborhood$ of $\FirstEmbedding(\ExtensionTarget)$ in $\BaseCarrier$ and a continuous map $\RetractionMap\colon\OpenNeighborhood\to\ExtensionTarget$ such that $\RetractionMap\circ\FirstEmbedding=\IdentityMap_\ExtensionTarget$.
\end{definition}

Unless otherwise stated, ANR means ANR for all metrizable spaces, and the same convention applies to AE, ANE, and AR.
Spaces described by these unqualified terms are assumed to be metrizable.

\begin{theorem}[{Dugundji~\cite[Theorem~12.1]{Dugundji1958}}]\label{thm:retract-extensor}
For a metrizable target in the category of metrizable spaces, the AE and AR properties are equivalent, and the ANE and ANR properties are equivalent.
\end{theorem}

\begin{definition}\label{def:local-contractibility}
A space $\ExtensionTarget $ is \emph{locally contractible at} $\BasePoint \in \ExtensionTarget $ if for every open neighborhood $\LargerNeighborhood $ of $\BasePoint $ there is an open neighborhood $\SmallerNeighborhood \subset \LargerNeighborhood $ of $\BasePoint $ such that the inclusion $\SmallerNeighborhood \to \LargerNeighborhood $ is null-homotopic.
\end{definition}

The homotopy may be chosen to end at the constant map with value $\BasePoint $.
Indeed, evaluating the homotopy at $\BasePoint $ defines a path from $\BasePoint $ to the constant value, and its reversed path can be appended to the homotopy path of every point of $\SmallerNeighborhood$.
\subsection{Residual sets in the Gromov--Hausdorff space}\label{subsec:measures-preliminaries-rouyer}
A subset of a topological space is \emph{residual} if its complement is contained in a countable union of closed sets with empty interior.
Since $\GHSpace$ is complete (\Cref{lem:gh-basics}), Baire's theorem implies that every residual subset of $\GHSpace$ is dense.

\begin{theorem}[{Rouyer~\cite[arXiv v1, Theorems~2 and~4, pp.~7--8]{Rouyer2011}}]\label{thm:rouyer-generic-metric-spaces}
In $\GHSpace$ the isometry classes of compact metric spaces $\MetricSpace{Y}{e}$ satisfying both of the following conditions form a residual set.
First, $Y$ has no isolated points.
Second, if $x,y,z,w\in Y$ satisfy $x\ne y$, $z\ne w$, and $e(x,y)=e(z,w)$, then $\{x,y\}=\{z,w\}$.
\end{theorem}

\clearpage
\part{The topology of Gromov--Hausdorff space I: Continuous invariant measures}\label{part:measures}
\begin{quote}
\small
\noindent\textbf{Abstract.}
We construct a continuous section of the map that forgets the measure on a compact metric measure space.
The probability measure selected on each nonempty compact metric space has full support and is invariant under every self-isometry.
When the underlying spaces converge in Hausdorff distance in a common compact metric space, the selected measures converge weakly there.
The section is the composition of a Milyutin selection, which assigns a probability law to each fiber, and the map taking the barycenter of that law.
Convex interpolation between a measure and the selected measure on the same space defines a strong deformation retraction onto the image of the section.

\par\smallskip
\noindent \textbf{Keywords.} Gromov--Hausdorff space, Gromov--Hausdorff--Prokhorov distance, invariant probability measure, full support, continuous selection.

\par\smallskip
\noindent \textbf{2020 Mathematics Subject Classification.} Primary 54E35; Secondary 54C65, 60B05, 60B10.
\end{quote}

\section{Introduction to Part I}\label{sec:intro-measures}
We ask whether one can choose a probability measure on each nonempty compact metric space, compatibly with isometries and continuously as the space varies.
Here continuity has a concrete meaning.
If compact subsets $X_n$ of a compact metric space $\MetricSpace{Z}{h}$ converge to a compact subset $X$ in Hausdorff distance, then the measures selected on $X_n$, considered as measures on $Z$, must converge weakly to the measure selected on $X$.
We require each measure to have full support.

Larrieu described a construction of a Borel probability measure on each nonempty compact metric space that is invariant under isometries between open subsets (\cite[arXiv v4, Theorem~3]{Larrieu2014InvariantProbabilities}).
The construction uses a fixed generalized limit and establishes neither full support nor continuity as the underlying space varies.
Leinster and Roff studied probability measures maximizing a family of metric entropies (\cite[Theorem~1.1 and Definition~9.2]{LeinsterRoff2021Entropy}).
When multiplication of the distance by every sufficiently large factor yields a unique maximizing measure and these measures have a weak limit as the factor tends to infinity, they call the limit the \emph{uniform measure}.
Our question concerns the ordinary Gromov--Hausdorff topology on the underlying compact spaces and weak convergence of measures in a common ambient space.

In our main result, \Cref{thm:invariant-measures}, we assign to every nonempty compact metric space $\MetricSpace{\BaseCarrier}{\MetricSymbol}$ a Borel probability measure $\SelectedMeasure{\BaseCarrier}{\MetricSymbol}$ with full support, so that
\[
 \supp(\SelectedMeasure{\BaseCarrier}{\MetricSymbol})=\BaseCarrier.
\]
For every isometry $\IdentifyingIsometry\colon\MetricSpace{\BaseCarrier}{\MetricSymbol}
\to\MetricSpace{\ComparisonCarrier}{\ComparisonMetric}$, we have
\[
 \IdentifyingIsometry\Pushforward\SelectedMeasure{\BaseCarrier}{\MetricSymbol}
 =\SelectedMeasure{\ComparisonCarrier}{\ComparisonMetric}.
\]
Thus the assignment is compatible with isometries.
This compatibility means that the resulting map from isometry classes of compact metric spaces to isometry classes of metric measure spaces is well defined.
Let compact subspaces $\BaseCarrier_\SequenceIndex$ and $\BaseCarrier$ of a compact metric space $\MetricSpace{\AmbientCarrier}{\AmbientMetric}$ carry the restricted metrics $\MetricSymbol_\SequenceIndex$ and $\MetricSymbol$.
Then
\[
 \HausdorffDistance{\AmbientMetric}(\BaseCarrier_\SequenceIndex,\BaseCarrier)\to0
 \quad\Longrightarrow\quad
 \SelectedMeasure{\BaseCarrier_\SequenceIndex}{\MetricSymbol_\SequenceIndex}
 \to\SelectedMeasure{\BaseCarrier}{\MetricSymbol}.
\]
The right-hand side denotes weak convergence of the selected measures as probability measures on the common space $\AmbientCarrier$.
Since $\SelectedMeasure{X}{d}$ has full support, the natural map $\ContinuousFunctions(X)\to\LebesgueSpace^2(X,\SelectedMeasure{X}{d})$ is injective.
Indeed, for any nonzero continuous function $f$, its absolute value $|f|$ is bounded below by a positive constant on some nonempty open set.
Hence $\int_X|f|^2\,\IntegrationDifferential\SelectedMeasure{X}{d}>0$.
Compatibility with isometries and weak continuity allow us to compare integrals as the compact space varies.

On a fixed compact metric space $\MetricSpace{X}{d}$, averaging a full-support probability measure against Haar probability on $\Isom(X,d)$ produces an invariant probability measure of full support.
Continuity as $X$ varies requires a further argument.
By Rouyer's theorem (\Cref{thm:rouyer-generic-metric-spaces}), spaces with trivial isometry group are dense in the Gromov--Hausdorff space $\GHSpace$.
A limit of such spaces can have a nontrivial isometry group.
We therefore consider all invariant probability measures of full support.
Write $\InvariantMeasuredSpaces$ for the space of isomorphism classes of spaces carrying such a measure and $\MeasureProjection:\InvariantMeasuredSpaces\to\GHSpace$ for the map that forgets the measure.

For $m=\IsometryClass{X}{d}\in\GHSpace$, a point of $\MeasureProjection^{-1}(m)$ specifies an admissible measure up to measure-preserving isometry.
We first prove that $\MeasureProjection$ is a continuous open surjection between Polish spaces.
A consequence of a theorem of Banakh and Radul (\cite[Theorem~1.2, p.~20]{BanakhRadul1999}, \Cref{thm:valov}) then lets us choose continuously a probability law $\LawSectionAt{X}{d}$ on each fiber.
Let $\FiberLawSpace$ be the subspace of probability laws on $\InvariantMeasuredSpaces$ concentrated on a single fiber, with the relative weak topology.
The desired section is the composition
\[
 \GHSpace\xrightarrow{\ \LawSection\ }\FiberLawSpace
 \xrightarrow{\ \FiberBarycenterMap\ }\InvariantMeasuredSpaces.
\]
The map $\FiberBarycenterMap$ takes the barycenter of the law on a representative of the underlying space.
Figure~\ref{fig:probability-on-fibers} depicts the fibers of the forgetful map $\MeasureProjection$.
The identity $\MeasureProjection_*\LawSectionAt{X}{d}=\DiracProbability_{\IsometryClass{X}{d}}$ means that $\LawSectionAt{X}{d}$ is concentrated on this fiber.
Figure~\ref{fig:fiber-probability-density} depicts the laws selected on the individual fibers.
For example, on a two-point space $\MetricSpace{X}{d}$ the uniform measure $u$ is the only invariant probability measure.
Its fiber is the singleton $\MeasuredIsometryClass{X}{d}{u}$, and $\LawSectionAt{X}{d}$ is the Dirac measure at that point.
In contrast, a three-point space of side lengths $1,3/2,2$ has trivial isometry group.
Once the vertices are ordered, its fiber is in bijection with the triples of positive masses whose sum is $1$.

\begin{figure}[H]
\centering
\begin{tikzpicture}[font=\small,>=Stealth,x=1cm,y=1cm]
\def\FigureFiberSegments{.61/1.16/2.68/.04,1.72/1.18/2.88/.03,3.10/1.09/3.07/-.01,4.49/1.05/2.86/-.04,5.52/1.10/2.68/0}
\path[draw=black!60,fill=black!4,line width=.45pt]
  (0,1.08) .. controls (1.65,1.4) and (3.8,.88) .. (6,1.15)
  -- (6,2.70) .. controls (4.65,2.52) and (4.19,3.16) .. (3.10,3.07)
  .. controls (1.82,3.02) and (.88,2.52) .. (0,2.70) -- cycle;
\draw[semithick] (-.12,0) .. controls (1.4,.19) and (3.66,-.23) .. (6.13,.04);
\node[left=5pt] at (0,1.90) {\large$\InvariantMeasuredSpaces$};
\node[left=5pt] at (0,0) {\large$\GHSpace$};
\foreach \fiberx/\fiberlo/\fiberhi/\basey [count=\fiberindex from 0] in \FigureFiberSegments {
  \ifnum\fiberindex=2\relax
    \draw[line width=1.05pt] (\fiberx,\fiberlo) -- (\fiberx,\fiberhi);
    \draw[->,semithick] (\fiberx,\fiberlo-.05) -- node[right=3pt] {$\MeasureProjection$} (\fiberx,\basey+.10);
    \fill (\fiberx,\basey) circle[radius=1.9pt] node[below=4pt] {$\IsometryClass{X}{d}$};
  \else
    \draw[black!60] (\fiberx,\fiberlo) -- (\fiberx,\fiberhi);
    \draw[->,black!50,thin] (\fiberx,\fiberlo-.04) -- (\fiberx,\basey+.10);
    \fill[black!65] (\fiberx,\basey) circle[radius=1.55pt];
  \fi
}
  \draw[thin] (3.13,2.85) -- (3.56,3.33) -- (3.90,3.33);
  \node[right=2pt] at (3.90,3.33) {$\MeasureProjection^{-1}(\IsometryClass{X}{d})$};
\end{tikzpicture}
\caption{Fibers of the forgetful map $\MeasureProjection\colon\InvariantMeasuredSpaces\to\GHSpace$. Each point on the bold vertical line represents an isomorphism class $\MeasuredIsometryClass{X}{d}{\mu}$ with an isometry-invariant measure $\mu$ of full support. All these points map to $\IsometryClass{X}{d}$.}
\label{fig:probability-on-fibers}
\end{figure}
\begin{figure}[H]
\centering
\begin{tikzpicture}[
  x={(1.36cm,-.12cm)},y={(.76cm,.58cm)},z={(0cm,.82cm)},
  font=\small,>=Stealth
]
\def\FigureFiberSegments{.61/1.16/2.68/.04,1.72/1.18/2.88/.03,3.10/1.09/3.07/-.01,4.49/1.05/2.86/-.04,5.52/1.10/2.68/0}
\path[draw=black!60,fill=black!4,line width=.45pt]
  (0,1.08) .. controls (1.65,1.4) and (3.8,.88) .. (6,1.15)
  -- (6,2.70) .. controls (4.65,2.52) and (4.19,3.16) .. (3.10,3.07)
  .. controls (1.82,3.02) and (.88,2.52) .. (0,2.70) -- cycle;
\draw[semithick] (-.12,0) .. controls (1.4,.19) and (3.66,-.23) .. (6.13,.04);
\node[left=5pt] at (0,1.90) {\large$\InvariantMeasuredSpaces$};
\node[left=5pt] at (0,0) {\large$\GHSpace$};
\foreach \fiberx/\fiberlo/\fiberhi/\basey [count=\fiberindex from 0] in \FigureFiberSegments {
  \ifnum\fiberindex=2\relax
    \draw[line width=1.05pt] (\fiberx,\fiberlo) -- (\fiberx,\fiberhi);
    \draw[->,semithick] (\fiberx,\fiberlo-.05) -- node[right=3pt] {$\MeasureProjection$} (\fiberx,\basey+.10);
    \fill (\fiberx,\basey) circle[radius=1.9pt] node[below=4pt] {$\IsometryClass{X}{d}$};
  \else
    \draw[black!60] (\fiberx,\fiberlo) -- (\fiberx,\fiberhi);
    \draw[->,black!50,thin] (\fiberx,\fiberlo-.04) -- (\fiberx,\basey+.10);
    \fill[black!65] (\fiberx,\basey) circle[radius=1.55pt];
  \fi
}
  \foreach \fiberx/\fiberlo/\fiberhi/\basey [count=\fiberindex from 0] in \FigureFiberSegments {
    \ifcase\fiberindex
      \def\leftpower{2}\def\rightpower{3}\def\normalizer{60}\or
      \def\leftpower{3}\def\rightpower{3}\def\normalizer{140}\or
      \def\leftpower{3}\def\rightpower{2}\def\normalizer{60}\or
      \def\leftpower{4}\def\rightpower{2}\def\normalizer{105}\or
      \def\leftpower{2}\def\rightpower{2}\def\normalizer{30}%
    \fi
    \ifnum\fiberindex=2\relax
      \tikzset{bell face/.style={fill=black!19},bell edge/.style={black,line width=.95pt}}
    \else
      \tikzset{bell face/.style={fill=black!7},bell edge/.style={black!70,line width=.60pt}}
    \fi
    \def\bellcurve{plot[domain=0:1,samples=48,variable=\u]
      ({\fiberx},{\fiberlo+(\fiberhi-\fiberlo)*\u},
       {\normalizer*(\u^\leftpower)*((1-\u)^\rightpower)})}
    \path[bell face] \bellcurve -- cycle;
    \draw[bell edge] \bellcurve;
    \draw[bell edge] (\fiberx,\fiberlo,0) -- (\fiberx,\fiberhi,0);
    \pgfmathsetmacro{\modeu}{\leftpower/(\leftpower+\rightpower)}
    \pgfmathsetmacro{\modey}{\fiberlo+(\fiberhi-\fiberlo)*\modeu}
    \pgfmathsetmacro{\modez}{\normalizer*(\modeu^\leftpower)*((1-\modeu)^\rightpower)}
    \draw[black!53,densely dashed,line width=.4pt] (\fiberx,\modey,0) -- (\fiberx,\modey,\modez);
    \ifnum\fiberindex=2\relax
      \node[above=3pt] at (\fiberx,\modey,\modez) {$\LawSectionAt{X}{d}$};
    \fi
  }
  \draw[->,thin] (6.65,1.12,0) -- (6.65,1.12,2.45) node[above] {$z$};
\end{tikzpicture}
\caption{A schematic view of the probability law chosen on each fiber in Figure~\ref{fig:probability-on-fibers}. Each bell lies in the vertical plane above its own fiber. Points on the central fiber represent classes $\MeasuredIsometryClass{X}{d}{\mu}$. The law $\LawSectionAt{X}{d}$ depicted above it assigns full mass to $\MeasureProjection^{-1}(\IsometryClass{X}{d})$.}
\label{fig:fiber-probability-density}
\end{figure}

To describe the second arrow, denote the fiber over $\IsometryClass{X}{d}$ by $\MeasureFiber{\IsometryClass{X}{d}}=\MeasureProjection^{-1}(\IsometryClass{X}{d})$.
Each of its points $\xi$ has a representative $\MeasuredIsometryClass{Y}{e}{\nu}$.
Fix a representative $\MetricSpace{X}{d}$ of the underlying space, choose an isometry $u:\MetricSpace{Y}{e}\to\MetricSpace{X}{d}$, and set $\FiberRealization{X}{d}(\xi)=u_*\nu$.
Invariance of $\nu$ makes the result independent of $u$.
Replacing the measured representative by a measure-preservingly isometric one also leaves the measure unchanged.
The integral of a real continuous function $\TestFunction$ against the barycenter is
\begin{equation}\label{eq:part-i-average-interface}
 \int_\BaseCarrier\TestFunction(\BasePoint)
 \,\IntegrationDifferential\SelectedMeasure{\BaseCarrier}{\MetricSymbol}(\BasePoint)
 =\int_{\MeasureFiber{\IsometryClass{\BaseCarrier}{\MetricSymbol}}}
 \left(\int_\BaseCarrier\TestFunction\,\IntegrationDifferential\bigl(\FiberRealization{\BaseCarrier}{\MetricSymbol}(\FiberPoint)\bigr)\right)
 \,\IntegrationDifferential\LawSectionAt{\BaseCarrier}{\MetricSymbol}(\FiberPoint).
\end{equation}
The inner integral is on $X$.
The outer integral averages its values over the fiber with the law $\LawSectionAt{X}{d}$.

Since $\MeasureProjection\circ s=\IdentityMap_{\GHSpace}$, the forgetful map restricted to the image of $s$ is the inverse of $s$.
Forgetting the measure on the whole space of measured spaces yields the retraction in \Cref{cor:ghp-retract}.

For metric measure spaces, Kazukawa, Nakajima, and Shioya proved contractibility and local path connectedness in both the box and concentration topologies (\cite[arXiv v1, Theorems~1.4--1.6]{KazukawaNakajimaShioya2024}).
Our section and retraction connect the ordinary Gromov--Hausdorff space and the Gromov--Hausdorff--Prokhorov space.

\medskip
\noindent \textbf{Dependence on the other parts.}
The proof uses the common preliminaries and the results of Section~\ref{sec:prelim-measures}.
In the local models of Part~\ref{part:spectral}, full support makes the denominator of each normalized bump positive.
Weak convergence of the selected measures allows us to pass to the limit in the integrals defining the coordinates.
Part~\ref{part:topology} uses these models to prove the absolute retract property.

\medskip
\noindent \textbf{Organization.}
In \Cref{subsec:measures-ghp} we recall the Gromov--Hausdorff--Prokhorov topology.
In \Cref{subsec:measures-preliminaries-laws,subsec:measures-open-projection} we obtain the selection from the softness theorem of Banakh and Radul for Polish spaces.
In \Cref{subsec:measures-laws-and-averages} we construct the barycenter in a prescribed common compact metric space and prove \Cref{thm:invariant-measures}.
In \Cref{subsec:measures-retract} we define the retraction onto the image of the section and describe the intermediate law as a measure on a compact metric space of probability measures.

\medskip
\noindent \textbf{Conventions and notation.}
All compact metric spaces are nonempty.
When a metric space represents a point of $\GHSpace$, we identify it with its isometry class.
An isometry is a surjective isometric embedding.
We denote by $\ProbabilityMeasures(\DomainSpace)$ the set of Radon probability measures on a metrizable space $\DomainSpace$, equipped with the topology of weak convergence against bounded continuous real functions.
Subscripts in $\SelectedMeasure{\BaseCarrier}{\MetricSymbol}$ refer to the specified metric space $\MetricSpace{\BaseCarrier}{\MetricSymbol}$.
For isometrically embedded subsets of a common ambient space, we use the same symbols for the subsets and their restricted metrics.
Sequences are indexed by $\NonnegativeIntegers=\{0,1,2,\ldots\}$.

\section{Preliminaries}\label{sec:prelim-measures}
The Gromov--Hausdorff--Prokhorov distance compares two compact metric spaces and their measures in one ambient space.
We recall this distance first.
We use Haar probability to average over a single isometry group.
To choose measures on a family of fibers, we use the softness theorem of Banakh and Radul, which Valov extended to completely metrizable spaces.
\subsection{The Gromov--Hausdorff--Prokhorov space}\label{subsec:measures-ghp}
To compare probabilities on two compact metric spaces, we embed both spaces isometrically in a common metric space.
The Gromov--Hausdorff--Prokhorov distance compares their images in Hausdorff distance and their pushed-forward measures in the L\'evy--Prokhorov metric of Subsection~\ref{subsec:preliminaries-probability}.

\begin{definition}\label{def:measured-space}
A \emph{compact metric measure space} is a triple $\MeasuredSpace{\BaseCarrier}{\MetricSymbol}{\FirstMeasure}$ consisting of a nonempty compact metric space $\MetricSpace{\BaseCarrier}{\MetricSymbol}$ and a probability measure $\FirstMeasure\in\ProbabilityMeasures(\BaseCarrier)$.
For two such triples $\MeasuredSpace{\BaseCarrier}{\MetricSymbol}{\FirstMeasure}$ and $\MeasuredSpace{\ComparisonCarrier}{\ComparisonMetric}{\SecondMeasure}$, an \emph{isomorphism} from $\MeasuredSpace{\BaseCarrier }{\MetricSymbol }{\FirstMeasure }$ to $\MeasuredSpace{\ComparisonCarrier }{\ComparisonMetric }{\SecondMeasure }$ is an isometry $\IdentifyingIsometry \colon\MetricSpace{\BaseCarrier }{\MetricSymbol }\to\MetricSpace{\ComparisonCarrier }{\ComparisonMetric }$ satisfying $\IdentifyingIsometry \Pushforward\FirstMeasure =\SecondMeasure $.
\end{definition}

Let $\MeasuredSpaces$ be the set of isomorphism classes of these triples.
The underlying space of $\MeasuredIsometryClass{X}{d}{\mu}$ is the whole of $X$, including points outside $\supp\mu$.

\begin{definition}\label{def:ghp-distance}
Let $\MeasuredSpace{\BaseCarrier}{\MetricSymbol}{\FirstMeasure}$ and $\MeasuredSpace{\ComparisonCarrier}{\ComparisonMetric}{\SecondMeasure}$ be representatives of two elements of $\MeasuredSpaces$.
Let $\MetricSpace{\AmbientCarrier}{\AmbientMetric}$ range over compact metric spaces and let $\FirstEmbedding\colon\MetricSpace{\BaseCarrier}{\MetricSymbol}\to\MetricSpace{\AmbientCarrier}{\AmbientMetric}$ and $\SecondEmbedding\colon\MetricSpace{\ComparisonCarrier}{\ComparisonMetric}\to\MetricSpace{\AmbientCarrier}{\AmbientMetric}$ range over isometric embeddings.
Let $\mathcal S$ be the collection of quadruples $(\AmbientCarrier,\AmbientMetric,\FirstEmbedding,\SecondEmbedding)$ consisting of a compact metric space and two isometric embeddings with these domains.
Taking the infimum over $\mathcal S$, we define the \emph{Gromov--Hausdorff--Prokhorov distance} by
\begin{equation}\label{eq:ghp-distance}
 \begin{split}
 &\GHPDistance(\MeasuredSpace{\BaseCarrier }{\MetricSymbol }{\FirstMeasure },\MeasuredSpace{\ComparisonCarrier }{\ComparisonMetric }{\SecondMeasure })\\
 &\qquad=\inf_{(\AmbientCarrier,\AmbientMetric,\FirstEmbedding,\SecondEmbedding)\in\mathcal S}
 \max\{\HausdorffDistance{\AmbientMetric }(\FirstEmbedding (\BaseCarrier ),\SecondEmbedding (\ComparisonCarrier )),
          \ProkhorovDistance{\AmbientMetric }(\FirstEmbedding \Pushforward\FirstMeasure ,\SecondEmbedding \Pushforward\SecondMeasure )\}.
 \end{split}
\end{equation}
For this unpointed distance we use the maximum convention.
These conventions agree with Khezeli's definition for the additional structure of a probability measure (\cite[arXiv v5, Example~2.1(iii), Definition~2.4 and equation~(8)]{Khezeli2023Framework}).
\end{definition}

If a measure on $\MetricSpace{X}{d}$ is invariant under $\Isom(X,d)$, pushing it forward by any isometry to another representative yields the same measure there.
We use this fact to realize all measures represented by points of one fiber on the same space before taking their barycenter.

\begin{definition}\label{def:invariant-measured-space}
Let $\InvariantMeasuredSpaces\subset\MeasuredSpaces$ consist of the classes represented by triples $\MeasuredSpace{\BaseCarrier}{\MetricSymbol}{\FirstMeasure}$ whose measures satisfy the following two conditions, \emph{full support} and \emph{isometry invariance}.
\begin{equation}\label{eq:invariant-full-support-1}
 \supp(\FirstMeasure )=\BaseCarrier,
\end{equation}
and for every $\GroupElement\in\Isom\MetricSpace{\BaseCarrier}{\MetricSymbol}$,
\begin{equation}\label{eq:invariant-full-support-2}
 \GroupElement \Pushforward\FirstMeasure =\FirstMeasure.
\end{equation}
\end{definition}

For full-support measures, isomorphism in \Cref{def:measured-space} agrees with the usual mm-isomorphism of Kazukawa, Nakajima, and Shioya (\cite[arXiv v1, Definition~2.1]{KazukawaNakajimaShioya2024}), who compare the supports of the measures.
Our definition compares the entire compact spaces even when the supports are smaller.
Thus the Hausdorff term in \eqref{eq:ghp-distance} retains points outside the support.

\begin{theorem}[{Khezeli~\cite[arXiv v5, Theorems~2.6 and~2.12 and Example~2.1(iii)]{Khezeli2023Framework}}]\label{thm:ghp-polish}
The metric $\GHPDistance$ makes $\MeasuredSpaces$ a complete separable metric space.
\end{theorem}

Convergence in $\MeasuredSpaces$ can be realized simultaneously as Hausdorff convergence and weak convergence in a common space.
Conversely, evaluating the upper bound in \eqref{eq:ghp-distance} in such a realization yields convergence in $\MeasuredSpaces$.

\begin{lemma}\label{lem:common-measured-embedding}
Suppose that
\[
 \GHPDistance(\MeasuredSpace{\BaseCarrier _\SequenceIndex }{\MetricSymbol _\SequenceIndex }{\FirstMeasure _\SequenceIndex },
                   \MeasuredSpace{\BaseCarrier }{\MetricSymbol }{\FirstMeasure })\longrightarrow0.
\]
Then there exist a compact metric space $\MetricSpace{\AmbientCarrier }{\AmbientMetric }$ and isometric embeddings
\[
 \FirstEmbedding _\SequenceIndex \colon\MetricSpace{\BaseCarrier _\SequenceIndex }{\MetricSymbol _\SequenceIndex }\to\MetricSpace{\AmbientCarrier }{\AmbientMetric },
\]
\[
 \FirstEmbedding \colon\MetricSpace{\BaseCarrier }{\MetricSymbol }\to\MetricSpace{\AmbientCarrier }{\AmbientMetric }
\]
such that
\begin{equation}\label{eq:common-measured-embedding-1}
 \HausdorffDistance{\AmbientMetric }(\FirstEmbedding _\SequenceIndex (\BaseCarrier _\SequenceIndex ),\FirstEmbedding (\BaseCarrier ))\longrightarrow0,
\end{equation}
\begin{equation}\label{eq:common-measured-embedding-2}
 (\FirstEmbedding _\SequenceIndex )\Pushforward\FirstMeasure _\SequenceIndex \longrightarrow\FirstEmbedding \Pushforward\FirstMeasure,
\end{equation}
where the convergence in \eqref{eq:common-measured-embedding-2} is weak convergence in $\ProbabilityMeasures(\AmbientCarrier)$.
\end{lemma}

\begin{proof}
Khezeli's theorem (\cite[arXiv v5, Lemma~2.5 and Example~2.1(iii)]{Khezeli2023Framework}) provides a common compact ambient space where the underlying spaces converge in Hausdorff distance and the pushforward probability measures converge in L\'evy--Prokhorov distance.
On a compact metric space the latter convergence is equivalent to weak convergence by Prokhorov's theorem (\cite[author's manuscript, Lemma~1.15]{ShioyaMetricMeasure}).
These are precisely \eqref{eq:common-measured-embedding-1} and \eqref{eq:common-measured-embedding-2}.
\end{proof}
\subsection{Haar probability and continuous laws}\label{subsec:measures-preliminaries-laws}
Two kinds of averaging enter the construction.
Haar probability averages the translates of a measure on a fixed space.
A choice map for an open projection later specifies how to average the invariant measures over each fiber.

We first equip the isometry group of a compact metric space with the uniform topology.
Its compactness allows us to use normalized Haar measure.

\begin{lemma}\label{lem:compact-isometry-group}
Let $\MetricSpace{\BaseCarrier}{\MetricSymbol}$ be a nonempty compact metric space.
Then $\Isom\MetricSpace{\BaseCarrier}{\MetricSymbol}$, equipped with the uniform topology, is a compact metrizable topological group.
\end{lemma}

\begin{proof}
Let $\IsometryGroup=\Isom\MetricSpace{\BaseCarrier}{\MetricSymbol}$.
For $\Mapping,\SecondMapping\in\ContinuousFunctions(\BaseCarrier,\BaseCarrier)$, define the uniform metric by
\[
 \UniformMetric(\Mapping,\SecondMapping)
 =\sup_{\BasePoint\in\BaseCarrier}
    \MetricSymbol(\Mapping(\BasePoint),\SecondMapping(\BasePoint)).
\]
Since every element of $\IsometryGroup$ is an isometry and $\BaseCarrier$ is compact, the Arzel\`a--Ascoli theorem (\Cref{thm:arzela-ascoli}) implies that $\IsometryGroup$ has compact closure in $\ContinuousFunctions(\BaseCarrier,\BaseCarrier)$.

A uniform limit $\Mapping\in\ContinuousFunctions(\BaseCarrier,\BaseCarrier)$ of elements of $\IsometryGroup$ preserves distances.
It is surjective by a theorem of Freudenthal and Hurewicz (\cite[Satz~VI]{FreudenthalHurewicz1936}, see also \cite[Theorem~1.6.14]{BuragoBuragoIvanov2001}).
Hence $\IsometryGroup$ is closed and compact.
For sequences $\GroupElement_\SequenceIndex\to\GroupElement$ and $\SecondGroupElement_\SequenceIndex\to\SecondGroupElement$ in $\IsometryGroup$, the following estimates prove continuity of composition and inversion.
\[
 \UniformMetric(\GroupElement_\SequenceIndex\circ\SecondGroupElement_\SequenceIndex,
               \GroupElement\circ\SecondGroupElement)
 \leq\UniformMetric(\GroupElement_\SequenceIndex,\GroupElement)
       +\UniformMetric(\SecondGroupElement_\SequenceIndex,\SecondGroupElement),
\]
\[
 \UniformMetric(\GroupElement_\SequenceIndex^{-1},\GroupElement^{-1})
 =\UniformMetric(\GroupElement_\SequenceIndex,\GroupElement).
\]
Thus $\IsometryGroup$ is a compact metrizable topological group.
\end{proof}

\begin{theorem}[{Existence by Haar~\cite{Haar1933Measure}, uniqueness by von Neumann~\cite{vonNeumann1936Haar}, see also \cite[Theorem~5.14 and Chapter~5, Section~3]{DiestelSpalsbury2014}}]\label{thm:haar}
Every compact metrizable topological group $\ActingGroup$ admits a unique left-invariant Radon probability measure $\HaarProbability_{\ActingGroup}\in\ProbabilityMeasures(\ActingGroup)$, its \emph{Haar probability measure}.
It satisfies
\[
 \HaarProbability_{\ActingGroup}(\ActingGroup)=1,
\]
and for every $\GroupElement\in\ActingGroup$ and every $\BorelSubset\in\BorelSets(\ActingGroup)$, we have
\[
 \HaarProbability_{\ActingGroup}(\GroupElement\BorelSubset)=\HaarProbability_{\ActingGroup}(\BorelSubset).
\]
This probability is also right invariant and has full support.
\end{theorem}

Fix $a\in G$ and put $L_a(g)=ag$ and $R_a(g)=ga$.
By \Cref{thm:haar} we have $(L_a)\Pushforward\HaarProbability_G=(R_a)\Pushforward\HaarProbability_G=\HaarProbability_G$.
Thus by \eqref{eq:pushforward-integration}, for every Borel function $F$ integrable with respect to $\HaarProbability_G$,
\begin{equation}\label{eq:haar-translation-integrals}
 \int_G F(ag)\,\IntegrationDifferential\HaarProbability_G(g)
 =\int_G F(g)\,\IntegrationDifferential\HaarProbability_G(g)
 =\int_G F(ga)\,\IntegrationDifferential\HaarProbability_G(g).
\end{equation}
In particular, if $G$ acts continuously on a compact metric space $Y$ from the left, take $F(g)=\varphi(gy)$ for $y\in Y$ and $\varphi\in\ContinuousFunctions(Y)$ to obtain
\begin{equation}\label{eq:haar-action-integrals}
 \begin{aligned}
 \int_G\varphi\bigl(a(gy)\bigr)\,\IntegrationDifferential\HaarProbability_G(g)
 &=\int_G\varphi(gy)\,\IntegrationDifferential\HaarProbability_G(g)\\
 &=\int_G\varphi\bigl(g(ay)\bigr)\,\IntegrationDifferential\HaarProbability_G(g).
 \end{aligned}
\end{equation}
The first equality uses left invariance and the next uses right invariance.

Following Valov (\cite[arXiv v2, pp.~1--2]{Valov2009}), consider a continuous surjection $\OpenSurjection\colon\LawTotalSpace\to\LawBaseSpace$ between metrizable spaces.
Suppose there exists a continuous map $\LawSection\colon\LawBaseSpace\to\ProbabilityMeasures(\LawTotalSpace)$ such that for every $\LawBasePoint\in\LawBaseSpace$,
\[
 \OpenSurjection\Pushforward\LawSection(\LawBasePoint)
 =\DiracProbability_{\LawBasePoint}.
\]
Then $\OpenSurjection$ is called a \emph{Milyutin map} and $\LawSection$ a \emph{choice map} for $\OpenSurjection$.

A continuous map $q\colon E\to B$ between metrizable spaces is \emph{soft} if for every metrizable space $T$, closed subset $A\subseteq T$, and continuous maps $f\colon A\to E$ and $g\colon T\to B$ with $q\circ f=g|_A$, there exists a continuous map $h\colon T\to E$ satisfying $h|_A=f$ and $q\circ h=g$.

Banakh and Radul proved that a continuous open surjection between Polish spaces induces a soft map between their spaces of Radon probability measures (\cite[Theorem~1.2, p.~20]{BanakhRadul1999}).
Valov extended the result to completely metrizable spaces without assuming separability (\cite[arXiv v2, Theorem~1.1]{Valov2009}).
The measures here may have noncompact support.
Part~\ref{part:measures} uses the following consequence.

\begin{theorem}[{Banakh and Radul~\cite[Theorem~1.2, p.~20]{BanakhRadul1999} for Polish spaces, Valov~\cite[arXiv v2, Theorem~1.1]{Valov2009} in general}]\label{thm:valov}
Let $\OpenSurjection \colon \LawTotalSpace \to \LawBaseSpace $ be an open continuous surjection between completely metrizable spaces.
Then $\OpenSurjection$ is a Milyutin map.
Namely, there is a continuous map $\LawSection\colon \LawBaseSpace \to\ProbabilityMeasures(\LawTotalSpace )$ such that $\OpenSurjection \Pushforward\LawSection(\LawBasePoint )=\DiracProbability _\LawBasePoint $ for every $\LawBasePoint \in \LawBaseSpace $.
\end{theorem}

\begin{proof}
Let $\delta:\LawBaseSpace\to\ProbabilityMeasures(\LawBaseSpace)$ take $t$ to the Dirac measure $\DiracProbability_t$.
Choose a metric $d_T$ inducing the topology of $\LawBaseSpace$ and replace it by the bounded metric $\min\{1,d_T\}$.
By Banakh's metrization lemma (\cite[Lemma~4.3]{Banakh1995ProbabilityII}), the weak topology on $\ProbabilityMeasures(\LawBaseSpace)$ is metrizable (see also Valov, \cite[arXiv v2, p.~3]{Valov2009}).
By Valov's softness theorem, $\OpenSurjection_*:\ProbabilityMeasures(\LawTotalSpace)\to\ProbabilityMeasures(\LawBaseSpace)$ is soft.
Take the empty set as a closed subset of $\ProbabilityMeasures(\LawBaseSpace)$ and lift the identity in the diagram
\[
 \begin{tikzpicture}[baseline=(current bounding box.center),>=Stealth]
 \node (a) at (0,1.3) {$\emptyset$};
 \node (b) at (4,1.3) {$\ProbabilityMeasures(\LawTotalSpace)$};
 \node (c) at (0,0) {$\ProbabilityMeasures(\LawBaseSpace)$};
 \node (d) at (4,0) {$\ProbabilityMeasures(\LawBaseSpace)$};
 \draw[->] (a)--(b); \draw[->] (a)--(c);
 \draw[->] (b)--node[right]{$\OpenSurjection_*$}(d);
 \draw[->] (c)--node[below]{$\IdentityMap$}(d);
 \draw[->,dashed] (c)--node[above left]{$S$}(b);
 \end{tikzpicture}
\]
The lift $S$ is continuous and $\OpenSurjection_*\circ S=\IdentityMap$.
The map $\delta$ is continuous because $\int f\,\IntegrationDifferential\DiracProbability_t=f(t)$ for every $f\in C_b(\LawBaseSpace)$.
Thus $\LawSection=S\circ\delta$ is continuous.
For each $t\in\LawBaseSpace$,
\[
 \OpenSurjection_*\LawSection(t)=\DiracProbability_t,\qquad
 \LawSection(t)(\OpenSurjection^{-1}(\{t\}))=\DiracProbability_t(\{t\})=1.
\]
Hence $\LawSection$ is a choice map for $\OpenSurjection$.
Alternatively, taking the test space $\LawBaseSpace$ and closed subset $\emptyset$ in the definition of softness lifts the Dirac map $\delta$ along $\OpenSurjection_*$ and directly yields a continuous $\LawSection$ with $\OpenSurjection_*\circ\LawSection=\delta$.
\end{proof}

\begin{remark}\label{rem:laws-averaging-operator}
The probability laws in Valov's theorem (\Cref{thm:valov}) also determine a linear averaging operator.
Assume $\LawBaseSpace\ne\emptyset$ and define the \emph{averaging operator} by
\[
 {\AveragingOperator\colon}
 \BoundedContinuousFunctions(\LawTotalSpace,\RealNumbers)
 \to\BoundedContinuousFunctions(\LawBaseSpace,\RealNumbers).
\]
For every $\TestFunction\in\BoundedContinuousFunctions(\LawTotalSpace,\RealNumbers)$ and $\LawBasePoint\in\LawBaseSpace$, put
\[
 (\AveragingOperator\TestFunction)(\LawBasePoint)
 =\int_{\LawTotalSpace}\TestFunction(\FiberPoint)
   \,\IntegrationDifferential\LawSection(\LawBasePoint)(\FiberPoint).
\]
Since $\LawSection$ is a weakly continuous family of probabilities, $\AveragingOperator$ is a positive linear operator satisfying
\begin{equation}\label{eq:law-averaging-identities-1}
 \AveragingOperator1=1.
\end{equation}
For every $\TestFunction\in\BoundedContinuousFunctions(\LawTotalSpace,\RealNumbers)$,
\begin{equation}\label{eq:law-averaging-identities-2}
 \|\AveragingOperator\TestFunction\|_\infty\leq\|\TestFunction\|_\infty.
\end{equation}
The fiber support condition implies that, for every $f\in\BoundedContinuousFunctions(\LawBaseSpace,\RealNumbers)$,
\begin{equation}\label{eq:law-averaging-identities-3}
 \AveragingOperator(f\circ\OpenSurjection)=f.
\end{equation}
Define $\OpenSurjection^*\colon
 \BoundedContinuousFunctions(\LawBaseSpace,\RealNumbers)
 \to\BoundedContinuousFunctions(\LawTotalSpace,\RealNumbers)$ by $\OpenSurjection^*f=f\circ\OpenSurjection$.
By \eqref{eq:law-averaging-identities-1}--\eqref{eq:law-averaging-identities-3}, $\OpenSurjection^*\AveragingOperator$ is a projection of norm one onto $\OpenSurjection^*\BoundedContinuousFunctions(\LawBaseSpace,\RealNumbers)$.
\end{remark}
\subsection{Convergence of graphs and pushforwards of measures}\label{subsec:measures-preliminaries-graphs}
Denote the \emph{graph} of a map $f\colon K\to W$ by $\FunctionGraph(f)=\{(x,f(x))\mid x\in K\}$.

\begin{lemma}\label{lem:graph-pushforward-convergence}
Let $Z,W$ be compact metric spaces and $K_j,K$ nonempty compact subsets of $Z$.
Suppose the graphs of continuous maps $f_j\colon K_j\to W$ and $f\colon K\to W$ converge in Hausdorff distance, $\FunctionGraph(f_j)\to\FunctionGraph(f)$, in $Z\times W$ with the maximum metric.
If $\mu_j\in\ProbabilityMeasures(K_j)$ converge weakly to $\mu\in\ProbabilityMeasures(K)$ on $Z$, then $(f_j)\Pushforward\mu_j$ converge weakly to $f\Pushforward\mu$ on $W$.
\end{lemma}

\begin{proof}
Fix $\varphi\in\ContinuousFunctions(W,\RealNumbers)$.
By the Tietze--Urysohn extension theorem (\Cref{thm:tietze}), $\varphi\circ f$ extends to a continuous function $F$ on $Z$.
Then
\[
 \begin{aligned}
 &\left|\int_W\varphi\,\IntegrationDifferential((f_j)\Pushforward\mu_j)
             -\int_W\varphi\,\IntegrationDifferential(f\Pushforward\mu)\right|\\
 &\qquad\leq \max_{(x,y)\in\FunctionGraph(f_j)}|\varphi(y)-F(x)|
      +\left|\int_Z F\,\IntegrationDifferential\mu_j-\int_Z F\,\IntegrationDifferential\mu\right|.
 \end{aligned}
\]
The function $(x,y)\mapsto|\varphi(y)-F(x)|$ is uniformly continuous on $Z\times W$ and vanishes on $\FunctionGraph(f)$.
Thus the first term on the right tends to $0$ by Hausdorff convergence of the graphs, and the second does so by weak convergence of $\mu_j$.
\end{proof}

\begin{theorem}[{Khezeli~\cite[arXiv v5, Example~2.1(iii), Lemma~2.13, (10)]{Khezeli2023Framework}}]\label{thm:khezeli-probability-hausdorff}
Fix a compact metric space $\MetricSpace{Z}{h}$.
For a nonempty compact subset $K\subset Z$, write $\AmbientProbabilitySet{Z}{K}=\{(j_K)\Pushforward\mu\mid\mu\in\ProbabilityMeasures(K)\}$ for the set of pushforwards of probability measures by the inclusion $j_K\colon K\hookrightarrow Z$.
For any nonempty compact subsets $A,B\subset Z$,
\[
 \HausdorffDistance{\ProkhorovDistance{h}}\bigl(\AmbientProbabilitySet{Z}{A},\AmbientProbabilitySet{Z}{B}\bigr)
 \leq \HausdorffDistance{h}(A,B).
\]
Consequently the probability measure functor is Hausdorff continuous.
If nonempty compact subsets $K_j$ of $Z$ converge to $K$ in Hausdorff distance, then $\AmbientProbabilitySet{Z}{K_j}$ converges to $\AmbientProbabilitySet{Z}{K}$ in Hausdorff distance in $\ProbabilityMeasures(Z)$.
In particular, for each $\mu\in\ProbabilityMeasures(K)$ we can choose $\mu_j\in\ProbabilityMeasures(K_j)$ so that $(j_{K_j})\Pushforward\mu_j$ converges weakly to $(j_K)\Pushforward\mu$.
\end{theorem}

\begin{lemma}[{Khezeli~\cite[arXiv v5, Example~2.1(iii) and Lemma~2.8]{Khezeli2023Framework}}]\label{lem:khezeli-uniform-pushforward}
Let $\MetricSpace{X}{d}$ and $\MetricSpace{Y}{e}$ be nonempty compact metric spaces and fix $\mu\in\ProbabilityMeasures(X)$.
For all isometric embeddings $f,g\colon X\to Y$,
\[
 \ProkhorovDistance{e}(f\Pushforward\mu,g\Pushforward\mu)
 \leq \sup_{x\in X}e(f(x),g(x)).
\]
In particular, with the uniform topology on $\Isom(X,d)$ and the weak topology on $\ProbabilityMeasures(X)$, the map $g\mapsto g\Pushforward\mu$ is continuous.
\end{lemma}

\begin{lemma}[{Weak convergence in a prescribed ambient space}]\label{lem:ghp-invariant-limit}
Let $\MetricSpace{Z}{h}$ be a compact metric space and suppose $K_j,K\in\CompactHyperspace(Z)$ satisfy $\HausdorffDistance{h}(K_j,K)\to0$.
Denote the restricted metrics by $h_j,h_K$ and the inclusions by $\iota_j\colon K_j\hookrightarrow Z$ and $\iota\colon K\hookrightarrow Z$.
Let $\mu_j\in\ProbabilityMeasures(K_j)$ and $\mu\in\ProbabilityMeasures(K)$, and suppose $\MeasuredIsometryClass{K_j}{h_j}{\mu_j}$ converges to $\MeasuredIsometryClass{K}{h_K}{\mu}$ in $\GHPDistance$.
If $\mu$ is invariant under $\Isom(K,h_K)$, then $(\iota_j)\Pushforward\mu_j$ converges weakly to $\iota\Pushforward\mu$ on $Z$.
\end{lemma}
\begin{proof}
Set $\bar\mu_j=(\iota_j)\Pushforward\mu_j$.
By the theorems of Varadarajan and Prokhorov (\Cref{thm:probability-topology}), $\ProbabilityMeasures(Z)$ is compact and metrizable in the weak topology.
Let $\bar\mu_{j_k}$ be any weakly convergent subsequence with limit $\nu$.
By \Cref{lem:limit-support}, $\nu(K)=1$.
Setting $\nu_K=\nu|_K$, we have $\iota\Pushforward\nu_K=\nu$.
Using the two inclusions in \eqref{eq:ghp-distance} yields
\[
 \begin{aligned}
 &\GHPDistance\bigl(\MeasuredIsometryClass{K_{j_k}}{h_{j_k}}{\mu_{j_k}},
                    \MeasuredIsometryClass{K}{h_K}{\nu_K}\bigr)\\
 &\qquad\leq\max\{\HausdorffDistance{h}(K_{j_k},K),\ProkhorovDistance{h}(\bar\mu_{j_k},\nu)\}
 \longrightarrow0.
 \end{aligned}
\]
Comparing with the assumed limit yields $\MeasuredIsometryClass{K}{h_K}{\nu_K}=\MeasuredIsometryClass{K}{h_K}{\mu}$.
Thus some $u\in\Isom(K,h_K)$ satisfies $\nu_K=u\Pushforward\mu$.
Invariance implies $\nu_K=\mu$, so every subsequential limit is $\iota\Pushforward\mu$.
By compactness of $\ProbabilityMeasures(Z)$ the full sequence converges weakly to $\iota\Pushforward\mu$.
\end{proof}

Khezeli showed that convergence in $\GHPDistance$ admits a common ambient realization, and that for a fixed underlying space some subsequence of the measures converges to an isometric pushforward of the limit measure (\cite[arXiv v5, Lemmas~2.5 and~2.9]{Khezeli2023Framework}).
In \Cref{lem:ghp-invariant-limit} invariance of the limit measure allows the ambient space and inclusions to be prescribed.
For example, take the two-point space $K_j=K=Z=\{a,b\}$ and alternate the masses at $a,b$ between $(3/4,1/4)$ and $(1/4,3/4)$.
Interchanging the points keeps the measured isometry class constant, but the integrals of the indicator of $\{a\}$ do not converge.
More generally, if $u\in\Isom(K,h_K)$ satisfies $u\Pushforward\mu\ne\mu$, alternating $\mu$ and $u\Pushforward\mu$ has the same effect.

\section{Continuous invariant probability measures}\label{sec:measures}
The fiber of $\MeasureProjection:\InvariantMeasuredSpaces\to\GHSpace$ over $\IsometryClass{X}{d}$ consists of the measured isomorphism classes of invariant full-support probabilities on $\MetricSpace{X}{d}$.
We seek a section whose values converge weakly whenever the underlying compact sets converge in a prescribed realization.

Openness of $\MeasureProjection$ allows us to select a law on each fiber.
To recover a measure from this law, we transport the points of the fiber to probability measures on a common representative.
The effect of moving the representative is the subject of \Cref{prop:ambient-fiber-laws}.
\subsection{The open projection to unmeasured spaces}\label{subsec:measures-open-projection}
We verify the hypotheses of the theorem of Banakh and Radul (\Cref{thm:valov}) for $\MeasureProjection$.
We first express full support and isometry invariance as countable intersections of open sets in $\MeasuredSpaces$, proving that $\InvariantMeasuredSpaces$ is Polish.
We also show that $\InvariantMeasuredSpaces$ is dense in $\MeasuredSpaces$.
For openness of the projection we approximate a full-support invariant measure on a prescribed sequence of underlying spaces.

In the next proposition we extract a subsequence along which the graphs of isometries converge.
By \Cref{lem:graph-pushforward-convergence} the pushforward measures converge weakly along the same subsequence.

\begin{proposition}\label{lem:isometry-limits}
Let $\BaseCarrier _\SequenceIndex ,\BaseCarrier \subset \AmbientCarrier $ be nonempty compact subsets of a compact metric space $\MetricSpace{\AmbientCarrier }{\AmbientMetric }$ with $\HausdorffDistance{\AmbientMetric }(\BaseCarrier _\SequenceIndex ,\BaseCarrier )\to0$.
Give each subset of $\AmbientCarrier$ the restricted metric and the product $\AmbientCarrier\times\AmbientCarrier$ the maximum metric.
For each $\SequenceIndex\in\NonnegativeIntegers$ let $\Isometry_\SequenceIndex\colon\BaseCarrier_\SequenceIndex\to\BaseCarrier_\SequenceIndex$ be an isometry.
Then there exist a strictly increasing sequence $(j_k)_{k\in\NonnegativeIntegers}$ of nonnegative integers and an isometry $\Isometry\colon\BaseCarrier\to\BaseCarrier$ such that, as $k\to\infty$, the following convergence holds in $\CompactHyperspace(\AmbientCarrier\times\AmbientCarrier)$.
\[
 \FunctionGraph(\Isometry_{j_k})\longrightarrow\FunctionGraph(\Isometry).
\]
Furthermore, if $\GraphProbability_{j_k}\in\ProbabilityMeasures(\BaseCarrier_{j_k})$ converge weakly on $\AmbientCarrier$ to $\GraphProbability\in\ProbabilityMeasures(\BaseCarrier)$, then as $k\to\infty$ the following weak convergence holds in $\ProbabilityMeasures(\AmbientCarrier)$.
\[
 (\Isometry_{j_k})\Pushforward\GraphProbability_{j_k}
 \longrightarrow\Isometry\Pushforward\GraphProbability.
\]
\end{proposition}

\begin{proof}
We first identify the limit graph and then compare test functions along it.
For $x,y,z,w\in\AmbientCarrier$ write
\[
 \AmbientMetric_\times((x,y),(z,w))
 =\max\{\AmbientMetric(x,z),
          \AmbientMetric(y,w)\}.
\]
For each $\SequenceIndex$, put $\Correspondence_\SequenceIndex=\FunctionGraph(\Isometry_\SequenceIndex)$.
By compactness of $\CompactHyperspace(\AmbientCarrier\times\AmbientCarrier)$, choose strictly increasing indices $j_0<j_1<\cdots$ and a nonempty compact set $\Correspondence\subset\AmbientCarrier\times\AmbientCarrier$ such that
\begin{equation}\label{eq:isometry-graph-convergence}
 \HausdorffDistance{\AmbientMetric_\times}
 (\Correspondence_{j_k},\Correspondence)\longrightarrow0.
\end{equation}
\textbf{Step 1.
The limiting graph.}
Let $\CoordinateProjection_1,\CoordinateProjection_2\colon
 \AmbientCarrier\times\AmbientCarrier\to\AmbientCarrier$ be the coordinate projections.
They are $1$-Lipschitz for $\AmbientMetric_\times$, and $\CoordinateProjection_i(\Correspondence_{j_k})=\BaseCarrier_{j_k}$ for $i=1,2$, because each $\Isometry_{j_k}$ is surjective.
For each $i=1,2$, we apply \Cref{lem:hausdorff-lipschitz-images} to $\CoordinateProjection_i$ with $L=1$.
Together with the triangle inequality, this yields
\[
 \begin{aligned}
 \HausdorffDistance{\AmbientMetric}(\CoordinateProjection_i(\Correspondence),\BaseCarrier)
 &\leq\HausdorffDistance{\AmbientMetric}
       (\CoordinateProjection_i(\Correspondence),\CoordinateProjection_i(\Correspondence_{j_k}))
       +\HausdorffDistance{\AmbientMetric}(\BaseCarrier_{j_k},\BaseCarrier)\\
 &\leq\HausdorffDistance{\AmbientMetric_\times}(\Correspondence,\Correspondence_{j_k})
       +\HausdorffDistance{\AmbientMetric}(\BaseCarrier_{j_k},\BaseCarrier)
       \longrightarrow0.
 \end{aligned}
\]
Thus $\CoordinateProjection_1(\Correspondence)=\CoordinateProjection_2(\Correspondence)=\BaseCarrier$.
In particular, $\Correspondence\subset\BaseCarrier\times\BaseCarrier$.

Fix $(x,y),(z,w)\in\Correspondence$.
For each $k\in\NonnegativeIntegers$ choose points $(x_k,y_k)$ and $(z_k,w_k)$ of the nonempty compact set $\Correspondence_{j_k}$ nearest to $(x,y)$ and $(z,w)$, respectively.
Then
\[
 \max\{\AmbientMetric_\times((x_k,y_k),(x,y)),\,
         \AmbientMetric_\times((z_k,w_k),(z,w))\}
 \leq\HausdorffDistance{\AmbientMetric_\times}(\Correspondence_{j_k},\Correspondence)
 \longrightarrow0.
\]
For each $k$ we have $y_k=\Isometry_{j_k}(x_k)$ and $w_k=\Isometry_{j_k}(z_k)$.
Since $\Isometry_{j_k}$ is an isometry,
\begin{equation}\label{eq:limit-graph-distance}
 \AmbientMetric(x,z)=\lim_{k\to\infty}\AmbientMetric(x_k,z_k)
 =\lim_{k\to\infty}\AmbientMetric(y_k,w_k)=\AmbientMetric(y,w).
\end{equation}
In particular, if $(x,y),(x,w)\in\Correspondence$ then $\AmbientMetric(y,w)=0$, so $y=w$.
Surjectivity of the first projection makes $\Correspondence$ the graph of a unique map $\Isometry:\BaseCarrier\to\BaseCarrier$.
Surjectivity of the second projection shows that $\Isometry$ is onto.
For $(x,\Isometry(x)),(z,\Isometry(z))\in\Correspondence$, identity \eqref{eq:limit-graph-distance} implies $\AmbientMetric(\Isometry(x),\Isometry(z))=\AmbientMetric(x,z)$.
Thus $\Isometry$ is an isometry.

\textbf{Step 2.
Weak convergence of the pushforwards.}
The limit in \eqref{eq:isometry-graph-convergence} is $\FunctionGraph(\Isometry)$.
Taking both ambient spaces in \Cref{lem:graph-pushforward-convergence} to be $\AmbientCarrier$ and setting $f_k=\Isometry_{j_k}$ yields $(\Isometry_{j_k})\Pushforward\GraphProbability_{j_k}\to\Isometry\Pushforward\GraphProbability$.
\end{proof}

\begin{proposition}\label{lem:full-support-gdelta}
For each open-ball radius $r>0$ the following function is lower semicontinuous for $\GHPDistance$.
\begin{equation}\label{eq:minimum-open-ball-mass}
 \LowerBallMass_r\colon\MeasuredSpaces\longrightarrow[0,1],\qquad
 \LowerBallMass_r\MeasuredIsometryClass{X}{d}{\mu}=\min_{x\in X}\mu\bigl(\MetricBall(x,r;d)\bigr).
\end{equation}
The minimum is over the whole underlying space $X$.
The isometry classes of full-support measured spaces form the following $G_\delta$ set.
\[
 \FullSupportMeasuredClasses
 =\bigcap_{k\in\NonnegativeIntegers}\{\xi\in\MeasuredSpaces:\LowerBallMass_{2^{-k}}(\xi)>0\}.
\]
\end{proposition}

\begin{proof}
On a compact metric space $\MetricSpace{Z}{h}$ put $\BallEvaluation_r(z,\nu)=\nu(\MetricBall(z,r;h))$.
Suppose $z_j\to z$ and $\nu_j\Rightarrow\nu$, and fix $0<s<r$.
For large $j$ we have $\MetricBall(z,s;h)\subset\MetricBall(z_j,r;h)$.
By Alexandroff's theorem (\Cref{thm:portmanteau}),
\[
 \liminf_j \BallEvaluation_r(z_j,\nu_j)\geq
 \liminf_j\nu_j(\MetricBall(z,s;h))\geq\nu(\MetricBall(z,s;h)).
\]
Letting $s\uparrow r$ proves lower semicontinuity of $\BallEvaluation_r$.
In particular $\BallEvaluation_r(\cdot,\nu)$ attains its minimum on each nonempty compact set.
Athreya, L\"ohr, and Winter also use lower semicontinuity in the center for a fixed measure (\cite[proof of Lemma~3.4]{AthreyaLoehrWinter2016Gap}).

Suppose this minimum on $K\in\CompactHyperspace(Z)$ exceeds $a$.
The open set $\{\BallEvaluation_r>a\}$ contains $K\times\{\nu\}$.
By compactness of $K$ there are open neighborhoods $U$ of $K$ and $V$ of $\nu$ such that $U\times V\subset\{\BallEvaluation_r>a\}$.
All compact sets sufficiently close to $K$ in Hausdorff distance are contained in $U$.
Hence $(K,\nu)\mapsto\min_{z\in K}\BallEvaluation_r(z,\nu)$ is lower semicontinuous.
By \Cref{lem:common-measured-embedding}, any convergent sequence in $\MeasuredSpaces$ is realized in a single compact space with Hausdorff convergence of the underlying spaces and weak convergence of the measures.
This proves lower semicontinuity of $\LowerBallMass_r$.

If $\mu$ has full support, every open ball has positive mass, including one at a minimizing center.
Conversely, if $\LowerBallMass_{2^{-k}}\MeasuredIsometryClass{X}{d}{\mu}>0$ for every $k$, then every nonempty open set contains a ball of one of these radii and has positive mass.
Each condition is open, proving the stated representation.
\end{proof}

This function need not be continuous.
Give both points of a two-point space mass $1/2$ and let their distance $1-2^{-(j+1)}$ tend to $1$.
In $\RealNumbers$ keep the left point at $0$.
The underlying spaces converge in Hausdorff distance and the measures converge weakly, whereas $\LowerBallMass_1$ equals $1$ along the sequence and $1/2$ at the limit.

For a measured isometry class $\FiberPoint=\MeasuredIsometryClass{X}{d}{\mu}\in\MeasuredSpaces$, define the \emph{maximum displacement} $\InvarianceDefect(\FiberPoint)$ of $\mu$ under self-isometries by
\[
 \InvarianceDefect(\FiberPoint)=\InvarianceDefect\MeasuredIsometryClass{X}{d}{\mu}
 =\max_{g\in\Isom(X,d)}\ProkhorovDistance{d}(\mu,g_*\mu).
\]
\begin{proposition}\label{lem:invariant-measures-gdelta}
The function $\InvarianceDefect$ is upper semicontinuous on $\MeasuredSpaces$.
Its zero set is the set $\IsometryInvariantMeasuredClasses$ of classes carrying an isometry-invariant measure and is a $\CountableIntersectionType$ subset of $\MeasuredSpaces$.
\end{proposition}

\begin{proof}
Fix $\MetricSpace{X}{d}$ and $\mu$.
Khezeli's estimate (\Cref{lem:khezeli-uniform-pushforward}) implies that $g\mapsto g\Pushforward\mu$ is continuous.
The maximum therefore exists by compactness of $\Isom(X,d)$ (\Cref{lem:compact-isometry-group}).
Conjugation by an isometry $u\colon\MetricSpace{X}{d}\to\MetricSpace{Y}{e}$ is a bijection between the isometry groups, and $u_*$ preserves the L\'evy--Prokhorov distance (Khezeli, \Cref{lem:prokhorov-isometric-embedding}).
Thus $\InvarianceDefect(\FiberPoint)$ is independent of the representative.
On a fixed space $\MetricSpace{X}{d}$, each $g_*$ is isometric.
Applying the triangle inequality and maximizing over $g$ yields, for all $\mu,\nu\in\ProbabilityMeasures(X)$,
\[
 \bigl|\InvarianceDefect\MeasuredIsometryClass{X}{d}{\mu}-\InvarianceDefect\MeasuredIsometryClass{X}{d}{\nu}\bigr|
 \leq 2\ProkhorovDistance{d}(\mu,\nu).
\]
We next prove upper semicontinuity of $\InvarianceDefect$ on $\MeasuredSpaces$ as the underlying space varies.
Suppose $\FiberPoint_j=\MeasuredIsometryClass{X_j}{d_j}{\mu_j}\to\FiberPoint=\MeasuredIsometryClass{X}{d}{\mu}$ and choose a subsequence such that $\InvarianceDefect(\FiberPoint_{j_k})\to L=\limsup_{j\to\infty}\InvarianceDefect(\FiberPoint_j)$.
For each $k$ choose a maximizer $g_k\in\Isom(X_{j_k},d_{j_k})$.
By \Cref{lem:common-measured-embedding} we may realize $X_{j_k}\to X$ in Hausdorff distance and $\mu_{j_k}\to\mu$ weakly in a common compact metric space $\MetricSpace{Z}{h}$.
By \Cref{lem:isometry-limits} a further subsequence satisfies $(g_k)_*\mu_{j_k}\to g_*\mu$ weakly for some $g\in\Isom(X,d)$.
Since $\ProkhorovDistance{h}$ metrizes weak convergence and is invariant under isometric embeddings, along this subsequence we obtain
\[
 L=\lim_{k\to\infty}\ProkhorovDistance{h}(\mu_{j_k},(g_k)_*\mu_{j_k})
  =\ProkhorovDistance{d}(\mu,g_*\mu)\leq\InvarianceDefect(\FiberPoint).
\]
This proves upper semicontinuity.
Since $\mu$ is isometry invariant precisely when $\InvarianceDefect(\FiberPoint)=0$,
\[
 \IsometryInvariantMeasuredClasses
 =\bigcap_{n\in\NonnegativeIntegers}
   \{\,\FiberPoint\in\MeasuredSpaces\mid\InvarianceDefect(\FiberPoint)<2^{-n}\,\}.
\]
Each set on the right is open by upper semicontinuity.
\end{proof}

\begin{proposition}\label{prop:full-support-measured-approximation}
Let $\xi=\MeasuredIsometryClass{X}{d}{\mu}\in\MeasuredSpaces$ and let $\IsometryClass{Y_j}{e_j}\to\IsometryClass{X}{d}$ in $\GHSpace$.
On each $Y_j$ one can choose a full-support probability measure $\lambda_j$ such that the following convergence holds in $\MeasuredSpaces$.
\begin{equation}\label{eq:full-support-measured-approximation}
 \MeasuredIsometryClass{Y_j}{e_j}{\lambda_j}\longrightarrow\xi.
\end{equation}
\end{proposition}

\begin{proof}
By \Cref{lem:common-gh-embedding} we may regard $Y_j$ and $X$ as subsets of a common compact metric space $\MetricSpace{Z}{h}$ where $Y_j$ converges to $X$ in Hausdorff distance.
By Khezeli's Hausdorff continuity theorem for the probability measure functor (\Cref{thm:khezeli-probability-hausdorff}), we can choose $\sigma_j\in\ProbabilityMeasures(Y_j)$ converging weakly to $\mu$ on $Z$.
Assign positive weights summing to $1$ to a dense sequence in each $Y_j$ to obtain a full-support probability measure $\theta_j$.
Put $t_j=2^{-(j+1)}$ and $\lambda_j=(1-t_j)\sigma_j+t_j\theta_j$.
Then $\lambda_j$ also has full support.
For each $\varphi\in\ContinuousFunctions(Z)$,
\[
 \left|\int_Z\varphi\,\IntegrationDifferential\lambda_j-\int_Z\varphi\,\IntegrationDifferential\sigma_j\right|
 \leq 2t_j\norm{\varphi}_\infty\longrightarrow0,
\]
so $\lambda_j$ converges weakly to $\mu$ on $Z$.
Using these common embeddings in \eqref{eq:ghp-distance} proves \eqref{eq:full-support-measured-approximation}.
\end{proof}

We next average a probability measure on a compact metric space over the isometry group.
The maximum displacement under isometries bounds how far the average moves the measure.

\begin{lemma}\label{lem:invariant-full-support}
Let $\MetricSpace{Y}{e}$ be a nonempty compact metric space and $\lambda\in\ProbabilityMeasures(Y)$.
Let $\bar\lambda$ be the average of $\lambda$ over $G=\Isom(Y,e)$ with respect to normalized Haar probability.
Then $\bar\lambda$ is $G$-invariant and
\begin{equation}\label{eq:haar-average-displacement}
 \ProkhorovDistance{e}(\lambda,\bar\lambda)\leq\InvarianceDefect\MeasuredIsometryClass{Y}{e}{\lambda}.
\end{equation}
If $\lambda$ has full support, so does $\bar\lambda$.
In particular, $\MetricSpace{Y}{e}$ admits an invariant probability measure of full support.
\end{lemma}

\begin{proof}
By \Cref{lem:compact-isometry-group} the group $G$ is compact, so Haar's theorem (\Cref{thm:haar}) supplies a normalized left Haar probability $\HaarProbability_G$.
The pushforward of $\HaarProbability_G\otimes\lambda$ by the continuous action $(g,y)\mapsto gy$ is $\bar\lambda$.
By Tonelli's theorem, for every Borel set $B\subseteq Y$,
\[
 \bar\lambda(B)=\int_G (g_*\lambda)(B)\,\IntegrationDifferential\HaarProbability_G(g).
\]
For $a\in G$ and $\varphi\in\ContinuousFunctions(Y)$, the change of variables formula \eqref{eq:pushforward-integration} and Fubini's theorem yield
\[
 \begin{aligned}
 \int_Y\varphi\,\IntegrationDifferential(a\Pushforward\bar\lambda)
 &=\int_Y\int_G\varphi\bigl(a(gy)\bigr)\,
       \IntegrationDifferential\HaarProbability_G(g)\,\IntegrationDifferential\lambda(y)\\
 &=\int_Y\int_G\varphi(gy)\,
       \IntegrationDifferential\HaarProbability_G(g)\,\IntegrationDifferential\lambda(y)
 =\int_Y\varphi\,\IntegrationDifferential\bar\lambda.
 \end{aligned}
\]
The second equality uses left invariance in \eqref{eq:haar-action-integrals}.
Equality of the integrals of all continuous functions implies $a\Pushforward\bar\lambda=\bar\lambda$.

Let $r>\InvarianceDefect\MeasuredIsometryClass{Y}{e}{\lambda}$.
By definition of the L\'evy--Prokhorov distance, for every Borel set $B\subseteq Y$ and every $g\in G$,
\[
 \lambda(B)\leq (g_*\lambda)(B^r)+r,\qquad
 (g_*\lambda)(B)\leq\lambda(B^r)+r.
\]
Integrating both inequalities in $g$ replaces $g_*\lambda$ by $\bar\lambda$.
Thus $\ProkhorovDistance{e}(\lambda,\bar\lambda)\leq r$.
Taking the infimum over such $r$ proves \eqref{eq:haar-average-displacement}.

If $\lambda$ has full support and $U\subseteq Y$ is nonempty and open, then $\lambda(g^{-1}U)>0$ for every $g\in G$.
Hence $\bar\lambda(U)>0$.
Averaging the measure obtained by assigning positive weights of sum $1$ to a dense sequence in $Y$ yields an invariant probability measure of full support.
\end{proof}

\begin{theorem}\label{thm:dense-invariant-measured-spaces}
The subspace $\InvariantMeasuredSpaces$ of spaces carrying isometry-invariant full-support probability measures is a dense $\CountableIntersectionType$ subset of the Gromov--Hausdorff--Prokhorov space $\MeasuredSpaces$.
In particular, $\InvariantMeasuredSpaces$ is Polish in the relative topology.
\end{theorem}

\begin{proof}
Fix $\xi=\MeasuredIsometryClass{X}{d}{\mu}\in\MeasuredSpaces$.
The residual set in Rouyer's \Cref{thm:rouyer-generic-metric-spaces} is dense in $\GHSpace$.
Thus for every $j\in\NonnegativeIntegers$ we may choose $\MetricSpace{Y_j}{e_j}$ satisfying both conditions of that theorem and
\[
 \GHDistance(\IsometryClass{Y_j}{e_j},\IsometryClass{X}{d})<2^{-(j+1)}.
\]
We show that each $Y_j$ has trivial isometry group.
Take $g\in\Isom(Y_j,e_j)$ and $x\in Y_j$.
Since $Y_j$ is nonempty and has no isolated points, it is infinite.
Choose $y,z\in Y_j$ so that $x,y,z$ are distinct.
Since $e_j(g(x),g(y))=e_j(x,y)>0$, Rouyer's second condition implies $\{g(x),g(y)\}=\{x,y\}$.
Likewise $\{g(x),g(z)\}=\{x,z\}$.
Hence
\[
 g(x)\in\{x,y\}\cap\{x,z\}=\{x\}.
\]
Since $x$ was arbitrary, $g=\IdentityMap_{Y_j}$.

As $\IsometryClass{Y_j}{e_j}\to\IsometryClass{X}{d}$, by \Cref{prop:full-support-measured-approximation} we can choose a full-support probability measure $\lambda_j$ on each $Y_j$ such that
\[
 \GHPDistance\bigl(\MeasuredIsometryClass{Y_j}{e_j}{\lambda_j},\xi\bigr)\longrightarrow0.
\]
Every self-isometry of $Y_j$ is the identity and thus preserves $\lambda_j$.
Hence $\MeasuredIsometryClass{Y_j}{e_j}{\lambda_j}\in\InvariantMeasuredSpaces$.
This sequence converges to $\xi$, so $\InvariantMeasuredSpaces$ is dense in $\MeasuredSpaces$.

By \Cref{lem:full-support-gdelta} full support is equivalent to $\LowerBallMass_{2^{-k}}(\xi)>0$ for every $k\in\NonnegativeIntegers$.
By \Cref{lem:invariant-measures-gdelta} isometry invariance is equivalent to $\InvarianceDefect(\xi)=0$.
Since $\InvarianceDefect$ is nonnegative,
\[
 \InvariantMeasuredSpaces
 =\bigcap_{k\in\NonnegativeIntegers}
 \{\,\xi\in\MeasuredSpaces\mid
       \LowerBallMass_{2^{-k}}(\xi)>0,\ \InvarianceDefect(\xi)<2^{-k}\,\}.
\]
The functions $\LowerBallMass_{2^{-k}}$ are lower semicontinuous and $\InvarianceDefect$ is upper semicontinuous, so the right-hand side is a countable intersection of open sets.
Finally, by Khezeli's \Cref{thm:ghp-polish} the space $\MeasuredSpaces$ is Polish.
Every $\CountableIntersectionType$ subset of a Polish space is Polish in the relative topology.
Hence $\InvariantMeasuredSpaces$ is Polish.
\end{proof}

\begin{theorem}\label{prop:open-invariant-projection}
The map that forgets the measure,
\[
 \MeasureProjection\colon\InvariantMeasuredSpaces\longrightarrow\GHSpace,
 \qquad \MeasuredIsometryClass{X}{d}{\mu}\longmapsto\IsometryClass{X}{d}
\]
is a continuous open surjection.
\end{theorem}

\begin{proof}
The Hausdorff term in \eqref{eq:ghp-distance} implies $\GHDistance(\MeasureProjection(\xi),\MeasureProjection(\eta))\leq\GHPDistance(\xi,\eta)$ for all $\xi,\eta\in\InvariantMeasuredSpaces$.
Thus $\MeasureProjection$ is continuous and is surjective by \Cref{lem:invariant-full-support}.

For openness, fix $\xi=\MeasuredIsometryClass{X}{d}{\mu}\in\InvariantMeasuredSpaces$ and any sequence $m_j=\IsometryClass{X_j}{d_j}$ converging to $\MeasureProjection(\xi)$.
It suffices to construct $\xi_j\in\InvariantMeasuredSpaces$ with $\MeasureProjection(\xi_j)=m_j$ and $\xi_j\to\xi$.
By \Cref{prop:full-support-measured-approximation}, choose a full-support probability measure $\lambda_j$ on each $X_j$ such that $\zeta_j=\MeasuredIsometryClass{X_j}{d_j}{\lambda_j}\to\xi$.
Since $\InvarianceDefect$ is nonnegative and invariance of $\mu$ implies $\InvarianceDefect(\xi)=0$, upper semicontinuity in \Cref{lem:invariant-measures-gdelta} yields $\InvarianceDefect(\zeta_j)\to0$.
Let $\bar\lambda_j$ be the Haar average of $\lambda_j$.
By \Cref{lem:invariant-full-support} we have $\xi_j=\MeasuredIsometryClass{X_j}{d_j}{\bar\lambda_j}\in\InvariantMeasuredSpaces$.
Using the identity on $X_j$ and \eqref{eq:haar-average-displacement}, we obtain
\[
 \GHPDistance(\zeta_j,\xi_j)
 \leq\ProkhorovDistance{d_j}(\lambda_j,\bar\lambda_j)
 \leq\InvarianceDefect(\zeta_j)\longrightarrow0.
\]
Thus $\xi_j\to\xi$ and $\MeasureProjection(\xi_j)=m_j$ for every $j$.
For every open neighborhood $U$ of $\xi$, we have $m_j\in\MeasureProjection(U)$ for all sufficiently large $j$.
Since $\GHSpace$ is metrizable, $\MeasureProjection(U)$ is a neighborhood of $\MeasureProjection(\xi)$.
\end{proof}
\subsection{Continuous laws and their averages}\label{subsec:measures-laws-and-averages}
Each point of the fiber over $\IsometryClass{X}{d}$ represents a probability measure up to measure-preserving isometry.
To take the barycenter of a law on the fiber, first fix a representative $\MetricSpace{X}{d}$ and realize each point of the fiber as a measure on $X$.
We then prove continuity as the representatives vary in a common compact metric space.

For a compact metric space $\MetricSpace{X}{d}$, the space $\ProbabilityMeasures(X)$ is weakly compact.
Giry proved weak continuity of the operation of averaging probability measures on a Polish space against their law (\cite[Part~I, Theorem~1, pp.~69--70]{Giry1982Probability}).
On a compact metric space this is the \emph{barycenter map} characterized by \eqref{eq:continuous-barycenter} (see also \cite[proof of Theorem~8.10.5, p.~218]{Bogachev2007II}).
Its integral representation identifies the properties inherited by the barycenter.

\begin{proposition}\label{lem:continuous-barycenter}
Let $\MetricSpace{\BaseCarrier}{\MetricSymbol}$ be a nonempty compact metric space.
Given a law $\FirstProbabilityLaw\in\ProbabilityMeasures(\ProbabilityMeasures(\BaseCarrier))$, consider the identity, for all $\TestFunction\in\ContinuousFunctions(\BaseCarrier)$,
\begin{equation}\label{eq:continuous-barycenter}
 \int_\BaseCarrier\TestFunction\,\IntegrationDifferential\SecondMeasure
 =\int_{\ProbabilityMeasures(\BaseCarrier)}
       \left(\int_\BaseCarrier\TestFunction\,\IntegrationDifferential\FirstMeasure\right)
       \,\IntegrationDifferential\FirstProbabilityLaw(\FirstMeasure).
\end{equation}
There exists exactly one probability measure $\SecondMeasure\in\ProbabilityMeasures(\BaseCarrier)$ satisfying this identity.
We call it the \emph{barycenter} of $\FirstProbabilityLaw$ and write $\Barycenter{\BaseCarrier}(\FirstProbabilityLaw)$.
The barycenter map
\[
 \Barycenter{\BaseCarrier}\colon
 \ProbabilityMeasures(\ProbabilityMeasures(\BaseCarrier))\longrightarrow\ProbabilityMeasures(\BaseCarrier)
\]
is continuous for the weak topologies on both spaces and has the following properties.
\begin{enumerate}[label=\textup{(\roman*)},ref=\textup{(\roman*)}]
\item\label{item:barycenter-full-support}
Suppose that $\supp(\FirstMeasure)=\BaseCarrier$ for $\FirstProbabilityLaw$-almost every $\FirstMeasure\in\ProbabilityMeasures(\BaseCarrier)$. Then $\supp(\Barycenter{\BaseCarrier}(\FirstProbabilityLaw))=\BaseCarrier$.
\item\label{item:barycenter-invariance}
Fix $\GroupElement\in\Isom\MetricSpace{\BaseCarrier}{\MetricSymbol}$.
If
\begin{equation}\label{eq:barycenter-invariant-law}
 \FirstProbabilityLaw\bigl(\{\FirstMeasure\in\ProbabilityMeasures(\BaseCarrier)\mid
                      \GroupElement\Pushforward\FirstMeasure=\FirstMeasure\}\bigr)=1,
\end{equation}
then $\GroupElement\Pushforward\Barycenter{\BaseCarrier}(\FirstProbabilityLaw)
=\Barycenter{\BaseCarrier}(\FirstProbabilityLaw)$.
\end{enumerate}
\end{proposition}

\begin{proof}
For $\TestFunction\in\ContinuousFunctions(\BaseCarrier)$, define $\widehat\TestFunction\colon\ProbabilityMeasures(\BaseCarrier)\to\RealNumbers$ as follows.
For every $\FirstMeasure\in\ProbabilityMeasures(\BaseCarrier)$, set
\[
 \widehat\TestFunction(\FirstMeasure)
 =\int_\BaseCarrier\TestFunction\,\IntegrationDifferential\FirstMeasure.
\]
This function is continuous and bounded in absolute value by $\norm{\TestFunction}_\infty$.
In Bogachev's construction (\cite[proof of Theorem~8.10.5, p.~218]{Bogachev2007II}), take $M=\ProbabilityMeasures(\BaseCarrier)$.
The resulting measure satisfies \eqref{eq:continuous-barycenter}.
Taking $\TestFunction=1$ shows that its total mass is $1$.
Uniqueness follows from the Riesz--Markov--Kakutani representation theorem (\Cref{thm:riesz-markov-kakutani}).
If $\FirstProbabilityLaw_\SequenceIndex\to\FirstProbabilityLaw$ weakly in $\ProbabilityMeasures(\ProbabilityMeasures(\BaseCarrier))$, then for every $\TestFunction\in\ContinuousFunctions(\BaseCarrier)$ we have
\[
 \int_{\ProbabilityMeasures(\BaseCarrier)}\widehat\TestFunction
       \,\IntegrationDifferential\FirstProbabilityLaw_\SequenceIndex
 \longrightarrow
 \int_{\ProbabilityMeasures(\BaseCarrier)}\widehat\TestFunction
       \,\IntegrationDifferential\FirstProbabilityLaw.
\]
Thus the barycenters converge weakly.
Since both probability spaces are metrizable in their weak topologies, this proves continuity.

Assume that $\FirstProbabilityLaw$ is concentrated on a Borel set of full-support measures.
For a nonempty open subset $\OpenSubset\subset\BaseCarrier$, choose a continuous function $\TestFunction\colon\BaseCarrier\to[0,1]$ that is nonzero and vanishes outside $\OpenSubset$.
For every $\FirstMeasure\in\ProbabilityMeasures(\BaseCarrier)$ with $\supp(\FirstMeasure)=\BaseCarrier$, we have $\widehat\TestFunction(\FirstMeasure)>0$.
By the hypothesis and \eqref{eq:continuous-barycenter},
\[
 \Barycenter{\BaseCarrier}(\FirstProbabilityLaw)(\OpenSubset)
 \geq\int_{\ProbabilityMeasures(\BaseCarrier)}
       \widehat\TestFunction(\FirstMeasure)
       \,\IntegrationDifferential\FirstProbabilityLaw(\FirstMeasure)>0.
\]
This proves full support.

Fix an isometry $\GroupElement\in\Isom\MetricSpace{\BaseCarrier}{\MetricSymbol}$ for which \eqref{eq:barycenter-invariant-law} holds.
For every $\TestFunction\in\ContinuousFunctions(\BaseCarrier)$ and every $\FirstMeasure\in\ProbabilityMeasures(\BaseCarrier)$ with $\GroupElement\Pushforward\FirstMeasure=\FirstMeasure$, we have
\begin{equation}\label{eq:fiber-invariant-integrals}
 \int_\BaseCarrier\TestFunction\circ\GroupElement\,
       \IntegrationDifferential\FirstMeasure
 =\int_\BaseCarrier\TestFunction\,
       \IntegrationDifferential\FirstMeasure.
\end{equation}
By \eqref{eq:barycenter-invariant-law}, integrating \eqref{eq:fiber-invariant-integrals} against $\FirstProbabilityLaw$ and using \eqref{eq:continuous-barycenter} proves invariance of $\Barycenter{\BaseCarrier}(\FirstProbabilityLaw)$ under $\GroupElement$.
\end{proof}

\begin{theorem}\label{thm:invariant-measures}
To each nonempty compact metric space $\MetricSpace{\BaseCarrier}{\MetricSymbol}$ one can assign a probability measure $\SelectedMeasure{\BaseCarrier}{\MetricSymbol}\in\ProbabilityMeasures(\BaseCarrier)$ so that all three of the following conditions hold.
\begin{enumerate}[label=\textup{(M\arabic*)},ref=\textup{(M\arabic*)}]
\item
For every $\MetricSpace{\BaseCarrier}{\MetricSymbol}$, we have $\supp(\SelectedMeasure{\BaseCarrier}{\MetricSymbol})=\BaseCarrier$.
\item
For every surjective isometry $\IdentifyingIsometry\colon\MetricSpace{\BaseCarrier}{\MetricSymbol}\to\MetricSpace{\ComparisonCarrier}{\ComparisonMetric}$ between nonempty compact metric spaces,
\[
 \IdentifyingIsometry \Pushforward\SelectedMeasure{\BaseCarrier}{\MetricSymbol}=\SelectedMeasure{\ComparisonCarrier}{\ComparisonMetric}.
\]
\item\label{item:invariant-measures-continuity}
Let $\MetricSpace{\BaseCarrier _\SequenceIndex }{\MetricSymbol _\SequenceIndex }$ and $\MetricSpace{\BaseCarrier }{\MetricSymbol }$ be nonempty compact metric spaces.
Let $\MetricSpace{\AmbientCarrier}{\AmbientMetric}$ be any compact metric space and $\FirstEmbedding _\SequenceIndex \colon\MetricSpace{\BaseCarrier _\SequenceIndex }{\MetricSymbol _\SequenceIndex }\to\MetricSpace{\AmbientCarrier }{\AmbientMetric }$ and $\FirstEmbedding \colon\MetricSpace{\BaseCarrier }{\MetricSymbol }\to\MetricSpace{\AmbientCarrier }{\AmbientMetric }$ isometric embeddings.
If
\[
 \HausdorffDistance{\AmbientMetric}
   \bigl(\FirstEmbedding_\SequenceIndex(\BaseCarrier_\SequenceIndex),
         \FirstEmbedding(\BaseCarrier)\bigr)\longrightarrow0,
\]
then the pushforwards $(\FirstEmbedding_\SequenceIndex)_*\SelectedMeasure{\BaseCarrier_\SequenceIndex}{\MetricSymbol_\SequenceIndex}$ converge weakly in $\ProbabilityMeasures(\AmbientCarrier)$ to $\FirstEmbedding_*\SelectedMeasure{\BaseCarrier}{\MetricSymbol}$.
\end{enumerate}
\end{theorem}

\begin{lemma}[{Naturality of the barycenter, Giry~\cite[Part~I, Theorem~1, p.~70]{Giry1982Probability}}]\label{lem:barycenter-pushforward}
Let $f:\MetricSpace{X}{d}\to\MetricSpace{Y}{e}$ be a continuous map of nonempty compact metric spaces.
If $\sigma\in\ProbabilityMeasures(\ProbabilityMeasures(X))$, then
\[
 f_*\Barycenter{X}(\sigma)
    =\Barycenter{Y}\bigl((f_*)_*\sigma\bigr).
\]
\end{lemma}

\begin{proof}
For $\varphi\in C(Y)$ put $\widehat\varphi(\nu)=\int_Y\varphi\,\IntegrationDifferential\nu$.
For each $\mu\in\ProbabilityMeasures(X)$ the change of variables formula \eqref{eq:pushforward-integration} implies
\[
 \widehat\varphi(f_*\mu)=\int_X\varphi\circ f\,\IntegrationDifferential\mu
 =\widehat{\varphi\circ f}(\mu).
\]
Integrate this identity in $\mu$ against $\sigma$ and use \eqref{eq:continuous-barycenter} on both $X$ and $Y$ to conclude.
\end{proof}

Write $\MeasureFiber{\IsometryClass{X}{d}}=\MeasureProjection^{-1}(\IsometryClass{X}{d})$.
We equip
\[
 \FiberLawSpace
 =\{\lambda\in\ProbabilityMeasures(\InvariantMeasuredSpaces)\mid
           \exists m\in\GHSpace,\ \MeasureProjection_*\lambda=\DiracProbability_m\},
\]
with the relative weak topology.
Dirac measures at distinct points are distinct, so the point $m$ in the defining condition is unique.
Denote it by $\LawBaseMap(\lambda)$.
Then $\lambda$ assigns mass $1$ to $\MeasureProjection^{-1}(m)$.
This construction extracting the base point of a law concentrated on one fiber occurs in Giry (\cite[Part~III, Section~2, Proposition~1, p.~81]{Giry1982Probability}).
Fix a compact metric space $\MetricSpace{Z}{h}$.
For $K\in H=\CompactHyperspace(Z)$ let $j_K:K\hookrightarrow Z$ be the inclusion and $h_K=h|_{K\times K}$ the restricted metric.
Denote the isometry class of $K$ by $\AmbientClassMap{Z}{h}(K)=\IsometryClass{K}{h_K}$.
We use the following two subspaces.
\[
 \begin{aligned}
 \MeasureIncidenceSpace{Z}{h}&=\{(K,\xi)\in H\times\InvariantMeasuredSpaces\mid
                      \AmbientClassMap{Z}{h}(K)=\MeasureProjection(\xi)\},\\
 \LawIncidenceSpace{Z}{h}&=\{(K,\lambda)\in H\times\FiberLawSpace\mid
                      {\AmbientClassMap{Z}{h}(K)=\LawBaseMap(\lambda)\}}.
 \end{aligned}
\]
A point of $\MeasureIncidenceSpace{Z}{h}$ specifies a measured class and a copy of its underlying space in $Z$.
A point of $\LawIncidenceSpace{Z}{h}$ specifies such a copy for an entire law.

\begin{lemma}\label{lem:law-on-incidence}
For $(K,\lambda)\in \LawIncidenceSpace{Z}{h}$ the product measure $\DiracProbability_K\otimes\lambda$ assigns mass $1$ to $\MeasureIncidenceSpace{Z}{h}$, and restriction defines a continuous map
\[
 \IncidenceLaw{Z}{h}:\LawIncidenceSpace{Z}{h}\longrightarrow\ProbabilityMeasures(\MeasureIncidenceSpace{Z}{h}),\qquad
 (K,\lambda)\longmapsto(\DiracProbability_K\otimes\lambda)|_{\MeasureIncidenceSpace{Z}{h}}.
\]
\end{lemma}

\begin{proof}
The section of $\MeasureIncidenceSpace{Z}{h}$ at $K$ is $\MeasureProjection^{-1}(\AmbientClassMap{Z}{h}(K))$, which has $\lambda$-measure $1$.
The inequality $\GHDistance(\AmbientClassMap{Z}{h}(K),\AmbientClassMap{Z}{h}(L))\leq\HausdorffDistance{h}(K,L)$ shows that $\AmbientClassMap{Z}{h}$ is continuous.
Thus $\MeasureIncidenceSpace{Z}{h}$ is closed as the preimage of the diagonal under $\AmbientClassMap{Z}{h}\times\MeasureProjection$.
Suppose $(K_n,\lambda_n)\to(K,\lambda)$ in $\LawIncidenceSpace{Z}{h}$.
Apply Giry's product-measure result (\Cref{lem:dirac-product-convergence}) to the compact parameter space $\{K\}\cup\{K_n:n\in\NonnegativeIntegers\}$ and the Polish space $\InvariantMeasuredSpaces$.
After inclusion into $H\times\InvariantMeasuredSpaces$ the product measures converge weakly to $\DiracProbability_K\otimes\lambda$.
The Tietze--Urysohn extension theorem (\Cref{thm:tietze}) extends each bounded continuous real function on $\MeasureIncidenceSpace{Z}{h}$ to a bounded continuous function on the product.
Its integrals against the restricted measures therefore converge, proving continuity of $\IncidenceLaw{Z}{h}$.
\end{proof}

\begin{lemma}[{Realizing measures in a fixed ambient space}]\label{lem:continuous-measure-realization}
Let $(K,\xi)\in \MeasureIncidenceSpace{Z}{h}$.
Choose a representative $(Y,e,\nu)$ of $\xi=\MeasuredIsometryClass{Y}{e}{\nu}$ and a surjective isometry $u:\MetricSpace{Y}{e}\to\MetricSpace{K}{h_K}$.
For the inclusion $j_K:K\hookrightarrow Z$, the measure
\[
 {\AmbientMeasureRealization{Z}{h}(K,\xi)=(j_K\circ u)_*\nu}
\]
is independent of the representative and of $u$, and the map $\AmbientMeasureRealization{Z}{h}:\MeasureIncidenceSpace{Z}{h}\to\ProbabilityMeasures(Z)$ is continuous.
In particular, taking $(Z,h)=(X,d)$ and $K=X$, the realization map $\FiberRealization{X}{d}\colon\MeasureFiber{\IsometryClass{X}{d}}\to\ProbabilityMeasures(X)$ on a fixed representative $(X,d)$ is continuous.
\end{lemma}
\begin{proof}
Fix a representative $\MetricSpace{X}{d}$ of a class in $\GHSpace$.
For $\xi=\MeasuredIsometryClass{Y}{e}{\nu}\in\MeasureFiber{\IsometryClass{X}{d}}$ choose an isometry $u:\MetricSpace{Y}{e}\to\MetricSpace{X}{d}$ and set $\FiberRealization{X}{d}(\xi)=u_*\nu$.
If $v:\MetricSpace{Y}{e}\to\MetricSpace{X}{d}$ is another isometry, then $v^{-1}u$ preserves $\nu$, so $u_*\nu=v_*\nu$.
If a measured isometry $a:Y\to Y_0$ changes the representative, then $(u\circ a^{-1})_*a_*\nu=u_*\nu$.
Thus $\FiberRealization{X}{d}:\MeasureFiber{\IsometryClass{X}{d}}\to\ProbabilityMeasures(X)$ is well defined.
Each of its values has full support, and for $g\in\Isom(X,d)$ we have $g_*\FiberRealization{X}{d}(\xi)=(gu)_*\nu=\FiberRealization{X}{d}(\xi)$.

Suppose $(K_j,\xi_j)\to(K,\xi)$ in $\MeasureIncidenceSpace{Z}{h}$.
Put $h_j=h|_{K_j\times K_j}$ and take $\mu_j=\FiberRealization{K_j}{h_j}(\xi_j)$ and $\mu=\FiberRealization{K}{h_K}(\xi)$ as above.
Then $\xi_j=\MeasuredIsometryClass{K_j}{h_j}{\mu_j}$ and $\xi=\MeasuredIsometryClass{K}{h_K}{\mu}$, and $\mu$ is invariant under $\Isom(K,h_K)$.

The spaces $K_j$ converge to $K$ in Hausdorff distance and $\xi_j$ converges to $\xi$ in $\GHPDistance$.
Hence by \Cref{lem:ghp-invariant-limit}, the measures $\AmbientMeasureRealization{Z}{h}(K_j,\xi_j)=(j_{K_j})\Pushforward\mu_j$ converge weakly to $\AmbientMeasureRealization{Z}{h}(K,\xi)=(j_K)\Pushforward\mu$.
The domain is metrizable, so $\AmbientMeasureRealization{Z}{h}$ is continuous.
Taking $(Z,h)=(X,d)$ and restricting to $K=X$ proves continuity of $\FiberRealization{X}{d}$.
\end{proof}

For $\lambda\in\FiberLawSpace$ with $\LawBaseMap(\lambda)=\IsometryClass{X}{d}$, use the map in \Cref{lem:continuous-measure-realization} to define the \emph{law of measures on a fixed representative} and its barycenter by
\begin{equation}\label{eq:fixed-representative-law}
 \RealizedLaw{X}{d}^{\lambda}=(\FiberRealization{X}{d})_*(\lambda|_{\MeasureFiber{\IsometryClass{X}{d}}}),
 \qquad \RealizedBarycenter{X}{d}(\lambda)=\Barycenter{X}\bigl(\RealizedLaw{X}{d}^{\lambda}\bigr).
\end{equation}
The former is a probability measure on $\ProbabilityMeasures(X)$, the latter a probability measure on $X$.

\begin{proposition}\label{prop:ambient-fiber-laws}
For $\lambda\in\FiberLawSpace$ with $\LawBaseMap(\lambda)=\IsometryClass{X}{d}$ and a surjective isometry $v\colon(X,d)\to(Y,e)$, the law and its barycenter satisfy
\begin{equation}\label{eq:fiber-law-transport}
 \RealizedLaw{Y}{e}^{\lambda}=(v_*)_*\RealizedLaw{X}{d}^{\lambda}.
\end{equation}
\begin{equation}\label{eq:fiber-barycenter-transport}
 \RealizedBarycenter{Y}{e}(\lambda)=v_*\RealizedBarycenter{X}{d}(\lambda).
\end{equation}
For $(K,\lambda)\in \LawIncidenceSpace{Z}{h}$, restrict $\DiracProbability_K\otimes\lambda$ to $\MeasureIncidenceSpace{Z}{h}$ and push the resulting measure $\IncidenceLaw{Z}{h}(K,\lambda)$ forward by $\AmbientMeasureRealization{Z}{h}$.
This defines a continuous map
\[
 {\AmbientLawMap{Z}{h}=(\AmbientMeasureRealization{Z}{h})_*\circ\IncidenceLaw{Z}{h}\colon}
 \LawIncidenceSpace{Z}{h}\longrightarrow\ProbabilityMeasures(\ProbabilityMeasures(Z)).
\]
Its value $\AmbientLawMap{Z}{h}(K,\lambda)$ is the law of $\AmbientMeasureRealization{Z}{h}(K,\xi)$ when $\xi$ has distribution $\lambda|_{\MeasureFiber{\IsometryClass{K}{h_K}}}$.
\end{proposition}

\begin{proof}
For $\TestFunction\in\ContinuousFunctions(X)$ put
\begin{equation}\label{eq:average-functional}
 \RepresentingFunctional_{X,d}^{\lambda}(\TestFunction)
 =\int_{\MeasureFiber{\IsometryClass{X}{d}}}
   \left(\int_X\TestFunction\,\IntegrationDifferential\bigl(\FiberRealization{X}{d}(\xi)\bigr)\right)\IntegrationDifferential\lambda(\xi).
\end{equation}
By \eqref{eq:continuous-barycenter} and the definition of pushforward,
\begin{equation}\label{eq:measure-barycenter}
 \int_X\TestFunction\,db_{X,d}(\lambda)=\RepresentingFunctional_{X,d}^{\lambda}(\TestFunction).
\end{equation}
Figure~\ref{fig:measure-average} shows the construction for the particular law $\LawSectionAt{X}{d}$.
\begin{figure}[H]
\centering
\[
 \GHSpace \xrightarrow{\ \LawSection\ } \FiberLawSpace
 \xrightarrow{\ \FiberBarycenterMap\ } \InvariantMeasuredSpaces
 \xrightarrow{\ \MeasureProjection\ } \GHSpace
\]
\begin{tikzpicture}[>=Stealth,node distance=9mm,font=\small]
\node (law) {$\LawSectionAt{X}{d}|_{\MeasureFiber{\IsometryClass{X}{d}}}
  \in\ProbabilityMeasures(\MeasureFiber{\IsometryClass{X}{d}})$};
\node[below=of law] (transport) {$
  \RealizedLaw{X}{d}^{\LawSectionAt{X}{d}}\in\ProbabilityMeasures(\ProbabilityMeasures(X))$};
\node[below=of transport] (average) {$
  \SelectedMeasure{X}{d}\in\ProbabilityMeasures(X)$};
\draw[->] (law) -- node[right]{pushforward by $\FiberRealization{\BaseCarrier}{\MetricSymbol}$} (transport);
\draw[->] (transport) -- node[right]{barycenter} (average);
\end{tikzpicture}
\caption{The top row gives the section on isometry classes. The first vertical arrow pushes the law forward to $\ProbabilityMeasures(X)$, and the second takes its barycenter.}
\label{fig:measure-average}
\end{figure}

For an isometry $v:\MetricSpace{X}{d}\to\MetricSpace{Y}{e}$ the fibers $\MeasureFiber{\IsometryClass{X}{d}}$ and $\MeasureFiber{\IsometryClass{Y}{e}}$ coincide, and on this fiber $\FiberRealization{Y}{e}=v_*\circ\FiberRealization{X}{d}$.
Pushing $\lambda|_{\MeasureFiber{\IsometryClass{X}{d}}}$ forward by both sides yields \eqref{eq:fiber-law-transport}.
Naturality of the barycenter under pushforward (\Cref{lem:barycenter-pushforward}) then implies \eqref{eq:fiber-barycenter-transport}.

For general $\MetricSpace{Z}{h}$ the ambient map $\AmbientMeasureRealization{Z}{h}$ and the fiberwise map satisfy the following identity by definition.
\[
 \AmbientMeasureRealization{Z}{h}(K,\xi)=(j_K)_*\FiberRealization{K}{h_K}(\xi)
 \qquad ((K,\xi)\in \MeasureIncidenceSpace{Z}{h}).
\]
If $(K,\lambda)\in \LawIncidenceSpace{Z}{h}$, the map $\xi\mapsto(K,\xi)$ on $\MeasureFiber{\IsometryClass{K}{h_K}}$ carries $\lambda|_{\MeasureFiber{\IsometryClass{K}{h_K}}}$ to $\IncidenceLaw{Z}{h}(K,\lambda)$.
Since $\AmbientMeasureRealization{Z}{h}(K,\xi)=(j_K)_*\FiberRealization{K}{h_K}(\xi)$, functoriality of pushforward yields
\[
 (\AmbientMeasureRealization{Z}{h})_*\IncidenceLaw{Z}{h}(K,\lambda)=\bigl((j_K)_*\bigr)_*\RealizedLaw{K}{h_K}^{\lambda}.
\]
Thus its value is the law of the transported measures on $Z$.

For $\eta\in\ProbabilityMeasures(\MeasureIncidenceSpace{Z}{h})$ and $F\in C(\ProbabilityMeasures(Z))$ the pushforward by $\AmbientMeasureRealization{Z}{h}$ satisfies
\begin{equation}\label{eq:parameter-law-products}
 \int_{\ProbabilityMeasures(Z)}F\,\IntegrationDifferential\bigl((\AmbientMeasureRealization{Z}{h})_*\eta\bigr)
 =\int_{\MeasureIncidenceSpace{Z}{h}}(F\circ \AmbientMeasureRealization{Z}{h})\,\IntegrationDifferential\eta.
\end{equation}
Since $F\circ \AmbientMeasureRealization{Z}{h}$ is bounded and continuous on $\MeasureIncidenceSpace{Z}{h}$, \eqref{eq:parameter-law-products} implies weak continuity of $(\AmbientMeasureRealization{Z}{h})_*$.
By \Cref{lem:law-on-incidence}, $\IncidenceLaw{Z}{h}$ is also continuous, so the composition $\AmbientLawMap{Z}{h}=(\AmbientMeasureRealization{Z}{h})_*\circ\IncidenceLaw{Z}{h}$ is continuous.
\end{proof}

\begin{corollary}\label{cor:continuous-fiber-barycenter}
For $\lambda\in\FiberLawSpace$ choose a representative $\MetricSpace{X}{d}$ of $\LawBaseMap(\lambda)$.
Then the assignment
\[
 \FiberBarycenterMap:\FiberLawSpace\longrightarrow\InvariantMeasuredSpaces,
 \qquad \lambda\longmapsto\MeasuredIsometryClass{X}{d}{\RealizedBarycenter{X}{d}(\lambda)}
\]
defines a continuous map satisfying $\MeasureProjection\circ\FiberBarycenterMap=\LawBaseMap$.
\end{corollary}

\begin{proof}
By second countability of $X$ and Alexandroff's weak convergence theorem (\Cref{thm:portmanteau}), the full-support probability measures form a Borel subset of $\ProbabilityMeasures(X)$.
For every $g\in\Isom(X,d)$, the pushforward $g_*$ is weakly continuous, so its fixed point set is closed.
Each $\FiberRealization{X}{d}(\xi)$ has full support and belongs to this fixed point set.
Thus \Cref{lem:continuous-barycenter} implies that $\RealizedBarycenter{X}{d}(\lambda)$ is invariant and has full support.
By \eqref{eq:fiber-barycenter-transport} the measured isomorphism class is independent of the representative.
Forgetting the measure recovers $\LawBaseMap(\lambda)$.

For any compact metric space $\MetricSpace{Z}{h}$, \Cref{prop:ambient-fiber-laws,lem:barycenter-pushforward} identify the barycenter in $Z$ as the composite
\[
 \begin{gathered} \AmbientMeanMap{Z}{h}=\Barycenter{Z}\circ\AmbientLawMap{Z}{h}:\LawIncidenceSpace{Z}{h}\longrightarrow\ProbabilityMeasures(Z),\\
 \AmbientMeanMap{Z}{h}(K,\lambda)=(j_K)_*\RealizedBarycenter{K}{h_K}(\lambda). \end{gathered}
\]
Both maps are continuous, so $\AmbientMeanMap{Z}{h}$ is continuous.
The map extracting the base point of a law concentrated on a fiber is
\[
 \LawBaseMap=\delta_{\GHSpace}^{-1}\circ\MeasureProjection_*,
 \qquad \delta_{\GHSpace}(a)=\DiracProbability_a,
\]
with the inverse taken on the Dirac measures.
The map $\MeasureProjection$ is a continuous surjection between Polish spaces.
Apply Giry's proposition (\cite[Part~III, Section~2, Proposition~1, p.~81]{Giry1982Probability}) to the category on $\GHSpace$ with only identity arrows.
It follows that $\FiberLawSpace$ is closed in $\ProbabilityMeasures(\InvariantMeasuredSpaces)$ and that $\LawBaseMap$ is continuous.
Hence $\lambda_n\to\lambda$ in $\FiberLawSpace$ implies $\LawBaseMap(\lambda_n)\to\LawBaseMap(\lambda)$.
By \Cref{lem:common-gh-embedding} we can choose representatives $K_n,K$ in one compact metric space $\MetricSpace{Z}{h}$ satisfying $\HausdorffDistance{h}(K_n,K)\to0$.
By continuity of $\AmbientMeanMap{Z}{h}$ their barycenters converge weakly on $Z$.
Combining this weak convergence and Hausdorff convergence in \eqref{eq:ghp-distance} proves Gromov--Hausdorff--Prokhorov convergence.
\end{proof}

\begin{proof}[{Proof of \Cref{thm:invariant-measures}}]
\Cref{thm:dense-invariant-measured-spaces,prop:open-invariant-projection} show that $\MeasureProjection$ is an open continuous surjection between Polish spaces.
The theorem of Banakh and Radul (\Cref{thm:valov}) therefore gives a continuous choice map $\LawSection$ with $\LawBaseMap\circ\LawSection=\IdentityMap_{\GHSpace}$.
Define
\[
 s=\FiberBarycenterMap\circ\LawSection,\qquad
 \SelectedMeasure{X}{d}=\RealizedBarycenter{X}{d}\bigl(\LawSectionAt{X}{d}\bigr).
\]
By \Cref{cor:continuous-fiber-barycenter}, $s$ is a continuous section and the measures on representatives satisfy (M1) and (M2).
To prove \ref{item:invariant-measures-continuity}, fix the prescribed embeddings into $\MetricSpace{Z}{h}$.
The map taking each $K$ of $H=\CompactHyperspace(Z)$ to the pushforward of the measure selected on $K$ to $Z$ is the composition
\[
 \begin{aligned}
 H&\longrightarrow \LawIncidenceSpace{Z}{h}
   \xrightarrow{\ \IncidenceLaw{Z}{h}\ }\ProbabilityMeasures(\MeasureIncidenceSpace{Z}{h})
   \xrightarrow{\ (\AmbientMeasureRealization{Z}{h})_*\ }\ProbabilityMeasures(\ProbabilityMeasures(Z))
   \xrightarrow{\ \Barycenter{Z}\ }\ProbabilityMeasures(Z),\\
 K&\longmapsto(K,\LawSectionAt{K}{h_K}).
 \end{aligned}
\]
The first arrow is continuous by continuity of $\AmbientClassMap{Z}{h}$ and $\LawSection$.
Restriction of the product law is continuous by \Cref{lem:law-on-incidence}.
By \Cref{lem:continuous-measure-realization} the map $\AmbientMeasureRealization{Z}{h}$ is continuous, so its pushforward is weakly continuous.
The final barycenter map is continuous by \Cref{lem:continuous-barycenter}.
By Giry's identity relating barycenters and pushforwards (\Cref{lem:barycenter-pushforward}), the value of the composition at $K$ is $(j_K)_*\SelectedMeasure{K}{h_K}$.
By (M2) this is also the pushforward of the original measure onto the image of the prescribed embedding.
Hence Hausdorff convergence of the images implies the weak convergence in \ref{item:invariant-measures-continuity}.
\end{proof}

\begin{remark}[{Gluing by a partition of unity}]\label{rem:barycenter-partition-unity}
Taking the barycenter in \Cref{thm:invariant-measures} generalizes gluing probability measures by a partition of unity.
Fix a compact metric space $\MetricSpace{X}{d}$.
Let $\{U_i\}_{i\in I}$ be an open cover of a topological space $T$, let $\{\PartitionWeight_i\}_{i\in I}$ be a locally finite partition of unity with $\supp\PartitionWeight_i\subset U_i$, and let $\mu_i:U_i\to\ProbabilityMeasures(X)$ be continuous.
For $t\in T$ put $I(t)=\{i\in I\mid\PartitionWeight_i(t)>0\}$.
Use this finite set to form a law $\sigma_t$ on $\ProbabilityMeasures(X)$.
Then \eqref{eq:continuous-barycenter} yields
\[
 \sigma_t=\sum_{i\in I(t)}\PartitionWeight_i(t)\,\DiracProbability_{\mu_i(t)},
 \qquad
 \Barycenter{X}(\sigma_t)=\sum_{i\in I(t)}\PartitionWeight_i(t)\,\mu_i(t).
\]
The right-hand side is the usual gluing by a partition of unity.
For each $F\in C(\ProbabilityMeasures(X))$ the function $t\mapsto\PartitionWeight_i(t)F(\mu_i(t))$ extends continuously by zero outside $U_i$, and the sum of these functions is locally finite.
Thus $t\mapsto\sigma_t$ is weakly continuous, as is its composition with the barycenter map.
Here we push the selected law $\LawSectionAt{X}{d}$ on the fiber $\MeasureFiber{\IsometryClass{X}{d}}$ forward by $\FiberRealization{X}{d}$ and take its barycenter as $\SelectedMeasure{X}{d}$.
The same integral formula defines the barycenter even when the pushforward law has infinite support.
\end{remark}
\subsection{A retract of the measured space}\label{subsec:measures-retract}
The section $s=\FiberBarycenterMap\circ\LawSection$ takes values in $\InvariantMeasuredSpaces$.
We now regard its image as a subspace of the entire Gromov--Hausdorff--Prokhorov space $\MeasuredSpaces$.

\begin{corollary}\label{cor:ghp-retract}
The map defined by the selected measures,
\[
 s\colon\GHSpace\longrightarrow\MeasuredSpaces,\qquad
 \IsometryClass{X}{d}\longmapsto\MeasuredIsometryClass{X}{d}{\SelectedMeasure{X}{d}}
\]
is a topological embedding.
There is a strong deformation retraction $H$ onto its image $s(\GHSpace)\subseteq\InvariantMeasuredSpaces$ preserving each fiber of the forgetful map $p\colon\MeasuredSpaces\to\GHSpace$.
Restricting $H$ to the spaces carrying full-support measures, or to those also carrying isometry-invariant measures, yields a strong deformation retraction onto the same image.
\end{corollary}

\begin{proof}
Let $p:\MeasuredSpaces\to\GHSpace$ be the forgetful map $p(\MeasuredIsometryClass{X}{d}{\mu})=\IsometryClass{X}{d}$.
For $\xi,\eta\in\MeasuredSpaces$, \eqref{eq:ghp-distance} implies
\[
 \GHDistance(p(\xi),p(\eta))\leq\GHPDistance(\xi,\eta).
\]
Thus $p$ is continuous.
By \Cref{cor:continuous-fiber-barycenter,thm:invariant-measures}, the composite $s=\FiberBarycenterMap\circ\LawSection$, viewed as a map into $\MeasuredSpaces$, is a continuous section of $p$.

Since $p\circ s=\IdentityMap_{\GHSpace}$, $s$ is a topological embedding.
For $0\leq t\leq1$ define
\[
 H(\MeasuredIsometryClass{X}{d}{\mu},t)
 =\MeasuredIsometryClass{X}{d}{(1-t)\mu+t\,\SelectedMeasure{X}{d}}.
\]
The right-hand side is independent of the representative by (M2) and preservation of convex combinations under pushforward.\par
To prove continuity suppose $\zeta_j\to\zeta$ in $\MeasuredSpaces$ and $t_j\to t$.
By \Cref{lem:common-measured-embedding} we may embed the underlying spaces into a compact metric space $Z$ with Hausdorff convergence of the spaces and weak convergence of the given measures.
In the same realization the selected measures converge weakly by (M3).
Integrating each continuous function on $Z$ shows that the convex combinations with coefficients $1-t_j,t_j$ converge weakly to the corresponding limit.
Prokhorov's theorem (\Cref{thm:probability-topology}) and the definition \eqref{eq:ghp-distance} of the Gromov--Hausdorff--Prokhorov distance yield $H(\zeta_j,t_j)\to H(\zeta,t)$.
Also $H(-,0)=\IdentityMap_{\MeasuredSpaces}$ and $H(-,1)=s\circ p$, while $H(s(\xi),t)=s(\xi)$ and $p(H(\zeta,t))=p(\zeta)$ for all $t$.
Both full support and isometry invariance are preserved under convex combinations, proving the two restricted assertions.
\end{proof}

\begin{proposition}\label{prop:law-of-measures-ghp}
For a nonempty compact metric space $\MetricSpace{X}{d}$ abbreviate the law $\RealizedLaw{X}{d}^{\LawSectionAt{X}{d}}$ on $\ProbabilityMeasures(X)$ defined in \eqref{eq:fixed-representative-law} by $\RealizedLaw{X}{d}$.
Then the following map is continuous.
\[
 \MeasureLawAssignment:\GHSpace\longrightarrow\MeasuredSpaces,\qquad
 \IsometryClass{X}{d}\longmapsto
 \MeasuredIsometryClass{\ProbabilityMeasures(X)}{\ProkhorovDistance{d}}{\RealizedLaw{X}{d}}.
\]
\end{proposition}

\begin{proof}
An isometry $v:\MetricSpace{X}{d}\to\MetricSpace{Y}{e}$ induces an isometry $v_*$ of the probability measure spaces.
By \eqref{eq:fiber-law-transport} with $\lambda=\LawSectionAt{X}{d}$, we have $(v_*)_*\RealizedLaw{X}{d}=\RealizedLaw{Y}{e}$.
Thus $\MeasureLawAssignment$ is independent of the representative.
Use \Cref{lem:common-gh-embedding} to realize a convergent sequence as compact subsets $K_n\to K$ of one compact metric space $\MetricSpace{Z}{h}$.
By \Cref{lem:prokhorov-isometric-embedding} the pushforwards induced by the inclusions $j_{K_n}$ and $j_K$ are isometric embeddings.
Since $\HausdorffDistance{h}(K_n,K)\to0$, it follows from \Cref{thm:khezeli-probability-hausdorff} that
\[
 \HausdorffDistance{\ProkhorovDistance{h}}\bigl(\AmbientProbabilitySet{Z}{K_n},\AmbientProbabilitySet{Z}{K}\bigr)
 \longrightarrow0.
\]
By \Cref{prop:ambient-fiber-laws}, the laws $\AmbientLawMap{Z}{h}(K_n,\LawSectionAt{K_n}{h_{K_n}})$ on these images converge weakly to $\AmbientLawMap{Z}{h}(K,\LawSectionAt{K}{h_K})$.
This common realization proves convergence in $\MeasuredSpaces$.
\end{proof}

The image $s(\GHSpace)$ is closed as the fixed point set of the continuous map $s\circ p$.
The existence of a retraction onto this image thus also follows from the absolute retract property of $\GHSpace$ (\Cref{thm:absolute-extensor}).
The homotopy in \Cref{cor:ghp-retract} preserves each underlying isometry class at all times.

\clearpage
\part{The topology of Gromov--Hausdorff space II: {Local finite-dimensional approximations of compact metric spaces}}\label{part:spectral}
\begin{quote}
\small
\noindent\textbf{Abstract.}
We construct local finite-dimensional approximations of compact metric spaces that vary continuously in the Gromov--Hausdorff topology.
For every nonempty compact metric space and every positive tolerance, a neighborhood of that space admits a continuous assignment of finite-dimensional normed coordinate classes approximating the distances to the prescribed accuracy.
These classes are compatible with pullback by isometries, and the uniform error of their induced pseudometrics is continuous.
We normalize a bump at each point using the selected measure, project it orthogonally onto a finite-dimensional spectral subspace, and apply the distance operator.
The resulting function approximates the Kuratowski map $d(x,\cdot)$ of its center $x$ uniformly.
We choose one radius and one spectral parameter throughout the neighborhood and estimate the uniform error by a continuous function on isometry classes.
At a center with at least two points, the coefficients of the approximating functions in a real orthonormal basis serve as coordinates.
Continuity of the coordinate class follows from Anselone and Palmer's convergence theorem for spectral projections and uniform convergence of parameter-dependent integrals.
Near a one-point center we use the zero map into $\RealNumbers$ with the absolute value norm.

\par\smallskip
\noindent \textbf{Keywords.} Gromov--Hausdorff space, finite-dimensional approximation, distance operator, spectral subspace, spectral projection.

\par\smallskip
\noindent \textbf{2020 Mathematics Subject Classification.} Primary 54E35; Secondary 47A55, 47A58.
\end{quote}

\section{Introduction to Part II}\label{sec:intro-spectral}
The Kuratowski map represents each point of a compact metric space as a continuous function.
Using the continuously selected invariant probability measures from Part~\ref{part:measures}, we construct finite-dimensional approximations of this representation that are equivariant under self-isometries and vary continuously with the compact space.
Write $\GHSpace$ for the space of isometry classes of nonempty compact metric spaces.
The main result of this part is \Cref{thm:local-models}.

The symmetries of the input are part of the approximation problem.
By Rouyer's theorem (\Cref{thm:rouyer-generic-metric-spaces}), compact metric spaces with trivial isometry group are dense in $\GHSpace$.
The approximation must also be continuous at spaces with nontrivial isometries.
We retain whole spectral subspaces and allow changes of orthonormal basis.
The resulting orthogonal action accounts both for input symmetries and for changes of coordinates as the input varies.

Given $\IsometryClass{X}{d}\in\GHSpace$ and $\LocalTolerance>0$, we choose an open neighborhood $\ModelNeighborhood$ of $\IsometryClass{X}{d}$ and an integer $\ModelDimension\geq1$.
To each $\IsometryClass{Y}{e}\in\ModelNeighborhood$ we continuously assign an invariant $\ModelDimension$-dimensional normed coordinate class $\ModelValue{Y}{e}$.
If $\LocalPseudometric{Y}{e}$ is the pseudometric it induces on $\MetricSpace{Y}{e}$, then
\[
 \LocalError\IsometryClass{Y}{e}=\max_{y,z\in Y}|e(y,z)-\LocalPseudometric{Y}{e}(y,z)|<\LocalTolerance
\]
throughout the neighborhood, and $\LocalError$ is continuous.
The graph of the map sending a point to its class in the metric quotient is a correspondence of distortion $\LocalError\IsometryClass{Y}{e}$.
Thus the correspondence formula of Kalton and Ostrovskii (\Cref{thm:gh-correspondence}) yields
\[
 \GHDistance\bigl((Y,e),\MetricQuotient{Y}{\LocalPseudometric{Y}{e}}\bigr)
 \leq\tfrac12\LocalError\IsometryClass{Y}{e}<\frac{\LocalTolerance}{2}.
\]
For every representative $(\LocalCoordinates{Y}{e},\LocalNorm{Y}{e})$ of the coordinate class, the metric quotient is isometric to the image $\LocalCoordinates{Y}{e}(Y)$ with its norm metric.
If the spaces converge in Hausdorff distance in a common compact metric space, the orbits of the graph and norm pairs representing the coordinate classes also converge in Hausdorff distance.
This specifies the topology for the assignment $\ModelAssignment$.

The map $\LocalCoordinates{\ComparisonCarrier}{\ComparisonMetric}$ is equivariant under an orthogonal representation of the isometry group of the input space.
In \Cref{prop:local-isometry-representation} we choose a continuous homomorphism
\[
 \IsometryRepresentation{\ComparisonCarrier}{\ComparisonMetric}
 \colon\Isom\MetricSpace{\ComparisonCarrier}{\ComparisonMetric}
 \to\OrthogonalGroup(\ModelDimension)
\]
that realizes this equivariance and preserves the chosen norm.
For every $\FirstIsometry\in\Isom\MetricSpace{\ComparisonCarrier}{\ComparisonMetric}$ and every $\ComparisonPoint\in\ComparisonCarrier$, the equivariance identity is
\[
 \LocalCoordinates{\ComparisonCarrier}{\ComparisonMetric}(\FirstIsometry\ComparisonPoint)
 =\IsometryRepresentation{\ComparisonCarrier}{\ComparisonMetric}(\FirstIsometry)
   \LocalCoordinates{\ComparisonCarrier}{\ComparisonMetric}(\ComparisonPoint).
\]
For every $\FirstIsometry\in\Isom\MetricSpace{\ComparisonCarrier}{\ComparisonMetric}$ and every $\CoefficientVector\in\RealNumbers^\ModelDimension$, the norm-preservation identity is
\[
 \LocalNorm{\ComparisonCarrier}{\ComparisonMetric}
   \bigl(\IsometryRepresentation{\ComparisonCarrier}{\ComparisonMetric}(\FirstIsometry)
          \CoefficientVector\bigr)
 =\LocalNorm{\ComparisonCarrier}{\ComparisonMetric}(\CoefficientVector).
\]
Thus the varying isometry groups act through one fixed orthogonal group on the model neighborhood.
Orthogonal changes of basis remain relevant even when the input space has trivial isometry group.

Two applications of spectral methods in geometry provide context.
Bates proved that closed connected Riemannian manifolds of fixed dimension at least two and volume one, with a common lower Ricci curvature bound and a common positive lower bound for the injectivity radius, admit smooth embeddings by a uniformly bounded number of Laplacian eigenfunctions (\cite[arXiv v1, Theorem~2]{Bates2014Eigenfunctions}).
Kuwae and Shioya studied spectral structures associated with nonnegative self-adjoint operators on varying Hilbert spaces.
Under compact convergence of the spectral structures, the dimensions of the spectral subspaces for a bounded interval eventually agree if its endpoints lie outside the limit spectrum (\cite[Theorem~2.6]{KuwaeShioya2003Spectral}).
Here we use the distance kernel on arbitrary compact metric spaces to approximate distances uniformly.
We use the continuous assignment of invariant full-support probability measures constructed in \Cref{thm:invariant-measures} of Part~\ref{part:measures}.
The distance operator $\DistanceOperator{\BaseCarrier}{\MetricSymbol}
 \colon\LebesgueSpace^2(\BaseCarrier,\SelectedMeasure{\BaseCarrier}{\MetricSymbol})
 \to\LebesgueSpace^2(\BaseCarrier,\SelectedMeasure{\BaseCarrier}{\MetricSymbol})$ associated with $\MetricSpace{\BaseCarrier}{\MetricSymbol}$ is a self-adjoint Hilbert--Schmidt integral operator.
For $f\in\LebesgueSpace^2(\BaseCarrier,\SelectedMeasure{\BaseCarrier}{\MetricSymbol})$ and $\BasePoint\in\BaseCarrier$ it is defined by
\[
 (\DistanceOperator{\BaseCarrier}{\MetricSymbol}f)(\BasePoint)
 =\int_{\BaseCarrier}\MetricSymbol(\BasePoint,\IntegrationPoint)
 f(\IntegrationPoint)\,\IntegrationDifferential\SelectedMeasure{\BaseCarrier}{\MetricSymbol}(\IntegrationPoint).
\]
Maria, Oudot, and Solomon studied this operator and its spectral coordinates (\cite[Section~3]{MariaOudotSolomon2020}).
We apply the orthogonal projection onto a finite-dimensional spectral subspace and then the distance operator to a normalized bump at each point.
The uniform distance between the resulting functions approximates the distance between the two points.
Retaining whole eigenspaces accommodates eigenvalue multiplicities through orthogonal changes of basis.
We realize Gromov--Hausdorff convergence as Hausdorff convergence in a common compact space.
Extending the distance kernel to this space yields strongly convergent operators whose differences from the limit are collectively compact.
By the spectral projection theorem of Anselone and Palmer (\cite[Proposition~6.3]{AnselonePalmer1968}), the Riesz projections associated with a boundary disjoint from the limit spectrum converge strongly and eventually have the same rank as the limit projection.
Comparing the compositions of the ambient operators and projections on a compact family of bumps yields convergence of the approximating functions uniformly in both variables along correspondences (\Cref{prop:intrinsic-spectral-kernel}).
This controls the approximation error throughout a neighborhood.
Once the bases are aligned, convergence of coordinates follows from convergence of parameter-dependent integrals and division by a positive denominator.


The author's interpolation theorem (\cite[Theorem~1.1]{Ishiki2026Interpolation}) provides a uniform error estimate for extending metrics prescribed on a discrete family of closed subsets of a fixed metrizable space.
Here \Cref{thm:local-models} constructs maps and norms on varying compact spaces for use in the absolute retract argument.

\medskip
\noindent \textbf{Dependence on the other parts.}
The local models use the continuous assignment of full-support invariant measures in \Cref{thm:invariant-measures} of Part~\ref{part:measures}.
The quantitative spectral estimates also use \Cref{lem:full-support-gdelta} on the minimum mass of open balls.
In Part~\ref{part:topology} we use \Cref{thm:local-models} to represent the varying norms in a fixed Banach space and combine the local point maps.
The orthogonal actions determine the compact group in that construction.

\medskip
\noindent \textbf{Organization.}
Section~\ref{sec:prelim-spectral} recalls spectral decomposition and perturbation theory and the square root of a positive definite matrix.
In Section~\ref{subsec:spectral-distance-operator} we define the distance operator and study its finite-dimensional spectral subspaces.
Section~\ref{subsec:spectral-distance-approximation} constructs uniform approximations to the Kuratowski maps from normalized bumps and bounds the spectral parameter and the dimension in terms of the masses of open balls.
In Section~\ref{subsec:spectral-subspace-continuity} we establish continuity of the spectral subspaces and construct orthonormal bases with uniformly convergent continuous extensions to a common compact ambient space, as stated in \Cref{lem:spectral-continuity}\ref{item:spectral-basis-uniform-convergence}.
In Section~\ref{subsec:spectral-local-models} we combine these results to construct the invariant coordinate classes of \Cref{thm:local-models} and prove their continuity.

\medskip
\noindent \textbf{Conventions and notation.}
All compact metric spaces are nonempty, and finite-dimensional coordinate spaces are real.
Weak convergence of probabilities and the uniform integral estimates are recalled in Subsection~\ref{subsec:preliminaries-probability}.
Function spaces are real unless complex scalars are explicitly indicated for spectral perturbation.
An isometry is a surjective isometric embedding.
The pair $\MetricSpace{\BaseCarrier}{\MetricSymbol}$ specifies both the underlying set and its metric.
In particular, $\SelectedMeasure{\BaseCarrier}{\MetricSymbol}$ and $\DistanceOperator{\BaseCarrier}{\MetricSymbol}$ depend on that pair.
When compact spaces are embedded isometrically into a common ambient space, we use the same symbols for the embedded subsets and the restricted metrics.
The notation $\IdentityMap_\BaseCarrier$ denotes the identity map on $\BaseCarrier$.
We omit the subscript when its domain is specified.
Sequences are indexed by $\NonnegativeIntegers=\{0,1,2,\ldots\}$.

\section{Preliminaries}\label{sec:prelim-spectral}
We recall the spectral theorem of Hilbert and Schmidt for compact self-adjoint operators, the Riesz--Schauder theorem on the spectrum of a compact operator, and the perturbation theorem of Anselone and Palmer used to compare finite-dimensional spectral subspaces.
The last subsection records the matrix result used to normalize their bases.
\subsection{Compact operators and spectral decomposition}\label{subsec:spectral-preliminaries-operators}
A bounded linear operator $\CompactOperator \colon \HilbertDomain \to \HilbertDomain $ on a Hilbert space $\HilbertDomain $ is \emph{compact} if the image of the closed unit ball of $\HilbertDomain$ has compact closure.
If for all $\HilbertFunction,\HilbertTestFunction\in\HilbertDomain$
\[
 \langle \CompactOperator \HilbertFunction ,\HilbertTestFunction \rangle_\HilbertDomain =\langle \HilbertFunction ,\CompactOperator \HilbertTestFunction \rangle_\HilbertDomain,
\]
then $\CompactOperator$ is \emph{self-adjoint}.
A nonzero vector $\HilbertFunction $ satisfying $\CompactOperator \HilbertFunction =\Eigenvalue  \HilbertFunction $ is an \emph{eigenvector} for the \emph{eigenvalue} $\Eigenvalue $.
For operators on function spaces we also call it an \emph{eigenfunction}.
Hilbert (\cite{Hilbert1904Integral}) and Schmidt (\cite{Schmidt1907Integral}) studied eigenfunction expansions for integral operators with symmetric kernels.
We use their spectral theorem in the form for compact self-adjoint operators.

For a real Hilbert space $\HilbertDomain$, write $\HilbertDomain_{\ComplexNumbers}=\HilbertDomain\oplus\ImaginaryUnit\HilbertDomain$ for its complexification.
We denote the complex-linear extension of $\CompactOperator$ by the same symbol.
For $\HilbertFunction,\HilbertTestFunction\in\HilbertDomain$, set
\[
 \CompactOperator(\HilbertFunction+\ImaginaryUnit\HilbertTestFunction)
 =\CompactOperator\HilbertFunction+\ImaginaryUnit\CompactOperator\HilbertTestFunction.
\]
For a complex Banach space, write $\IdentityMap$ for its identity operator.
For a bounded operator $\CompactOperator$ on that space, its \emph{spectrum} $\spec(\CompactOperator)$ is the set of complex numbers $\Eigenvalue$ for which $\Eigenvalue\IdentityMap-\CompactOperator$ has no bounded inverse.
For an operator on a real Hilbert space, we use the spectrum of its complex-linear extension to $\HilbertDomain_{\ComplexNumbers}$.
For a self-adjoint operator this set is real.
In finite dimension it is exactly the set of eigenvalues.
In infinite dimension, the number $0$ can belong to the spectrum without being an eigenvalue.

\begin{theorem}[{Spectral theorem of Hilbert and Schmidt, \cite[Theorem~15.22, pp.~400--401]{Muscat2024}}]\label{prop:spectral-basics}
Let $\CompactOperator $ be a compact self-adjoint operator on a separable real Hilbert space $\HilbertDomain $.
Its spectrum is real.
By the Riesz--Schauder theorem (\Cref{thm:compact-spectrum}) for the complexified operator, every nonzero spectral value $\lambda$ is an eigenvalue with finite-dimensional eigenspace $\ker(\CompactOperator-\lambda \IdentityMap)$.
For each $\SpectralCutoff >0$ only finitely many eigenvalues, counted with multiplicity, have absolute value greater than $\SpectralCutoff $.
There exist a finite or countable orthonormal family $\{\Eigenfunction _\CoordinateIndex \}_{\CoordinateIndex \in \EigenIndexSet }$ indexed by $\EigenIndexSet\subset\NonnegativeIntegers$ and corresponding nonzero eigenvalues $\{\Eigenvalue _\CoordinateIndex \}_{\CoordinateIndex \in \EigenIndexSet }$ such that
\[
 \HilbertDomain =\ker \CompactOperator \ \mathbin{\oplus}\
       \overline{\LinearSpan }\{\Eigenfunction _\CoordinateIndex \mid \CoordinateIndex \in \EigenIndexSet \}.
\]
Every $\HilbertFunction \in \HilbertDomain $ admits the decomposition
\[
 \HilbertFunction =\HilbertFunction _0+\sum_{\CoordinateIndex \in \EigenIndexSet }\langle \HilbertFunction ,\Eigenfunction _\CoordinateIndex \rangle \Eigenfunction _\CoordinateIndex.
\]
The operator acts on this decomposition by
\[
 \CompactOperator \HilbertFunction =\sum_{\CoordinateIndex \in \EigenIndexSet }\Eigenvalue _\CoordinateIndex \langle \HilbertFunction ,\Eigenfunction _\CoordinateIndex \rangle \Eigenfunction _\CoordinateIndex.
\]
The remaining component satisfies
\[
 \HilbertFunction _0\in\ker \CompactOperator.
\]
Both sums converge in $\HilbertDomain $.
If $\EigenIndexSet $ is infinite, its eigenvalues tend to $0$ in any enumeration.
\end{theorem}

The orthogonal expansion implies that every nonzero compact self-adjoint operator has a nonzero eigenvalue.
\subsection{Spectral perturbation and generalized eigenspaces}\label{subsec:spectral-preliminaries-perturbation}
Let $\BanachAmbientSpace$ be a complex Banach space.
Write
\[
 {\ClosedUnitBall=\{\FirstVector\in\BanachAmbientSpace\mid\norm{\FirstVector}\leq1\}}
\]
for its closed unit ball.
A family of bounded linear operators $\{\AmbientOperator_\SequenceIndex\colon\BanachAmbientSpace\to\BanachAmbientSpace\}_{\SequenceIndex\in\NonnegativeIntegers}$ is \emph{collectively compact} if
\[
 \overline{\bigcup_{\SequenceIndex\in\NonnegativeIntegers}
 \AmbientOperator_\SequenceIndex(\ClosedUnitBall)}
 \text{ is compact}.
\]
The operators \emph{converge strongly} to a bounded linear operator $\AmbientOperator\colon\BanachAmbientSpace\to\BanachAmbientSpace$ if $\norm{\AmbientOperator_\SequenceIndex\FirstVector-\AmbientOperator\FirstVector}\to0$ for every $\FirstVector\in\BanachAmbientSpace$.

\begin{theorem}[{Theorem of Riesz~\cite{Riesz1916Linear} and Schauder~\cite{Schauder1930Linear}, see also \cite[Theorem~14.19, pp.~361--362]{Muscat2024}}]\label{thm:compact-spectrum}
Let $\BanachAmbientSpace$ be a complex Banach space and $\CompactOperator\colon\BanachAmbientSpace\to\BanachAmbientSpace$ a compact linear operator.
\begin{enumerate}[label=\textup{(\roman*)}]
\item The spectrum $\spec(\CompactOperator)$ is at most countable and its only possible accumulation point is $0$.
In particular, for each $a>0$ only finitely many spectral values $\Eigenvalue$ satisfy $|\Eigenvalue|\geq a$.
\item Every nonzero spectral value $\Eigenvalue$ is an eigenvalue and $0<\dim\ker(\CompactOperator-\Eigenvalue\IdentityMap)<\infty$.
\item Let $\BanachAmbientSpace^*$ denote the continuous complex-linear functionals on $\BanachAmbientSpace$, and define the \emph{dual operator} $\CompactOperator^{\TransposeSymbol}\colon\BanachAmbientSpace^*\to\BanachAmbientSpace^*$ by $\CompactOperator^{\TransposeSymbol}\ell=\ell\circ\CompactOperator$.
The operators $\CompactOperator$ and $\CompactOperator^{\TransposeSymbol}$ have the same nonzero eigenvalues, and for each such eigenvalue,
\[
 \dim\ker(\CompactOperator-\Eigenvalue\IdentityMap)
 =\dim\ker(\CompactOperator^{\TransposeSymbol}-\Eigenvalue\IdentityMap).
\]
\end{enumerate}
\end{theorem}

\begin{theorem}[{Spectral inclusion, Anselone and Palmer~\cite[Theorem~5.3]{AnselonePalmer1968}}]\label{thm:spectral-inclusion}
Let $\AmbientOperator_\SequenceIndex,\AmbientOperator\colon\BanachAmbientSpace\to\BanachAmbientSpace$ be bounded linear operators on a complex Banach space.
Assume that $\AmbientOperator_\SequenceIndex\to\AmbientOperator$ strongly and that $\{\AmbientOperator_\SequenceIndex-\AmbientOperator\}_{\SequenceIndex\in\NonnegativeIntegers}$ is collectively compact.
Then for every open set $\SpectralNeighborhood\subset\ComplexNumbers$ containing $\spec(\AmbientOperator)$, there is $\SequenceIndex_0\in\NonnegativeIntegers$ such that for every $\SequenceIndex\geq\SequenceIndex_0$ we have
\[
 \spec(\AmbientOperator_\SequenceIndex)\subset\SpectralNeighborhood.
\]
\end{theorem}

Under the hypotheses of the spectral inclusion theorem of Anselone and Palmer (\Cref{thm:spectral-inclusion}), fix $\ResolventParameter\in\ComplexNumbers\setminus\spec(\AmbientOperator)$ and define
\[
 {\ResolventErrorOperator_\SequenceIndex}
 =(\ResolventParameter\IdentityMap-\AmbientOperator)^{-1}
 (\AmbientOperator_\SequenceIndex-\AmbientOperator).
\]
The operators $\ResolventErrorOperator_\SequenceIndex$ converge strongly to $0$ and form a collectively compact family.
By Anselone and Palmer (\cite[Proposition~3.3]{AnselonePalmer1968}), $\norm{\ResolventErrorOperator_\SequenceIndex^2}\to0$.
For all sufficiently large $\SequenceIndex$ the Neumann series makes $\IdentityMap-\ResolventErrorOperator_\SequenceIndex^2$ invertible, and
\[
 (\IdentityMap-\ResolventErrorOperator_\SequenceIndex)^{-1}
 =(\IdentityMap+\ResolventErrorOperator_\SequenceIndex)
 (\IdentityMap-\ResolventErrorOperator_\SequenceIndex^2)^{-1}.
\]
Consequently,
\[
 \ResolventParameter\IdentityMap-\AmbientOperator_\SequenceIndex
 =(\ResolventParameter\IdentityMap-\AmbientOperator)
 (\IdentityMap-\ResolventErrorOperator_\SequenceIndex)
\]
is invertible, so $\ResolventParameter\notin\spec(\AmbientOperator_\SequenceIndex)$.
This is the argument of Anselone and Palmer (\cite[proof of Theorem~5.3]{AnselonePalmer1968}) for a fixed $\ResolventParameter$.

\begin{definition}[{Riesz projection, \cite[Definition~14.23 and Proposition~14.28]{Muscat2024}}]\label{def:riesz-projection}
Let $\CompactOperator\colon\BanachAmbientSpace\to\BanachAmbientSpace$ be a bounded linear operator on a complex Banach space.
Let $\SpectralContour=\sum_{k=0}^{m-1}\gamma_k$ be a finite sum of positively oriented rectifiable Jordan curves disjoint from $\spec(\CompactOperator)$.
Integration and winding number along $\SpectralContour$ are the sums of the corresponding integrals and winding numbers for the individual curves.
Suppose the winding number of $\SpectralContour$ at each point of $\spec(\CompactOperator)$ is either $0$ or $1$.
The \emph{Riesz projection} of $\CompactOperator$ is the bounded linear operator $\RieszProjection{\CompactOperator}{\SpectralContour}\colon\BanachAmbientSpace\to\BanachAmbientSpace$ defined by
\[
 \RieszProjection{\CompactOperator}{\SpectralContour}
 =\frac{1}{2\CircleConstant\ImaginaryUnit}
 \int_{\SpectralContour}(\ResolventParameter\IdentityMap-\CompactOperator)^{-1}
 \,\IntegrationDifferential\ResolventParameter,
\]
where the line integral is taken in the operator norm.
If $\Eigenvalue\in\spec(\CompactOperator)$ and $\SpectralContour$ has winding number $1$ at $\Eigenvalue$ and $0$ at every point of $\spec(\CompactOperator)\setminus\{\Eigenvalue\}$, then we call $\RieszProjection{\CompactOperator}{\SpectralContour}$ the \emph{Riesz projection associated with $\Eigenvalue$}.
\end{definition}

\begin{lemma}[{Consequence of the holomorphic functional calculus of Riesz and Dunford, \cite[Examples~14.24(1) and Proposition~14.28]{Muscat2024}}]\label{lem:riesz-projection-properties}
Let $\CompactOperator$ and $\SpectralContour$ be as in \Cref{def:riesz-projection}.
Then we have
\[
 \RieszProjection{\CompactOperator}{\SpectralContour}^2
 =\RieszProjection{\CompactOperator}{\SpectralContour},
\]
and
\[
 \CompactOperator\RieszProjection{\CompactOperator}{\SpectralContour}
 =\RieszProjection{\CompactOperator}{\SpectralContour}\CompactOperator.
\]
Moreover, $\RieszProjection{\CompactOperator}{\SpectralContour}$ depends only on the subset of $\spec(\CompactOperator)$ at which $\SpectralContour$ has winding number $1$.
\end{lemma}

\begin{theorem}[{Convergence of spectral projections, Anselone and Palmer~\cite[Proposition~6.3]{AnselonePalmer1968}}]\label{thm:riesz-convergence}
Let $\AmbientOperator_\SequenceIndex,\AmbientOperator\colon\BanachAmbientSpace\to\BanachAmbientSpace$ satisfy the operator hypotheses of the spectral inclusion theorem of Anselone and Palmer (\Cref{thm:spectral-inclusion}).
Let $\SpectralContour$ be as in \Cref{def:riesz-projection} for $\AmbientOperator$, and put
\[
 \AmbientProjection=\RieszProjection{\AmbientOperator}{\SpectralContour}.
\]
Then for all sufficiently large $\SequenceIndex$ the curves $\SpectralContour$ satisfy the conditions of \Cref{def:riesz-projection} for $\AmbientOperator_\SequenceIndex$. Define the corresponding projections by
\[
 \AmbientProjection_\SequenceIndex
 =\RieszProjection{\AmbientOperator_\SequenceIndex}{\SpectralContour}.
\]
The projections $\AmbientProjection_\SequenceIndex\colon\BanachAmbientSpace\to\BanachAmbientSpace$ satisfy
\[
 \dim\range\AmbientProjection_\SequenceIndex=\dim\range\AmbientProjection.
\]
For every $\FirstVector\in\BanachAmbientSpace$, we also have
\[
 \norm{\AmbientProjection_\SequenceIndex\FirstVector-\AmbientProjection\FirstVector}\to0.
\]
If $\AmbientOperator_\SequenceIndex$ and $\AmbientOperator$ are compact and $\SpectralContour$ has winding number $0$ at $0$, then $\dim\range\AmbientProjection_\SequenceIndex$ and $\dim\range\AmbientProjection$ are finite.
\end{theorem}

For a bounded linear operator $\CompactOperator\colon\BanachAmbientSpace\to\BanachAmbientSpace$ on a complex Banach space and $\Eigenvalue\in\ComplexNumbers$, denote the \emph{generalized eigenspace} at $\Eigenvalue$ by
\[
 \GeneralizedEigenspace{\CompactOperator}{\Eigenvalue}
 =\bigcup_{\RootOrder\geq1}
  \ker\bigl((\CompactOperator-\Eigenvalue\IdentityMap)^\RootOrder\bigr).
\]
\begin{theorem}[{Anselone and Palmer~\cite[Section~7, p.~429]{AnselonePalmer1968}, \cite[Proposition~14.28 and Examples~14.29]{Muscat2024}}]\label{lem:riesz-range-generalized-eigenspace}
Let $\CompactOperator\colon\BanachAmbientSpace\to\BanachAmbientSpace$ be a compact linear operator on a complex Banach space.
For each $\Eigenvalue\in\spec(\CompactOperator)\setminus\{0\}$ denote the Riesz projection associated with $\Eigenvalue$ by $\EigenvalueProjection{\CompactOperator}{\Eigenvalue}$.
Then
\[
 \range\EigenvalueProjection{\CompactOperator}{\Eigenvalue}
 =\GeneralizedEigenspace{\CompactOperator}{\Eigenvalue}.
\]
More generally, let $\SelectedEigenvalues\subset\spec(\CompactOperator)\setminus\{0\}$ be finite, and suppose $\SpectralContour$ satisfies the conditions of \Cref{def:riesz-projection} for $\CompactOperator$.
Assume the winding number of $\SpectralContour$ is $1$ at each point of $\SelectedEigenvalues$ and $0$ at each point of $\spec(\CompactOperator)\setminus\SelectedEigenvalues$.
Then
\[
 \RieszProjection{\CompactOperator}{\SpectralContour}
 =\sum_{\Eigenvalue\in\SelectedEigenvalues}
   \EigenvalueProjection{\CompactOperator}{\Eigenvalue},
\]
and
\[
 \range\RieszProjection{\CompactOperator}{\SpectralContour}
 =\bigoplus_{\Eigenvalue\in\SelectedEigenvalues}
   \GeneralizedEigenspace{\CompactOperator}{\Eigenvalue}.
\]
\end{theorem}

\begin{proof}
Since $\Eigenvalue\ne0$, the operator $-\CompactOperator/\Eigenvalue$ is compact, and $\CompactOperator-\Eigenvalue\IdentityMap=-\Eigenvalue(\IdentityMap-\CompactOperator/\Eigenvalue)$ has finite ascent and descent (\cite[Proposition~14.18]{Muscat2024}).
Also $\Eigenvalue$ is an isolated eigenvalue by the theorem of Riesz and Schauder (\Cref{thm:compact-spectrum}).
Hence by Anselone and Palmer (\cite[Section~7, p.~429]{AnselonePalmer1968}) the range of its Riesz projection equals its generalized eigenspace.
For a finite set $\SelectedEigenvalues$, spectral decomposition (\cite[Proposition~14.28]{Muscat2024}) identifies the corresponding projection with the sum of the individual projections and its range with the direct sum of their ranges.
\end{proof}

\begin{lemma}[{Generalized eigenspaces of a normal operator}]\label{lem:normal-power-kernel}
Suppose $T$ is a bounded normal operator on a complex Hilbert space, so $T^*T=TT^*$.
For all $\lambda\in\ComplexNumbers$ and integers $m\geq1$,
\[
 \ker\bigl((T-\lambda\IdentityMap)^m\bigr)
 =\ker(T-\lambda\IdentityMap)
 =\GeneralizedEigenspace{T}{\lambda}.
\]
\end{lemma}
\begin{proof}
Since $N=T-\lambda\IdentityMap$ is normal, $\norm{Nx}=\norm{N^*x}$ and hence $\ker N=\ker N^*$.
Thus $\range N$ is orthogonal to $\ker N$.
If $N^2x=0$ then $Nx\in\range N\cap\ker N=\{0\}$.
Applying this to $N^{m-2}x$ yields $\ker N^m=\ker N^{m-1}$ for $m\geq2$.
Taking the union over $m\geq1$ proves the assertion.
For normal operators on Banach spaces Fong also proved $\ker N^2=\ker N$ (\cite[Theorem~A, Corollary~3]{Fong1979Normal}).
\end{proof}

\begin{lemma}[{Reversing products and generalized eigenspaces}]\label{lem:product-power-kernel}
Let $E,H$ be complex Banach spaces, $R\colon E\to H$ and $S\colon H\to E$ bounded linear maps, and put $B=SR$ and $T=RS$.
If $\GeneralizedEigenspace{T}{\lambda}=\ker(T-\lambda\IdentityMap)$ at $\lambda\in\ComplexNumbers\setminus\{0\}$, then $\GeneralizedEigenspace{B}{\lambda}=\ker(B-\lambda\IdentityMap)$.
Moreover, $R$ induces a linear isomorphism
\[
 R\colon\ker(B-\lambda\IdentityMap)\longrightarrow\ker(T-\lambda\IdentityMap),
\]
whose inverse is the restriction of $\lambda^{-1}S$.
\end{lemma}
\begin{proof}
Let $f\in\GeneralizedEigenspace{B}{\lambda}$ and choose an integer $m\geq1$ with $(B-\lambda\IdentityMap)^m f=0$.
Since $RB=TR$,
\[
 (T-\lambda\IdentityMap)^m Rf=R(B-\lambda\IdentityMap)^m f=0,
\]
so by hypothesis $Rf\in\ker(T-\lambda\IdentityMap)$.

Set $g=(B-\lambda\IdentityMap)f$.
Then $Rg=(T-\lambda\IdentityMap)Rf=0$ and $Bg=SRg=0$.
Thus
\[
 0=(B-\lambda\IdentityMap)^m f
   =(B-\lambda\IdentityMap)^{m-1}g
   =(-\lambda)^{m-1}g.
\]
Since $\lambda\neq0$, we obtain $g=0$, and hence $f\in\ker(B-\lambda\IdentityMap)$.
The reverse inclusion follows by taking $m=1$ in the definition of the generalized eigenspace.

By $RB=TR$ and $BS=ST$, the map $R$ sends the $\lambda$-eigenspace of $B$ into the $\lambda$-eigenspace of $T$, and $S$ maps in the opposite direction.
For $f\in\ker(B-\lambda\IdentityMap)$ and $u\in\ker(T-\lambda\IdentityMap)$ we have
\[
 (\lambda^{-1}S)Rf=\lambda^{-1}Bf=f,\qquad
 R(\lambda^{-1}S)u=\lambda^{-1}Tu=u.
\]
Hence the restrictions of $R$ and $\lambda^{-1}S$ to these eigenspaces are inverse to one another.
\end{proof}
\subsection{Square roots of positive definite matrices}\label{subsec:spectral-preliminaries-matrices}
Sylvester proved orthogonal diagonalizability of real quadratic forms.
For a positive definite Gram matrix, taking the positive square roots of the eigenvalues in this diagonalization produces its positive definite square root.

\begin{theorem}[{Consequence of Sylvester's orthogonal diagonalization~\cite{Sylvester1852Orthogonal}, see also \cite[Theorems~7.29 and~7.39]{Axler2024}}]\label{thm:matrix-square-root}
Every real symmetric matrix has an orthonormal basis of eigenvectors.
If it is positive definite, then it has a unique symmetric positive definite square root.
In an orthonormal eigenbasis, this square root is the diagonal matrix whose entries are the positive square roots of the original eigenvalues.
\end{theorem}

Fix a nonempty compact metric space $Z$.
Let $n\geq1$ and arrange $n$ continuous real functions as a row vector $\ProjectedBasisRow\in C(Z,\RealNumbers^{1\times n})$.
For $\mu\in\ProbabilityMeasures(Z)$ define the \emph{Gram matrix} and, when it is positive definite, the \emph{orthonormalization} by
\[
 \GramMatrix(\mu,\ProjectedBasisRow)=\int_Z\ProjectedBasisRow(z)^{\TransposeSymbol}\ProjectedBasisRow(z)\,\IntegrationDifferential\mu(z),
 \qquad
 \GramNormalization(\mu,\ProjectedBasisRow)=\ProjectedBasisRow\,\GramMatrix(\mu,\ProjectedBasisRow)^{-1/2}.
\]
\begin{lemma}\label{lem:continuous-gram-normalization}
Equip $\ProbabilityMeasures(Z)$ with the weak topology and $C(Z,\RealNumbers^{1\times n})$ with the uniform norm.
Then $\GramMatrix$ is continuous, and $\GramNormalization$ is continuous on the open set where $\GramMatrix>0$.
Its Gram matrix is the $n$-dimensional identity matrix $\IdentityMatrix_n$, and for any $U\in O(n)$,
\[
 \GramNormalization(\mu,\ProjectedBasisRow U)=\GramNormalization(\mu,\ProjectedBasisRow)U.
\]
\end{lemma}
\begin{proof}
Uniform convergence of the rows implies uniform convergence of products of their entries.
Combining this with weak convergence against each limit product proves convergence of the integrals.
Near a positive definite matrix, all eigenvalues lie in a common interval $[a,b]\subset(0,\infty)$.
Uniformly approximate $t^{-1/2}$ there by polynomials and use Sylvester's orthogonal diagonalization (\Cref{thm:matrix-square-root}).
The inverse square root is then a uniform limit of matrix polynomials on that neighborhood.
This proves continuity of $\GramNormalization$.
If $\SampleGramMatrix=\GramMatrix(\mu,\ProjectedBasisRow)$, the normalized Gram matrix is $\SampleGramMatrix^{-1/2}\SampleGramMatrix\SampleGramMatrix^{-1/2}=\IdentityMatrix_n$.
Moreover, $\GramMatrix(\mu,\ProjectedBasisRow U)=U^{\TransposeSymbol}\SampleGramMatrix U$ and uniqueness of the positive square root implies $(U^{\TransposeSymbol}\SampleGramMatrix U)^{-1/2}=U^{\TransposeSymbol}\SampleGramMatrix^{-1/2}U$.
\end{proof}

\section{Finite-dimensional approximation of distances}\label{sec:spectral}
We approximate the Kuratowski map of each point of a compact metric space by a function in a finite-dimensional vector space.
The selected probability measure makes this space invariant under isometries.
Using the spectral results of Section~\ref{sec:prelim-spectral}, we first establish approximation on a single space.
We then compare the functions and coordinates on nearby spaces.

A simple construction explains how finite-dimensional coordinates approximate distances.
Let $\MetricSpace{\BaseCarrier}{\MetricSymbol}$ be a nonempty compact metric space.
Fix $\ErrorTolerance>0$ and choose a finite $\ErrorTolerance$-net $\{\NetPoint_0,\ldots,\NetPoint_{\NetSize-1}\}$.
Define the map
\[
 \NetCoordinateMap\colon\BaseCarrier\to\RealNumbers^{\NetSize}
\]
\[
 \NetCoordinateMap(\BasePoint)
 =\bigl(\MetricSymbol(\BasePoint,\NetPoint_0),\ldots,
        \MetricSymbol(\BasePoint,\NetPoint_{\NetSize-1})\bigr).
\]
For all $\BasePoint,\ComparisonPoint\in\BaseCarrier$ we have
\begin{equation}\label{eq:finite-net-distance-approximation}
 \MetricSymbol(\BasePoint,\ComparisonPoint)-2\ErrorTolerance
 \leq\norm{\NetCoordinateMap(\BasePoint)-\NetCoordinateMap(\ComparisonPoint)}_\infty
 \leq\MetricSymbol(\BasePoint,\ComparisonPoint).
\end{equation}
A construction with finitely many functions also underlies Shioya's approximation of metric measure spaces (\cite[author's manuscript, Section~4.4, Definition~4.37, Theorem~4.43 and Corollary~4.44]{ShioyaMetricMeasure}).
More precisely, for any mm-space $\mathbf X=(X,d_X,\mu_X)$ and $\ErrorTolerance>0$ there are an integer $N\geq1$ and real $1$-Lipschitz functions $f_1,\ldots,f_N$ on $X$ such that the map $\Phi=(f_1,\ldots,f_N)$ satisfies
\[
 \BoxDistance\bigl(\mathbf X,(\RealNumbers^N,\norm{\cdot}_\infty,
                         \Phi_*\mu_X)\bigr)<\ErrorTolerance.
\]
Here $\BoxDistance$ is the distance in \Cref{def:box-distance}.
Estimate \eqref{eq:finite-net-distance-approximation} controls the uniform error for every pair of points in a single compact space.
Our construction must also be continuous on a Gromov--Hausdorff neighborhood and compatible with isometries.
Finite nets on individual spaces do not specify such a continuous choice.
At a one-point center we use the zero map into $\RealNumbers$ and the absolute value norm throughout the neighborhood in \Cref{thm:local-models}.
When the center has at least two points, a finite spectral subspace determines a coordinate class that varies continuously with the space and is compatible with pullback by isometries.
In an orthonormal basis the coordinates are the coefficients of the function obtained by applying the spectral projection and then the distance operator to the normalized bump.
Its uniform norm defines the norm on the coefficient space.
An orthogonal change of basis replaces each coordinate vector $v$ by $\OrthogonalChange v$.
In the new coordinates the vector $v$ represents $\OrthogonalChange^{-1}v$ in the old coordinates, so the new norm is $\LocalNorm{\ComparisonCarrier}{\ComparisonMetric}\circ\OrthogonalChange^{-1}$.
Thus distances between coordinate vectors are independent of the basis.

Throughout this section, $\SelectedMeasure{\BaseCarrier}{\MetricSymbol}$ denotes the assignment of \Cref{thm:invariant-measures}.
The two subscripts denote the carrier $\BaseCarrier$ and its metric $\MetricSymbol$.
\subsection{The distance operator and its spectral subspaces}\label{subsec:spectral-distance-operator}
Fix a nonempty compact metric space $\MetricSpace{\BaseCarrier}{\MetricSymbol}$.
Equip the space
\[
 \HilbertDomain_{\BaseCarrier,\MetricSymbol} =\LebesgueSpace^2(\BaseCarrier ,\SelectedMeasure{\BaseCarrier}{\MetricSymbol}),
\]
with the inner product
\[
 \langle \HilbertFunction ,\HilbertTestFunction \rangle_{\HilbertDomain_{\BaseCarrier,\MetricSymbol} }=\int_\BaseCarrier  \HilbertFunction  \HilbertTestFunction \,\IntegrationDifferential \SelectedMeasure{\BaseCarrier}{\MetricSymbol}.
\]
This is a separable real Hilbert space.
Separability follows from density of continuous functions in $\LebesgueSpace^2$ and separability of $\ContinuousFunctions(\BaseCarrier )$ for compact metric spaces (\cite[Section~6.4]{Muscat2024}).
Define the \emph{distance operator} $\DistanceOperator{\BaseCarrier}{\MetricSymbol}\colon\HilbertDomain_{\BaseCarrier,\MetricSymbol}\to\HilbertDomain_{\BaseCarrier,\MetricSymbol}$ for $\HilbertFunction\in\HilbertDomain_{\BaseCarrier,\MetricSymbol}$ and $\BasePoint\in\BaseCarrier$ by
\[
 \DistanceOperator{\BaseCarrier}{\MetricSymbol}\HilbertFunction (\BasePoint )
 =\int_\BaseCarrier  \MetricSymbol (\BasePoint ,\IntegrationPoint )\HilbertFunction (\IntegrationPoint )\,\IntegrationDifferential \SelectedMeasure{\BaseCarrier}{\MetricSymbol}(\IntegrationPoint ).
\]
The integral is independent of the representative of $\HilbertFunction $.
This is the distance kernel operator of Maria, Oudot, and Solomon (\cite[Definition~11]{MariaOudotSolomon2020}).

\begin{lemma}\label{lem:distance-operator}
The operator $\DistanceOperator{\BaseCarrier}{\MetricSymbol}\colon \HilbertDomain_{\BaseCarrier,\MetricSymbol} \to \HilbertDomain_{\BaseCarrier,\MetricSymbol} $ is compact and self-adjoint.
It maps $\HilbertDomain_{\BaseCarrier,\MetricSymbol} $ into $\ContinuousFunctions(\BaseCarrier )$, and for every $\HilbertFunction \in \HilbertDomain_{\BaseCarrier,\MetricSymbol} $ and $\BasePoint ,\ComparisonPoint \in \BaseCarrier $,
\begin{equation}\label{eq:operator-uniform-bound-1}
 \norm{\DistanceOperator{\BaseCarrier}{\MetricSymbol}\HilbertFunction }_\infty
 \leq\diam\MetricSpace{\BaseCarrier }{\MetricSymbol }\norm{\HilbertFunction }_2,
\end{equation}
and
\begin{equation}\label{eq:operator-uniform-bound-2}
 |\DistanceOperator{\BaseCarrier}{\MetricSymbol}\HilbertFunction (\BasePoint )-\DistanceOperator{\BaseCarrier}{\MetricSymbol}\HilbertFunction (\ComparisonPoint )|
 \leq \MetricSymbol (\BasePoint ,\ComparisonPoint )\norm{\HilbertFunction }_2.
\end{equation}
Every eigenfunction with nonzero eigenvalue has a unique continuous representative.
\end{lemma}

Maria, Oudot, and Solomon also proved compactness and self-adjointness of the distance operator (\cite[Propositions~12 and~13]{MariaOudotSolomon2020}).

\begin{proof}
Let $\HilbertFunction\in\HilbertDomain_{\BaseCarrier,\MetricSymbol}$ and $\BasePoint,\ComparisonPoint\in\BaseCarrier$.
The Cauchy--Schwarz inequality (\cite[Lemma~3.2.27]{BreitGmeineder2026}) and $\SelectedMeasure{\BaseCarrier}{\MetricSymbol}(\BaseCarrier)=1$ imply
\[
 \int_\BaseCarrier |\HilbertFunction |\,\IntegrationDifferential \SelectedMeasure{\BaseCarrier}{\MetricSymbol}\leq\norm{\HilbertFunction }_2.
\]
Consequently,
\begin{equation}\label{eq:distance-operator-pointwise-bounds-1}
|\DistanceOperator{\BaseCarrier}{\MetricSymbol}\HilbertFunction (\BasePoint )|
 \leq\int_\BaseCarrier  \MetricSymbol (\BasePoint ,\IntegrationPoint )|\HilbertFunction (\IntegrationPoint )|\,\IntegrationDifferential \SelectedMeasure{\BaseCarrier}{\MetricSymbol}(\IntegrationPoint )
 \leq\diam\MetricSpace{\BaseCarrier }{\MetricSymbol }\norm{\HilbertFunction }_2,
\end{equation}
and
\begin{equation}\label{eq:distance-operator-pointwise-bounds-2}
\begin{aligned}
|\DistanceOperator{\BaseCarrier}{\MetricSymbol}\HilbertFunction (\BasePoint )-\DistanceOperator{\BaseCarrier}{\MetricSymbol}\HilbertFunction (\ComparisonPoint )|
 &\leq\int_\BaseCarrier |\MetricSymbol (\BasePoint ,\IntegrationPoint )-\MetricSymbol (\ComparisonPoint ,\IntegrationPoint )|\,|\HilbertFunction (\IntegrationPoint )|\,\IntegrationDifferential \SelectedMeasure{\BaseCarrier}{\MetricSymbol}(\IntegrationPoint )\\
 &\leq \MetricSymbol (\BasePoint ,\ComparisonPoint )\norm{\HilbertFunction }_2.
\end{aligned}
\end{equation}
By \eqref{eq:distance-operator-pointwise-bounds-1}, \eqref{eq:distance-operator-pointwise-bounds-2}, and the Arzel\`a--Ascoli theorem (\Cref{thm:arzela-ascoli}), the operator $\DistanceOperator{\BaseCarrier}{\MetricSymbol}\colon
\HilbertDomain_{\BaseCarrier,\MetricSymbol}\to\ContinuousFunctions(\BaseCarrier)$ is compact.
The inclusion $\ContinuousFunctions(\BaseCarrier)\to\HilbertDomain_{\BaseCarrier,\MetricSymbol}$ is bounded, so the operator on $\HilbertDomain_{\BaseCarrier,\MetricSymbol}$ is also compact.
Self-adjointness follows from symmetry of the kernel and Fubini's theorem (\cite[Theorem~2.8.20(a)]{BreitGmeineder2026}).

If $\DistanceOperator{\BaseCarrier}{\MetricSymbol}\HilbertFunction =\Eigenvalue  \HilbertFunction $ in $\HilbertDomain_{\BaseCarrier,\MetricSymbol} $ and $\Eigenvalue \ne0$, then $\Eigenvalue ^{-1}\DistanceOperator{\BaseCarrier}{\MetricSymbol}\HilbertFunction $ is a continuous representative of $\HilbertFunction $.
The representative is unique because $\SelectedMeasure{\BaseCarrier}{\MetricSymbol}$ has full support.
\end{proof}

For each eigenfunction with nonzero eigenvalue we use its unique continuous representative.
Apply the spectral theorem of Hilbert and Schmidt (\Cref{prop:spectral-basics}) to $\DistanceOperator{\BaseCarrier}{\MetricSymbol}$ on $\HilbertDomain_{\BaseCarrier,\MetricSymbol}$ and choose an orthonormal family $\{\Eigenfunction_\CoordinateIndex\}_{\CoordinateIndex\in\EigenIndexSet}$ and the corresponding nonzero eigenvalues $\{\Eigenvalue_\CoordinateIndex\}_{\CoordinateIndex\in\EigenIndexSet}$ as in that theorem.
For $\SpectralCutoff >0$ put
\[
 {\EigenIndexSet_{\BaseCarrier,\MetricSymbol,\SpectralCutoff}=\{\CoordinateIndex \in \EigenIndexSet \mid|\Eigenvalue _\CoordinateIndex |>\SpectralCutoff \},}
\]
and
\[
 \SpectralSubspace{\BaseCarrier,\MetricSymbol}{\SpectralCutoff}
 =\LinearSpan \{\Eigenfunction _\CoordinateIndex \mid \CoordinateIndex \in \EigenIndexSet_{\BaseCarrier,\MetricSymbol,\SpectralCutoff}\}.
\]
The operator $\DistanceOperatorSymbol=\DistanceOperator{X}{d}$ is compact.
By the Riesz--Schauder theorem (\Cref{thm:compact-spectrum}), only finitely many eigenvalues $\lambda$ satisfy $|\lambda|>a$, and each eigenspace $\ker(\DistanceOperatorSymbol-\lambda \IdentityMap)$ is finite-dimensional.
Thus the \emph{spectral subspace} $\SpectralSubspace{X,d}{a}$, their direct sum, is finite-dimensional.
The orthogonal projection $\SpectralProjection{\BaseCarrier,\MetricSymbol}{\SpectralCutoff}
 \colon\HilbertDomain_{\BaseCarrier,\MetricSymbol}\to\SpectralSubspace{\BaseCarrier,\MetricSymbol}{\SpectralCutoff}$ is
\begin{equation}\label{eq:spectral-projection-definition}
 \SpectralProjection{\BaseCarrier,\MetricSymbol}{\SpectralCutoff}\HilbertFunction
 =\sum_{\CoordinateIndex \in \EigenIndexSet_{\BaseCarrier,\MetricSymbol,\SpectralCutoff}}\langle \HilbertFunction ,\Eigenfunction _\CoordinateIndex \rangle_{\HilbertDomain_{\BaseCarrier,\MetricSymbol} }\Eigenfunction _\CoordinateIndex.
\end{equation}
Equivalently,
\[
 \begin{aligned}
 \SpectralSubspace{\BaseCarrier,\MetricSymbol}{\SpectralCutoff}
 =\bigoplus_{\substack{\Eigenvalue \in\spec(\DistanceOperator{\BaseCarrier}{\MetricSymbol})\\|\Eigenvalue |>\SpectralCutoff }}
       \ker(\DistanceOperator{\BaseCarrier}{\MetricSymbol}-\Eigenvalue \IdentityMap),
 \end{aligned}
\]
and
\[
 \range\SpectralProjection{\BaseCarrier,\MetricSymbol}{\SpectralCutoff}=\SpectralSubspace{\BaseCarrier,\MetricSymbol}{\SpectralCutoff}.
\]
Write $\SelectedMeasureSymbol=\SelectedMeasure{X}{d}$ and $D=\diam(X,d)$.
Counting eigenvalues with multiplicity, self-adjointness and the integral kernel formula for the Hilbert--Schmidt norm (\cite[Lemma~6.10 and (6.16)]{Teschl2014}) imply
\[
 a^2\dim \SpectralSubspace{X,d}{a}\leq\sum_{|\lambda_k|>a}|\lambda_k|^2
 \leq\|\DistanceOperatorSymbol\|_{\mathrm{HS}}^2
 =\int_X\!\int_X d(x,y)^2\,\IntegrationDifferential\SelectedMeasureSymbol(x)\,\IntegrationDifferential\SelectedMeasureSymbol(y)\leq D^2.
\]
The last inequality uses $\SelectedMeasureSymbol(X)=1$.
The spectral subspace $\SpectralSubspace{X,d}{a}$ and its orthogonal projection $\SpectralProjection{X,d}{a}$ are independent of the orthonormal system $\{\Eigenfunction_\CoordinateIndex\}_{\CoordinateIndex\in\EigenIndexSet}$.
For every $\HilbertFunction\in\HilbertDomain_{\BaseCarrier,\MetricSymbol}$, as $\SpectralCutoff\downarrow0$ the vectors $\SpectralProjection{\BaseCarrier,\MetricSymbol}{\SpectralCutoff}\HilbertFunction$ converge in the norm of $\HilbertDomain_{\BaseCarrier,\MetricSymbol}$ to the orthogonal projection of $\HilbertFunction$ onto $(\ker\DistanceOperator{\BaseCarrier}{\MetricSymbol})^\perp$.
\subsection{Approximation of Kuratowski maps}\label{subsec:spectral-distance-approximation}
We approximate all the Kuratowski maps in the uniform norm by functions in a single finite-dimensional spectral subspace.
For $\BasePoint\in\BaseCarrier$, we call the function
\[
 \MetricSymbol_\BasePoint=\MetricSymbol(\BasePoint,\cdot)
 \in\ContinuousFunctions(\BaseCarrier)
\]
the \emph{Kuratowski map} of $\BasePoint$.
Its values are the distances from $\BasePoint$ to the points of $\BaseCarrier$.
We follow the terminology of Dress, Huber, Koolen, Moulton, and Spillner for finite metric spaces (\cite[arXiv v1, Section~2, p.~5]{DressHuberKoolenMoultonSpillner2010}).
The map
\[
 \BaseCarrier\to\ContinuousFunctions(\BaseCarrier),
 \qquad\BasePoint\longmapsto\MetricSymbol_\BasePoint
\]
is an isometric embedding for the uniform norm, because
\[
 \norm{\MetricSymbol_\BasePoint-\MetricSymbol_\ComparisonPoint}_\infty
 =\MetricSymbol(\BasePoint,\ComparisonPoint).
\]
This is the distance-function construction of the Kuratowski embedding.
Since $\MetricSpace{\BaseCarrier}{\MetricSymbol}$ is compact, $\MetricSymbol_\BasePoint$ is bounded.
The following theorem approximates all the Kuratowski maps in one spectral subspace.

\begin{theorem}\label{lem:spectral-density}
Let $\MetricSpace{X}{d}$ be a nonempty compact metric space and let $\SmallError>0$.
There exist a radius $r>0$ and a spectral parameter $\SpectralCutoff>0$ with $\pm\SpectralCutoff\notin\spec(\DistanceOperator{X}{d})$ such that the following holds.
For each $x\in X$, set $b_x=\NormalizedBump{X}{d}{x}{r}{\SelectedMeasure{X}{d}}$.
The function
\[
 u_x=\DistanceOperator{X}{d}\SpectralProjection{X,d}{\SpectralCutoff}b_x
 \in\SpectralSubspace{X,d}{\SpectralCutoff}
\]
satisfies the uniform estimate
\begin{equation}\label{eq:projected-bump-error}
 \sup_{x\in X}\|d_x-u_x\|_{\infty,X}<\SmallError.
\end{equation}
If $\Card(\BaseCarrier)\geq2$, we may also require $\dim\SpectralSubspace{\BaseCarrier,\MetricSymbol}{\SpectralCutoff}\geq1$.
\end{theorem}

\begin{proof}
Write $\SelectedMeasureSymbol=\SelectedMeasure{X}{d}$, $H=L^2(X,\SelectedMeasureSymbol;\RealNumbers)$, $\DistanceOperatorSymbol=\DistanceOperator{X}{d}$, and $\SpectralProjectionSymbol_a=\SpectralProjection{X,d}{a}$.
For $f\in H$ and $z\in X$, the continuous representative of $\DistanceOperatorSymbol f$ satisfies $(\DistanceOperatorSymbol f)(z)=\langle f,d_z\rangle_H$.
If $f\in\ker\DistanceOperatorSymbol$, full support implies that this representative vanishes everywhere.
Hence $d_z\in(\ker\DistanceOperatorSymbol)^\perp$.
Since $\|d_z-d_w\|_2\leq d(z,w)$, the family $\{d_z:z\in X\}$ is compact in $H$.

Fix $0<r<\SmallError$ and set $b_x=\NormalizedBump{X}{d}{x}{r}{\SelectedMeasureSymbol}$.
By \eqref{eq:normalized-bump-definition}, the function $b_x$ is nonnegative, has integral $1$, and satisfies $b_x(w)=0$ whenever $d(x,w)\geq r/2$.
The triangle inequality bounds $|(\DistanceOperatorSymbol b_x)(z)-d(x,z)|$ by the mean of $d(x,w)$ with respect to $b_x(w)\,\IntegrationDifferential\SelectedMeasureSymbol(w)$.
Equality holds at $z=x$, so
\begin{equation}\label{eq:bump-distance-approximation}
 {\begin{aligned}
 \|d_x-\DistanceOperatorSymbol b_x\|_\infty
 &=(\DistanceOperatorSymbol b_x)(x)\\
 &=\int_X d(x,w)b_x(w)\,\IntegrationDifferential\SelectedMeasureSymbol(w)\leq r/2.
 \end{aligned}}
\end{equation}
The numerator of the bump is continuous in $(x,w)$, and its denominator is positive and continuous in $x$.
Thus $x\mapsto b_x$ is continuous into $C(X)$, and $M=\max_{x\in X}\|b_x\|_2$ is finite.

By the spectral theorem of Hilbert and Schmidt (\Cref{prop:spectral-basics}), as $a\downarrow0$ the projections $\IdentityMap-\SpectralProjectionSymbol_a$ converge to zero at each point of $(\ker\DistanceOperatorSymbol)^\perp$.
Their norms are at most $1$.
A finite $\varepsilon$-net therefore makes the convergence uniform on the compact set $\{d_z:z\in X\}$ (Anselone and Palmer, \cite[Section~2, p.~424, last paragraph]{AnselonePalmer1968}).
Since $\SpectralProjectionSymbol_a$ is self-adjoint, the dual formula for the Hilbert space norm yields
\[
 \begin{aligned}
 \|\DistanceOperatorSymbol(\IdentityMap-\SpectralProjectionSymbol_a)\|_{H\to C(X)}
 &=\sup_{z\in X}\sup_{\|f\|_2\leq1}
       |\langle f,(\IdentityMap-\SpectralProjectionSymbol_a)d_z\rangle_H|\\
 &=\max_{z\in X}\|(\IdentityMap-\SpectralProjectionSymbol_a)d_z\|_2
 \ \longrightarrow\ 0.
 \end{aligned}
\]
For $u_x=\DistanceOperatorSymbol\SpectralProjectionSymbol_a b_x\in\SpectralSubspace{X,d}{a}$ we consequently obtain from \eqref{eq:bump-distance-approximation} that
\[
 \sup_{x\in X}\|d_x-u_x\|_\infty
 \leq r/2+M\|\DistanceOperatorSymbol(\IdentityMap-\SpectralProjectionSymbol_a)\|_{H\to C(X)}.
\]
The right-hand side is less than $\SmallError$ for all sufficiently small $a>0$.
By the Riesz--Schauder theorem (\Cref{thm:compact-spectrum}), the spectrum is at most countable, so we may choose $a$ with $\pm a\notin\spec(\DistanceOperatorSymbol)$.

Suppose that $X$ has at least two points. Full support then gives $(\DistanceOperatorSymbol1)(x)=\int_X d(x,y)\,\IntegrationDifferential\SelectedMeasureSymbol(y)>0$.
Hence $\DistanceOperatorSymbol\ne0$, and the spectral theorem of Hilbert and Schmidt (\Cref{prop:spectral-basics}) yields a nonzero eigenvalue of $\DistanceOperatorSymbol$.
Choosing $a$ less than its absolute value ensures $\dim \SpectralSubspace{X,d}{a}\geq1$.
For a one-point space the Kuratowski map vanishes.
\end{proof}

At a center with at least two points, \Cref{thm:local-models} uses the coefficients of the functions $u_x$ from \Cref{lem:spectral-density} in an orthonormal basis as local coordinates.
Near a one-point center we use the zero map into $\RealNumbers$ with the absolute value norm.
We next estimate the parameter and the dimension of the spectral subspace in this construction.

\begin{corollary}[{Spectral estimates in terms of the masses of balls}]\label{cor:quantitative-spectral-cutoff}
Let $(X,d)$ be a compact metric space with at least two points.
The spectral parameter and dimension in \Cref{lem:spectral-density} satisfy the following bounds.
Set $D=\diam(X,d)$ and $\SelectedMeasureSymbol=\SelectedMeasure{X}{d}$.
For $t>0$, write
\[
 \BallMassBound_t=\LowerBallMass_t\MeasuredIsometryClass{X}{d}{\SelectedMeasureSymbol}
 =\min_{x\in X}\SelectedMeasureSymbol\bigl(\MetricBall(x,t;d)\bigr).
\]
Since $\SelectedMeasureSymbol$ has full support, \Cref{lem:full-support-gdelta} implies $\BallMassBound_t>0$.
Fix $\eta>0$ and set
\[
 r=\frac{\eta}{4},\qquad M_r=\sqrt{\frac{2}{\BallMassBound_{r/4}}},\qquad
 s=\frac{\eta}{4M_r},\qquad
 a_*=\min\left\{\frac{\eta\sqrt{\BallMassBound_s}}{4M_r},\,{\frac D2\sqrt{\frac{\BallMassBound_{D/4}}{2}}}\right\}.
\]
There exists $a$ with $a_*/2<a<a_*$ and $\pm a\notin\spec(\DistanceOperator{X}{d})$.
For each such $a$, the functions $u_x$ in \Cref{lem:spectral-density} with radius $r$ and parameter $a$ satisfy
\[
 {\sup_{x\in X}\|d_x-u_x\|_{\infty,X}<\frac{\eta}{2}},
\]
and the dimension satisfies
\[
 1\leq\dim \SpectralSubspace{X,d}{a}\leq\frac{4D^2}{a_*^2}.
\]
\end{corollary}

\begin{proof}
By the Riesz--Schauder theorem (\Cref{thm:compact-spectrum}), only finitely many absolute values of eigenvalues lie in $[a_*/2,a_*]$, so such an $a$ exists.
Write $\DistanceOperatorSymbol=\DistanceOperator{X}{d}$ and $\SpectralProjectionSymbol_a=\SpectralProjection{X,d}{a}$.
For a Lipschitz function $g$, average $|g(x)|\leq s\operatorname{Lip}(g)+|g(y)|$ over $y\in B(x,s;d)$ and use the Cauchy--Schwarz inequality to obtain $\|g\|_\infty\leq s\operatorname{Lip}(g)+\|g\|_2/\sqrt{\BallMassBound_s}$.
If $\|f\|_2\leq1$ and $g=\DistanceOperatorSymbol(\IdentityMap-\SpectralProjectionSymbol_a)f$, then \eqref{eq:operator-uniform-bound-2} implies $\operatorname{Lip}(g)\leq1$, and the spectral theorem implies $\|g\|_2\leq a$.
Thus the operator norm from $L^2(X,\SelectedMeasureSymbol)$ to $C(X)$ satisfies
\[
 \|\DistanceOperatorSymbol(\IdentityMap-\SpectralProjectionSymbol_a)\|_{L^2\to C(X)}
 \leq s+\frac{a}{\sqrt{\BallMassBound_s}}.
\]
Fix $x\in X$ and put $v=\BumpFunction{X}{d}{x}{r/2}$ and $E=\MetricBall(x,r/4;d)$.
For a Borel set $F\subseteq X$ write $S_F=\int_F v\,\IntegrationDifferential\SelectedMeasureSymbol$ and $Q_F=\int_F v^2\,\IntegrationDifferential\SelectedMeasureSymbol$.
On $E$ we have $1/2\leq v\leq1$, and on $X\setminus E$ we have $0\leq v\leq1/2$.
Thus $Q_{X\setminus E}\leq S_{X\setminus E}/2$ and $Q_E\geq S_E/2$, whence $S_EQ_X\leq Q_ES_X$.
Since $S_X\geq S_E\geq\SelectedMeasureSymbol(E)/2>0$,
\[
 \frac{Q_X}{S_X^2}\leq\frac{Q_E}{S_E^2}.
\]
Put $t=S_E/\SelectedMeasureSymbol(E)$.
On $E$ we have $(v-1/2)(v-1)\leq0$, so $Q_E\leq\SelectedMeasureSymbol(E)(3t-1)/2$.
The inequality $4(3t-1)\leq9t^2$ follows from $(3t-2)^2\geq0$.
Hence the normalized bump $b_x=v/S_X$ in \eqref{eq:normalized-bump-definition} satisfies
\[
 \|b_x\|_2^2=\frac{Q_X}{S_X^2}
 \leq\frac{3t-1}{2\SelectedMeasureSymbol(E)t^2}
 \leq\frac{9}{8\SelectedMeasureSymbol(E)}
 \leq\frac{9}{8\BallMassBound_{r/4}}.
\]
By \eqref{eq:bump-distance-approximation} in the proof of \Cref{lem:spectral-density}, $\|d_x-\DistanceOperatorSymbol b_x\|_\infty\leq r/2$.
Using $\|b_x\|_2\leq3M_r/4$ and $s+a/\sqrt{\BallMassBound_s}<\eta/(2M_r)$, we obtain
\[
 \sup_{x\in X}\|d_x-\DistanceOperatorSymbol\SpectralProjectionSymbol_a b_x\|_\infty
 {\leq\frac r2+\frac{3M_r}{4}\left(s+\frac{a}{\sqrt{\BallMassBound_s}}\right)<\frac{\eta}{8}+\frac{3\eta}{8}=\frac{\eta}{2}}.
\]
Let $A,B$ be the open balls of radius $D/4$ centered at the endpoints of a diameter.
They are disjoint, and their masses $\alpha=\SelectedMeasureSymbol(A)$ and $\beta=\SelectedMeasureSymbol(B)$ are both at least $\BallMassBound_{D/4}$.
Set $f=\Indicator_A/\alpha+\Indicator_B/\beta$.
For $z\in X$ and $(x,y)\in A\times B$ the triangle inequality implies $d(z,x)+d(z,y)\geq d(x,y)\geq D/2$.
Integration against $(\SelectedMeasureSymbol|_A/\alpha)\otimes(\SelectedMeasureSymbol|_B/\beta)$ yields $(\DistanceOperatorSymbol f)(z)\geq D/2$.
Since $\|f\|_2^2=1/\alpha+1/\beta\leq2/\BallMassBound_{D/4}$,
\[
 \|\DistanceOperatorSymbol\|_{L^2\to L^2}
 \geq\frac{\|\DistanceOperatorSymbol f\|_2}{\|f\|_2}
 \geq\frac D2\sqrt{\frac{\alpha\beta}{\alpha+\beta}}
 \geq\frac D2\sqrt{\frac{\BallMassBound_{D/4}}{2}}
 \geq a_*>a.
\]
As the operator is compact and self-adjoint, it has an eigenvalue of absolute value $\|\DistanceOperatorSymbol\|$.
The dimension bound follows from $a^2\dim \SpectralSubspace{X,d}{a}\leq D^2$ and $a>a_*/2$.
\end{proof}

Henceforth fix a radius $r>0$ and a parameter $a>0$.
On $\MetricSpace{Y}{e}$ write $\DistanceOperatorSymbol =\DistanceOperator{Y}{e}$, $\SpectralProjectionSymbol =\SpectralProjection{Y,e}{a}$, and $V=\SpectralSubspace{Y,e}{a}$.
For $y\in Y$ define the approximating function and the \emph{uniform approximation error of the Kuratowski maps} by
\begin{equation}\label{eq:projected-bump-profile}
 \ModelBump{Y}{e}{y}=\NormalizedBump{Y}{e}{y}{r}{\SelectedMeasure{Y}{e}},
 \qquad \ProjectedProfile{Y}{e}{y}=\DistanceOperatorSymbol\SpectralProjectionSymbol\ModelBump{Y}{e}{y}\in V,
\end{equation}
\begin{equation}\label{eq:kuratowski-uniform-error}
 \KuratowskiError{r}{a}\IsometryClass{Y}{e}
   =\max_{y,z\in Y}|e(y,z)-\ProjectedProfile{Y}{e}{y}(z)|.
\end{equation}
The numerator of the bump is continuous in $(y,z)$ and its denominator is positive and continuous in $y$.
Thus $y\mapsto\ModelBump{Y}{e}{y}$ is continuous into $C(Y)$.
Since $\DistanceOperatorSymbol\SpectralProjectionSymbol$ is bounded from $L^2$ to $C(Y)$, the function $(y,z)\mapsto\ProjectedProfile{Y}{e}{y}(z)$ is continuous and the maximum exists.

A surjective isometry $\iota\colon(Y,e)\to(Y_0,e_0)$ transports the selected measure by condition (M2).
The pullback $Uf=f\circ\iota$ is an isometric isomorphism on $L^2$.
The integral formula for the distance kernel yields $U\DistanceOperatorSymbol_0=\DistanceOperatorSymbol U$.
It therefore identifies the eigenspaces and their orthogonal projections, so $U\SpectralProjectionSymbol_0=\SpectralProjectionSymbol U$.
Since $U\ModelBump{Y_0}{e_0}{\iota y}=\ModelBump{Y}{e}{y}$,
\[
 U\ProjectedProfile{Y_0}{e_0}{\iota y}=\ProjectedProfile{Y}{e}{y}.
\]
Both sides are continuous and the selected measure has full support, so the equality holds everywhere.
Thus \eqref{eq:kuratowski-uniform-error} depends only on the isometry class.
\subsection{Continuity of finite spectral subspaces}\label{subsec:spectral-subspace-continuity}
Distance operators for different compact spaces act on different Hilbert spaces.
Extending the distance kernels to a single compact metric space $\MetricSpace{\AmbientCarrier}{\AmbientMetric}$ produces operators on the fixed Banach space $\ContinuousFunctions(\AmbientCarrier ,\ComplexNumbers)$.
Choose one path of integration avoiding the limit spectrum.
By the spectral projection theorem of Anselone and Palmer (\Cref{thm:riesz-convergence}), this path also works at all sufficiently large indices, the corresponding projections converge strongly, and their ranks eventually agree.

By \Cref{lem:parameter-integrals}, the parameter-dependent integrals used below converge uniformly.

We first compare the generalized eigenspaces of the ambient kernel operator with those of the self-adjoint distance operator.

\begin{proposition}\label{lem:ambient-eigenspaces}
For a nonempty compact subset $Y$ of a compact metric space $\MetricSpace{Z}{h}$ put $e=h|_{Y\times Y}$ and $D=\diam\MetricSpace{Z}{h}$.
Let $\DistanceOperator{Y}{e}$ denote the distance operator on the complex Hilbert space $L^2(Y,\SelectedMeasure{Y}{e};\ComplexNumbers)$.
Define the restriction operator $R=\RestrictionOperator_{Y,e}$ and the integral operator $S=\KernelExtensionOperator_{Y,e}$ by
\[
 R\colon\ContinuousFunctions(Z,\ComplexNumbers)\longrightarrow
 \LebesgueSpace^2(Y,\SelectedMeasure{Y}{e};\ComplexNumbers),\qquad Rf=[f|_Y],
\]
\[
 S\colon\LebesgueSpace^2(Y,\SelectedMeasure{Y}{e};\ComplexNumbers)\longrightarrow\ContinuousFunctions(Z,\ComplexNumbers),
 \qquad(Su)(z)=\int_Y h(z,y)u(y)\,\IntegrationDifferential\SelectedMeasure{Y}{e}(y).
\]
Denote their products in opposite orders by
\[
 B=\CarrierAmbientOperator_{Y,e}=SR,\qquad \DistanceOperatorSymbol=\DistanceOperator{Y}{e}=RS.
\]
Then $B$ is compact and
\[
 \spec(B)\setminus\{0\}=\spec(\DistanceOperatorSymbol)\setminus\{0\}\subset[-D,D].
\]
For all $\lambda\in\ComplexNumbers\setminus\{0\}$, $A\in\{B,\DistanceOperatorSymbol\}$, and integers $m\geq1$,
\[
 \ker\bigl((A-\lambda\IdentityMap)^m\bigr)=\ker(A-\lambda\IdentityMap).
\]
Hence the generalized eigenspace at $\lambda$ equals the eigenspace.
Moreover the restriction
\[
 R:\ker(B-\lambda\IdentityMap)\longrightarrow\ker(\DistanceOperatorSymbol-\lambda\IdentityMap)
\]
is a linear isomorphism whose inverse is the restriction of $\lambda^{-1}S$.
\end{proposition}

\begin{proof}
\textbf{Step 1.
Compactness and factorization.}
Since $\SelectedMeasure{\ComparisonCarrier}{\ComparisonMetric}$ has total mass one, for every $\HilbertFunction\in\ContinuousFunctions(\AmbientCarrier,\ComplexNumbers)$ we have
\[
 \norm{\RestrictionOperator_{\ComparisonCarrier,\ComparisonMetric}\HilbertFunction}_2
 \leq\norm{\HilbertFunction}_\infty.
\]
For all $\HilbertTestFunction\in\LebesgueSpace^2(\ComparisonCarrier,\SelectedMeasure{\ComparisonCarrier}{\ComparisonMetric};\ComplexNumbers)$ and $\IntegrationPoint,\IntegrationPoint_0\in\AmbientCarrier$ the Cauchy--Schwarz inequality implies
\[
 \norm{\KernelExtensionOperator_{\ComparisonCarrier,\ComparisonMetric}\HilbertTestFunction}_\infty
 \leq\DiameterBound\norm{\HilbertTestFunction}_2,
\]
and
\[
 |\KernelExtensionOperator_{\ComparisonCarrier,\ComparisonMetric}\HilbertTestFunction(\IntegrationPoint)
  -\KernelExtensionOperator_{\ComparisonCarrier,\ComparisonMetric}\HilbertTestFunction(\IntegrationPoint_0)|
 \leq\AmbientMetric(\IntegrationPoint,\IntegrationPoint_0)
       \norm{\HilbertTestFunction}_2.
\]
Thus $\RestrictionOperator_{\ComparisonCarrier,\ComparisonMetric}$ and $\KernelExtensionOperator_{\ComparisonCarrier,\ComparisonMetric}$ are bounded, and $\KernelExtensionOperator_{\ComparisonCarrier,\ComparisonMetric}$ is compact by the Arzel\`a--Ascoli theorem (\Cref{thm:arzela-ascoli}).
Hence $\CarrierAmbientOperator_{\ComparisonCarrier,\ComparisonMetric}$ is compact.
By definition,
\begin{equation}\label{eq:operator-factorization-1}
 \CarrierAmbientOperator_{\ComparisonCarrier,\ComparisonMetric} =\KernelExtensionOperator_{\ComparisonCarrier,\ComparisonMetric} \RestrictionOperator_{\ComparisonCarrier,\ComparisonMetric} ,
\end{equation}
and
\begin{equation}\label{eq:operator-factorization-2}
 \DistanceOperator{\ComparisonCarrier}{\ComparisonMetric}=\RestrictionOperator_{\ComparisonCarrier,\ComparisonMetric} \KernelExtensionOperator_{\ComparisonCarrier,\ComparisonMetric} ,
\end{equation}
and
\begin{equation}\label{eq:operator-factorization-3}
 \RestrictionOperator_{\ComparisonCarrier,\ComparisonMetric}  \CarrierAmbientOperator_{\ComparisonCarrier,\ComparisonMetric} =\DistanceOperator{\ComparisonCarrier}{\ComparisonMetric}\RestrictionOperator_{\ComparisonCarrier,\ComparisonMetric}.
\end{equation}
The intertwining identity is displayed in Figure~\ref{fig:operator-comparison}.
\begin{figure}[H]
\centering
\begin{tikzcd}[column sep=large,row sep=large]
 \ContinuousFunctions(\AmbientCarrier,\ComplexNumbers)
 \arrow[r,"\CarrierAmbientOperator_{\ComparisonCarrier,\ComparisonMetric}"]
 \arrow[d,"\RestrictionOperator_{\ComparisonCarrier,\ComparisonMetric}"'] &
 \ContinuousFunctions(\AmbientCarrier,\ComplexNumbers)
 \arrow[d,"\RestrictionOperator_{\ComparisonCarrier,\ComparisonMetric}"] \\
 \LebesgueSpace^2(\ComparisonCarrier,\SelectedMeasure{\ComparisonCarrier}{\ComparisonMetric};\ComplexNumbers)
 \arrow[r,"\DistanceOperator{\ComparisonCarrier}{\ComparisonMetric}"']
 \arrow[ur,"\KernelExtensionOperator_{\ComparisonCarrier,\ComparisonMetric}" description] &
 \LebesgueSpace^2(\ComparisonCarrier,\SelectedMeasure{\ComparisonCarrier}{\ComparisonMetric};\ComplexNumbers)
\end{tikzcd}
\caption{Restriction intertwines the ambient operator $\CarrierAmbientOperator_{\ComparisonCarrier,\ComparisonMetric}$ and the distance operator $\DistanceOperator{\ComparisonCarrier}{\ComparisonMetric}$. The factorizations $\CarrierAmbientOperator_{\ComparisonCarrier,\ComparisonMetric}=\KernelExtensionOperator_{\ComparisonCarrier,\ComparisonMetric}\RestrictionOperator_{\ComparisonCarrier,\ComparisonMetric}$ and $\DistanceOperator{\ComparisonCarrier}{\ComparisonMetric}=\RestrictionOperator_{\ComparisonCarrier,\ComparisonMetric}\KernelExtensionOperator_{\ComparisonCarrier,\ComparisonMetric}$ allow us to compare the subspaces corresponding to their nonzero spectral values.}
\label{fig:operator-comparison}
\end{figure}

\textbf{Step 2.
Generalized eigenspaces and restriction.}
The operator $\DistanceOperatorSymbol=\DistanceOperator{Y}{e}$ is self-adjoint on the complex Hilbert space $L^2(Y,\SelectedMeasure{Y}{e};\ComplexNumbers)$ (\Cref{lem:distance-operator}).
It is therefore normal, so \Cref{lem:normal-power-kernel} identifies each generalized eigenspace with the eigenspace, even at $\lambda=0$.
For $B=\CarrierAmbientOperator_{Y,e}$, identities \eqref{eq:operator-factorization-1} and \eqref{eq:operator-factorization-2} imply $B=SR$ and $\DistanceOperatorSymbol=RS$.
If $\lambda\ne0$, \Cref{lem:product-power-kernel} identifies the generalized eigenspace of $B$ with its eigenspace.
On the two eigenspaces the inverse of $R$ is $\lambda^{-1}S$.

\textbf{Step 3.
Equality of the nonzero spectra.}
Put $E=C(Z,\mathbb C)$ and $H=L^2(Y,\SelectedMeasure{Y}{e};\mathbb C)$.
Since $B=SR$ and $\DistanceOperatorSymbol=RS$, if $\zeta\ne0$ and $\zeta \IdentityMap_H-\DistanceOperatorSymbol$ is invertible, then
\[
 (\zeta \IdentityMap_E-B)^{-1}
 =\zeta^{-1}\bigl(\IdentityMap_E+S(\zeta \IdentityMap_H-\DistanceOperatorSymbol)^{-1}R\bigr).
\]
Multiplying the right-hand side by $\zeta\IdentityMap_E-SR$ on either side yields $\IdentityMap_E$.
This is the operator form of Jacobson's lemma.
In a unital algebra, if $1-ab$ has inverse $c$, then $1-ba$ has inverse $1+bca$ (\cite[Remark~1.2.1, p.~6]{Murphy1990}).
Interchanging $R,S$ yields the reverse implication, so $B$ and $\DistanceOperatorSymbol$ have the same nonzero spectrum.
Since $\DistanceOperatorSymbol$ is self-adjoint and $\|\DistanceOperatorSymbol\|\leq\diam(Y,e)\leq D$, their common nonzero spectrum is contained in $[-D,D]$.
\end{proof}

\textbf{An example at eigenvalue zero.} Take $Z=\{y,z\}$, $h(y,z)=1$, and $Y=\{y\}$.
The distance operator on $Y$ vanishes, whereas the ambient operator on $C(Z)$ is $(f(y),f(z))\mapsto(0,f(y))$.
This operator is nonzero but its square vanishes.
Thus at eigenvalue $0$ its eigenspace is a proper subspace of its generalized eigenspace.

We first compare the ambient Riesz projection with the orthogonal projection for a fixed compact carrier.
We then prove convergence of the ambient projections as the carrier varies.

\begin{proposition}\label{lem:projection-restriction}
For a nonempty compact subset $Y$ of a compact metric space $\MetricSpace{Z}{h}$ put $e=h|_{Y\times Y}$.
Fix $a>0$ with $\pm a\notin\spec(\DistanceOperator{Y}{e})$.
Let $\Pi=\AmbientProjection_{Y,e,a}$ be the Riesz projection of $\CarrierAmbientOperator_{Y,e}$ associated with the eigenvalues of absolute value greater than $a$.
This projection preserves real-valued functions, and the restriction operator $R=\RestrictionOperator_{Y,e}$ of \Cref{lem:ambient-eigenspaces} induces a complex-linear isomorphism
\[
 R:\range\Pi\longrightarrow\SpectralSubspace{Y,e}{a}+\ImaginaryUnit\SpectralSubspace{Y,e}{a}.
\]
Write $\SpectralProjectionSymbol$ for the complex-linear extension of the orthogonal projection $\SpectralProjection{Y,e}{a}$.
As operators from $\ContinuousFunctions(Z,\ComplexNumbers)$ to $\LebesgueSpace^2(Y,\SelectedMeasure{Y}{e};\ComplexNumbers)$,
\[
 R\Pi=\SpectralProjectionSymbol R.
\]
For each $f\in C(Z,\ComplexNumbers)$ the function $(\Pi f)|_Y$ is the unique continuous representative of $\SpectralProjectionSymbol(Rf)$.
\end{proposition}

\begin{proof}
By \Cref{lem:ambient-eigenspaces},
\begin{equation}\label{eq:fixed-carrier-nonzero-spectra}
 \spec(\CarrierAmbientOperator_{\ComparisonCarrier,\ComparisonMetric})\setminus\{0\}
 =\spec(\DistanceOperator{\ComparisonCarrier}{\ComparisonMetric})\setminus\{0\}
 \subset\RealNumbers.
\end{equation}
Choose a bounded open set $\SpectralRegion\subset\ComplexNumbers$ consisting of two rectangles symmetric about the real axis with $0\notin\overline{\SpectralRegion}$.
By \eqref{eq:fixed-carrier-nonzero-spectra} and the assumption $\pm\SpectralCutoff\notin\spec(\DistanceOperator{\ComparisonCarrier}{\ComparisonMetric})$ we may choose $\SpectralRegion$ so that, for each operator $\CompactOperator\in\{\CarrierAmbientOperator_{\ComparisonCarrier,\ComparisonMetric},
\DistanceOperator{\ComparisonCarrier}{\ComparisonMetric}\}$, we have
\begin{equation}\label{eq:fixed-carrier-selected-spectrum}
 \SpectralRegion\cap\spec(\CompactOperator)
 =\{\Eigenvalue\in\spec(\CompactOperator)\mid|\Eigenvalue|>\SpectralCutoff\},
\end{equation}
and
\[
 \partial\SpectralRegion\cap\spec(\CompactOperator)=\emptyset.
\]
Give $\SpectralContour=\partial\SpectralRegion$ the positive orientation.
By \Cref{def:riesz-projection} and the holomorphic functional calculus of Riesz and Dunford (\Cref{lem:riesz-projection-properties}),
\[
 \AmbientProjection_{\ComparisonCarrier,\ComparisonMetric,\SpectralCutoff}
 =\RieszProjection{\CarrierAmbientOperator_{\ComparisonCarrier,\ComparisonMetric}}{\SpectralContour}.
\]
By the description of the range by Anselone and Palmer (\Cref{lem:riesz-range-generalized-eigenspace}),
\begin{equation}\label{eq:fixed-carrier-riesz-range}
 \range\AmbientProjection_{\ComparisonCarrier,\ComparisonMetric,\SpectralCutoff}
 =\bigoplus_{\Eigenvalue\in\SpectralRegion\cap\spec(\CarrierAmbientOperator_{\ComparisonCarrier,\ComparisonMetric})}
 \GeneralizedEigenspace{\CarrierAmbientOperator_{\ComparisonCarrier,\ComparisonMetric}}{\Eigenvalue}.
\end{equation}
By the Riesz--Schauder theorem (\Cref{thm:compact-spectrum}), only finitely many eigenvalues have absolute value greater than $a$.
Applying \Cref{lem:ambient-eigenspaces} to each summand in \eqref{eq:fixed-carrier-riesz-range} shows that restriction $R$ maps $\range\Pi$ isomorphically onto $\SpectralSubspace{Y,e}{a}+\ImaginaryUnit\SpectralSubspace{Y,e}{a}$.
With the integral operator $S=\KernelExtensionOperator_{Y,e}$ from the same proposition, write $B=SR$ and $\DistanceOperatorSymbol=RS$.
Then $RB=\DistanceOperatorSymbol R$, so for $\zeta\in\SpectralContour$,
\[
 R(\zeta\IdentityMap-B)^{-1}=(\zeta\IdentityMap-\DistanceOperatorSymbol)^{-1}R.
\]
For the self-adjoint operator $\DistanceOperatorSymbol$, the spectral theorem of Hilbert and Schmidt (\Cref{prop:spectral-basics}) identifies the Riesz projection for this path with the orthogonal projection $\SpectralProjectionSymbol$ onto the sum of the selected eigenspaces (see \cite[Chapter~III, Sections~6.4--6.5]{Kato1995}).
Integrating both sides along $\SpectralContour$ and dividing by $2\pi\ImaginaryUnit$ yields $R\Pi=\SpectralProjectionSymbol R$.

For each $f\in\ContinuousFunctions(Z,\ComplexNumbers)$ the left-hand side of this identity is the $L^2$ class of the continuous function $(\Pi f)|_Y$.
Thus $(\Pi f)|_Y$ is a continuous representative of $\SpectralProjectionSymbol(Rf)$.
Full support of $\SelectedMeasure{Y}{e}$ makes it unique, and at every point of $Y$ we have
\begin{equation}\label{eq:projection-restriction}
 (\AmbientProjection_{\ComparisonCarrier,\ComparisonMetric,\SpectralCutoff}\HilbertFunction)|_\ComparisonCarrier
 =\SpectralProjection{\ComparisonCarrier,\ComparisonMetric}{\SpectralCutoff}
   (\HilbertFunction|_\ComparisonCarrier).
\end{equation}
On the left we first project in the ambient space and then restrict.
On the right we first restrict and then take the $L^2$-orthogonal projection.

Finally, the distance kernels are real and $\SpectralRegion$ is invariant under complex conjugation.
Conjugating the line integrals and reversing the orientation of $\SpectralContour$ therefore yields
\[
 \overline{\AmbientProjection_{\ComparisonCarrier,\ComparisonMetric,\SpectralCutoff}\HilbertFunction}
 =\AmbientProjection_{\ComparisonCarrier,\ComparisonMetric,\SpectralCutoff}\overline{\HilbertFunction}.
\]
Thus $\AmbientProjection_{\ComparisonCarrier,\ComparisonMetric,\SpectralCutoff}$ preserves real functions.
\end{proof}

\begin{proposition}\label{lem:spectral-projections}
Let $X_j,X$ be nonempty compact subsets of a compact metric space $\MetricSpace{Z}{h}$ satisfying $\HausdorffDistance{h}(X_j,X)\to0$, and denote their restricted metrics by $d_j,d$.
Fix $a>0$ with $\pm a\notin\spec(\DistanceOperator{X}{d})$ and put
\[
 A_j=\CarrierAmbientOperator_{X_j,d_j},\quad A=\CarrierAmbientOperator_{X,d},\quad
 \SpectralRank=\dim\SpectralSubspace{X,d}{a}.
\]
There exist $j_0\in\NonnegativeIntegers$ and a common positively oriented path of integration $\SpectralContour$ such that for all $j\geq j_0$ the following holds.
The numbers $\pm a$ lie outside the spectrum of $\DistanceOperator{X_j}{d_j}$, and $\SpectralContour$ lies in the resolvent sets of both $A$ and $A_j$.
For each operator the spectrum enclosed by $\SpectralContour$ consists precisely of the eigenvalues of absolute value greater than $a$.
The Riesz projections
\[
 \Pi=\RieszProjection{A}{\SpectralContour},\qquad
 \Pi_j=\RieszProjection{A_j}{\SpectralContour}\quad(j\geq j_0)
\]
defined by this path have finite rank and preserve real-valued functions.
Each range is the direct sum of the generalized eigenspaces associated with eigenvalues enclosed by $\SpectralContour$.
Moreover,
\[
 \dim\SpectralSubspace{X_j,d_j}{a}=\SpectralRank\quad(j\geq j_0),
\]
and for all $f\in C(Z,\ComplexNumbers)$,
\[
 \norm{\Pi_j f-\Pi f}_\infty\longrightarrow0.
\]
\end{proposition}

\begin{proof}
\textbf{Step 1.
Strong convergence and collective compactness.}
For the inclusions $\iota_j:X_j\hookrightarrow Z$ and $\iota:X\hookrightarrow Z$ put $\bar\mu_j=(\iota_j)_*\SelectedMeasure{X_j}{d_j}$ and $\bar\mu=\iota_*\SelectedMeasure{X}{d}$.
By \Cref{thm:invariant-measures}\ref{item:invariant-measures-continuity} for $\iota_j$ and $\iota$, we have $\bar\mu_\SequenceIndex\to\bar\mu$ weakly.
For all $\SequenceIndex\in\NonnegativeIntegers$, $\HilbertFunction\in\ContinuousFunctions(\AmbientCarrier,\ComplexNumbers)$, and $\IntegrationPoint\in\AmbientCarrier$, by definition
\[
\AmbientOperator _\SequenceIndex \HilbertFunction (\IntegrationPoint )=\int_\AmbientCarrier  \AmbientMetric (\IntegrationPoint ,\SecondVector )\HilbertFunction (\SecondVector )\,\IntegrationDifferential \bar\mu_\SequenceIndex(\SecondVector ),
\]
and
\[
\AmbientOperator \HilbertFunction (\IntegrationPoint )=\int_\AmbientCarrier  \AmbientMetric (\IntegrationPoint ,\SecondVector )\HilbertFunction (\SecondVector )\,\IntegrationDifferential \bar\mu(\SecondVector ).
\]
For each fixed $\HilbertFunction $, the family $\{\AmbientMetric (\IntegrationPoint ,\cdot)\HilbertFunction (\cdot)\mid \IntegrationPoint \in \AmbientCarrier \}$ is compact in $\ContinuousFunctions(\AmbientCarrier ,\ComplexNumbers)$.
By \Cref{lem:parameter-integrals} with parameter $\IntegrationPoint\in\AmbientCarrier$,
\begin{equation}\label{eq:ambient-strong-convergence}
 \norm{\AmbientOperator _\SequenceIndex \HilbertFunction -\AmbientOperator \HilbertFunction }_\infty\to0.
\end{equation}
Thus $\AmbientOperator_\SequenceIndex\to\AmbientOperator$ strongly.
Writing $\DiameterBound =\diam\MetricSpace{\AmbientCarrier }{\AmbientMetric }$, for every $\SequenceIndex\in\NonnegativeIntegers$, $\HilbertFunction\in\ContinuousFunctions(\AmbientCarrier,\ComplexNumbers)$ with $\norm{\HilbertFunction}_\infty\leq1$, and $\IntegrationPoint,\IntegrationPoint_0\in\AmbientCarrier$, we obtain
\begin{equation}\label{eq:ambient-operator-bounds-1}
 \norm{\AmbientOperator _\SequenceIndex \HilbertFunction }_\infty\leq \DiameterBound ,
\end{equation}
and
\begin{equation}\label{eq:ambient-operator-bounds-2}
 |\AmbientOperator _\SequenceIndex \HilbertFunction (\IntegrationPoint )-\AmbientOperator _\SequenceIndex \HilbertFunction (\IntegrationPoint _0)|\leq \AmbientMetric (\IntegrationPoint ,\IntegrationPoint _0).
\end{equation}
Hence $\bigcup_\SequenceIndex  \AmbientOperator _\SequenceIndex \{\HilbertFunction \mid\norm{\HilbertFunction }_\infty\leq1\}$ has compact closure by the Arzel\`a--Ascoli theorem (\Cref{thm:arzela-ascoli}).
Thus $\{\AmbientOperator_\SequenceIndex\}_{\SequenceIndex\in\NonnegativeIntegers}$ is collectively compact.
Since $\AmbientOperator$ is compact by \Cref{lem:ambient-eigenspaces}, the family $\{\AmbientOperator_\SequenceIndex-\AmbientOperator\}_{\SequenceIndex\in\NonnegativeIntegers}$ is also collectively compact.

\textbf{Step 2.
A common boundary and spectral projections.}
Set $\DiameterBound=\diam(Z,h)$ and choose $\SpectralOuterBound>\max\{a,\DiameterBound\}+1$.
Put
\[
 \SpectralRegion=\{\zeta\in\ComplexNumbers\mid
        a<|\RealPart\zeta|<\SpectralOuterBound,\ |\ImaginaryPart\zeta|<1\},
\]
and let $\SpectralContour=\partial\SpectralRegion$ be the sum of the positively oriented boundaries of the two rectangles.
By \Cref{lem:ambient-eigenspaces}, the spectra of $A$ and of every $A_j$ lie in $[-\DiameterBound,\DiameterBound]$.
The only possible intersections of $\SpectralContour$ with this interval are $\pm a$, and by hypothesis and the same proposition, $\pm a\notin\spec(A)$.
Thus $\SpectralContour$ lies in the resolvent set of $A$.

Strong convergence and collective compactness of the difference operators from Step~1 allow us to use the spectral projection theorem of Anselone and Palmer (\Cref{thm:riesz-convergence}) on $\SpectralContour$.
Hence for all sufficiently large $j$, the path $\SpectralContour$ also lies in the resolvent set of $A_j$, the projections defined by this boundary converge strongly, and their ranks agree.
Since $\pm a$ lie on the boundary, they belong neither to $\spec(A_j)$ nor, by \Cref{lem:ambient-eigenspaces}, to $\spec(\DistanceOperator{X_j}{d_j})$.
For a real spectral value the interiors of the rectangles are characterized by $|\lambda|>a$ (\Cref{fig:spectral-contour}).
\begin{figure}[H]
\centering
\begin{tikzpicture}[>=Stealth,x=1cm,y=1cm,font=\small]
\fill[black!6] (-4,-.8) rectangle (-1,.8);
\fill[black!6] (1,-.8) rectangle (4,.8);
\draw[->] (-4.5,0)--(4.7,0) node[right] {$\RealPart\ResolventParameter$};
\draw[->] (0,-1.05)--(0,1.3) node[above] {$\ImaginaryPart\ResolventParameter$};
\draw[thick] (-4,-.8) rectangle (-1,.8);
\draw[thick] (1,-.8) rectangle (4,.8);
\draw[->,thick] (-3.2,-.8)--(-2.5,-.8);
\draw[->,thick] (2,-.8)--(2.7,-.8);
\foreach \x in {-3.5,-2.8,-1.7,1.6,2.5,3.4} \fill (\x,0) circle (1.5pt);
\foreach \x in {-.55,-.3,.25,.5} \fill[black!45] (\x,0) circle (1.2pt);
\node[below left] at (0,0) {$0$};
\node[below] at (-4,-.85) {$-\SpectralOuterBound$};
\node[below] at (-1,-.85) {$-\SpectralCutoff$};
\node[below] at (1,-.85) {$\SpectralCutoff$};
\node[below] at (4,-.85) {$\SpectralOuterBound$};
\node[above] at (-2.5,.85) {$\SpectralContour$};
\node[above] at (2.5,.85) {$\SpectralContour$};
\end{tikzpicture}
\caption{The positively oriented boundary $\SpectralContour=\partial\SpectralRegion$ encloses the positive and negative spectral values retained by the cutoff. The dots are schematic. The gaps at $\pm\SpectralCutoff$ allow the same contour to be used for $\AmbientOperator_\SequenceIndex$ for all sufficiently large $\SequenceIndex$.}
\label{fig:spectral-contour}
\end{figure}

Denote the Riesz projections of $A$ and $A_j$ defined by this common path $\SpectralContour$ by $\AmbientProjection$ and $\AmbientProjection_j$, respectively.

\textbf{Step 3.
Identification with the real spectral subspaces.}
The set $\overline{\SpectralRegion}$ omits zero.
By the Riesz--Schauder theorem (\Cref{thm:compact-spectrum}) and the range formula of Anselone and Palmer (\Cref{lem:riesz-range-generalized-eigenspace}), the range of each projection is a direct sum of finite-dimensional generalized eigenspaces for finitely many nonzero eigenvalues.
By \Cref{lem:ambient-eigenspaces} these are ordinary eigenspaces.
The rank equality and strong convergence from Step~2 read, for large $j$,
\begin{equation}\label{eq:spectral-projection-convergence-1}
 \dim\range \AmbientProjection _\SequenceIndex =\dim\range \AmbientProjection,
\end{equation}
and for each $\HilbertFunction\in\ContinuousFunctions(\AmbientCarrier,\ComplexNumbers)$ the following convergence is uniform.
\begin{equation}\label{eq:spectral-projection-convergence-2}
 \norm{\AmbientProjection _\SequenceIndex \HilbertFunction -\AmbientProjection \HilbertFunction }_\infty\to0.
\end{equation}
By \Cref{lem:projection-restriction} for $(X_j,d_j)$ and $(X,d)$, the ranges of the Riesz projections are isomorphic to the complexifications of the corresponding real spectral subspaces.
Thus \eqref{eq:spectral-projection-convergence-1} yields the dimension identity
\begin{equation}\label{eq:spectral-dimension-identification}
 \begin{aligned}
 \dim_\RealNumbers\SpectralSubspace{\BaseCarrier_\SequenceIndex,\MetricSymbol_\SequenceIndex}{\SpectralCutoff}
 &=\dim_\ComplexNumbers\range\AmbientProjection_\SequenceIndex
 =\dim_\ComplexNumbers\range\AmbientProjection\\
 &=\dim_\RealNumbers\SpectralSubspace{\BaseCarrier,\MetricSymbol}{\SpectralCutoff}
 =\SpectralRank.
 \end{aligned}
\end{equation}
Finally, applying \Cref{lem:projection-restriction} to $\BaseCarrier$ and to $\BaseCarrier_\SequenceIndex$ for all sufficiently large $\SequenceIndex$ shows that $\AmbientProjection$ and $\AmbientProjection_\SequenceIndex$ preserve real functions.
\end{proof}

\textbf{An example with two-point spaces.} In the common space $Z=[0,1]$ set $X=\{0,1\}$ and $X_j=\{t_j,1\}$, where $t_j=2^{-j-1}$ for $j\geq0$.
On every two-point space the transposition is an isometry, so the invariant probability assigns mass $1/2$ to each point.
The spaces $X_j$ converge to $X$ in Hausdorff distance and the measures converge weakly on $Z$.
The corresponding operators on $C(Z)$ are
\[
 (A_j f)(z)=\tfrac12|z-t_j|f(t_j)+\tfrac12|z-1|f(1),\qquad
 (Af)(z)=\tfrac12z f(0)+\tfrac12|z-1|f(1).
\]
For fixed $f$ we have $\|(A_j-A)f\|_\infty\leq(t_j\|f\|_\infty+|f(t_j)-f(0)|)/2\to0$.
On the other hand let $f_j(0)=-1$ and $f_j(t_j)=1$, interpolate linearly on $[0,t_j]$, and set the function equal to $1$ on $[t_j,1]$.
Then $\|f_j\|_\infty=1$ and
\[
 \|A_j-A\|\geq |((A_j-A)f_j)(1)|=1-\tfrac12t_j\longrightarrow1.
\]
Every $A_j$ and the limit $A$ take the unit ball into a uniformly bounded family of $1$-Lipschitz functions.
By Arzel\`a--Ascoli the differences are collectively compact.
Thus the spectral projection theorem of Anselone and Palmer (\Cref{thm:riesz-convergence}) applies.
Kato's perturbation theorem under generalized convergence of bounded operators (\cite[Chapter~IV, Theorems~2.23 and~3.16]{Kato1995}) would require operator-norm convergence, which fails here.

In the same example take $a=1/8$.
The eigenvalues on $X_j$ are $\pm(1-t_j)/2$, so the spectral projections on $X_j$ and $X$ are both the identity.
By \eqref{eq:projection-restriction}, $(\Pi_j f)(t_j)=f(t_j)$.
At the limit, $\range A$ is spanned by $z\mapsto z$ and $z\mapsto 1-z$, and $\range\Pi$ is a two-dimensional subspace of this range.
They are therefore equal, so $\Pi f$ linearly interpolates the endpoint values of $f$.
Define a piecewise linear function by $g_j(0)=1$, $g_j(t_j)=-1$, and $g_j(1)=0$.
Then $\|g_j\|_\infty=1$, but
\[
 \|\Pi_j-\Pi\|\geq|(\Pi_jg_j)(t_j)-(\Pi g_j)(t_j)|=2-t_j.
\]
Thus strong convergence of the projections need not imply operator-norm convergence.

\begin{theorem}\label{lem:spectral-continuity}
Let $\MetricSpace{\AmbientCarrier}{\AmbientMetric}$ be a compact metric space and suppose $\BaseCarrier_\SequenceIndex,\BaseCarrier\in\CompactHyperspace(\AmbientCarrier)$ satisfy $\HausdorffDistance{\AmbientMetric}(\BaseCarrier_\SequenceIndex,\BaseCarrier)\to0$, with respective restricted metrics $\MetricSymbol_\SequenceIndex,\MetricSymbol$ on $\BaseCarrier_\SequenceIndex$ and $\BaseCarrier$.
Fix $\SpectralCutoff>0$ with $\pm\SpectralCutoff\notin\spec(\DistanceOperator{\BaseCarrier}{\MetricSymbol})$ and put $\SpectralRank=\dim\SpectralSubspace{\BaseCarrier,\MetricSymbol}{\SpectralCutoff}$.
There exists $\SequenceIndex_0\in\NonnegativeIntegers$ such that for all $\SequenceIndex\geq\SequenceIndex_0$ the following two assertions hold.
\begin{equation}\label{eq:stable-spectral-endpoints}
 \SpectralCutoff,-\SpectralCutoff
 \notin\spec(\DistanceOperator{\BaseCarrier_\SequenceIndex}{\MetricSymbol_\SequenceIndex}),
\end{equation}
and
\begin{equation}\label{eq:stable-spectral-rank}
 \dim\SpectralSubspace{\BaseCarrier_\SequenceIndex,\MetricSymbol_\SequenceIndex}{\SpectralCutoff}
 =\SpectralRank.
\end{equation}
In particular, if $\SpectralRank=0$ then $\SpectralSubspace{\BaseCarrier_\SequenceIndex,\MetricSymbol_\SequenceIndex}{\SpectralCutoff}=\{0\}$ for all $\SequenceIndex\geq\SequenceIndex_0$.
Now suppose $\SpectralRank>0$ and fix any real orthonormal basis $\{\Eigenfunction_\BasisIndex\}_{0\leq\BasisIndex<\SpectralRank}$ of $\SpectralSubspace{\BaseCarrier,\MetricSymbol}{\SpectralCutoff}$.
Increasing $\SequenceIndex_0$ if necessary for this basis, we can find functions
\[
 \begin{aligned}
 &\BasisExtension_\BasisIndex\colon\AmbientCarrier\to\RealNumbers
     &&(0\leq\BasisIndex<\SpectralRank),\\
 &\Eigenfunction_{\SequenceIndex,\BasisIndex}\colon\BaseCarrier_\SequenceIndex\to\RealNumbers,\qquad
   \BasisExtension_{\SequenceIndex,\BasisIndex}\colon\AmbientCarrier\to\RealNumbers
     &&(\SequenceIndex\geq\SequenceIndex_0,\ 0\leq\BasisIndex<\SpectralRank)
 \end{aligned}
\]
satisfying the following four conditions.
\begin{enumerate}[label=\textup{(B\arabic*)},ref=\textup{(B\arabic*)}]
\item\label{item:spectral-basis-orthonormality}
For every $\SequenceIndex\geq\SequenceIndex_0$, the family $\{\Eigenfunction_{\SequenceIndex,\BasisIndex}\}_{0\leq\BasisIndex<\SpectralRank}$ is a real orthonormal basis of $\SpectralSubspace{\BaseCarrier_\SequenceIndex,\MetricSymbol_\SequenceIndex}{\SpectralCutoff}$.
\item\label{item:spectral-basis-extensions}
For every $\SequenceIndex\geq\SequenceIndex_0$ and $0\leq\BasisIndex<\SpectralRank$, the functions $\BasisExtension_{\SequenceIndex,\BasisIndex},\BasisExtension_\BasisIndex
 \colon\AmbientCarrier\to\RealNumbers$ are continuous and satisfy
\[
 \BasisExtension_{\SequenceIndex,\BasisIndex}|_{\BaseCarrier_\SequenceIndex}
 =\Eigenfunction_{\SequenceIndex,\BasisIndex},
\]
and
\[
 \BasisExtension_\BasisIndex|_\BaseCarrier=\Eigenfunction_\BasisIndex.
\]
\item\label{item:spectral-basis-uniform-convergence}
These extensions also satisfy
\begin{equation}\label{eq:aligned-basis-extensions}
 \max_{0\leq\BasisIndex<\SpectralRank}
 \norm{\BasisExtension_{\SequenceIndex,\BasisIndex}-\BasisExtension_\BasisIndex}_\infty
 \longrightarrow0.
\end{equation}
\item\label{item:spectral-basis-correspondence-convergence}
Let $\{\Correspondence_\SequenceIndex\subset\BaseCarrier_\SequenceIndex\times\BaseCarrier\}_{\SequenceIndex\geq\SequenceIndex_0}$ be a sequence of correspondences such that
\[
 \sup_{(\BasePoint,\ComparisonPoint)\in\Correspondence_\SequenceIndex}
 \AmbientMetric(\BasePoint,\ComparisonPoint)\longrightarrow0.
\]
Then the bases $\{\Eigenfunction_{\SequenceIndex,\BasisIndex}\}_{0\leq\BasisIndex<\SpectralRank}$ and $\{\Eigenfunction_\BasisIndex\}_{0\leq\BasisIndex<\SpectralRank}$ satisfy
\begin{equation}\label{eq:aligned-bases}
 \max_{0\leq\BasisIndex<\SpectralRank}
 \sup_{(\BasePoint,\ComparisonPoint)\in\Correspondence_\SequenceIndex}
 |\Eigenfunction_{\SequenceIndex,\BasisIndex}(\BasePoint)
       -\Eigenfunction_\BasisIndex(\ComparisonPoint)|\longrightarrow0.
\end{equation}
\end{enumerate}
\end{theorem}

Figure~\ref{fig:aligned-spectral-bases} illustrates the choice of complete orthonormal systems in the spectral subspaces aligned with the prescribed limit basis.
\begin{figure}[H]
\centering
\begin{tikzpicture}[font=\small,>=Stealth]
  \foreach \xx/\turnangle in {0/26,3.15/11,7.25/0} {
    \begin{scope}[shift={(\xx,0)}]
      \draw[fill=black!4,draw=black!58] (0,0) circle[radius=.89cm];
      \begin{scope}[rotate=\turnangle]
        \draw[->,line width=.85pt] (0,0)--(.89,0);
        \draw[->,densely dashed,line width=.85pt] (0,0)--(0,.89);
        \draw[thin] (.16,0)--(.16,.16)--(0,.16);
      \end{scope}
      \fill (0,0) circle[radius=1.25pt];
    \end{scope}
  }
  \node at (0,1.28) {$\SpectralSubspace{X_j,d_j}{\SpectralCutoff}$};
  \node at (3.15,1.28) {$\SpectralSubspace{X_{j+1},d_{j+1}}{\SpectralCutoff}$};
  \node at (7.25,1.28) {$\SpectralSubspace{X,d}{\SpectralCutoff}$};
  \node at (5.20,0) {$\cdots$};
  \node at (0,-1.23) {$(\Eigenfunction_{j,0},\Eigenfunction_{j,1})$};
  \node at (3.15,-1.23) {$(\Eigenfunction_{j+1,0},\Eigenfunction_{j+1,1})$};
  \node at (7.25,-1.23) {$(\Eigenfunction_0,\Eigenfunction_1)$};
  \node at (3.55,-2.07)
    {$\displaystyle
       (\BasisExtension_{j,0},\BasisExtension_{j,1})
       \xrightarrow[\ j\to\infty\ ]{\ \ContinuousFunctions(Z)^2\ }
       (\BasisExtension_0,\BasisExtension_1)$};
\end{tikzpicture}
\caption{A two-dimensional spectral subspace ($n=2$) in Theorem~\ref{lem:spectral-continuity}. The circles are the $L^2$ unit circles of the spectral subspaces. Fix the basis on the right. For all sufficiently large $j$, one can choose an orthonormal basis on the left in which the solid vectors correspond and the dashed vectors correspond. The lower diagram shows uniform convergence of the respective extensions to the common space $Z$.}
\label{fig:aligned-spectral-bases}
\end{figure}

\begin{proof}
Let $\iota_\SequenceIndex\colon\BaseCarrier_\SequenceIndex\hookrightarrow\AmbientCarrier$ and $\iota\colon\BaseCarrier\hookrightarrow\AmbientCarrier$ be the inclusions and put $\bar\mu_\SequenceIndex=(\iota_\SequenceIndex)_*\SelectedMeasure{\BaseCarrier_\SequenceIndex}{\MetricSymbol_\SequenceIndex}$ and $\bar\mu=\iota_*\SelectedMeasure{\BaseCarrier}{\MetricSymbol}$.
By \Cref{thm:invariant-measures}\ref{item:invariant-measures-continuity}, $\bar\mu_\SequenceIndex$ converges weakly to $\bar\mu$ in $\ProbabilityMeasures(\AmbientCarrier)$.
By \Cref{lem:spectral-projections} we obtain \eqref{eq:stable-spectral-endpoints} and \eqref{eq:stable-spectral-rank}.
Denote the projections constructed there by $\AmbientProjection_\SequenceIndex,\AmbientProjection$.
If $\SpectralRank=0$, the spectral subspaces eventually vanish and no basis or extension data are needed.
Assume $\SpectralRank>0$.
To construct the real orthonormal basis $\{\Eigenfunction_{\SequenceIndex,\BasisIndex}\}_{0\leq\BasisIndex<\SpectralRank}$ of $\SpectralSubspace{\BaseCarrier_\SequenceIndex,\MetricSymbol_\SequenceIndex}{\SpectralCutoff}$ in \ref{item:spectral-basis-orthonormality}, first project continuous extensions of the prescribed basis $\{\Eigenfunction_\BasisIndex\}_{0\leq\BasisIndex<\SpectralRank}$.
The resulting functions need not be orthonormal for the varying measures $\SelectedMeasure{\BaseCarrier_\SequenceIndex}{\MetricSymbol_\SequenceIndex}$.
We normalize them by their Gram matrices while preserving the limit.
We then verify \ref{item:spectral-basis-extensions} and \ref{item:spectral-basis-uniform-convergence} for the normalized extensions and deduce \ref{item:spectral-basis-correspondence-convergence} by uniform continuity.

\ProofStep{Projected extensions and convergence of the Gram matrices.}\label{step:spectral-basis-projection}
For every integer $0\leq\BasisIndex<\SpectralRank$, use the Tietze--Urysohn theorem (\Cref{thm:tietze}) to choose a continuous real extension $\widetilde\Eigenfunction_\BasisIndex\colon\AmbientCarrier\to\RealNumbers$ of the prescribed function $\Eigenfunction_\BasisIndex\colon\BaseCarrier\to\RealNumbers$.
For every sufficiently large $\SequenceIndex\in\NonnegativeIntegers$ and every integer $0\leq\BasisIndex<\SpectralRank$, put
\[
 \ProjectedBasisFunction _{\SequenceIndex ,\BasisIndex }=\AmbientProjection _\SequenceIndex \widetilde \Eigenfunction _\BasisIndex.
\]
For every integer $0\leq\BasisIndex<\SpectralRank$, put
\[
 \ProjectedBasisFunction _\BasisIndex =\AmbientProjection \widetilde \Eigenfunction _\BasisIndex.
\]
Applying \eqref{eq:spectral-projection-convergence-2} to each fixed function $\widetilde\Eigenfunction_\BasisIndex$, we obtain $\norm{\ProjectedBasisFunction _{\SequenceIndex ,\BasisIndex }-\ProjectedBasisFunction _\BasisIndex }_\infty\to0$ for every $0\leq\BasisIndex<\SpectralRank$.
Moreover, \eqref{eq:projection-restriction} implies
\[
 \RestrictionOperator_{\BaseCarrier,\MetricSymbol}\ProjectedBasisFunction_\BasisIndex
 =\SpectralProjection{\BaseCarrier,\MetricSymbol}{\SpectralCutoff}
   \RestrictionOperator_{\BaseCarrier,\MetricSymbol}\widetilde\Eigenfunction_\BasisIndex
 =\SpectralProjection{\BaseCarrier,\MetricSymbol}{\SpectralCutoff}\Eigenfunction_\BasisIndex
 =\Eigenfunction_\BasisIndex.
\]
Thus $\ProjectedBasisFunction_\BasisIndex|_\BaseCarrier=\Eigenfunction_\BasisIndex$.
Equation \eqref{eq:projection-restriction} on $\BaseCarrier_\SequenceIndex$ shows that
\[
 \ProjectedBasisFunction_{\SequenceIndex,\BasisIndex}|_{\BaseCarrier_\SequenceIndex}
 =\SpectralProjection{\BaseCarrier_\SequenceIndex,\MetricSymbol_\SequenceIndex}{\SpectralCutoff}
   (\widetilde\Eigenfunction_\BasisIndex|_{\BaseCarrier_\SequenceIndex})
 \in\SpectralSubspace{\BaseCarrier_\SequenceIndex,\MetricSymbol_\SequenceIndex}{\SpectralCutoff}.
\]
Arranging these functions in a row, we obtain in the uniform norm of $C(Z,\RealNumbers^{1\times n})$
\begin{equation}\label{eq:projected-row-convergence}
 \ProjectedBasisRow_j=(u_{j,0},\ldots,u_{j,n-1})
 \longrightarrow\ProjectedBasisRow=(u_0,\ldots,u_{n-1}).
\end{equation}
For each $0\leq k,l<n$ orthonormality of the prescribed basis implies
\begin{equation}\label{eq:gram-weak-integration-error}
 \GramMatrix(\bar\mu,\ProjectedBasisRow)_{kl}
 =\int_Z u_k u_l\,\IntegrationDifferential\bar\mu
 =\int_X \Eigenfunction_k\Eigenfunction_l\,\IntegrationDifferential\SelectedMeasure{X}{d}
 =\KroneckerSymbol_{kl}.
\end{equation}
Thus \Cref{lem:continuous-gram-normalization} implies $\SampleGramMatrix_j:=\GramMatrix(\bar\mu_j,\ProjectedBasisRow_j)\to \IdentityMatrix_n$.
In particular all $\SampleGramMatrix_j$ are positive definite from some index $j_0$ onward.

\ProofStep{Normalization by inverse square roots.}\label{step:spectral-basis-normalization}
For $j\geq j_0$ define the normalized extensions and their restrictions by
\begin{equation}\label{eq:normalized-basis-row}
 (\widetilde\Eigenfunction_{j,0},\ldots,\widetilde\Eigenfunction_{j,n-1})
  =\GramNormalization(\bar\mu_j,\ProjectedBasisRow_j),\qquad
 \Eigenfunction_{j,k}=\widetilde\Eigenfunction_{j,k}|_{X_j}.
\end{equation}
Their Gram matrix is $\IdentityMatrix_n$.
Each $\Eigenfunction_{j,k}$ is a linear combination of the restrictions $u_{j,l}|_{X_j}$ and hence belongs to $\SpectralSubspace{X_j,d_j}{\SpectralCutoff}$.
This space has dimension $n$, so these functions form the real orthonormal basis in \ref{item:spectral-basis-orthonormality}.
This normalization also takes the isometric factor $F_j\SampleGramMatrix_j^{-1/2}$ of the operator $F_j$ defined by $F_jc=\sum_{k=0}^{n-1}c_k u_{j,k}|_{X_j}$.
After choosing an orthonormal basis in the range of $F_j$, this is Autonne's polar decomposition (\cite{Autonne1902Polar}, see also \cite[Theorem~8.1]{Higham2008}).

\ProofStep{Continuous extensions and convergence.}\label{step:spectral-basis-extensions}
Put $\BasisExtension_{j,k}=\widetilde\Eigenfunction_{j,k}$ and $\BasisExtension_k=u_k$.
These are the continuous real-valued extensions in \ref{item:spectral-basis-extensions}.
Since $\GramNormalization(\bar\mu,\ProjectedBasisRow)=\ProjectedBasisRow$, by \eqref{eq:projected-row-convergence} and continuity in \Cref{lem:continuous-gram-normalization},
\[
 \max_{0\leq k<n}\norm{\BasisExtension_{j,k}-\BasisExtension_k}_\infty\longrightarrow0.
\]
This proves \ref{item:spectral-basis-uniform-convergence}.
Finally, for the correspondences $R_j$ in \ref{item:spectral-basis-correspondence-convergence} put $\CorrespondenceDisplacement_j=\sup_{(x,y)\in R_j}h(x,y)$.
Let $\CommonModulus$ be a common modulus of continuity for the finitely many functions $u_k$, namely $\CommonModulus(t)=\max_{k<n}\sup_{h(z,w)\leq t}|u_k(z)-u_k(w)|$.
Then
\begin{equation}\label{eq:basis-alignment-error}
 \max_{0\leq k<n}\sup_{(x,y)\in R_j}|\Eigenfunction_{j,k}(x)-\Eigenfunction_k(y)|
 \leq \max_{0\leq k<n}\norm{\BasisExtension_{j,k}-u_k}_\infty+\CommonModulus(\CorrespondenceDisplacement_j)
 \longrightarrow0.
\end{equation}
Indeed, uniform continuity on $Z$ implies $\CommonModulus(t)\to0$ as $t\downarrow0$, and $\CorrespondenceDisplacement_j\to0$ by hypothesis.
\end{proof}
\subsection{Local maps to normed spaces}\label{subsec:spectral-local-models}
We now turn the spectral subspace into coordinates in $\RealNumbers^\ModelDimension $.
A choice of orthonormal basis determines the coordinates, but the resulting pseudometric is independent of this choice.

On column vectors in $\RealNumbers^\ModelDimension $ we use the orthogonal group
\[
 \OrthogonalGroup(\ModelDimension )=\{\OrthogonalChange \in\RealNumbers^{\ModelDimension \times \ModelDimension }\mid \OrthogonalChange ^{\TransposeSymbol }\OrthogonalChange =\IdentityMatrix _\ModelDimension \}.
\]
Thus $\OrthogonalChange $ preserves the Euclidean norm and $\OrthogonalChange ^{-1}=\OrthogonalChange ^{\TransposeSymbol }$.

We consider the equivalence class of a coordinate map and its norm, so that a change of orthonormal basis does not change the model.

\begin{definition}\label{def:normed-coordinate-class}
Let $\ComparisonCarrier$ be a nonempty compact space and let $\ModelDimension\geq1$.
Write
\[
 \FiniteNormSpace{\ModelDimension}
 =\{\NormSymbol\colon\RealNumbers^\ModelDimension\to[0,\infty)
       \mid\NormSymbol\text{ is a norm}\}.
\]
The group $\OrthogonalGroup(\ModelDimension)$ acts on $\ContinuousFunctions(\ComparisonCarrier,\RealNumbers^\ModelDimension)
 \times\FiniteNormSpace{\ModelDimension}$ by
\begin{equation}\label{eq:frame-action}
 \OrthogonalChange\cdot(\FeatureCoordinates,\NormSymbol)
 =(\OrthogonalChange\circ\FeatureCoordinates,
   \NormSymbol\circ\OrthogonalChange^{-1}).
\end{equation}
We call an element of the quotient set
\begin{equation}\label{eq:normed-coordinate-class-space}
 \CoordinateClassSpace{\ModelDimension}{\ComparisonCarrier}
 =\bigl(\ContinuousFunctions(\ComparisonCarrier,\RealNumbers^\ModelDimension)
       \times\FiniteNormSpace{\ModelDimension}\bigr)
       /\OrthogonalGroup(\ModelDimension)
\end{equation}
an \emph{$\ModelDimension$-dimensional normed coordinate class on $\ComparisonCarrier$}.
The class of a pair $(\FeatureCoordinates,\NormSymbol)$ is denoted by $\CoordinateClassOf{\ModelDimension}{\FeatureCoordinates}{\NormSymbol}$.
\end{definition}

A coordinate class ${\CoordinateClassElement}=\CoordinateClassOf{\ModelDimension}{\FeatureCoordinates}{\NormSymbol}\in\CoordinateClassSpace{\ModelDimension}{Y}$ determines a continuous pseudometric $\ClassPseudometric{\CoordinateClassElement}$ on $Y$ and the maximum distance from the origin to the coordinate image, its \emph{coordinate radius} $\CoordinateRadius(\CoordinateClassElement)$.
They are defined by
\[
 \ClassPseudometric{\CoordinateClassElement}(y,z)=\NormSymbol\bigl(\FeatureCoordinates(y)-\FeatureCoordinates(z)\bigr),
 \qquad \CoordinateRadius(\CoordinateClassElement)=\max_{y\in Y}\NormSymbol\bigl(\FeatureCoordinates(y)\bigr).
\]
Both are independent of the representative since $(\NormSymbol\circ \OrthogonalChange^{-1})(\OrthogonalChange\xi)=\NormSymbol(\xi)$ (Figure~\ref{fig:orthogonal-coordinate-change}).

\begin{figure}[H]
\centering
\begin{tikzpicture}[>=Stealth,font=\small]
\node at (0,1.65)
 {$(\FeatureCoordinates,\NormSymbol)$};
\node at (5.8,1.65)
 {$(\OrthogonalChange\circ\FeatureCoordinates,
     \NormSymbol\circ\OrthogonalChange^{-1})$};
\begin{scope}
 \filldraw[fill=black!4] (0,0) ellipse (1.45 and 0.88);
 \draw[->,black!35] (-1.8,0) -- (1.8,0);
 \draw[->,black!35] (0,-1.15) -- (0,1.2);
 \foreach \x/\y in {-0.8/-0.3,-0.55/0.42,0.12/-0.55,0.48/0.3,0.92/-0.15}
  \fill (\x,\y) circle (1.4pt);
 \draw[black!55] (0.48,0.3) -- (0.85,1.03);
 \node at (0.85,1.25)
  {$\FeatureCoordinates(\ComparisonPoint)$};
 \draw[black!55] (-0.8,-0.3) -- (-0.95,-0.9);
 \node at (-0.95,-1.12)
  {$\FeatureCoordinates(\IntegrationPoint)$};
\end{scope}
\begin{scope}[xshift=5.8cm]
 \begin{scope}[rotate=30]
  \filldraw[fill=black!4] (0,0) ellipse (1.45 and 0.88);
 \end{scope}
 \draw[->,black!35] (-1.8,0) -- (1.8,0);
 \draw[->,black!35] (0,-1.15) -- (0,1.2);
 \begin{scope}[rotate=30]
  \foreach \x/\y in {-0.8/-0.3,-0.55/0.42,0.12/-0.55,0.48/0.3,0.92/-0.15}
   \fill (\x,\y) circle (1.4pt);
  \coordinate (rotatedy) at (0.48,0.3);
  \coordinate (rotatedz) at (-0.8,-0.3);
 \end{scope}
 \draw[black!55] (rotatedy) -- (0.4,1.03);
 \node at (0.4,1.25)
  {$\OrthogonalChange\FeatureCoordinates(\ComparisonPoint)$};
 \draw[black!55] (rotatedz) -- (-1,-0.92);
 \node at (-1,-1.14)
  {$\OrthogonalChange\FeatureCoordinates(\IntegrationPoint)$};
\end{scope}
\draw[->] (2.05,0.15) -- node[above]
 {$\OrthogonalChange\in\OrthogonalGroup(2)$} (3.75,0.15);
\node at (0,-1.55)
 {$\{\FirstVector\mid\NormSymbol(\FirstVector)\leq1\}$};
\node at (5.8,-1.55)
 {$\{\FirstVector\mid
     (\NormSymbol\circ\OrthogonalChange^{-1})(\FirstVector)\leq1\}$};
\node at (2.9,-2.35)
 {$(\NormSymbol\circ\OrthogonalChange^{-1})
     \bigl(\OrthogonalChange\FeatureCoordinates(\ComparisonPoint)
           -\OrthogonalChange\FeatureCoordinates(\IntegrationPoint)\bigr)
   =\NormSymbol\bigl(\FeatureCoordinates(\ComparisonPoint)
                     -\FeatureCoordinates(\IntegrationPoint)\bigr)$};
\end{tikzpicture}
\caption{An orthogonal change of coordinates acts on both the coordinate image and the unit ball of the norm. The dots represent a finite coordinate image in two dimensions, and the shaded regions are the two unit balls. Rotating the image and the ball together preserves the induced pseudometric.}
\label{fig:orthogonal-coordinate-change}
\end{figure}

\begin{lemma}\label{lem:coordinate-pseudometric-estimate}
Let $\BaseCarrier,\ComparisonCarrier$ be nonempty compact spaces, let $\Correspondence\subset\BaseCarrier\times\ComparisonCarrier$ be a correspondence, and let $\ModelDimension\geq1$.
Let $\FeatureCoordinates_0\colon\BaseCarrier\to\RealNumbers^\ModelDimension$ and $\FeatureCoordinates_1\colon\ComparisonCarrier\to\RealNumbers^\ModelDimension$ be continuous, and let $\NormSymbol_0,\NormSymbol_1$ be norms on $\RealNumbers^\ModelDimension$.
Denote the induced pseudometrics by $\Pseudometric_0$ and $\Pseudometric_1$, respectively.
Put
\[
 \CoefficientBound=\sup_{\BasePoint\in\BaseCarrier}
       \norm{\FeatureCoordinates_0(\BasePoint)}_2,
\]
\[
 {\CoordinateComparisonError}
 =\sup_{(\BasePoint,\ComparisonPoint)\in\Correspondence}
       \norm{\FeatureCoordinates_0(\BasePoint)-\FeatureCoordinates_1(\ComparisonPoint)}_2,
\]
\[
 {\NormComparisonError}
 =\sup_{\substack{\CoefficientVector\in\RealNumbers^\ModelDimension\\\norm{\CoefficientVector}_2=1}}
       |\NormSymbol_0(\CoefficientVector)-\NormSymbol_1(\CoefficientVector)|,
\]
and
\[
 {\NormComparisonBound}
 =\max_{\substack{\CoefficientVector\in\RealNumbers^\ModelDimension\\\norm{\CoefficientVector}_2=1}}\NormSymbol_1(\CoefficientVector).
\]
Then
\begin{equation}\label{eq:coordinate-pseudometric-estimate}
 \PseudometricError_\Correspondence(\Pseudometric_0,\Pseudometric_1)
 \leq2\CoefficientBound\NormComparisonError
       +2\NormComparisonBound\CoordinateComparisonError.
\end{equation}
\end{lemma}

\begin{proof}
Fix $(x,y),(z,w)\in\Correspondence$ and set
\[
 a=\FeatureCoordinates_0(x)-\FeatureCoordinates_0(z),\qquad
 b=\FeatureCoordinates_1(y)-\FeatureCoordinates_1(w).
\]
Then $\norm{a}_2\leq2\CoefficientBound$ and $\norm{a-b}_2\leq2\CoordinateComparisonError$.
By homogeneity and the triangle inequality,
\[
 \begin{aligned}
 |\NormSymbol_0(a)-\NormSymbol_1(b)|
 &\leq|\NormSymbol_0(a)-\NormSymbol_1(a)|+\NormSymbol_1(a-b)\\
 &\leq\NormComparisonError\norm{a}_2+\NormComparisonBound\norm{a-b}_2\\
 &\leq2\CoefficientBound\NormComparisonError+2\NormComparisonBound\CoordinateComparisonError.
 \end{aligned}
\]
Taking the supremum over pairs $(x,y),(z,w)$ proves \eqref{eq:coordinate-pseudometric-estimate}.
\end{proof}

For an isometry $u\colon(Y_0,e_0)\to(Y_1,e_1)$ define \emph{pullback of coordinate classes} by
\[
 u^*\CoordinateClassOf{\ModelDimension}{\FeatureCoordinates}{\NormSymbol}
       =\CoordinateClassOf{\ModelDimension}{\FeatureCoordinates\circ u}{\NormSymbol}.
\]
Pullback commutes with the orthogonal group action in \eqref{eq:frame-action}.

Denote the set of \emph{invariant coordinate classes} on $\MetricSpace{Y}{e}$, those fixed by every self-isometry, by
\[
 \InvariantCoordinateClasses{\ModelDimension}{Y}{e}
   =\{\CoordinateClassElement\in\CoordinateClassSpace{\ModelDimension}{Y}
           \mid \FirstIsometry^*\CoordinateClassElement=\CoordinateClassElement\quad(\FirstIsometry\in\Isom(Y,e))\}.
\]
Pullback identifies these sets on isometric representatives.
Indeed, every isometry with the same domain and codomain as $u$ has the form $\FirstIsometry\circ u$ for some $\FirstIsometry\in\Isom(Y_1,e_1)$.
For an invariant class $\CoordinateClassElement$ we have $(\FirstIsometry\circ u)^*\CoordinateClassElement=u^*\FirstIsometry^*\CoordinateClassElement=u^*\CoordinateClassElement$, so the identification is independent of $u$.

To compare coordinate classes on different underlying sets, fix a common compact metric space $\MetricSpace{Z}{h}$.
Equip $Z\times\RealNumbers^\ModelDimension$ with the maximum of $h$ and the Euclidean distance.
On $\FiniteNormSpace{\ModelDimension}$ measure distances between norms uniformly on the Euclidean unit sphere $\UnitSphere^{\ModelDimension-1}$.
On the product
\[
 \CoordinateGraphSpace{\ModelDimension}{Z}{h}=\CompactHyperspace(Z\times\RealNumbers^\ModelDimension)\times\FiniteNormSpace{\ModelDimension}
\]
take the maximum $\CoordinateGraphMetric$ of the Hausdorff distance in the first component and uniform distance in the second.
More explicitly, write the product metric on $Z\times\RealNumbers^\ModelDimension$ as
\[
 \CoordinateProductMetric{h}\bigl((x,u),(y,v)\bigr)
 =\max\{h(x,y),\norm{u-v}_2\}.
\]
Then for any $(C_0,N_0),(C_1,N_1)\in\CoordinateGraphSpace{\ModelDimension}{Z}{h}$,
\[
 \CoordinateGraphMetric\bigl((C_0,N_0),(C_1,N_1)\bigr)
 =
 \max\left\{
     \HausdorffDistance{\CoordinateProductMetric{h}}(C_0,C_1),
     \max_{u\in\UnitSphere^{\ModelDimension-1}}|N_0(u)-N_1(u)|
 \right\},
\]
where $\HausdorffDistance{\CoordinateProductMetric{h}}$ is the Hausdorff distance induced by $\CoordinateProductMetric{h}$.

Let $K\in\CompactHyperspace(Z)$ and let $h_K$ be the restricted metric on $K$.
To $\CoordinateClassElement\in\CoordinateClassSpace{\ModelDimension}{K}$ assign the set of graph and norm pairs representing it,
\[
 \CoordinateOrbit{Z}{h}{\CoordinateClassElement}
  =\bigl\{(\FunctionGraph(\FeatureCoordinates),\NormSymbol)\mid
                    (\FeatureCoordinates,\NormSymbol)\in\CoordinateClassElement\bigr\}
       \subset\CoordinateGraphSpace{\ModelDimension}{Z}{h},
\]
where $\FunctionGraph(\FeatureCoordinates)=\{(x,\FeatureCoordinates(x))\mid x\in K\}$.
This is an orbit of the continuous isometric action of $\OrthogonalGroup(\ModelDimension)$ defined by $\OrthogonalChange(C,\NormSymbol)=((\IdentityMap_Z\times \OrthogonalChange)C,\NormSymbol\circ \OrthogonalChange^{-1})$, so it is compact.

Suppose $\ModelAssignment$ assigns an invariant coordinate class to each $\IsometryClass{Y}{e}$ in $\ModelNeighborhood\subset\GHSpace$.
Write $\LocalModelClass{Y}{e}$ for the class representing $\ModelValue{Y}{e}$ on $\MetricSpace{Y}{e}$.
The assignment is \emph{continuous} if for every compact metric space $\MetricSpace{Z}{h}$ the map
\[
 \{K\in\CompactHyperspace(Z)\mid \IsometryClass{K}{h_K}\in\ModelNeighborhood\}
 \ \longrightarrow\ \CompactHyperspace(\CoordinateGraphSpace{\ModelDimension}{Z}{h}),
 \qquad K\longmapsto\CoordinateOrbit{Z}{h}{\LocalModelClass{K}{h_K}}
\]
is continuous for the Hausdorff distances on both sides.

\begin{theorem}[{Approximation of distances by finite-dimensional coordinate classes}]\label{thm:local-models}
Given $\IsometryClass{X}{d}\in\GHSpace$ and $\LocalTolerance>0$, there exist a radius $r>0$, a spectral parameter $a>0$, an open neighborhood $\ModelNeighborhood$ of $\IsometryClass{X}{d}$, an integer $\ModelDimension\geq1$, and a continuous assignment
\[
 \IsometryClass{Y}{e}\longmapsto\ModelValue{Y}{e}\in\InvariantCoordinateClasses{\ModelDimension}{Y}{e}
 \qquad (\IsometryClass{Y}{e}\in\ModelNeighborhood),
\]
with the following properties.
At every $\IsometryClass{Y}{e}$ in $\ModelNeighborhood$ we have $\pm a\notin\spec(\DistanceOperator{Y}{e})$, and the approximating functions of \eqref{eq:projected-bump-profile} satisfy
\[
 \KuratowskiError{r}{a}\IsometryClass{Y}{e}
 =\max_{y,z\in Y}|e(y,z)-\ProjectedProfile{Y}{e}{y}(z)|<\LocalTolerance/2.
\]
The pseudometric induced by the coordinate class $\LocalModelClass{Y}{e}$ is $\LocalPseudometric{Y}{e}(y,z)=\|\ProjectedProfile{Y}{e}{y}-\ProjectedProfile{Y}{e}{z}\|_{\infty,Y}$.
Its distance error and the coordinate radius satisfy
\[
 \LocalError\IsometryClass{Y}{e}
 :=\max_{y,z\in Y}|e(y,z)-\LocalPseudometric{Y}{e}(y,z)|<\LocalTolerance,
 \qquad \CoordinateRadius(\LocalModelClass{Y}{e})\leq2\diam(Y,e).
\]
The functions $\KuratowskiError{r}{a}$, $\LocalError$, and the coordinate radius are continuous on $\ModelNeighborhood$.
\end{theorem}
Before proving the theorem, we fix terminology for local models and express pullback and convergence of coordinate classes in terms of representatives.
For every isometry $u\colon\MetricSpace{Y}{e}\to(Y_0,e_0)$ we have $u^*\LocalModelClass{Y_0}{e_0}=\LocalModelClass{Y}{e}$.
In particular, $\LocalPseudometric{Y_0}{e_0}(u(y),u(z))=\LocalPseudometric{Y}{e}(y,z)$.
For arbitrary representative pairs $(\LocalCoordinates{Y}{e},\LocalNorm{Y}{e})$ and $(\LocalCoordinates{Y_0}{e_0},\LocalNorm{Y_0}{e_0})$ on the two spaces, this pullback identity is equivalent to the existence of a single $\OrthogonalChange\in \OrthogonalGroup(\ModelDimension)$ satisfying both
\begin{equation}\label{eq:input-isometry-1}
\LocalCoordinates{Y_0}{e_0}\circ u=\OrthogonalChange\circ\LocalCoordinates{Y}{e},
\end{equation}
\begin{equation}\label{eq:input-isometry-2}
\LocalNorm{Y_0}{e_0}=\LocalNorm{Y}{e}\circ \OrthogonalChange^{-1}.
\end{equation}
The orthogonal group $G=\OrthogonalGroup(\ModelDimension)$ acts continuously by isometries on $\CoordinateGraphSpace{\ModelDimension}{Z}{h}$.
By \Cref{lem:compact-isometric-orbit-distance} the distance between the two orbits representing coordinate classes is
\[
 \HausdorffDistance{\CoordinateGraphMetric}(G\eta_0,G\eta_1)
 =\min_{\OrthogonalChange\in G}\CoordinateGraphMetric(\eta_0,\OrthogonalChange\,\eta_1)
 \qquad(\eta_0,\eta_1\in\CoordinateGraphSpace{\ModelDimension}{Z}{h}).
\]
Suppose $K_j,K\subset Z$ represent isometry classes in $\ModelNeighborhood$ and $K_j$ converges to $K$ in Hausdorff distance.
Fix a representative $(\FeatureCoordinates,\NormSymbol)$ at the limit.
If the orbits converge, the preceding formula allows us to choose representatives $(\FeatureCoordinates_j,\NormSymbol_j)$ on each $K_j$ whose graphs converge in Hausdorff distance and whose norms $\NormSymbol_j$ converge uniformly to $\NormSymbol$ on $\UnitSphere^{\ModelDimension-1}$.
For these same representatives and for every sequence of correspondences $R_j\subset K_j\times K$ satisfying $\sup_{(x,y)\in R_j}h(x,y)\to0$,
\begin{equation}\label{eq:local-coordinate-convergence}
\sup_{(x,y)\in R_j}\|\FeatureCoordinates_j(x)-\FeatureCoordinates(y)\|_2\longrightarrow0.
\end{equation}
Indeed, if $s_j$ is the distance between the graphs, then for each $(x,y)\in R_j$ there exists $z\in K$ with $h(x,z)\leq s_j$ and $\|\FeatureCoordinates_j(x)-\FeatureCoordinates(z)\|_2\leq s_j$.
Since $h(z,y)\leq s_j+h(x,y)$, uniform continuity of $\FeatureCoordinates$ proves the assertion.

Conversely, suppose there are representatives for which the norms converge uniformly on the sphere and \eqref{eq:local-coordinate-convergence} holds.
Put $t_j=\HausdorffDistance{h}(K_j,K)$.
The set of all $(x,y)\in K_j\times K$ satisfying $h(x,y)\leq t_j$ is a correspondence.
The graph distance is at most the maximum of $t_j$ and the supremum of the coordinate differences along this correspondence.
Hence the graph and norm pairs converge in $\CoordinateGraphMetric$, and the orbit distance tends to $0$.

In terms of representatives the induced pseudometric is
\begin{equation}\label{eq:local-model-pseudometric}
\LocalPseudometric{Y}{e}(y,z)=\LocalNorm{Y}{e}\bigl(\LocalCoordinates{Y}{e}(y)-\LocalCoordinates{Y}{e}(z)\bigr).
\end{equation}
Thus the metric quotient $\MetricQuotient{Y}{\LocalPseudometric{Y}{e}}$ is isometric to the coordinate image $\LocalCoordinates{Y}{e}(Y)$ with the distance induced by $\LocalNorm{Y}{e}$.

Fix a compact metric space $\MetricSpace{Z}{h}$.
To a representative graph $C\subset Z\times\RealNumbers^n$ and norm $N$ assign the maximum discrepancy between the two distances.
Put
\begin{equation}\label{eq:coordinate-graph-distortion}
 (C,N)\longmapsto\max_{(x,u),(y,v)\in C}|h(x,y)-N(u-v)|.
\end{equation}
For a model representative on $K\subset Z$ this value is $\LocalError\IsometryClass{K}{h_K}$.
Simultaneous orthogonal change of the graph and norm preserves the value.

This function is continuous for Hausdorff convergence of graphs and uniform convergence of norms on the sphere.
Indeed, the coordinate components of a convergent sequence of graphs lie in a common bounded ball.
By homogeneity the norms converge uniformly on its difference set.
For the fixed limit norm the function being maximized is uniformly continuous on a common compact set.
Thus the maxima converge.
The function is constant on orbits and hence induces a continuous function in the quotient topology.
For a quantitative comparison of the pseudometrics along a correspondence, see \Cref{lem:coordinate-pseudometric-estimate}.

We call $\LocalModel=(\ModelNeighborhood,\ModelDimension,\ModelAssignment,\LocalError)$ a \emph{local model}.
By the pullback identity the isometry class of the metric quotient depends only on $\IsometryClass{Y}{e}$ and defines a map from $\ModelNeighborhood$ to $\GHSpace$.\par
For fixed $r,a>0$ use the notation $\DistanceOperatorSymbol,\SpectralProjectionSymbol,V$ and the approximating functions $\ProjectedProfile{Y}{e}{y}$ of \eqref{eq:projected-bump-profile}.
Assume $n=\dim V\geq1$ and let $\BasisFunction_0,\ldots,\BasisFunction_{n-1}$ be a real $L^2$-orthonormal basis of $V$.
Choose the unique continuous representative of each basis function, using full support of the measure.
For $c=(c_0,\ldots,c_{n-1})\in\RealNumbers^n$ define the coordinates and norm by
\begin{equation}\label{eq:projected-profile-coordinates}
 \begin{aligned}
 (\LocalCoordinates{Y}{e}(y))_k
   &=\langle\ProjectedProfile{Y}{e}{y},\BasisFunction_k\rangle_2
     =\langle\ModelBump{Y}{e}{y},T\BasisFunction_k\rangle_2,\\
 \LocalNorm{Y}{e}(c)&=\left\|\sum_{k=0}^{n-1}c_k\BasisFunction_k\right\|_{\infty,Y}.
 \end{aligned}
\end{equation}
The inner products are in $L^2(Y,\SelectedMeasure{Y}{e})$.
The two inner products agree since $T$ and $P$ are self-adjoint and $T\BasisFunction_k\in V$.

Any real orthonormal basis of $V$ will do.
For example, on the two-point space $\{y_0,y_1\}$ with interpoint distance $\ell>0$, if $0<a<\ell/2$ then $V=L^2$ and the functions $\BasisFunction_k=\sqrt{2}\,\mathbf 1_{\{y_k\}}$ for $k=0,1$ form an orthonormal basis.
Yet $\DistanceOperatorSymbol \BasisFunction_0=(\ell/2)\BasisFunction_1$ and $\DistanceOperatorSymbol \BasisFunction_1=(\ell/2)\BasisFunction_0$, so neither is an eigenfunction.

\begin{lemma}[{Convergence of spectral coordinates and norms}]\label{lem:local-basis-limits}
Suppose nonempty compact subsets $Y_j$ of a compact metric space $\MetricSpace{Z}{h}$ converge in Hausdorff distance to $Y$, with restricted metrics $e_j,e$.
Fix $r,a>0$ and suppose $\pm a\notin\spec(\DistanceOperator{Y}{e})$ and $n=\dim\SpectralSubspace{Y,e}{a}\geq1$.
For every real $L^2$-orthonormal basis at the limit, for all sufficiently large $j$ there are real $L^2$-orthonormal bases of $\SpectralSubspace{Y_j,e_j}{a}$ for which the norms in \eqref{eq:projected-profile-coordinates} converge uniformly on the Euclidean unit sphere to the limit norm.
With these same bases, for any sequence of correspondences $R_j\subseteq Y_j\times Y$ with $\sup_{(x,y)\in R_j}h(x,y)\to0$, we have $\sup_{(x,y)\in R_j}\|\LocalCoordinates{Y_j}{e_j}(x)-\LocalCoordinates{Y}{e}(y)\|_2\to0$.
\end{lemma}

\begin{proof}
Let $\mu_j,\mu$ be the pushforwards of the selected measures to $Z$.
By (M3), $\mu_j\Rightarrow\mu$.
Choose bases aligned with the limit basis by \Cref{lem:spectral-continuity} and denote the column vectors of their continuous extensions by $\widetilde{\EigenbasisColumn}_j,\widetilde{\EigenbasisColumn}$.
They converge uniformly on $Z$.
Apply the ambient operators of \Cref{lem:spectral-projections} componentwise and put $F_j=A_j\widetilde{\EigenbasisColumn}_j$ and $F=A\widetilde{\EigenbasisColumn}$.
With $D=\diam(Z,h)$, each component satisfies
\[
 \|A_j\widetilde\BasisFunction_{j,k}-A\widetilde\BasisFunction_k\|_\infty
 \leq D\|\widetilde\BasisFunction_{j,k}-\widetilde\BasisFunction_k\|_\infty
      +\|(A_j-A)\widetilde\BasisFunction_k\|_\infty\longrightarrow0.
\]
The second term tends to zero by \eqref{eq:ambient-strong-convergence}.
As the number $n$ of components is fixed, $F_j\to F$ uniformly.

For $p,z\in Z$ put $\kappa_p(z)=\max\{1-2h(p,z)/r,0\}$ and define
\begin{equation}\label{eq:local-bump-parameter-integrals}
 s_j(p)=\int_Z\kappa_p(z)\,\IntegrationDifferential\mu_j(z),\qquad
 w_j(p)=\int_Z\kappa_p(z)F_j(z)\,\IntegrationDifferential\mu_j(z).
\end{equation}
Define $s,w$ similarly using $\mu,F$.
The integrands are continuous on $Z\times Z$, and $\kappa_p(z)F_j(z)$ converges uniformly to $\kappa_p(z)F(z)$.
By \Cref{lem:parameter-integrals} we have $s_j\to s$ and $w_j\to w$ uniformly on $Z$.
By \eqref{eq:projected-profile-coordinates} the quotient $w_j/s_j$ agrees with the coordinate map $\LocalCoordinates{Y_j}{e_j}$ on $Y_j$.

Since the support of $\mu$ is $Y$, we have $m=\min_{p\in Y}s(p)>0$.
The closed set $C=\{p\in Z:s(p)\geq m/2\}$ contains $Y$ in its interior, so Hausdorff convergence implies $Y_j\subset C$ for large $j$.
Also $s_j\geq m/4$ on $C$ eventually.
Hence $w_j/s_j$ converges uniformly to $v=w/s$ on $C$.
For a correspondence $R_j$ put $\CorrespondenceDisplacement_j=\sup_{(x,y)\in R_j}h(x,y)$, and let $\omega_v$ be the modulus of continuity of $v$ on $C$.
Then
\[
 \sup_{(x,y)\in R_j}\|\LocalCoordinates{Y_j}{e_j}(x)-\LocalCoordinates{Y}{e}(y)\|_2
 \leq\|w_j/s_j-v\|_{\infty,C}+\omega_v(\CorrespondenceDisplacement_j).
\]
Thus for the same bases the coordinates converge along every sequence of correspondences with $\CorrespondenceDisplacement_j\to0$.
Since $R_j$ covers both spaces, for every $\|c\|_2=1$ we also have
\[
 |\LocalNorm{Y_j}{e_j}(c)-\LocalNorm{Y}{e}(c)|
 \leq\sup_{(x,y)\in R_j}\|\EigenbasisColumn_j(x)-\EigenbasisColumn(y)\|_2.
\]
The right-hand side tends to zero by uniform convergence of the extensions and uniform continuity of the limit extension.
This bound is independent of $c$ and proves uniform convergence of the norms on the unit sphere.
\end{proof}

\begin{proposition}[{Continuity of the approximating functions for Kuratowski maps}]\label{prop:intrinsic-spectral-kernel}
Let nonempty compact sets $Y_j$ converge to $Y$ in Hausdorff distance in a compact metric space $(Z,h)$, with restricted metrics $e_j,e$.
Fix $r,a>0$ and assume $\pm a\notin\spec(\DistanceOperator{Y}{e})$.
For every sequence of correspondences $R_j\subset Y_j\times Y$ satisfying $\CorrespondenceDisplacement_j=\sup_{(x,y)\in R_j}h(x,y)\to0$,
\begin{equation}\label{eq:intrinsic-spectral-kernel-convergence}
 \KernelComparisonError_j=\sup_{(x,y),(z,w)\in R_j}
  |\ProjectedProfile{Y_j}{e_j}{x}(z)-\ProjectedProfile{Y}{e}{y}(w)|\longrightarrow0.
\end{equation}
In particular, $\KuratowskiError{r}{a}$ is continuous for the Gromov--Hausdorff metric at every class $\IsometryClass{Y}{e}$ satisfying $\pm a\notin\spec(\DistanceOperator{Y}{e})$.
\end{proposition}

\begin{proof}
Use the ambient operators $A_j,A$ and Riesz projections $\Pi_j,\Pi$ of \Cref{lem:spectral-projections}.
Both sequences converge strongly, and $\|A_j\|\leq\diam(Z,h)$.
Hence for each $f\in C(Z)$,
\[
 \|A_j\Pi_j f-A\Pi f\|_\infty
 \leq\diam(Z,h)\|(\Pi_j-\Pi)f\|_\infty+\|(A_j-A)\Pi f\|_\infty\longrightarrow0.
\]
The map $p\mapsto\kappa_p$ from $Z$ to $C(Z)$ is $(2/r)$-Lipschitz and thus has compact image.
By uniform boundedness, the operator norms of $A_j\Pi_j$ are uniformly bounded.
Applying strong convergence on a finite net of this image (Anselone and Palmer, \cite[Section~2, p.~424, last paragraph]{AnselonePalmer1968}) yields
\[
 \sup_{p\in Z}\|(A_j\Pi_j-A\Pi)\kappa_p\|_{\infty,Z}\longrightarrow0.
\]
The denominators $s_j,s$ in \eqref{eq:local-bump-parameter-integrals} are independent of the choice of basis.
By \Cref{lem:parameter-integrals} and full support of the selected measure, the argument in the proof of \Cref{lem:local-basis-limits} applies.
Set $m=\min_Y s>0$ and $C=\{s\geq m/2\}$.
For large $j$ we have $Y_j\subset C$ and $s_j\geq m/4$ on this set.
Since $s_j\to s$ uniformly on $C$,
\[
 U_j(p,z)=\frac{(A_j\Pi_j\kappa_p)(z)}{s_j(p)}
 \ \longrightarrow\
 U(p,z)=\frac{(A\Pi\kappa_p)(z)}{s(p)}
 \qquad ((p,z)\in C\times Z),
\]
uniformly.
The map $p\mapsto A\Pi\kappa_p$ is continuous into $C(Z)$, so $U$ is continuous on $C\times Z$.

Denote restriction to $Y_j$ by $\RestrictionOperator_j$ and set $\DistanceOperatorSymbol_j=\DistanceOperator{Y_j}{e_j}$ and $\SpectralProjectionSymbol_j=\SpectralProjection{Y_j,e_j}{a}$.
By \eqref{eq:projection-restriction} we have $\RestrictionOperator_j\Pi_j=\SpectralProjectionSymbol_j\RestrictionOperator_j$.
By definition of the distance kernel, $\RestrictionOperator_j A_j=\DistanceOperatorSymbol_j\RestrictionOperator_j$.
Full support ensures uniqueness of continuous representatives, so $U_j(x,z)=\ProjectedProfile{Y_j}{e_j}{x}(z)$ holds pointwise for $x,z\in Y_j$.
Likewise $U(y,w)=\ProjectedProfile{Y}{e}{y}(w)$.
Give $C\times Z$ the maximum metric and let $\omega_U$ be the modulus of continuity of $U$.
Then $\KernelComparisonError_j\leq\|U_j-U\|_\infty+\omega_U(\CorrespondenceDisplacement_j)\to0$.
Since the correspondences cover both spaces and $|e_j(x,z)-e(y,w)|\leq2\CorrespondenceDisplacement_j$,
\[
 |\KuratowskiError{r}{a}\IsometryClass{Y_j}{e_j}
  -\KuratowskiError{r}{a}\IsometryClass{Y}{e}|\leq2\CorrespondenceDisplacement_j+\KernelComparisonError_j.
\]
By \Cref{lem:common-gh-embedding} every Gromov--Hausdorff convergent sequence has such a realization.
Pairs at distance at most the Hausdorff distance form a correspondence with $\CorrespondenceDisplacement_j\to0$.
This proves continuity, including at rank zero.
\end{proof}

\textbf{A spectral parameter at an eigenvalue.} Two-point spaces of interpoint distance $2a+2^{-j}$ converge in Gromov--Hausdorff distance to the two-point space of interpoint distance $2a$.
The selected measure assigns mass $1/2$ to each point.
The eigenvalues of the approximating spaces are $\pm(a+2^{-j-1})$ and those at the limit are $\pm a$, so the spectral projections are respectively the identity and zero.
If $0<r<4a$, then $b_y=2\mathbf 1_{\{y\}}$ on every space.
Thus for distance $2a+2^{-j}$ the approximating functions coincide with the corresponding Kuratowski maps and the error is zero.
At the limit all approximating functions vanish and the error is $2a$.

\begin{proof}[{Proof of \Cref{thm:local-models}}]
If the center is a one-point space, fix any $r,a>0$, take $\ModelDimension=1$, and set
\[
 \ModelNeighborhood=\{\IsometryClass{Y}{e}\in\GHSpace:
       \diam(Y,e)<\min\{a,\LocalTolerance/2\}\},
 \qquad
 \LocalModelClass{Y}{e}=\CoordinateClassOf{1}{0}{|\cdot|}.
\]
Throughout this neighborhood $\|\DistanceOperator{Y}{e}\|\leq\diam(Y,e)<a$, so the spectral projection and all approximating functions vanish.
In a common ambient space the graph of the zero map is $Y\times\{0\}$, which varies continuously in Hausdorff distance with $Y$.
The norm is constant and every isometry preserves the coordinate class.
The pseudometric and coordinate radius vanish, and both errors equal the continuous function $\diam(Y,e)$.

Now suppose that the center has at least two points.
Fix the radius and the spectral parameter at the center and use the approximation error of the Kuratowski maps to choose a neighborhood for the coordinates.

\ProofStep{A neighborhood with constant spectral dimension.}\label{step:local-model-spectral-dimension}
Put $D=\diam(X,d)>0$ and $\eta=\min\{\LocalTolerance/4,D/2\}$.
By \Cref{lem:spectral-density} we may choose a radius $r>0$ and spectral parameter $a>0$ such that $\pm a\notin\spec(\DistanceOperator{X}{d})$, $n=\dim\SpectralSubspace{X,d}{a}\geq1$, and
\[
 \sup_{x\in X}\|d(x,\cdot)-\ProjectedProfile{X}{d}{x}\|_{\infty,X}<\eta.
\]
By \Cref{lem:common-gh-embedding,lem:spectral-continuity} there is an open neighborhood $\SpectralRankDomain$ of the center on which $\pm a$ are outside the spectrum and the dimension equals $n$.
Otherwise a sequence violating one of these conditions would converge to the center, contradicting spectral stability in a common compact realization.

By \Cref{prop:intrinsic-spectral-kernel} the function $\KuratowskiError{r}{a}$ is continuous on $\SpectralRankDomain$.
The diameter is also Gromov--Hausdorff continuous, so
\begin{equation}\label{eq:local-neighbourhood-kuratowski-error}
 \ModelNeighborhood=\{\IsometryClass{Y}{e}\in\SpectralRankDomain:
 \KuratowskiError{r}{a}\IsometryClass{Y}{e}<\min\{\LocalTolerance/2,\diam(Y,e)\}\}
\end{equation}
is open.
The approximation error at the center is less than $\eta$, so $\IsometryClass{X}{d}\in\ModelNeighborhood$.

\ProofStep{Norms and continuous coordinate maps.}\label{step:local-model-coordinates}
For each $\IsometryClass{Y}{e}\in\ModelNeighborhood$ choose a real $L^2$-orthonormal basis of $V=\SpectralSubspace{Y,e}{a}$ and define $\LocalCoordinates{Y}{e},\LocalNorm{Y}{e}$ by \eqref{eq:projected-profile-coordinates}.
The uniform norm of $\sum_k c_k\BasisFunction_k$ dominates its $L^2$ norm with respect to the probability $\SelectedMeasure{Y}{e}$, which equals $\|c\|_2$.
Thus $\LocalNorm{Y}{e}$ is a norm on $\RealNumbers^n$.

Since $\LocalCoordinates{Y}{e}(y)$ consists of the basis coefficients of $\ProjectedProfile{Y}{e}{y}$,
\[
 \LocalNorm{Y}{e}(\LocalCoordinates{Y}{e}(y)-\LocalCoordinates{Y}{e}(z))
 =\|\ProjectedProfile{Y}{e}{y}-\ProjectedProfile{Y}{e}{z}\|_{\infty,Y}
 \qquad(y,z\in Y).
\]
Hence the local model pseudometric is the uniform distance between the approximating functions (\Cref{fig:local-coordinates}).
\begin{figure}[H]
\centering
\begin{tikzpicture}[>=Stealth,node distance=12mm,font=\small]
\node (profile) {${\ModelBump{Y}{e}{y}
 \in\LebesgueSpace^2}
 (\ComparisonCarrier,\SelectedMeasure{\ComparisonCarrier}{\ComparisonMetric})$};
\node[below=of profile] (approximation) {$
 \displaystyle\sum_{\BasisIndex=0}^{\ModelDimension-1}
 (\LocalCoordinates{\ComparisonCarrier}{\ComparisonMetric}(\ComparisonPoint))_\BasisIndex
 \BasisFunction_\BasisIndex
 \in\SpectralSubspace{\ComparisonCarrier,\ComparisonMetric}{\SpectralCutoff}$};
\node[below=of approximation] (coordinates) {$
 \LocalCoordinates{\ComparisonCarrier}{\ComparisonMetric}(\ComparisonPoint)
 \in(\RealNumbers^\ModelDimension,\LocalNorm{\ComparisonCarrier}{\ComparisonMetric})$};
\draw[->] (profile) -- node[right]
 {{$\DistanceOperator{Y}{e}\SpectralProjection{Y,e}{a}$}} (approximation);
\draw[->] (approximation) -- node[right]
 {{$\BasisCoordinateMap$}} (coordinates);
\end{tikzpicture}
\caption{The two maps that give the coordinates of $y$ from the normalized bump $\ModelBump{Y}{e}{y}$. The first arrow gives the approximating function in \eqref{eq:projected-bump-profile}. Here the spectral subspace $V$ has positive dimension $\ModelDimension$. The second arrow $\BasisCoordinateMap\colon(V,\|\cdot\|_{\infty,Y})\to(\RealNumbers^\ModelDimension,\LocalNorm{Y}{e})$ is the linear isometric isomorphism that takes coefficients in the chosen real $L^2$-orthonormal basis (\eqref{eq:projected-profile-coordinates}). The coordinates of $y$ are $\LocalCoordinates{Y}{e}(y)=\BasisCoordinateMap(\ProjectedProfile{Y}{e}{y})$.}
\label{fig:local-coordinates}
\end{figure}

Write $E=\KuratowskiError{r}{a}\IsometryClass{Y}{e}$.
Since the uniform distance between $e(y,\cdot)$ and $e(z,\cdot)$ is $e(y,z)$, the triangle inequality implies
\[
 |e(y,z)-\LocalPseudometric{Y}{e}(y,z)|
 \leq\|e(y,\cdot)-\ProjectedProfile{Y}{e}{y}\|_\infty
     +\|e(z,\cdot)-\ProjectedProfile{Y}{e}{z}\|_\infty\leq2E.
\]
Moreover, $\|\ProjectedProfile{Y}{e}{y}\|_\infty\leq\diam(Y,e)+E$, so
\[
 \CoordinateRadius(\LocalModelClass{Y}{e})
 =\max_{y\in Y}\|\ProjectedProfile{Y}{e}{y}\|_\infty\leq\diam(Y,e)+E.
\]
By \eqref{eq:local-neighbourhood-kuratowski-error} the distance error is less than $\LocalTolerance$ and the coordinate radius is at most $2\diam(Y,e)$ throughout the neighborhood.

On each $Y$ the coordinate map is the composition of $y\mapsto\ModelBump{Y}{e}{y}$, $TP$, and the basis coordinate map $\BasisCoordinateMap$ on $V$.
The first map is continuous in uniform norm by \eqref{eq:normalized-bump-definition}, and the inclusion from $C(Y)$ to $L^2(Y,\SelectedMeasure{Y}{e})$ is continuous.
The operator $P$ is an orthogonal projection, and $T$ is bounded from $L^2$ to $C(Y)$ by \eqref{eq:operator-uniform-bound-1}.
The finite-dimensional coordinate map is also continuous, so $\LocalCoordinates{Y}{e}$ is continuous.

\ProofStep{Changes of basis and isometries.}\label{step:local-model-isometries}
Let $\iota\colon(Y,e)\to(Y_0,e_0)$ be a surjective isometry.
As shown after \eqref{eq:kuratowski-uniform-error}, pullback $Uf=f\circ\iota$ is an isometric isomorphism of the spectral subspaces and satisfies $U\ProjectedProfile{Y_0}{e_0}{\iota y}=\ProjectedProfile{Y}{e}{y}$.
Denote the respective columns of basis functions by $\EigenbasisColumn,\EigenbasisColumn_0$.
Since $U$ takes one orthonormal basis to another, there exists $\OrthogonalChange\in\OrthogonalGroup(n)$ such that
\begin{equation}\label{eq:projected-profile-basis-transport}
 \EigenbasisColumn_0\circ\iota=\OrthogonalChange\EigenbasisColumn.
\end{equation}
By definition of the coefficients and the uniform norm, both $\LocalCoordinates{Y_0}{e_0}\circ\iota=\OrthogonalChange\LocalCoordinates{Y}{e}$ and $\LocalNorm{Y_0}{e_0}=\LocalNorm{Y}{e}\circ\OrthogonalChange^{-1}$ hold.
With $\iota=\IdentityMap_Y$ this proves independence of the coordinate class from the basis.
In general it proves compatibility with pullback by isometries and hence invariance under self-isometries.

\ProofStep{Continuity of the approximation error.}\label{step:local-model-error-continuity}
Let $Y_j\to Y$ in a common compact space $(Z,h)$ and take correspondences with $\CorrespondenceDisplacement_j=\sup_{(x,y)\in R_j}h(x,y)\to0$.
With $\KernelComparisonError_j$ as in \Cref{prop:intrinsic-spectral-kernel}, for any $(x,y),(z,w)\in R_j$ we have
\[
 |\LocalPseudometric{Y_j}{e_j}(x,z)-\LocalPseudometric{Y}{e}(y,w)|\leq2\KernelComparisonError_j.
\]
Indeed, the pseudometric is the uniform norm of the difference between two approximating functions.
Since $R_j$ covers both spaces, the difference of the two maxima is bounded by the sum of the errors at corresponding evaluation points.
Likewise the coordinate radius is the maximum of $|\ProjectedProfile{Y}{e}{y}(z)|$ over $(y,z)\in Y^2$.
Thus
\[
 \begin{aligned}
 |\LocalError\IsometryClass{Y_j}{e_j}-\LocalError\IsometryClass{Y}{e}|&\leq2\CorrespondenceDisplacement_j+2\KernelComparisonError_j,\\
 |\CoordinateRadius(\LocalModelClass{Y_j}{e_j})-\CoordinateRadius(\LocalModelClass{Y}{e})|&\leq\KernelComparisonError_j.
 \end{aligned}
\]
Both right-hand sides tend to zero.
By \Cref{lem:common-gh-embedding}, the distance error and coordinate radius are Gromov--Hausdorff continuous.

Finally use the same bases aligned with the limit basis in \Cref{lem:local-basis-limits}.
Both coordinates and norms then converge.
The criterion following \eqref{eq:local-coordinate-convergence} proves continuity of the coordinate class assignment on $\ModelNeighborhood$.
In these bases $\ProjectedProfile{Y}{e}{y}(z)=\LocalCoordinates{Y}{e}(y)^{\mathsf T}\EigenbasisColumn(z)$, so the simultaneous convergence of coordinates and basis values also verifies \eqref{eq:intrinsic-spectral-kernel-convergence}.
\end{proof}

Identities \eqref{eq:input-isometry-1} and \eqref{eq:input-isometry-2} associate coordinate changes with isometries.
For the self-isometries of a single space $\MetricSpace{Y}{e}$ we can choose these coordinate changes to form a continuous representation, taking composition to matrix multiplication.

\begin{proposition}[{Orthogonal representation of the isometry group in a local model}]\label{prop:local-isometry-representation}\label{rem:local-isometry-representation}
Fix a local model constructed in the proof of \Cref{thm:local-models} and let $\IsometryClass{Y}{e}\in\ModelNeighborhood$.
In the spectral case take the representative $(\LocalCoordinates{\ComparisonCarrier}{\ComparisonMetric},\LocalNorm{\ComparisonCarrier}{\ComparisonMetric})$ determined by a real orthonormal basis $\{\BasisFunction_\BasisIndex\}_{0\leq\BasisIndex<\ModelDimension}$ of $\SpectralSubspace{\ComparisonCarrier,\ComparisonMetric}{\SpectralCutoff}$.
In the constant model centered at a one-point space take the zero map and the absolute value norm.
There exists a group homomorphism, continuous for the uniform topology on the isometry group,
\begin{equation}\label{eq:local-isometry-representation-1}
 \IsometryRepresentation{\ComparisonCarrier}{\ComparisonMetric}
 \colon\Isom\MetricSpace{\ComparisonCarrier}{\ComparisonMetric}
       \to\OrthogonalGroup(\ModelDimension).
\end{equation}
In the spectral case, for all $\FirstIsometry\in\Isom(Y,e)$ and $0\leq\BasisIndex,\SecondBasisIndex<\ModelDimension$ its matrix entries are
\begin{equation}\label{eq:local-isometry-representation-2}
 \bigl(\IsometryRepresentation{\ComparisonCarrier}{\ComparisonMetric}(\FirstIsometry)\bigr)_{\BasisIndex\SecondBasisIndex}
 =\int_\ComparisonCarrier
   \BasisFunction_\BasisIndex(\FirstIsometry\ComparisonPoint)
   \BasisFunction_\SecondBasisIndex(\ComparisonPoint)
   \,\IntegrationDifferential\SelectedMeasure{\ComparisonCarrier}{\ComparisonMetric}(\ComparisonPoint),
\end{equation}
and the basis functions transform according to
\begin{equation}\label{eq:isometry-basis-action}
 \BasisFunction_\BasisIndex\circ\FirstIsometry
 =\sum_{\SecondBasisIndex=0}^{\ModelDimension-1}
  \bigl(\IsometryRepresentation{\ComparisonCarrier}{\ComparisonMetric}(\FirstIsometry)\bigr)_{\BasisIndex\SecondBasisIndex}
  \BasisFunction_\SecondBasisIndex.
\end{equation}
In the constant model put $\IsometryRepresentation{Y}{e}(\FirstIsometry)=1\in\OrthogonalGroup(1)$ for every $\FirstIsometry$.
In either model the coordinate map is equivariant for this representation and the representation preserves the norm.
Namely, for all $\FirstIsometry\in\Isom(Y,e)$ and $y\in Y$,
\begin{equation}\label{eq:local-model-equivariance-1}
 \LocalCoordinates{\ComparisonCarrier}{\ComparisonMetric}(\FirstIsometry\ComparisonPoint)
 =\IsometryRepresentation{\ComparisonCarrier}{\ComparisonMetric}(\FirstIsometry)
   \LocalCoordinates{\ComparisonCarrier}{\ComparisonMetric}(\ComparisonPoint).
\end{equation}
For every $\FirstIsometry\in\Isom\MetricSpace{\ComparisonCarrier}{\ComparisonMetric}$ and $\CoefficientVector\in\RealNumbers^\ModelDimension$,
\begin{equation}\label{eq:local-model-equivariance-2}
 \LocalNorm{\ComparisonCarrier}{\ComparisonMetric}
   \bigl(\IsometryRepresentation{\ComparisonCarrier}{\ComparisonMetric}(\FirstIsometry)\CoefficientVector\bigr)
 =\LocalNorm{\ComparisonCarrier}{\ComparisonMetric}(\CoefficientVector).
\end{equation}
Thus $\LocalPseudometric{Y}{e}(\FirstIsometry y,\FirstIsometry z)=\LocalPseudometric{Y}{e}(y,z)$ for all $\FirstIsometry\in\Isom(Y,e)$ and $y,z\in Y$.
The surjective isometry $[y]\mapsto\LocalCoordinates{Y}{e}(y)$ from the metric quotient $\MetricQuotient{Y}{\LocalPseudometric{Y}{e}}$ to the coordinate image $\LocalCoordinates{Y}{e}(Y)$ intertwines the action $[y]\mapsto[\FirstIsometry y]$ on the quotient with $\IsometryRepresentation{Y}{e}(\FirstIsometry)$ on the coordinate image.
Furthermore, represent the same local model on isometric spaces $\MetricSpace{Y}{e}$ and $(Y_0,e_0)$ and fix coordinates for each.
For any surjective isometry $u\colon\MetricSpace{Y}{e}\to(Y_0,e_0)$ there is a single $\OrthogonalChange\in\OrthogonalGroup(\ModelDimension)$ satisfying \eqref{eq:input-isometry-1} and \eqref{eq:input-isometry-2} such that for every $\FirstIsometry\in\Isom(Y,e)$,
\[
 \IsometryRepresentation{Y_0}{e_0}(u\circ\FirstIsometry\circ u^{-1})
 =\OrthogonalChange\,\IsometryRepresentation{Y}{e}(\FirstIsometry)\,\OrthogonalChange^{-1}.
\]
In particular, changing the orthonormal basis on a fixed $\MetricSpace{Y}{e}$ conjugates the orthogonal representation by the change of basis.
\end{proposition}

\begin{proof}
In the constant model $\ModelDimension=1$, $\LocalCoordinates{Y}{e}=0$, and $\LocalNorm{Y}{e}=|\cdot|$.
Take the trivial representation.
The quotient and coordinate image are singletons, and for any $u$ we can choose $\OrthogonalChange=1$.
In the spectral case put $V=\SpectralSubspace{Y,e}{\SpectralCutoff}$, $\SelectedMeasureSymbol=\SelectedMeasure{Y}{e}$, and $\DistanceOperatorSymbol=\DistanceOperator{Y}{e}$.
By \Cref{thm:invariant-measures} we have $\FirstIsometry_*\SelectedMeasureSymbol=\SelectedMeasureSymbol$ for every $\FirstIsometry\in\Isom(Y,e)$.
Thus $f\mapsto f\circ\FirstIsometry$ preserves the real inner product on $L^2(Y,\SelectedMeasureSymbol)$.
For $f\in L^2(Y,\SelectedMeasureSymbol)$ and $y\in Y$, substituting $w=\FirstIsometry z$ yields
\[
 \begin{aligned}
 \bigl(\DistanceOperatorSymbol(f\circ\FirstIsometry)\bigr)(y)
 &=\int_Y e(y,z)f(\FirstIsometry z)\,\IntegrationDifferential\SelectedMeasureSymbol(z)\\
 &=\int_Y e(\FirstIsometry y,w)f(w)\,\IntegrationDifferential\SelectedMeasureSymbol(w)
   =(\DistanceOperatorSymbol f)(\FirstIsometry y).
 \end{aligned}
\]
Hence pullback preserves each eigenspace, and its inverse is pullback by $\FirstIsometry^{-1}$.
Thus $\{\BasisFunction_\BasisIndex\circ\FirstIsometry\}_{0\leq\BasisIndex<\ModelDimension}$ is another orthonormal basis of $V$.

Write $\rho=\IsometryRepresentation{Y}{e}$.
The change of basis between two orthonormal bases is an orthogonal matrix with entries given by the inner products \eqref{eq:local-isometry-representation-2}.
By \Cref{lem:distance-operator}, the elements of $V$ have continuous representatives.
Since $\SelectedMeasureSymbol$ has full support, the basis expansion \eqref{eq:isometry-basis-action} holds at every point of $Y$.
Arranging the basis functions as the column $\EigenbasisColumn=(\BasisFunction_0,\ldots,\BasisFunction_{\ModelDimension-1})^{\TransposeSymbol}$ gives the identity $\EigenbasisColumn(\FirstIsometry y)=\rho(\FirstIsometry)\EigenbasisColumn(y)$.
Thus for $\FirstIsometry,\SecondIsometry\in\Isom(Y,e)$,
\[
 {\rho(\FirstIsometry\circ\SecondIsometry)\EigenbasisColumn(y)}
 =\EigenbasisColumn(\FirstIsometry(\SecondIsometry y))
 =\rho(\FirstIsometry)\rho(\SecondIsometry)\EigenbasisColumn(y)
 \qquad(y\in Y).
\]
Since the components of $\EigenbasisColumn$ are linearly independent functions, the coefficients in corresponding rows agree.
Hence $\rho(\FirstIsometry\circ\SecondIsometry)=\rho(\FirstIsometry)\rho(\SecondIsometry)$.
Orthonormality also implies $\rho(\IdentityMap_Y)=\IdentityMatrix_{\ModelDimension}$.

To prove continuity, note that $\SelectedMeasureSymbol(Y)=1$ and $\|\BasisFunction_\SecondBasisIndex\|_2=1$ imply $\|\BasisFunction_\SecondBasisIndex\|_1\leq1$.
Thus \eqref{eq:local-isometry-representation-2} yields
\[
 \bigl|\rho(\FirstIsometry)_{\BasisIndex\SecondBasisIndex}
       -\rho(\SecondIsometry)_{\BasisIndex\SecondBasisIndex}\bigr|
 \leq \|\BasisFunction_\BasisIndex\circ\FirstIsometry
         -\BasisFunction_\BasisIndex\circ\SecondIsometry\|_\infty.
\]
The function $\BasisFunction_\BasisIndex$ is uniformly continuous on the compact space $Y$, so as $\SecondIsometry$ tends uniformly to $\FirstIsometry$ the right-hand side tends to $0$.
Every matrix entry is continuous, hence so is $\rho$.

Write $R=\rho(\FirstIsometry)$ and for $c\in\RealNumbers^\ModelDimension$ put $F_c(z)=c^{\TransposeSymbol}\EigenbasisColumn(z)$.
By $\EigenbasisColumn(\FirstIsometry z)=R\EigenbasisColumn(z)$ and orthogonality of $R$, we have $F_{Rc}(\FirstIsometry z)=F_c(z)$.
Taking uniform norms yields $\LocalNorm{Y}{e}(Rc)=\LocalNorm{Y}{e}(c)$.
Moreover, the compatibility of the operator and bumps with isometries in the proof of \Cref{thm:local-models} yields
\[
 \ProjectedProfile{Y}{e}{\FirstIsometry y}(\FirstIsometry z)
 =\ProjectedProfile{Y}{e}{y}(z)\qquad(y,z\in Y).
\]
The right-hand side is $F_{\LocalCoordinates{Y}{e}(y)}(z)$.
Hence the left-hand side evaluated at $\FirstIsometry z$ equals $F_{R\LocalCoordinates{Y}{e}(y)}(\FirstIsometry z)$.
Uniqueness of basis coefficients implies $\LocalCoordinates{Y}{e}(\FirstIsometry y)=R\LocalCoordinates{Y}{e}(y)$, proving \eqref{eq:local-model-equivariance-1} and \eqref{eq:local-model-equivariance-2}.

By \eqref{eq:local-model-equivariance-1} and \eqref{eq:local-model-equivariance-2} in the definition \eqref{eq:local-model-pseudometric} of the pseudometric, for all $y,z\in Y$ we have
\[
 \LocalPseudometric{Y}{e}(\FirstIsometry y,\FirstIsometry z)
 =\LocalNorm{Y}{e}\bigl(R(\LocalCoordinates{Y}{e}(y)-\LocalCoordinates{Y}{e}(z))\bigr)
 =\LocalPseudometric{Y}{e}(y,z).
\]
Hence $[y]\mapsto[\FirstIsometry y]$ defines a self-isometry of the metric quotient.
As $\LocalNorm{Y}{e}$ is a norm, $\LocalPseudometric{Y}{e}(y,z)=0$ is equivalent to $\LocalCoordinates{Y}{e}(y)=\LocalCoordinates{Y}{e}(z)$.
Thus $[y]\mapsto\LocalCoordinates{Y}{e}(y)$ is a bijection onto the coordinate image and preserves distances by \eqref{eq:local-model-pseudometric}.
It sends $[\FirstIsometry y]$ to $R\LocalCoordinates{Y}{e}(y)$ by \eqref{eq:local-model-equivariance-1}.

Finally, fix a surjective isometry $u\colon\MetricSpace{Y}{e}\to(Y_0,e_0)$ and let $\EigenbasisColumn_0$ be the column of basis functions chosen on $(Y_0,e_0)$.
The change of basis $\OrthogonalChange$ in \eqref{eq:projected-profile-basis-transport} satisfies $\EigenbasisColumn_0\circ u=\OrthogonalChange\EigenbasisColumn$ and, as shown in the proof of \Cref{thm:local-models}, satisfies \eqref{eq:input-isometry-1} and \eqref{eq:input-isometry-2}.
Set $\rho_0=\IsometryRepresentation{Y_0}{e_0}$.
For every $\FirstIsometry\in\Isom(Y,e)$,
\[
 \rho_0(u\circ\FirstIsometry\circ u^{-1})\,\OrthogonalChange\EigenbasisColumn(y)
 =\EigenbasisColumn_0(u(\FirstIsometry y))
 =\OrthogonalChange\,\rho(\FirstIsometry)\EigenbasisColumn(y)
 \qquad(y\in Y).
\]
Linear independence of the components of $\EigenbasisColumn$ implies equality of the matrices.
Multiplication on the right by $\OrthogonalChange^{-1}$ proves the conjugacy formula.
\end{proof}

For a compact group $\ActingGroup$, the Peter--Weyl theorem states that matrix coefficients of finite-dimensional continuous unitary representations span a dense subspace of $\ContinuousFunctions(\ActingGroup,\ComplexNumbers)$ (\cite[Chapter~6, Section~2, pp.~165--167]{DiestelSpalsbury2014}).
Its proof uses convolution operators commuting with left translations, whose nonzero eigenspaces carry finite-dimensional unitary representations.
Here the distance operator commutes with the action of $\Isom\MetricSpace{\ComparisonCarrier}{\ComparisonMetric}$, and its invariant spectral subspace determines $\IsometryRepresentation{\ComparisonCarrier}{\ComparisonMetric}$.
By \Cref{lem:spectral-continuity} we can compare these subspaces as the compact metric space varies.

\begin{remark}\label{rem:mds-comparison}
Kroshnin, Stepanov, and Trevisan proved stability of the first $m$ coordinates in multidimensional scaling of compact metric measure spaces (\cite[arXiv v2, Theorem~5.4 and Remark~5.5]{KroshninStepanovTrevisan2022}).
They assume convergence in the Gromov--Wasserstein distance of order four and the gap $\lambda_m^+>\lambda_{m+1}^+$ between positive eigenvalues of the limit operator.
Under these hypotheses the $\LebesgueSpace^2$ coordinate error, measured using an optimal coupling and minimized over orthogonal transformations, tends to $0$.
Our \Cref{thm:local-models} constructs local models for the Gromov--Hausdorff distance.
When the center $\MetricSpace{X}{d}$ has at least two points, we choose a spectral parameter $\SpectralCutoff>0$ satisfying $\pm\SpectralCutoff\notin\spec(\DistanceOperator{X}{d})$.
The uniform distances between the functions obtained by applying the distance kernel operator and its spectral projection to normalized bumps then approximate the original distances.
\end{remark}

\clearpage
\part{The topology of Gromov--Hausdorff space III: The Gromov--Hausdorff space is an absolute retract}\label{part:topology}
\begin{quote}
\small
\noindent\textbf{Abstract.}
We prove that the Gromov--Hausdorff space is an absolute retract for metrizable spaces.
Starting with continuous finite-dimensional local approximations, we represent the varying norms in a fixed Banach space and combine the local maps using a partition of unity.
We take the quotient of the hyperspace of this Banach space by a compact group action.
Antonyan's equivariant hyperspace theorem and orbit space theorem imply that the quotient is an absolute retract.
Keeping the quotient and the realization map that assigns to a compact set its induced metric fixed, we obtain approximations to the identity with arbitrarily small homotopy tracks.
Hanner's domination theorem then yields the absolute neighborhood retract property.
Since the Gromov--Hausdorff space is contractible, it is an absolute retract.
\par\smallskip
\noindent \textbf{Keywords.} Gromov--Hausdorff space, absolute retract, absolute extensor, hyperspace, approximation.

\par\smallskip
\noindent \textbf{2020 Mathematics Subject Classification.} Primary 54C55; Secondary 54E35, 54B20.
\end{quote}

\section{Introduction to Part III}\label{sec:intro-topology}
The Gromov--Hausdorff space is a global analogue of a hyperspace.
For a fixed metric space, the hyperspace of its nonempty compact subsets is equipped with the Hausdorff distance.
In the Gromov--Hausdorff space $\GHSpace$, we consider the isometry classes of all nonempty compact metric spaces and allow the ambient metric space to vary.
The ordinary unpointed Gromov--Hausdorff distance $\GHDistance$ is the infimum of the Hausdorff distances between isometric copies in all common ambient metric spaces.

Curtis proved, in particular, that the hyperspace of nonempty compact subsets of a Banach space with the Hausdorff metric induced by the norm is an absolute retract and hence an absolute neighborhood retract (\cite[Theorem~1.6]{Curtis1980}).
This result for a fixed ambient space motivates the analogous question for $\GHSpace$.
Antonyan asked whether $\GHSpace$ is an absolute retract (\cite[p.~2]{Antonyan2020Euclidean}).
The following theorem answers this question affirmatively.

\begin{theorem}\label{thm:part-iii-main}
The space $\GHSpace$ is an absolute retract for all metrizable spaces.
\end{theorem}

Antonyan also proved that for each positive integer $\FiniteSize$, the subspace represented by compact subsets of $\RealNumbers^\FiniteSize$ is homeomorphic to the Hilbert cube with one point removed (\cite[Corollary~5.3]{Antonyan2021Euclidean}).
Hyperspaces occur in the study of Gromov--Hausdorff geodesics as well.
M\'emoli and Wan proved that every such geodesic can be realized as a Hausdorff geodesic of compact subsets of one compact metric space (\cite[arXiv v2, Theorem~1]{MemoliWan2023}).

For metrics on a fixed metrizable space, the author proved an interpolation theorem with a uniform error estimate for metrics prescribed on a discrete family of closed sets (\cite[Theorem~1.1]{Ishiki2026Interpolation}).
Here the domain of an extension problem maps into $\GHSpace$, so the underlying compact space varies with the parameter.
The coordinate classes in \Cref{thm:local-models} vary continuously.
Along each convergent sequence one can therefore choose representatives for which both norms and coordinate maps converge (\eqref{eq:local-coordinate-convergence}).
Our embeddings into a Banach space preserve the model norms, allowing the partition of unity to control the error in every pairwise distance.

To prove that $\GHSpace$ is an ANR we use Hanner's domination theorem (\Cref{thm:hanner-domination}).
A retract of an ANR is an ANR.
Hanner's theorem extends this fact by allowing the composition of two maps through an ANR to be homotopic to the identity with every track of diameter less than any prescribed $\ErrorTolerance>0$.
For $\GHSpace$ we construct the auxiliary ANR as the quotient of a Banach hyperspace by a compact group action.
Fix a sequence $(\ModelDimension_\CoverIndex)_{\CoverIndex\in\NonnegativeIntegers}$ in which every positive integer occurs infinitely often.
It determines both the Banach space and the group.
The absolute retract property of the quotient follows from Antonyan's results on equivariant hyperspaces and orbit spaces.

Fix a continuous error function $\ErrorControl\colon\GHSpace\to(0,\infty)$.
Choose local models whose errors are less than $\ErrorControl$ on their domains and a countable locally finite partition of unity subordinate to these domains.
For each partition function choose a model whose domain contains its support.
Assign distinct indices $\CoverIndex$ to these pairs so that the dimension of model number $\CoverIndex$ is $\ModelDimension_\CoverIndex$.
Denote its domain by $\ModelDomain_\CoverIndex$ and its weight by $\PartitionWeight_\CoverIndex$, and put $\PartitionWeight_\CoverIndex\equiv0$ at unused indices.
For each used index $\CoverIndex$ and each compact metric representative $\MetricSpace{\ComparisonCarrier}{\ComparisonMetric}$ with isometry class in $\ModelDomain_\CoverIndex$, choose a representative of the coordinate class $\BlockModelValue{i}{Y}{e}$ on $\MetricSpace{Y}{e}$.
Denote its map and norm by
\[
 \BlockCoordinates{\CoverIndex}{\ComparisonCarrier}{\ComparisonMetric}\colon
 \ComparisonCarrier\to\RealNumbers^{\ModelDimension_\CoverIndex},
 \qquad
 \BlockNorm{\CoverIndex}{\ComparisonCarrier}{\ComparisonMetric}\colon
 \RealNumbers^{\ModelDimension_\CoverIndex}\to[0,\infty).
\]
We first represent the varying local norms in fixed Banach spaces.
For each integer $\ModelDimension\geq1$, let $\UnitSphere^{\ModelDimension-1}$ be the Euclidean unit sphere and set
\[
 \SphereFunctionSpace{\ModelDimension}
 =\ContinuousFunctions(\UnitSphere^{\ModelDimension-1},\RealNumbers),
\]
with its uniform norm.
For an integer $\ModelDimension\geq1$ and a norm $\NormSymbol$ on $\RealNumbers^\ModelDimension$, write $\NormSymbol\DualNorm$ for the dual norm with respect to the Euclidean inner product.
For $\FirstVector\in\RealNumbers^\ModelDimension$ and $\SpherePoint\in\UnitSphere^{\ModelDimension-1}$, \Cref{lem:dual-representation} uses the formula
\[
 \NormEmbedding_\NormSymbol(\FirstVector)(\SpherePoint)
 =\frac{\langle\FirstVector,\SpherePoint\rangle}
       {\NormSymbol\DualNorm(\SpherePoint)}.
\]
to define a linear isometric embedding
\[
 \NormEmbedding_\NormSymbol\colon
 (\RealNumbers^\ModelDimension,\NormSymbol)
 \to\SphereFunctionSpace{\ModelDimension}.
\]
Along the fixed sequence $(\ModelDimension_\CoverIndex)_{\CoverIndex\in\NonnegativeIntegers}$ we take the $\SummableNorm$-sum of these function spaces on spheres and the product of the corresponding orthogonal groups.
Put
\[
 \BanachAmbientSpace
 =\left(\bigoplus_{\CoverIndex\in\NonnegativeIntegers}
        \SphereFunctionSpace{\ModelDimension_\CoverIndex}\right)_{\SummableNorm},
\]
and
\[
 \ActingGroup
 =\prod_{\CoverIndex\in\NonnegativeIntegers}
   \OrthogonalGroup(\ModelDimension_\CoverIndex).
\]
Each orthogonal group acts on the corresponding space of functions by precomposition with inverse orthogonal transformations of the sphere.
The resulting action of the compact metrizable group $\ActingGroup$ on the separable Banach space $\BanachAmbientSpace$ is continuous and consists of linear isometries.
Equation~\eqref{eq:dual-equivariance} makes the embeddings compatible with orthogonal changes of the local norms and coordinates.

Let $\CompactHyperspace(\BanachAmbientSpace)$ be the space of nonempty compact subsets of $\BanachAmbientSpace$, with the Hausdorff metric induced by its norm, and set
\[
 \AuxiliaryAR=\CompactHyperspace(\BanachAmbientSpace)/\ActingGroup.
\]
Since $\BanachAmbientSpace$ is connected and locally continuum-connected, Antonyan's equivariant hyperspace theorem (\cite[Proposition~3.1]{Antonyan2003West}) and its correction (\cite{Antonyan2006Correction}) imply that $\CompactHyperspace(\BanachAmbientSpace)$ is an equivariant absolute retract, or $\ActingGroup$-AR (\Cref{thm:equivariant-hyperspace}).
Since $\ActingGroup$ is a compact group with a countable base, Antonyan's orbit space theorem (\cite[Theorem~8 and Corollary~1]{Antonyan1990}) implies that $\AuxiliaryAR$ is an AR (\Cref{thm:orbit-retract}).
In \Cref{lem:hyperspace-orbit-ar} we also prove that $\AuxiliaryAR$ is separable and metrizable.
In particular, $\AuxiliaryAR$ is an ANR through which we can approximate the identity on $\GHSpace$.

We combine the local point maps using the partition of unity.
For each compact metric representative $\MetricSpace{\ComparisonCarrier}{\ComparisonMetric}$, define $\JointFeatureAt{\ComparisonCarrier}{\ComparisonMetric}\colon
\ComparisonCarrier\to\BanachAmbientSpace$ as follows.
For every $\CoverIndex\in\NonnegativeIntegers$ and $\ComparisonPoint\in\ComparisonCarrier$, set
\[
 \bigl(\JointFeatureAt{\ComparisonCarrier}{\ComparisonMetric}(\ComparisonPoint)\bigr)_\CoverIndex
 =
 \begin{cases}
 \ModelWeightAt{\CoverIndex}{\ComparisonCarrier}{\ComparisonMetric}
 \NormEmbedding_{\BlockNorm{\CoverIndex}{\ComparisonCarrier}{\ComparisonMetric}}
   (\BlockCoordinates{\CoverIndex}{\ComparisonCarrier}{\ComparisonMetric}(\ComparisonPoint))
 &\ModelWeightAt{\CoverIndex}{\ComparisonCarrier}{\ComparisonMetric}>0,\\
 0&\text{otherwise}
 \end{cases}.
\]
For fixed $\MetricSpace{\ComparisonCarrier}{\ComparisonMetric}$, only finitely many coordinate maps are nonzero.
Thus $\JointFeatureAt{\ComparisonCarrier}{\ComparisonMetric}$ is continuous and its image $\FeatureImageAt{\ComparisonCarrier}{\ComparisonMetric}
 =\JointFeatureAt{\ComparisonCarrier}{\ComparisonMetric}(\ComparisonCarrier)$ is nonempty and compact.
Changing local bases or isometric representatives moves $\FeatureImageAt{\ComparisonCarrier}{\ComparisonMetric}$ by an element of $\ActingGroup$, so its orbit is well defined on the input isometry class.
This construction works for both trivial and nontrivial isometry groups.
By Rouyer's theorem (\Cref{thm:rouyer-generic-metric-spaces}), spaces with trivial isometry group are dense in $\GHSpace$.
Write $\OrbitProjection\colon\CompactHyperspace(\BanachAmbientSpace)\to\AuxiliaryAR$ for the orbit map.
In \Cref{thm:global-domination} we construct the continuous map
\[
 \GHSpace\xrightarrow{\ \ApproximationMap\ }
 \AuxiliaryAR
 \xrightarrow{\ \RealizationMap\ }\GHSpace.
\]
For every compact metric representative $\MetricSpace{\ComparisonCarrier}{\ComparisonMetric}$, the first map satisfies
\[
 \ApproximationAt{\ComparisonCarrier}{\ComparisonMetric}
 =\OrbitProjection(\FeatureImageAt{\ComparisonCarrier}{\ComparisonMetric}).
\]
The second map sends the orbit of a nonempty compact subset $K\subset\BanachAmbientSpace$ to its isometry class with the metric induced by the norm of $\BanachAmbientSpace$.

\Cref{thm:global-domination} also constructs a continuous homotopy
\[
 \ApproximationHomotopy\colon\GHSpace\times\UnitInterval\to\GHSpace.
\]
For each $\IsometryClass{Y}{e}$, weight the pseudometric $\BlockPseudometric{i}{Y}{e}$ of local model number $i$ by $\ModelWeightAt{i}{Y}{e}$ and denote the sum over positive weights by $\GluedPseudometric{Y}{e}$.
At time $t$ set ${\HomotopyPseudometric{Y}{e}{t}}=(1-t)e+t\GluedPseudometric{Y}{e}$.
The value $\ApproximationHomotopyAt{Y}{e}{t}$ is the isometry class of the metric quotient of $Y$ by this pseudometric.
The endpoints satisfy
\[
 \ApproximationHomotopyAt{\ComparisonCarrier}{\ComparisonMetric}{0}
 =\IsometryClass{\ComparisonCarrier}{\ComparisonMetric},
 \qquad
 \ApproximationHomotopyAt{\ComparisonCarrier}{\ComparisonMetric}{1}
 =\RealizationMap\ApproximationAt{\ComparisonCarrier}{\ComparisonMetric}.
\]
For every compact metric representative $\MetricSpace{\ComparisonCarrier}{\ComparisonMetric}$, the point track satisfies
\begin{equation}\label{eq:part-iii-domination-3}
 \diam_{\GHDistance}
 \bigl\{\ApproximationHomotopyAt{\ComparisonCarrier}{\ComparisonMetric}{t}
       \bigm| t\in\UnitInterval\bigr\}
 <\ErrorControl\IsometryClass{\ComparisonCarrier}{\ComparisonMetric}.
\end{equation}
For each $\ErrorTolerance>0$, we apply this construction with $\ErrorControl\equiv\ErrorTolerance$.
Hanner's domination theorem (\Cref{thm:hanner-domination}) and \eqref{eq:part-iii-domination-3} then imply that $\GHSpace$ is an ANR for all separable metrizable spaces.
Hanner's theorem (\Cref{thm:anr-category}) extends this conclusion to all metrizable spaces.
Finally, scaling the metrics contracts $\GHSpace$ to the one-point space.
The ANR/ANE and AR/AE equivalences in Dugundji's theorem (\Cref{thm:retract-extensor}), together with Hanner's theorem (\Cref{thm:contractible-ane}), therefore imply the absolute retract property in \Cref{thm:part-iii-main}.

\medskip
\noindent \textbf{Dependence on the other parts.}
The result we need from Part~\ref{part:spectral} is \Cref{thm:local-models}.
The global construction uses the convergence \eqref{eq:local-coordinate-convergence} following from continuity of the coordinate classes and the coordinate identities \eqref{eq:input-isometry-1} and \eqref{eq:input-isometry-2} for pullback by isometries.
The proof of \Cref{thm:part-iii-main} uses Part~\ref{part:measures} through \Cref{thm:local-models}.
At the end of this part, \Cref{thm:invariant-measured-ar} uses openness of the forgetful map and continuity of fiberwise barycenters from Part~\ref{part:measures} directly.
Part~\ref{part:hilbert} combines the absolute retract property proved here with an independent discrete approximation to show that $\GHSpace$ is homeomorphic to Hilbert space.

\medskip
\noindent \textbf{Organization.}
In Section~\ref{sec:prelim-topology}, we recall the equivariant hyperspace results and the domination and extension theorems used in the proof.
In Subsection~\ref{subsec:ar-norm-representations}, we represent varying finite-dimensional norms isometrically in a fixed Banach space, with compatibility under orthogonal changes of coordinates.
In Subsection~\ref{subsec:ar-hyperspace-orbits}, we construct absolute retracts from hyperspaces modulo compact group actions and compare their orbit metrics with the Gromov--Hausdorff distance.
In Subsection~\ref{subsec:ar-global-approximation}, we use a partition of unity to combine local models into global approximations with small homotopy tracks.
In Subsection~\ref{subsec:ar-extension}, we apply Hanner's \Cref{thm:hanner-domination,thm:anr-category} and use contractibility to obtain extension from closed subsets of arbitrary metrizable spaces, proving \Cref{thm:absolute-extensor}.

\medskip
\noindent \textbf{Conventions and notation.}
All compact metric spaces are nonempty, and Banach and Hilbert spaces are real.
We identify compact metric spaces with their isometry classes in $\GHSpace$.
An isometry is a surjective isometric embedding.
For a metric space $\MetricSpace{\BaseCarrier}{\MetricSymbol}$, we write $\CompactHyperspace(\BaseCarrier)$ for its nonempty compact subsets, equipped with the Hausdorff metric $\HausdorffDistance{\MetricSymbol}$.
For a compact space with an isometric embedding into an ambient metric space, we identify the compact space with its image and use the restricted ambient metric.
The notation $\IdentityMap$ denotes the identity map on its specified domain.
Sequences are indexed by $\NonnegativeIntegers=\{0,1,2,\ldots\}$.

\section{Preliminaries}\label{sec:prelim-topology}
We state Antonyan's theorems on equivariant hyperspaces and orbit spaces and the domination and extension theorems used in this part.
We use AR, ANR, AE, and ANE in the sense of \Cref{def:extension}.
\subsection{Equivariant hyperspaces and orbit spaces}\label{subsec:ar-preliminaries-hyperspaces}
A metrizable \emph{$\ActingGroup$-space} is a metrizable space with a continuous action of the topological group $\ActingGroup$.
A continuous map $\Mapping:\BaseCarrier\to\ComparisonCarrier$ between metrizable $\ActingGroup$-spaces $\BaseCarrier$ and $\ComparisonCarrier$ is \emph{equivariant} if for all $\GroupElement\in\ActingGroup$ and $\BasePoint\in\BaseCarrier$,
\[
 \Mapping(\GroupElement\cdot\BasePoint)
 =\GroupElement\cdot\Mapping(\BasePoint).
\]
A metrizable $\ActingGroup$-space $\BaseCarrier$ is a \emph{$\ActingGroup$-AR} if for every closed equivariant embedding $\FirstEmbedding:\BaseCarrier\to\ComparisonCarrier$ into a metrizable $\ActingGroup$-space $\ComparisonCarrier$ there exists a continuous equivariant retraction $\EquivariantRetraction:\ComparisonCarrier\to\BaseCarrier$ satisfying
\[
 \EquivariantRetraction\circ\FirstEmbedding=\IdentityMap_\BaseCarrier.
\]
A space $\BanachAmbientSpace$ is \emph{locally continuum-connected} if for every $\BasePoint\in\BanachAmbientSpace$ and every open set $\OpenSubset\subset\BanachAmbientSpace$ with $\BasePoint\in\OpenSubset$ there is an open set $\InnerNeighborhood$ such that
\[
 \BasePoint\in\InnerNeighborhood\subset\OpenSubset,
 \qquad
 \forall\ComparisonPoint\in\InnerNeighborhood\quad
 \exists\JoiningContinuum\subset\BanachAmbientSpace\quad
 \begin{cases}
 \JoiningContinuum\text{ is compact and connected},\\
 \{\BasePoint,\ComparisonPoint\}\subset\JoiningContinuum\subset\OpenSubset
 \end{cases}.
\]
In a Banach space, take a sufficiently small open ball for $\InnerNeighborhood$ and use the line segment $\JoiningContinuum=\{(1-\HomotopyTime)\BasePoint+\HomotopyTime\ComparisonPoint
 \mid\HomotopyTime\in\UnitInterval\}$.

For an action of a group $\ActingGroup$ on a space $\BanachAmbientSpace$, the induced action on $\CompactHyperspace(\BanachAmbientSpace)$ is defined as follows.
For every $\GroupElement\in\ActingGroup$ and $\FirstCompactSet\in\CompactHyperspace(\BanachAmbientSpace)$, set
\[
 \GroupElement\FirstCompactSet
 =\{\GroupElement\FirstVector\mid\FirstVector\in\FirstCompactSet\}.
\]
The following hyperspace theorem is due to Antonyan (\cite[Proposition~3.1]{Antonyan2003West}), with the correction in \cite{Antonyan2006Correction}.
See also Antonyan (\cite[Theorem~2.3]{Antonyan2020Euclidean}).

\begin{theorem}[{Antonyan~\cite[Proposition~3.1]{Antonyan2003West}}]\label{thm:equivariant-hyperspace}
Let a compact group $\ActingGroup$ act continuously on a nonempty connected, locally continuum-connected metrizable space $\BanachAmbientSpace$.
Then $\CompactHyperspace(\BanachAmbientSpace)$, with its Vietoris topology and this induced action, is a $\ActingGroup$-AR.
\end{theorem}

\begin{theorem}[{Antonyan~\cite[Theorem~8 and Corollary~1]{Antonyan1990}}]\label{thm:orbit-retract}
Let $\ActingGroup$ be a compact group with a countable basis and let $\BanachAmbientSpace$ be a metrizable $\ActingGroup$-AR.
Then its orbit space $\BanachAmbientSpace/\ActingGroup$ is an AR for all metrizable spaces.
\end{theorem}
\subsection{Domination and absolute extension}\label{subsec:ar-preliminaries-domination}
We first use domination by ANRs to obtain the ANR property for all separable metrizable spaces.
For a contractible space, we obtain the AE property.
A retract of an ANR is an ANR.
For a retract, the inclusion followed by the retraction is the identity.
The next theorem of Hanner extends this fact.
The retraction identity may be replaced by a homotopy to the identity, provided that, for every $\ErrorTolerance>0$, one can choose an ANR, maps through it, and such a homotopy whose point tracks all have diameter less than $\ErrorTolerance$.

\begin{theorem}[{Hanner~\cite[Theorem~7.2(b)]{Hanner1951}}]\label{thm:hanner-domination}
Let $\MetricSpace{\RecognitionTarget}{\TargetMetric}$ be a separable metric space.
Assume that for every $\ErrorTolerance>0$ there exist a separable metrizable ANR $\AuxiliaryAR$ and continuous maps
\[
 \ApproximationMap\colon\RecognitionTarget\to\AuxiliaryAR,
 \quad\RealizationMap\colon\AuxiliaryAR\to\RecognitionTarget,
 \quad\ApproximationHomotopy\colon\RecognitionTarget\times\UnitInterval\to\RecognitionTarget
\]
such that for every $\DomainPoint\in\RecognitionTarget$,
\[
 \ApproximationHomotopy(\DomainPoint,0)=\DomainPoint,
 \quad\ApproximationHomotopy(\DomainPoint,1)=\RealizationMap\ApproximationMap(\DomainPoint),
 \quad
 \diam_{\TargetMetric}\{\ApproximationHomotopy(\DomainPoint,\HomotopyTime)
 \mid\HomotopyTime\in\UnitInterval\}<\ErrorTolerance.
\]
Then $\RecognitionTarget$ is an ANR for all separable metrizable spaces.
\end{theorem}

\begin{theorem}[{Hanner~\cite[Theorem~13.4]{Hanner1952}}]\label{thm:anr-category}
A separable metrizable space which is an ANR for all separable metrizable spaces is an ANR for all metrizable spaces.
\end{theorem}

\begin{theorem}[{Hanner~\cite[Theorem~12.3]{Hanner1952}}]\label{thm:contractible-ane}
A nonempty contractible metrizable ANE for all metrizable spaces is an AE for all metrizable spaces.
\end{theorem}

\subsection{Interpolation of compact correspondences}\label{subsec:compact-correspondence-interpolation}
Equip the product of metric spaces $\MetricSpace{Z}{h}$ and $(E,\rho)$ with the maximum metric and write $H_\times$ for its Hausdorff distance.
Every nonempty compact subset $C\subset Z\times E$ is a correspondence between its two projections.
For $t\in[0,1]$ consider the pseudometric
\[
 D_t((x,u),(y,v))=(1-t)h(x,y)+t\rho(u,v)
\]
on the product.
Denote the isometry class of the metric quotient of its restriction to $C$ by $\GraphInterpolation(C,t)$.
Chowdhury and M\'emoli used this interpolation on a fixed optimal correspondence to construct a Gromov--Hausdorff geodesic (\cite[Theorem~1.2]{ChowdhuryMemoli2018Geodesics}).
Here both the correspondence and time vary.

\begin{lemma}\label{lem:graph-quotient-interpolation}
The map $\GraphInterpolation:\CompactHyperspace(Z\times E)\times[0,1]\to\GHSpace$ is continuous.
Its values at $t=0,1$ are the isometry classes of the first and second projections, respectively.
For $C,D\in\CompactHyperspace(Z\times E)$ and $s,t\in[0,1]$ the following estimates hold.
\begin{align}
 \GHDistance(\GraphInterpolation(C,t),\GraphInterpolation(D,t))
   &\leq H_\times(C,D),\label{eq:graph-quotient-bounds}\\
 \GHDistance(\GraphInterpolation(C,s),\GraphInterpolation(C,t))
   &\leq\tfrac12|s-t|\,\dis(C),\label{eq:graph-quotient-time-bound}\\
 |\dis(C)-\dis(D)|
   &\leq4H_\times(C,D).\label{eq:graph-distortion-bound}.
\end{align}
Here $\dis(C)=\max_{(x,u),(y,v)\in C}|h(x,y)-\rho(u,v)|$.
\end{lemma}
\begin{proof}
The projection from the product with the maximum metric to the metric quotient for $D_t$ is $1$-Lipschitz.
The images of $C,D$ there are isometric realizations of their metric quotients, so their Hausdorff distance yields \eqref{eq:graph-quotient-bounds}.
At each endpoint $D_t$ is the pullback of the metric on one factor, and the quotient is isometric to that projection.
Relate points of the two quotients at $s,t$ if they have a common representative in $C$.
This correspondence has distortion $|s-t|\dis(C)$, so the correspondence formula of Kalton and Ostrovskii (\Cref{thm:gh-correspondence}) yields \eqref{eq:graph-quotient-time-bound}.
Varying the set in \eqref{eq:graph-quotient-bounds} and then the time on a fixed set in \eqref{eq:graph-quotient-time-bound} proves joint continuity.
Finally, replace two points of $C$ by nearest points in $D$.
Each distance in $h$ or $\rho$ changes by at most $2H_\times(C,D)$.
Comparing maxima and interchanging $C,D$ proves \eqref{eq:graph-distortion-bound}.
\end{proof}
\subsection{Convergence of graphs of weighted coordinate maps}\label{subsec:weighted-coordinate-graphs}
Let $(\PartitionWeight_i)_{i\geq0}$ be a countable locally finite partition of unity on $\GHSpace$, and assign to each $i$ a positive integer $n_i$.
For every $i$ with $\PartitionWeight_i\not\equiv0$ fix a local model of dimension $n_i$ in the sense of \Cref{thm:local-models} whose domain contains $\supp\PartitionWeight_i$.
Let $E_i$ be a Banach space.
Suppose that for each norm $N$ on $\RealNumbers^{n_i}$ there is a linear isometric embedding $\BlockNormEmbedding{i}{N}\colon(\RealNumbers^{n_i},N)\to E_i$, and that $\BlockNormEmbedding{i}{N}$ depends continuously on $N$ in operator norm for the Euclidean domain norm when the norms converge uniformly on the Euclidean unit sphere.
Put $\BanachAmbientSpace=(\bigoplus_{i\geq0}E_i)_{\SummableNorm}$.
At $\xi=\IsometryClass{Y}{e}$, for each $\PartitionWeight_i(\xi)>0$ choose a pair $(\BlockCoordinates{i}{Y}{e},\BlockNorm{i}{Y}{e})$ representing the coordinate class.
For $y\in Y$ set
\[
 \FeatureBlock{i}{Y}{e}(y)=
 \begin{cases}
 \PartitionWeight_i(\xi)\,\BlockNormEmbedding{i}{\BlockNorm{i}{Y}{e}}
            \bigl(\BlockCoordinates{i}{Y}{e}(y)\bigr),&\PartitionWeight_i(\xi)>0,\\
 0,&\PartitionWeight_i(\xi)=0
 \end{cases}.
\]
Define the \emph{weighted coordinate map} by
\begin{equation}\label{eq:joint-feature-map-1}
 \JointFeatureAt{Y}{e}\colon Y\longrightarrow\BanachAmbientSpace,\qquad
 y\longmapsto\bigl(\FeatureBlock{i}{Y}{e}(y)\bigr)_{i\geq0}.
\end{equation}
At each $\xi$ only finitely many weights are positive.
Thus $\JointFeatureAt{Y}{e}$ is continuous and its graph is compact.

\begin{lemma}\label{lem:weighted-coordinate-graphs}
Let nonempty compact subsets $Y_j$ converge to $Y$ in Hausdorff distance in a compact metric space $\MetricSpace{Z}{h}$, and denote by $e_j,e$ the induced metrics on $Y_j,Y$.
For each $j$ choose a correspondence $R_j\subset Y_j\times Y$ with $\AmbientDisplacement_j=\sup_{(x,y)\in R_j}h(x,y)\to0$.
Arbitrarily fix representatives of the coordinate classes with positive weights on $Y$.
Then one can choose representatives on $Y_j$ such that
\[
 \FeatureDisplacement_j
 =\sup_{(x,y)\in R_j}\norm{\JointFeatureAt{Y_j}{e_j}(x)-\JointFeatureAt{Y}{e}(y)}_{\BanachAmbientSpace}
 \longrightarrow0.
\]
For these choices, the graphs $C_j=\FunctionGraph(\JointFeatureAt{Y_j}{e_j})$ and $C=\FunctionGraph(\JointFeatureAt{Y}{e})$ satisfy the following bound in the Hausdorff distance of the maximum metric on $Z\times\BanachAmbientSpace$.
\begin{equation}\label{eq:global-graph-convergence}
 H_\times(C_j,C)\leq\max\{\AmbientDisplacement_j,\FeatureDisplacement_j\}\longrightarrow0.
\end{equation}
\end{lemma}
\begin{proof}
Put $\xi_j=\IsometryClass{Y_j}{e_j}$ and $\xi=\IsometryClass{Y}{e}$.
The common realization implies $\GHDistance(\xi_j,\xi)\leq\HausdorffDistance{h}(Y_j,Y)\to0$.
Since the supports are locally finite, a neighborhood of $\xi$ meets only finitely many of them.
For sufficiently large $j$ only these finitely many indices enter the sum.

First suppose $\PartitionWeight_i(\xi)>0$.
For all sufficiently large $j$ we have $\PartitionWeight_i(\xi_j)>0$.
Write $(v,N)$ for the fixed representative on $Y$.
By continuity of coordinate classes in \Cref{thm:local-models} and by \eqref{eq:local-coordinate-convergence}, we can choose representatives $(v_j,N_j)$ on $Y_j$ such that $N_j$ converges uniformly on the sphere to $N$ and simultaneously $\sup_{(x,y)\in R_j}\norm{v_j(x)-v(y)}_2\to0$.
Set $A_j=\PartitionWeight_i(\xi_j)\BlockNormEmbedding{i}{N_j}$ and $A=\PartitionWeight_i(\xi)\BlockNormEmbedding{i}{N}$.
By hypothesis $\norm{A_j-A}_{\OperatorNorm}\to0$.
Thus
\[
 \begin{aligned}
 \sup_{(x,y)\in R_j}\norm{A_jv_j(x)-Av(y)}_{E_i}
 &\leq\norm{A_j}_{\OperatorNorm}\sup_{(x,y)\in R_j}\norm{v_j(x)-v(y)}_2\\
 &\quad+\norm{A_j-A}_{\OperatorNorm}\max_{y\in Y}\norm{v(y)}_2
 \longrightarrow0.
 \end{aligned}
\]
Now suppose $\PartitionWeight_i(\xi)=0$.
Whenever $\PartitionWeight_i(\xi_j)>0$, we may choose any representative.
The isometry of $\BlockNormEmbedding{i}{N}$ and the coordinate radius estimate in \Cref{thm:local-models} imply
\[
 \sup_{x\in Y_j}\norm{\FeatureBlock{i}{Y_j}{e_j}(x)}_{E_i}
 \leq 2\PartitionWeight_i(\xi_j)\diam(Y_j,e_j)
 \longrightarrow0.
\]
Indeed, $\PartitionWeight_i(\xi_j)\to0$ and the diameters are at most $\diam(Z,h)$.
The same estimate holds for a zero weight.

The norm on $\BanachAmbientSpace$ is the sum of the component norms, so summing the finitely many estimates yields $\FeatureDisplacement_j\to0$.
The correspondence $R_j$ is surjective onto both spaces.
Thus for every point of either graph, the other graph contains a point at maximum distance at most $\max\{\AmbientDisplacement_j,\FeatureDisplacement_j\}$.
This proves \eqref{eq:global-graph-convergence}.
\end{proof}

\section{Absolute extension by global approximation}\label{sec:ar}
We now combine the local maps of \Cref{thm:local-models} into a global approximation through an absolute retract.
A partition of unity permits arbitrarily many local models to be active, with only local finiteness required.
Using this construction, we prove the absolute extension property without a dimension restriction on the extension domain.
\subsection{Isometric representations of varying norms}\label{subsec:ar-norm-representations}
The local model norm $\LocalNorm{\ComparisonCarrier}{\ComparisonMetric}$ varies with the compact metric space $\MetricSpace{\ComparisonCarrier}{\ComparisonMetric}$.
Fix the input $\IsometryClass{Y}{e}$ and take the $\SummableNorm$-sum of the summands $(\RealNumbers^{\ModelDimension_i},\BlockNorm{i}{Y}{e})$ with positive weights.
Weighting each coordinate makes the distance in the image the weighted sum of the local pseudometrics (\eqref{eq:global-pseudometric-sum}).
The norm of this sum depends on the input.
Moreover, a change of basis $\OrthogonalChange\in\OrthogonalGroup(\ModelDimension)$ in \eqref{eq:frame-action} defines an isometry
\[
 \OrthogonalChange\colon(\RealNumbers^\ModelDimension,\NormSymbol)\longrightarrow
   (\RealNumbers^\ModelDimension,\NormSymbol\circ\OrthogonalChange^{-1}),
 \qquad
 (\NormSymbol\circ\OrthogonalChange^{-1})(\OrthogonalChange\FirstVector)=\NormSymbol(\FirstVector).
\]
It is a self-isometry of the original normed space precisely when $\NormSymbol\circ\OrthogonalChange^{-1}=\NormSymbol$.
We therefore embed each normed space isometrically into the continuous functions on the Euclidean unit sphere.
The embedding depends continuously on $\NormSymbol$, and for every $\NormSymbol$ the change of basis $\OrthogonalChange$ is realized by the same isometric action $f\mapsto f\circ\OrthogonalChange^{-1}$ (\Cref{lem:dual-representation-continuity,lem:dual-representation-equivariance}).
In \Cref{thm:global-domination} the group $\ActingGroup=\prod_i\OrthogonalGroup(\ModelDimension_i)$ acts isometrically on the $\SummableNorm$-sum $\BanachAmbientSpace$ of these function spaces.
The orbit of the coordinate image is the value of the approximation map in the absolute retract $\CompactHyperspace(\BanachAmbientSpace)/\ActingGroup$ of \Cref{lem:hyperspace-orbit-ar}.

For $\ModelDimension \geq1$, write
\[
 {\UnitSphere^{\ModelDimension -1}=\{\SpherePoint \in\RealNumbers^\ModelDimension \mid\norm{\SpherePoint }_2=1\},}
\]
and
\[
 {\SphereFunctionSpace{\ModelDimension} =\ContinuousFunctions(\UnitSphere^{\ModelDimension -1},\RealNumbers)},
\]
and equip $\SphereFunctionSpace{\ModelDimension} $ with the uniform norm.
Given a norm $\NormSymbol $ on $\RealNumbers^\ModelDimension $, its \emph{dual norm} on $\RealNumbers^\ModelDimension$ is defined, for each $\SecondVector\in\RealNumbers^\ModelDimension$, by
\[
 \NormSymbol\DualNorm(\SecondVector)
 =\sup\{|\langle\FirstVector,\SecondVector\rangle|
       \mid\FirstVector\in\RealNumbers^\ModelDimension,
            \NormSymbol(\FirstVector)\leq1\}.
\]
The pairing is the Euclidean inner product.
For every nonzero $\SecondVector\in\RealNumbers^\ModelDimension$, we have $\NormSymbol\DualNorm(\SecondVector)>0$.
In particular, the restriction of the dual norm to $\UnitSphere^{\ModelDimension-1}$ is positive.
Figure~\ref{fig:dual-norm-embedding} illustrates the construction in \Cref{lem:dual-representation} in dimension $2$.

\begin{figure}[H]
\centering
\begin{tikzpicture}[>=Stealth,font=\small]
\node at (0,1.65) {Radial normalization};
\filldraw[fill=black!4] (0,0) ellipse (1.55 and 0.85);
\draw[dashed] (0,0) circle (1);
\draw[->,black!35] (-1.85,0) -- (1.9,0);
\draw[->,black!35] (0,-1.1) -- (0,1.25);
\pgfmathsetmacro{\radialfactor}{1/sqrt((cos(10)/1.55)^2+(sin(10)/0.85)^2)}
\coordinate (spherepoint) at ({cos(10)},{sin(10)});
\coordinate (dualpoint) at ({\radialfactor*cos(10)},{\radialfactor*sin(10)});
\draw[->] (0,0) -- (dualpoint);
\fill (spherepoint) circle (1.3pt);
\fill (dualpoint) circle (1.3pt);
\draw[black!55] (spherepoint) -- (0.65,0.3);
\node at (0.65,0.46) {$\SpherePoint$};
\draw[black!55] (dualpoint) -- (1.5,0.58);
\node at (1.5,0.98)
 {$\dfrac{\SpherePoint}{\NormSymbol\DualNorm(\SpherePoint)}$};
\node[below left,inner sep=2pt] at (0,0) {$0$};
\node at (-0.85,0.95) {$\UnitSphere^1$};
\node at (0,-1.4)
 {$\{\SecondVector\mid\NormSymbol\DualNorm(\SecondVector)=1\}$};
\draw[->] (2.45,0.05) -- node[above,align=center]
 {pair with\\$\FirstVector$} (3.45,0.05);
\begin{scope}[xshift=4.25cm]
 \node at (2.1,1.65)
  {$\NormEmbedding_\NormSymbol(\FirstVector)(\SpherePoint(\theta))$};
 \draw[->,black!50] (-0.15,0) -- (4.5,0) node[right] {$\theta$};
 \draw[->,black!50] (0,-1.35) -- (0,1.35);
 \draw[densely dotted,black!45] (0,1.1) -- (4.2,1.1);
 \draw[densely dotted,black!45] (0,-1.1) -- (4.2,-1.1);
 \draw[thick,smooth,samples=121,domain=0:360,variable=\angle]
  plot ({4.2*\angle/360},
        {1.1*cos(\angle)/(1.55*sqrt((cos(\angle)/1.55)^2+(sin(\angle)/0.85)^2))});
 \node[left,inner sep=3pt] at (0,1.1) {$\NormSymbol(\FirstVector)$};
 \node[left,inner sep=3pt] at (0,-1.1) {$-\NormSymbol(\FirstVector)$};
 \foreach \x/\tick in {0/0,2.1/\pi,4.2/2\pi}
 {
  \draw (\x,0.05) -- (\x,-0.05);
  \node[below,inner sep=3pt] at (\x,0) {$\tick$};
 }
\end{scope}
\node at (3.25,-2.15)
 {$\NormEmbedding_\NormSymbol(\FirstVector)(\SpherePoint)
   =\left\langle\FirstVector,
      \dfrac{\SpherePoint}{\NormSymbol\DualNorm(\SpherePoint)}\right\rangle,
   \qquad
   \norm{\NormEmbedding_\NormSymbol(\FirstVector)}_\infty
   =\NormSymbol(\FirstVector)$};
\end{tikzpicture}
\caption{The fixed Euclidean sphere parametrizes the dual unit sphere by radial projection. In this example the dual unit sphere is an ellipse, $\FirstVector=(1,0)$, and $\SpherePoint(\theta)=(\cos\theta,\sin\theta)$. Pairing with $\FirstVector$ gives a continuous function on $\UnitSphere^1$ with supremum norm $\NormSymbol(\FirstVector)$. Varying the norm changes the function while keeping its domain fixed.}
\label{fig:dual-norm-embedding}
\end{figure}

In their linear isometric embedding of separable real Banach spaces into $\ContinuousFunctions(\UnitInterval,\RealNumbers)$, Banach and Mazur represented each vector $\FirstVector$ by the evaluation function $f\mapsto f(\FirstVector)$ on the closed dual unit ball (\cite[Chapter~XI, Section~8, proof of Theorem~9, pp.~185--187]{Banach1932}).
By the dual norm formula, a consequence of the Hahn--Banach theorem (\cite[Proposition~11.21]{Muscat2024}), the supremum norm of this function equals the norm of $\FirstVector$.
Restricting the evaluation function to the dual unit sphere and pulling it back to the Euclidean unit sphere by $\SpherePoint\mapsto\SpherePoint/\NormSymbol\DualNorm(\SpherePoint)$ yields \eqref{eq:dual-embedding} below.

\begin{proposition}\label{lem:dual-representation}
For a norm $\NormSymbol$ on $\RealNumbers^\ModelDimension$ define the map
\[
 \NormEmbedding_\NormSymbol:(\RealNumbers^\ModelDimension,\NormSymbol)\longrightarrow
                              (\SphereFunctionSpace{\ModelDimension},\norm{\cdot}_\infty)
\]
for $\FirstVector\in\RealNumbers^\ModelDimension$ and $\SpherePoint\in\UnitSphere^{\ModelDimension-1}$ by the formula
\begin{equation}\label{eq:dual-embedding}
 (\NormEmbedding_\NormSymbol(\FirstVector))(\SpherePoint)
 =\frac{\langle\FirstVector,\SpherePoint\rangle}{\NormSymbol\DualNorm(\SpherePoint)}.
\end{equation}
This map is a linear isometric embedding.
In particular, for $\FirstVector,\SecondVector\in\RealNumbers^\ModelDimension$,
\begin{equation}\label{eq:dual-isometry}
 \norm{\NormEmbedding_\NormSymbol(\FirstVector)-\NormEmbedding_\NormSymbol(\SecondVector)}_\infty
 =\NormSymbol(\FirstVector-\SecondVector).
\end{equation}
\end{proposition}

\begin{proof}
The formula in \eqref{eq:dual-embedding} defines a continuous function because the dual norm is continuous and positive on the compact sphere.
It is linear in $\FirstVector$.
For every $\SpherePoint\in\UnitSphere^{\ModelDimension-1}$, the definition of the dual norm yields
\[
 \frac{|\langle\FirstVector,\SpherePoint\rangle|}
      {\NormSymbol\DualNorm(\SpherePoint)}
 \leq\NormSymbol(\FirstVector).
\]
Taking the supremum over $\UnitSphere^{\ModelDimension-1}$ shows that
\begin{equation}\label{eq:dual-embedding-upper-bound}
 \sup_{\SpherePoint\in\UnitSphere^{\ModelDimension-1}}
 \frac{|\langle\FirstVector,\SpherePoint\rangle|}
      {\NormSymbol\DualNorm(\SpherePoint)}
 \leq\NormSymbol(\FirstVector).
\end{equation}
If $\FirstVector\ne0$, define $\SubspaceFunctional\colon\LinearSpan\{\FirstVector\}\to\RealNumbers$ by $\SubspaceFunctional(\SpanCoefficient\FirstVector)=\SpanCoefficient\NormSymbol(\FirstVector)$.
Its operator norm is
\[
 \sup_{\substack{\HahnBanachTestVector\in\LinearSpan\{\FirstVector\}\\\HahnBanachTestVector\ne0}}
 \frac{|\SubspaceFunctional(\HahnBanachTestVector)|}{\NormSymbol(\HahnBanachTestVector)}=1.
\]
By the Hahn--Banach theorem (\cite[Theorem~11.19 and Proposition~11.20]{Muscat2024}), the functional $\SubspaceFunctional$ extends to a functional $\HahnBanachExtension\colon\RealNumbers^\ModelDimension\to\RealNumbers$ satisfying
\[
 \sup_{\HahnBanachTestVector\in\RealNumbers^\ModelDimension\setminus\{0\}}
 \frac{|\HahnBanachExtension(\HahnBanachTestVector)|}{\NormSymbol(\HahnBanachTestVector)}=1.
\]
Using the Euclidean inner product to identify the dual, choose $\NormingVector\in\RealNumbers^\ModelDimension$ such that $\HahnBanachExtension(\HahnBanachTestVector)=\langle\HahnBanachTestVector,\NormingVector\rangle$ for every $\HahnBanachTestVector\in\RealNumbers^\ModelDimension$.
Then $\NormSymbol\DualNorm(\NormingVector)=1$ and $\langle\FirstVector,\NormingVector\rangle=\NormSymbol(\FirstVector)$.
For any $\SpherePoint\in\RealNumbers^\ModelDimension\setminus\{0\}$ and $\HomogeneityScalar\in\RealNumbers\setminus\{0\}$ we have
\[
 \frac{|\langle\FirstVector,\HomogeneityScalar\SpherePoint\rangle|}
      {\NormSymbol\DualNorm(\HomogeneityScalar\SpherePoint)}
 =\frac{|\HomogeneityScalar|\,|\langle\FirstVector,\SpherePoint\rangle|}
        {|\HomogeneityScalar|\,\NormSymbol\DualNorm(\SpherePoint)}
 =\frac{|\langle\FirstVector,\SpherePoint\rangle|}
        {\NormSymbol\DualNorm(\SpherePoint)}.
\]
Set
\[
 \NormalizedNormingVector
 =\frac{\NormingVector}{\norm{\NormingVector}_2}
 \in\UnitSphere^{\ModelDimension-1}.
\]
Then
\begin{equation}\label{eq:dual-embedding-lower-bound}
 \sup_{\SpherePoint\in\UnitSphere^{\ModelDimension-1}}
 \frac{|\langle\FirstVector,\SpherePoint\rangle|}
      {\NormSymbol\DualNorm(\SpherePoint)}
 \geq
 \frac{|\langle\FirstVector,\NormalizedNormingVector\rangle|}
      {\NormSymbol\DualNorm(\NormalizedNormingVector)}
 =\frac{|\langle\FirstVector,\NormingVector\rangle|}
        {\NormSymbol\DualNorm(\NormingVector)}
 =\NormSymbol(\FirstVector).
\end{equation}
Combining the bounds \eqref{eq:dual-embedding-upper-bound} and \eqref{eq:dual-embedding-lower-bound} with \eqref{eq:dual-embedding}, we obtain, for $\FirstVector\ne0$,
\[
 \norm{\NormEmbedding_\NormSymbol(\FirstVector)}_\infty
 =\sup_{\SpherePoint\in\UnitSphere^{\ModelDimension-1}}
       \frac{|\langle\FirstVector,\SpherePoint\rangle|}
            {\NormSymbol\DualNorm(\SpherePoint)}
 =\NormSymbol(\FirstVector).
\]
This identity is immediate when $\FirstVector=0$.
Linearity now proves \eqref{eq:dual-isometry}.
\end{proof}

\begin{proposition}\label{lem:dual-representation-continuity}
Let $n\geq1$.
Equip the norms $N$ on $\RealNumbers^n$ with the topology of uniform convergence on the Euclidean unit sphere $S=\UnitSphere^{n-1}$.
For normed spaces $V_1,V_2$, denote by $\BoundedLinearMaps(V_1,V_2)$ the space of bounded linear maps from $V_1$ to $V_2$.
Then the map
\[
 N\longmapsto\NormEmbedding_N\in\BoundedLinearMaps((\RealNumbers^n,\norm{\cdot}_2),C(S,\RealNumbers))
\]
is continuous in operator norm.
Consequently, if $N_j\to N$, then $\NormEmbedding_{N_j}$ converges uniformly to $\NormEmbedding_N$ on bounded sets, and $x_j\to x$ implies $\NormEmbedding_{N_j}x_j\to\NormEmbedding_N x$.
\end{proposition}

\begin{proof}
Set $\SphereMultiplier{N}=1/(N^*|_S)$.
Define $\SphereMaxTransform$ on $C(S,\RealNumbers)$ by $(\SphereMaxTransform f)(\SpherePoint)=\max_{v\in S}|\langle \SpherePoint,v\rangle|f(v)$.
The function being maximized is uniformly continuous on $S\times S$, so $\SphereMaxTransform f$ is continuous.
Since $|\langle \SpherePoint,v\rangle|\leq1$, we have $\norm{\SphereMaxTransform f-\SphereMaxTransform g}_\infty\leq\norm{f-g}_\infty$.
By homogeneity,
\begin{equation}\label{eq:dual-embedding-continuity}
 N^*(\SpherePoint)=\max_{v\in S}\frac{|\langle \SpherePoint,v\rangle|}{N(v)},\qquad
 N^*|_S=\SphereMaxTransform\bigl(1/(N|_S)\bigr).
\end{equation}
Taking reciprocals is continuous on the positive continuous functions on the compact space $S$.
Hence $N\mapsto\SphereMultiplier{N}$ is continuous.
By \eqref{eq:dual-embedding} we have $(\NormEmbedding_N x)(\SpherePoint)=\SphereMultiplier{N}(\SpherePoint)\langle x,\SpherePoint\rangle$.
Consequently for any two norms $N,M$,
\begin{equation}\label{eq:dual-embedding-uniform-on-bounded}
 \norm{\NormEmbedding_N-\NormEmbedding_M}_{\OperatorNorm}
 =\max_{\SpherePoint\in S}|\SphereMultiplier{N}(\SpherePoint)-\SphereMultiplier{M}(\SpherePoint)|
 =\norm{\SphereMultiplier{N}-\SphereMultiplier{M}}_\infty.
\end{equation}
Indeed, if $\norm{x}_2\leq1$, the Cauchy--Schwarz inequality implies $|\langle x,\SpherePoint\rangle|\leq1$ for each $\SpherePoint\in S$, with equality at $x=\SpherePoint$.
Maximizing over $\SpherePoint$ yields \eqref{eq:dual-embedding-uniform-on-bounded}, which proves operator-norm continuity.
On the ball of radius $R$ the difference is at most $R$ times this operator norm.
If $N_j\to N$, the norms $\norm{\NormEmbedding_{N_j}}_{\OperatorNorm}$ are bounded, so the assertion for $x_j\to x$ follows from the triangle inequality.
\end{proof}

For $\OrthogonalChange\in\OrthogonalGroup(\ModelDimension)$ and a function $\SphereFunction\in\SphereFunctionSpace{\ModelDimension}$ on the sphere, define $\OrthogonalChange\SphereFunction$ by
\[
 (\OrthogonalChange\SphereFunction)(\SpherePoint)=\SphereFunction(\OrthogonalChange^{-1}\SpherePoint)
 \qquad(\SpherePoint\in\UnitSphere^{\ModelDimension-1}).
\]
\begin{proposition}\label{lem:dual-representation-equivariance}
With this formula, $\OrthogonalGroup(\ModelDimension)$ acts continuously on $\SphereFunctionSpace{\ModelDimension}$ by linear isometries.
For every norm $\NormSymbol$ on $\RealNumbers^\ModelDimension$, every $\OrthogonalChange\in\OrthogonalGroup(\ModelDimension)$, and every $\FirstVector\in\RealNumbers^\ModelDimension$,
\begin{equation}\label{eq:dual-equivariance}
 \NormEmbedding_{\NormSymbol\circ\OrthogonalChange^{-1}}(\OrthogonalChange\FirstVector)
 =\OrthogonalChange\NormEmbedding_\NormSymbol(\FirstVector).
\end{equation}
\end{proposition}

\begin{proof}
Orthogonal transformations preserve the Euclidean unit sphere, so $\norm{\OrthogonalChange \SphereFunction }_\infty=\norm{\SphereFunction }_\infty$.
If $\OrthogonalChange _\SequenceIndex \to \OrthogonalChange $ and $\SphereFunction _\SequenceIndex \to \SphereFunction $ uniformly, uniform continuity of $\SphereFunction $ implies
\[
 \norm{\OrthogonalChange _\SequenceIndex \SphereFunction _\SequenceIndex -\OrthogonalChange \SphereFunction }_\infty
 \leq\norm{\SphereFunction _\SequenceIndex -\SphereFunction }_\infty+
      \sup_{\SpherePoint \in\UnitSphere^{\ModelDimension -1}}|\SphereFunction (\OrthogonalChange _\SequenceIndex ^{-1}\SpherePoint )-\SphereFunction (\OrthogonalChange ^{-1}\SpherePoint )|\to0.
\]
Finally, substituting $\FirstVector =\OrthogonalChange \SecondVector $ in the definition of the dual norm $\NormSymbol\DualNorm$ yields
\[
 (\NormSymbol \circ \OrthogonalChange ^{-1})\DualNorm(\SpherePoint )=\NormSymbol \DualNorm(\OrthogonalChange ^{-1}\SpherePoint ).
\]
Substituting this identity into \eqref{eq:dual-embedding}, we obtain
\[
 \NormEmbedding _{\NormSymbol \circ \OrthogonalChange ^{-1}}(\OrthogonalChange \FirstVector )(\SpherePoint )
 =\frac{\langle \OrthogonalChange \FirstVector ,\SpherePoint \rangle}{\NormSymbol \DualNorm(\OrthogonalChange ^{-1}\SpherePoint )}
 =\NormEmbedding _\NormSymbol (\FirstVector )(\OrthogonalChange ^{-1}\SpherePoint ),
\]
which proves \eqref{eq:dual-equivariance}.
\end{proof}
\subsection{Absolute retracts from hyperspace orbit spaces}\label{subsec:ar-hyperspace-orbits}
Equation \eqref{eq:dual-equivariance} converts a change of coordinates into an isometry of a fixed Banach space.
The global approximation will use compact subsets of such a space modulo this action.
We apply Antonyan's two theorems (\Cref{thm:equivariant-hyperspace,thm:orbit-retract}) to prove that the resulting hyperspace orbit space is an absolute retract.

\begin{theorem}\label{lem:hyperspace-orbit-ar}
Let $(\BanachAmbientSpace ,\norm{\cdot}_\BanachAmbientSpace )$ be a separable Banach space, and let a compact metrizable group $\ActingGroup $ act continuously on $\BanachAmbientSpace $ by linear isometries.
Equip $\BanachAmbientSpace$ with its induced metric $\BanachMetric\colon\BanachAmbientSpace\times\BanachAmbientSpace\to[0,\infty)$.
For every $\FirstVector,\SecondVector\in\BanachAmbientSpace$, it satisfies $\BanachMetric(\FirstVector,\SecondVector)=\norm{\FirstVector-\SecondVector}_\BanachAmbientSpace$.
Then the orbit space $\CompactHyperspace(\BanachAmbientSpace )/\ActingGroup $ is a separable metrizable AR.
For nonempty compact subsets $\FirstCompactSet,\SecondCompactSet\subset\BanachAmbientSpace$, the formula
\begin{equation}\label{eq:orbit-metric}
  \OrbitMetric([\FirstCompactSet ],[\SecondCompactSet ])=\min_{\GroupElement \in \ActingGroup }\HausdorffDistance{\BanachMetric }(\FirstCompactSet ,\GroupElement \SecondCompactSet )
\end{equation}
defines a metric that induces the quotient topology.
\end{theorem}

\begin{proof}
\textbf{Step 1.
The absolute retract property.}
The induced action on $\CompactHyperspace(\BanachAmbientSpace )$ is continuous.
To see this, assume that $\GroupElement _\SequenceIndex \to \GroupElement $ and $\HausdorffDistance{\BanachMetric }(\FirstCompactSet _\SequenceIndex ,\FirstCompactSet )\to0$.
Then
\[
 \HausdorffDistance{\BanachMetric }(\GroupElement _\SequenceIndex \FirstCompactSet _\SequenceIndex ,\GroupElement \FirstCompactSet )
 \leq\HausdorffDistance{\BanachMetric }(\FirstCompactSet _\SequenceIndex ,\FirstCompactSet )+\sup_{\CompactSetPoint \in \FirstCompactSet }\norm{\GroupElement _\SequenceIndex \CompactSetPoint -\GroupElement \CompactSetPoint }_\BanachAmbientSpace \to0.
\]
The action is continuous and $\FirstCompactSet $ is compact, so the last supremum tends to $0$.
Each group element acts isometrically on the hyperspace $\CompactHyperspace(\BanachAmbientSpace)$.
By Michael's hyperspace theorem (\Cref{thm:compact-hyperspace-topology}), the Hausdorff metric induces the Vietoris topology on $\CompactHyperspace(\BanachAmbientSpace)$.
Since $\BanachAmbientSpace $ is connected and locally continuum-connected, Antonyan's hyperspace theorem (\Cref{thm:equivariant-hyperspace}) implies that $\CompactHyperspace(\BanachAmbientSpace )$ is a $\ActingGroup $-AR.
The compact metrizable group $\ActingGroup $ has a countable base.
By Antonyan's orbit space theorem (\Cref{thm:orbit-retract}), $\CompactHyperspace(\BanachAmbientSpace )/\ActingGroup $ is an AR for metrizable spaces.

\textbf{Step 2.
The orbit metric and quotient topology.}
As shown in Step~1, the compact group $\ActingGroup$ acts continuously on $\CompactHyperspace(\BanachAmbientSpace)$ by isometries of the Hausdorff metric $\HausdorffDistance{\BanachMetric}$.
Hence \Cref{lem:compact-isometric-orbit-distance} applies to this hyperspace.

The minimum formula in that lemma is \eqref{eq:orbit-metric}.
Thus $\OrbitMetric$ induces the quotient topology, and the orbit map $\OrbitProjection\colon\CompactHyperspace(\BanachAmbientSpace)\to\CompactHyperspace(\BanachAmbientSpace)/\ActingGroup$ is continuous.

Since $\BanachAmbientSpace$ is separable, Michael's hyperspace theorem (\Cref{thm:compact-hyperspace-topology}) implies that $\CompactHyperspace(\BanachAmbientSpace)$ is separable.
Its continuous image under $\OrbitProjection$ is therefore separable.
\end{proof}

\begin{lemma}\label{lem:hyperspace-orbit-realization}
Let $(\BanachAmbientSpace ,\norm{\cdot}_\BanachAmbientSpace )$ be a separable Banach space, and let a compact metrizable group $\ActingGroup $ act continuously on $\BanachAmbientSpace $ by linear isometries.
Equip $\BanachAmbientSpace$ with its induced metric $\BanachMetric\colon\BanachAmbientSpace\times\BanachAmbientSpace\to[0,\infty)$.
For every $\FirstVector,\SecondVector\in\BanachAmbientSpace$, it satisfies $\BanachMetric(\FirstVector,\SecondVector)=\norm{\FirstVector-\SecondVector}_\BanachAmbientSpace$.
Equip $\CompactHyperspace(\BanachAmbientSpace)/\ActingGroup$ with the metric in \eqref{eq:orbit-metric}.
Then the map
\[
 \RealizationMap([\FirstCompactSet ])=\MetricSpace{\FirstCompactSet }{\RestrictedMetric{\BanachMetric }{\FirstCompactSet }},
\]
defined by
\[
 \RealizationMap\colon\CompactHyperspace(\BanachAmbientSpace )/\ActingGroup \to\GHSpace
\]
is well-defined and $1$-Lipschitz.
\end{lemma}

\begin{proof}
For every $\GroupElement \in \ActingGroup $, restricting $\GroupElement$ defines an isometry $\SecondCompactSet \to \GroupElement \SecondCompactSet $.
Thus $\RealizationMap$ is independent of the orbit representative.
For each $\GroupElement\in\ActingGroup$ we also have
\[
 \begin{aligned}
 \GHDistance(\RealizationMap([\FirstCompactSet ]),\RealizationMap([\SecondCompactSet ]))
 &=\GHDistance(\MetricSpace{\FirstCompactSet }{\BanachMetric |_{\FirstCompactSet \times \FirstCompactSet }},
         \MetricSpace{\GroupElement \SecondCompactSet }{\BanachMetric |_{\GroupElement \SecondCompactSet \times \GroupElement \SecondCompactSet }})\\
 &\leq\HausdorffDistance{\BanachMetric }(\FirstCompactSet ,\GroupElement \SecondCompactSet ).
 \end{aligned}
\]
Taking the minimum over $\GroupElement\in\ActingGroup$ proves that $\RealizationMap$ is $1$-Lipschitz.
\end{proof}

\begin{remark}\label{rem:peano-hyperspace}
A classical hyperspace theorem connects compact subsets with the Hilbert cube.
Let $\MetricSpace{X}{d}$ be a nondegenerate Peano continuum, that is, a compact, connected, locally connected metric space with at least two points.
Curtis and Schori proved that the hyperspace $\bigl(\CompactHyperspace(X),\HausdorffDistance{d}\bigr)$ is homeomorphic to the Hilbert cube $\HilbertCube=[0,1]^{\NonnegativeIntegers}$ (\cite[Theorem~3.2]{CurtisSchori1978}).
This hyperspace consists of all nonempty compact subsets of a fixed space $X$.
In \Cref{lem:hyperspace-orbit-ar} the points are compact subsets of a Banach space up to a compact group action.
\end{remark}
\subsection{Global approximation with small homotopy tracks}\label{subsec:ar-global-approximation}
Fix a sequence $(\ModelDimension_\CoverIndex)_{\CoverIndex\in\NonnegativeIntegers}$ of positive integers in which each positive integer occurs infinitely often.
Using the spaces $\SphereFunctionSpace{\ModelDimension_\CoverIndex}$ of real continuous functions on the spheres and the orthogonal groups $\OrthogonalGroup(\ModelDimension_\CoverIndex)$, put
\begin{equation}\label{eq:fixed-domination-space}
 \BanachAmbientSpace=\left(\bigoplus_{\CoverIndex\in\NonnegativeIntegers}\SphereFunctionSpace{\ModelDimension_\CoverIndex}\right)_{\SummableNorm},
 \qquad \ActingGroup=\prod_{\CoverIndex\in\NonnegativeIntegers}\OrthogonalGroup(\ModelDimension_\CoverIndex),
 \qquad \AuxiliaryAR=\CompactHyperspace(\BanachAmbientSpace)/\ActingGroup.
\end{equation}
The orthogonal group of summand $\CoverIndex$ acts on functions on its sphere by $(g f)(u)=f(g^{-1}u)$.
Denote the norm metric on $\BanachAmbientSpace$ by $\BanachMetric$ and define the realization map $\RealizationMap$ of \Cref{lem:hyperspace-orbit-realization} by
\[
 \RealizationMap\colon\AuxiliaryAR\longrightarrow\GHSpace,\qquad
 [K]\longmapsto\IsometryClass{K}{\RestrictedMetric{\BanachMetric}{K}}.
\]
Below the values of $\ModelError$ and $\ErrorControl$ at $\MetricSpace{Y}{e}$ mean their values at the isometry class $\IsometryClass{Y}{e}$.

\begin{theorem}\label{thm:global-domination}
The space $\AuxiliaryAR$ defined above is a separable metrizable AR, and $\RealizationMap$ is $1$-Lipschitz.
For every positive continuous function $\ErrorControl\colon\GHSpace\to(0,\infty)$ there exist a continuous map $\ApproximationMap\colon\GHSpace\to\AuxiliaryAR$, a continuous function $\ModelError\colon\GHSpace\to[0,\infty)$, and a continuous homotopy $\ApproximationHomotopy\colon\GHSpace\times\UnitInterval\to\GHSpace$ with the following properties.
For every $\IsometryClass{Y}{e}\in\GHSpace$,
\[
 \ModelError\IsometryClass{\ComparisonCarrier}{\ComparisonMetric}
 <\ErrorControl\IsometryClass{\ComparisonCarrier}{\ComparisonMetric}.
\]
The endpoints of the homotopy satisfy
\begin{equation}\label{eq:global-domination-1}
 \ApproximationHomotopyAt{\ComparisonCarrier}{\ComparisonMetric}{0}=\IsometryClass{\ComparisonCarrier }{\ComparisonMetric },
\end{equation}
\begin{equation}\label{eq:global-domination-2}
 \ApproximationHomotopyAt{\ComparisonCarrier}{\ComparisonMetric}{1}=\RealizationMap\ApproximationAt{\ComparisonCarrier}{\ComparisonMetric}.
\end{equation}
For every $\IsometryClass{Y}{e}\in\GHSpace$ and all $\FirstTime,\HomotopyTime\in\UnitInterval$, the tracks of the homotopy satisfy
\begin{equation}\label{eq:global-track-control}
  \GHDistance(\ApproximationHomotopyAt{\ComparisonCarrier}{\ComparisonMetric}{\FirstTime},\ApproximationHomotopyAt{\ComparisonCarrier}{\ComparisonMetric}{\HomotopyTime})\leq\frac{\abs{\FirstTime -\HomotopyTime }}2 \ModelError\IsometryClass{\ComparisonCarrier }{\ComparisonMetric }.
\end{equation}
\end{theorem}

\begin{proof}
Each $\SphereFunctionSpace{n}$ is a separable Banach space, so its countable $\SummableNorm$-sum $\BanachAmbientSpace$ is a separable Banach space.
The group $\ActingGroup$ is compact and metrizable.
By \Cref{lem:dual-representation-equivariance} its action on each summand is continuous and isometric.
Finitely supported vectors are dense, so $g_j\to g$ implies $g_j u\to gu$ for every $u\in\BanachAmbientSpace$.
The inequality $\|g_j u_j-gu\|\leq\|u_j-u\|+\|g_j u-gu\|$ yields joint continuity of the action.
By \Cref{lem:hyperspace-orbit-ar}, which uses Antonyan's results, $\AuxiliaryAR$ is a separable metrizable AR.
By \Cref{lem:hyperspace-orbit-realization} the map $\RealizationMap$ taking an orbit to the induced metric on a representative is $1$-Lipschitz.
Fix $\ErrorControl$.
We now construct $\ApproximationMap,\ModelError,\ApproximationHomotopy$.

\ProofStep{Local models, a partition of unity, and the ambient space.}\label{step:global-local-models}
For each center $\MetricSpace{\BaseCarrier}{\MetricSymbol}\in\GHSpace$, apply \Cref{thm:local-models} with a positive tolerance smaller than $\ErrorControl\IsometryClass{\BaseCarrier}{\MetricSymbol}$.
By continuity of the model error and of $\ErrorControl$, shrink the model domain so that for every space $\MetricSpace{\ComparisonCarrier}{\ComparisonMetric}$ in that domain the model error is less than $\ErrorControl\IsometryClass{\ComparisonCarrier}{\ComparisonMetric}$.
Let $\OpenCover$ be the family of these domains.
Then $\OpenCover$ is an open cover of $\GHSpace$.

Since $\GHSpace$ is separable and metrizable, it admits a countable locally finite partition of unity subordinate to $\OpenCover$ (Stone and Michael, \Cref{thm:metric-partition-unity}).
For each partition function choose a model whose domain contains its support.
Several functions may use the same model.
Enumerate these pairs.
For a model of dimension $\ModelDimension$, assign an unused index $\CoverIndex$ satisfying $\ModelDimension_\CoverIndex=\ModelDimension$.
This is possible at every stage because $\ModelDimension$ occurs infinitely often.
Denote the resulting partition functions by $\PartitionWeight_\CoverIndex$ and put $\PartitionWeight_\CoverIndex\equiv0$ at unused indices.
When $\PartitionWeight_\CoverIndex\not\equiv0$ write the associated local model as
\[
 \LocalModel_\CoverIndex
 {=(\ModelDomain_\CoverIndex,\ModelDimension_\CoverIndex,}
   {\ModelAssignment_\CoverIndex,\BlockModelError{\CoverIndex})}.
\]
Here $\ModelDimension_\CoverIndex\geq1$ is its dimension, $\ModelAssignment_\CoverIndex$ the assignment of coordinate classes, and $\BlockModelError{\CoverIndex}\colon\ModelDomain_\CoverIndex\to[0,\infty)$ the error function.
For each compact metric space $\MetricSpace{\ComparisonCarrier}{\ComparisonMetric}$ with class in $\ModelDomain_\CoverIndex$, represent $\BlockModelValue{i}{Y}{e}$ by a pair on $\MetricSpace{Y}{e}$.
Denote its map by $\BlockCoordinates{\CoverIndex}{\ComparisonCarrier}{\ComparisonMetric}\colon\ComparisonCarrier\to\RealNumbers^{\ModelDimension_\CoverIndex}$, its norm by $\BlockNorm{\CoverIndex}{\ComparisonCarrier}{\ComparisonMetric}\colon\RealNumbers^{\ModelDimension_\CoverIndex}\to[0,\infty)$, and its induced local pseudometric by $\BlockPseudometric{i}{Y}{e}$.
For all $y,z\in Y$,
\[
 \BlockPseudometric{\CoverIndex}{\ComparisonCarrier}{\ComparisonMetric}(\ComparisonPoint,\IntegrationPoint)
 =\BlockNorm{\CoverIndex}{\ComparisonCarrier}{\ComparisonMetric}
   (\BlockCoordinates{\CoverIndex}{\ComparisonCarrier}{\ComparisonMetric}(\ComparisonPoint)
    -\BlockCoordinates{\CoverIndex}{\ComparisonCarrier}{\ComparisonMetric}(\IntegrationPoint)).
\]
The approximation error of this pseudometric satisfies
\[
 \BlockModelError{\CoverIndex}\IsometryClass{\ComparisonCarrier}{\ComparisonMetric}
 =\norm{\ComparisonMetric-\BlockPseudometric{\CoverIndex}{\ComparisonCarrier}{\ComparisonMetric}}_\infty
 <\ErrorControl\IsometryClass{\ComparisonCarrier}{\ComparisonMetric}.
\]
Let $\OrthogonalChange_\CoverIndex\in\OrthogonalGroup(\ModelDimension_\CoverIndex)$ be the orthogonal change of basis for each local model.
Setting $\OrthogonalChange_\CoverIndex=\IdentityMap_{\RealNumbers^{\ModelDimension_\CoverIndex}}$ at unused indices makes $(\OrthogonalChange_\CoverIndex)_{\CoverIndex\in\NonnegativeIntegers}$ an element of $\ActingGroup$.

\ProofStep{The global map and its error.}\label{step:global-map}
For $\xi=\IsometryClass{Y}{e}$ put $\ActiveBlocksAt{Y}{e}=\{i\geq0\mid\PartitionWeight_i(\xi)>0\}$.
By \Cref{lem:dual-representation,lem:dual-representation-continuity} the spaces $E_i=\SphereFunctionSpace{n_i}$ and embeddings $\BlockNormEmbedding{i}{N}=\NormEmbedding_N$ satisfy the hypotheses of \Cref{lem:weighted-coordinate-graphs}.
Thus \eqref{eq:joint-feature-map-1} defines a continuous map $\JointFeatureAt{Y}{e}\colon Y\to\BanachAmbientSpace$.
Denote its image by $\FeatureImageAt{Y}{e}=\JointFeatureAt{Y}{e}(Y)$.
Since each $\NormEmbedding_N$ is isometric and the weights sum to $1$, the coordinate radius bound in \Cref{thm:local-models} implies
\[
 \norm{\JointFeatureAt{Y}{e}(y)}_{\BanachAmbientSpace}\leq 2\diam(Y,e)
 \qquad(y\in Y).
\]
If the local coordinate pairs are changed as in \eqref{eq:frame-action} by $\GroupElement=\{\OrthogonalChange_\CoverIndex\}_\CoverIndex\in\ActingGroup$, then \eqref{eq:dual-equivariance} transforms the joint point map into $\GroupElement\circ\JointFeatureAt{\ComparisonCarrier}{\ComparisonMetric}$ and its image into $\GroupElement\FeatureImageAt{\ComparisonCarrier}{\ComparisonMetric}$.
If $\InputIsometry\colon\MetricSpace{\ComparisonCarrier}{\ComparisonMetric}
 \to\MetricSpace{\ComparisonCarrier_0}{\ComparisonMetric_0}$ is an isometry, then \eqref{eq:input-isometry-1}, \eqref{eq:input-isometry-2}, and \eqref{eq:dual-equivariance} imply that for some $\GroupElement\in\ActingGroup$,
\[
 \JointFeatureAt{\ComparisonCarrier_0}{\ComparisonMetric_0}\circ\InputIsometry
 =\GroupElement\circ\JointFeatureAt{\ComparisonCarrier}{\ComparisonMetric}.
\]
Fix $\xi=\IsometryClass{Y}{e}$ and regard $C=\FunctionGraph(\JointFeatureAt{Y}{e})$ as a correspondence between $Y$ and its coordinate image.
This graph determines the approximation map, the error, and the homotopy.
\[
 \begin{aligned}
 \ApproximationMap(\xi)&=[\CoordinateProjection_2 C],\\
 \ModelError(\xi)&=\dis(C),\\
 \ApproximationHomotopy(\xi,t)&=\GraphInterpolation(C,t)\qquad(0\leq t\leq1).
 \end{aligned}
\]
The brackets in the first formula denote the orbit in $\AuxiliaryAR=\CompactHyperspace(\BanachAmbientSpace)/\ActingGroup$.
By the preceding identity, changing representatives takes $C$ to the new graph by the componentwise isometry $(y,b)\mapsto(\InputIsometry(y),\GroupElement b)$.
Hence the orbit of the second projection, the distortion of the correspondence, and the isometry class of the interpolated metric quotient are independent of representatives.
In particular, the map $\RealizationMap$ fixed above sends $\ApproximationMap(\xi)$ to the isometry class of the coordinate image $\FeatureImageAt{Y}{e}$.

For $\ComparisonPoint,\IntegrationPoint\in\ComparisonCarrier$, define the pseudometric induced by $\JointFeatureAt{\ComparisonCarrier}{\ComparisonMetric}$ by
\begin{equation}\label{eq:global-pseudometric}
  \GluedPseudometric{\ComparisonCarrier}{\ComparisonMetric} (\ComparisonPoint ,\IntegrationPoint )=\norm{\JointFeatureAt{\ComparisonCarrier}{\ComparisonMetric} (\ComparisonPoint )-\JointFeatureAt{\ComparisonCarrier}{\ComparisonMetric} (\IntegrationPoint )}_\BanachAmbientSpace.
\end{equation}
By \eqref{eq:dual-isometry}, we have
\begin{equation}\label{eq:global-pseudometric-sum}
 \GluedPseudometric{\ComparisonCarrier}{\ComparisonMetric}(\ComparisonPoint,\IntegrationPoint)
 =\sum_{\substack{\CoverIndex\in\NonnegativeIntegers\\\ModelWeightAt{\CoverIndex}{\ComparisonCarrier}{\ComparisonMetric}>0}} \ModelWeightAt{\CoverIndex}{\ComparisonCarrier}{\ComparisonMetric}\BlockPseudometric{\CoverIndex}{\ComparisonCarrier}{\ComparisonMetric}(\ComparisonPoint ,\IntegrationPoint ).
\end{equation}
Every block in \eqref{eq:joint-feature-map-1} is evaluated at the same point $\ComparisonPoint\in\ComparisonCarrier$.
Comparing the distances between the projections of two graph points identifies the distortion with the following uniform error.
\begin{equation}\label{eq:global-error}
  \ModelError\IsometryClass{\ComparisonCarrier }{\ComparisonMetric }=\norm{\ComparisonMetric -\GluedPseudometric{\ComparisonCarrier}{\ComparisonMetric} }_\infty.
\end{equation}
Since the weights sum to one, \eqref{eq:global-pseudometric-sum} implies
\begin{equation}\label{eq:global-error-bound}
 \ModelError\IsometryClass{\ComparisonCarrier}{\ComparisonMetric}
 \leq
 \sum_{\CoverIndex\in \ActiveBlocksAt{\ComparisonCarrier}{\ComparisonMetric}}\ModelWeightAt{\CoverIndex}{\ComparisonCarrier}{\ComparisonMetric}\BlockModelError{\CoverIndex}\IsometryClass{\ComparisonCarrier }{\ComparisonMetric }
 <
 \sum_{\CoverIndex\in \ActiveBlocksAt{\ComparisonCarrier}{\ComparisonMetric}}\ModelWeightAt{\CoverIndex}{\ComparisonCarrier}{\ComparisonMetric}
 \ErrorControl\IsometryClass{\ComparisonCarrier }{\ComparisonMetric }
=\ErrorControl\IsometryClass{\ComparisonCarrier }{\ComparisonMetric }.
\end{equation}
\ProofStep{Convergence and interpolation of the graphs.}\label{step:global-homotopy}
Suppose $\xi_j=\IsometryClass{Y_j}{e_j}\to\xi=\IsometryClass{Y}{e}$.
By \Cref{lem:common-gh-embedding} we may realize $Y_j,Y$ isometrically in a single compact metric space $\MetricSpace{Z}{h}$ so that $Y_j$ converges to $Y$ in Hausdorff distance.
All pairs $(x,y)\in Y_j\times Y$ satisfying $h(x,y)\leq\HausdorffDistance{h}(Y_j,Y)$ form a correspondence.
By \Cref{lem:weighted-coordinate-graphs} we can choose representatives so that $C_j=\FunctionGraph(\JointFeatureAt{Y_j}{e_j})$ converges to $C=\FunctionGraph(\JointFeatureAt{Y}{e})$ in Hausdorff distance in $Z\times\BanachAmbientSpace$.

Both the second projection and the orbit map are $1$-Lipschitz (\eqref{eq:orbit-metric}), so $C_j\to C$ implies $\ApproximationMap(\xi_j)\to\ApproximationMap(\xi)$.
The same graph convergence and the distortion estimate \eqref{eq:graph-distortion-bound} imply $\ModelError(\xi_j)\to\ModelError(\xi)$.

To describe the metric quotient $\GraphInterpolation(C,t)$, identify the graph with $Y$ by $y\mapsto(y,\JointFeatureAt{Y}{e}(y))$.
The interpolation is then the metric quotient for the pseudometric
\begin{equation}\label{eq:interpolating-pseudometrics-1}
 \HomotopyPseudometric{Y}{e}{t}=(1-t)e+t\GluedPseudometric{Y}{e}.
\end{equation}
This quotient is compact.
For $t<1$ the continuous pseudometric $\HomotopyPseudometric{Y}{e}{t}$ dominates $(1-t)e$ and hence induces the original topology of $Y$.

At $t=1$ the quotient identifies exactly the points with equal values under $\JointFeatureAt{Y}{e}$.
Thus the quotient at this endpoint is isometric to $\FeatureImageAt{Y}{e}$ with the restriction of $\BanachMetric$.

By \Cref{lem:graph-quotient-interpolation} the endpoint values are the isometry classes of the two projections.
Hence $\ApproximationHomotopy(\xi,0)=\xi$ and, by definition of $\RealizationMap$, $\ApproximationHomotopy(\xi,1)=\RealizationMap(\ApproximationMap(\xi))$.
These are \eqref{eq:global-domination-1} and \eqref{eq:global-domination-2}.

Finally, combining \eqref{eq:graph-quotient-bounds} and \eqref{eq:graph-quotient-time-bound} with the triangle inequality, we obtain for $s,t\in[0,1]$
\[
 \GHDistance\bigl(\GraphInterpolation(C_j,s),\GraphInterpolation(C,t)\bigr)
 \leq H_\times(C_j,C)+\tfrac12|s-t|\,\dis(C).
\]
Taking $s=t_j\to t$ makes the right-hand side tend to $0$.
Thus $\ApproximationHomotopy$ is jointly continuous in the space and time variables.

Taking $C_j=C$ and using $\dis(C)=\ModelError(\xi)$ yields the track estimate \eqref{eq:global-track-control}.
\end{proof}
\subsection{Absolute retracts}\label{subsec:ar-extension}
We now apply the domination and extension theorems from Subsection~\ref{subsec:ar-preliminaries-domination}.

\begin{theorem}\label{thm:absolute-extensor}
The Gromov--Hausdorff space $\GHSpace$ is an ANR for metrizable spaces.
Since $\GHSpace$ is contractible by \Cref{lem:scaling-contraction}, it is an AR in the same category.
\end{theorem}

\begin{proof}
Apply \Cref{thm:global-domination} with $\ErrorControl \equiv\ErrorTolerance $ for each $\ErrorTolerance >0$.
We obtain continuous maps $\ApproximationMap\colon\GHSpace\to\AuxiliaryAR$ and $\RealizationMap\colon\AuxiliaryAR\to\GHSpace$ and a homotopy $\ApproximationHomotopy$ from $\IdentityMap_{\GHSpace}$ to $\RealizationMap\circ\ApproximationMap$.
By \eqref{eq:global-track-control}, for every $\MetricSpace{\ComparisonCarrier}{\ComparisonMetric}\in\GHSpace$, we have
\[
 \diam_{\GHDistance}
 \{\ApproximationHomotopyAt{\ComparisonCarrier}{\ComparisonMetric}{\HomotopyTime}
   \mid\HomotopyTime\in\UnitInterval\}
 \leq\tfrac12\ModelError\IsometryClass{\ComparisonCarrier}{\ComparisonMetric}
 <\tfrac12\ErrorTolerance.
\]
By Hanner's domination theorem (\Cref{thm:hanner-domination}), $\GHSpace$ is an ANR for all separable metrizable spaces.

By Hanner's theorem (\Cref{thm:anr-category}), the space $\GHSpace$ is an ANR for all metrizable spaces.
Thus Dugundji's theorem (\Cref{thm:retract-extensor}) implies that $\GHSpace$ is an ANE for all metrizable spaces.

\Cref{lem:scaling-contraction} proves that $\GHSpace$ is contractible by scaling every distance to $0$.
Thus Hanner's theorem (\Cref{thm:contractible-ane}) shows that the contractible ANE $\GHSpace$ is an AE for all metrizable spaces.
Finally, using Dugundji's theorem (\Cref{thm:retract-extensor}), we conclude that $\GHSpace$ is an AR for all metrizable spaces.
\end{proof}

\begin{proof}[{Proof of \Cref{thm:part-iii-main}}]
This is a direct consequence of \Cref{thm:absolute-extensor}.
\end{proof}

\begin{theorem}\label{thm:invariant-measured-ar}\label{thm:invariant-measured-anr}\label{rem:invariant-measured-anr}
Equip the space $\InvariantMeasuredSpaces$ of isomorphism classes of compact metric spaces carrying an invariant probability measure of full support with the Gromov--Hausdorff--Prokhorov topology.
The forgetful map $\MeasureProjection\colon\InvariantMeasuredSpaces\to\GHSpace$ is soft, and $\InvariantMeasuredSpaces$ is an AR for metrizable spaces.
\end{theorem}

\begin{proof}
Let $T$ be metrizable, $A\subseteq T$ closed, and $f\colon A\to\InvariantMeasuredSpaces$ and $g\colon T\to\GHSpace$ continuous with $\MeasureProjection\circ f=g|_A$.
For every $a\in A$ we have $\MeasureProjection_*\DiracProbability_{f(a)}=\DiracProbability_{g(a)}$.
By \Cref{thm:dense-invariant-measured-spaces,prop:open-invariant-projection}, $\MeasureProjection$ is a continuous open surjection between Polish spaces.
The softness theorem of Banakh and Radul (\cite[Theorem~1.2]{BanakhRadul1999}, see also Valov, \cite[Theorem~1.1]{Valov2009}) yields a continuous map $\widetilde g\colon T\to\ProbabilityMeasures(\InvariantMeasuredSpaces)$ such that
\[
 \MeasureProjection_*\widetilde g(t)=\DiracProbability_{g(t)}\quad(t\in T),\qquad
 \widetilde g(a)=\DiracProbability_{f(a)}\quad(a\in A).
\]
The first identity means that $\widetilde g(t)$ is concentrated on the fiber over $g(t)$.
By \Cref{cor:continuous-fiber-barycenter} the map $h=\FiberBarycenterMap\circ\widetilde g$ is continuous and satisfies $\MeasureProjection\circ h=g$.
The barycenter of a Dirac law is its point of concentration, so the second identity implies $h|_A=f$.
This proves softness of $\MeasureProjection$.

Now take any continuous map $f\colon A\to\InvariantMeasuredSpaces$ on a closed subset $A\subseteq T$.
By \Cref{thm:absolute-extensor} and Dugundji's theorem (\Cref{thm:retract-extensor}), $\GHSpace$ is an AE, so $\MeasureProjection\circ f$ extends to a continuous map $g$ on $T$.
The softness just proved lifts $g$ to an extension of $f$ over $T$.
Thus $\InvariantMeasuredSpaces$ is an AE and hence an AR by Dugundji's theorem.
\end{proof}

\clearpage
\part{The topology of Gromov--Hausdorff space IV: {The Gromov--Hausdorff space is homeomorphic to \texorpdfstring{$\SeparableHilbertSpace$}{ℓ²}}}\label{part:hilbert}
\begin{quote}
\small
\noindent\textbf{Abstract.}
Combining the discrete approximation theorem with the absolute retract property from Part~\ref{part:topology}, completeness, separability, and Toru\'nczyk's characterization, we prove that the Gromov--Hausdorff space is homeomorphic to separable infinite-dimensional Hilbert space.
We also construct an explicit strong deformation retraction onto the subspace of spaces of diameter one and show that this subspace is homeomorphic to the same Hilbert space.

\par\smallskip
\noindent \textbf{Keywords.} Gromov--Hausdorff space, Hilbert space, discrete approximation property, metric product.

\par\smallskip
\noindent \textbf{2020 Mathematics Subject Classification.} Primary 57N20; Secondary 54E35, 54F45.
\end{quote}

\section{Introduction to Part IV}\label{sec:intro-hilbert}
In Part~\ref{part:topology} we proved that the Gromov--Hausdorff space is an absolute retract.
Here we establish the discrete approximation property needed to identify its homeomorphism type.
Write $\GHSpace$ for the space of isometry classes of nonempty compact metric spaces with the usual unpointed Gromov--Hausdorff distance $\GHDistance$.
Antonyan asked whether $\GHSpace$ is homeomorphic to separable infinite-dimensional Hilbert space (\cite[p.~2]{Antonyan2020Euclidean}).
For each positive integer $\FiniteSize $, Antonyan proved that the subspace represented by compact subsets of $\RealNumbers^\FiniteSize $ is homeomorphic to the Hilbert cube with one point removed (\cite[Corollary~5.3]{Antonyan2021Euclidean}).

Earlier work studied the topology of $\GHSpace$ through explicit families and embeddings.
The author constructed families of branching Gromov--Hausdorff geodesics continuously parametrized by the Hilbert cube (\cite[Theorem~1.3]{Ishiki2022Branching}).
The author also embedded each compact metrizable space into the subspaces of $\GHSpace$ consisting of continua (\cite[Theorem~1.1]{Ishiki2022Continua}), of spaces of prescribed dimensions (\cite[Theorem~1.3]{Ishiki2023FractalDimensions}), and of compact metric trees (\cite[Theorem~1.1]{Ishiki2023MetricTrees}).
Each construction allows one to prescribe finitely many distinct values of the embedding.
Byakuno proved that every countable product of closed intervals of positive summable lengths, with the supremum metric, embeds isometrically into $\GHSpace$ (\cite[arXiv v1, Corollary~1.4]{Byakuno2026Embedding}).

For finite metric spaces, Nakajima, Yamauchi, and Zava proved that the subspace of spaces with at most $\FiniteSize$ points has topological dimension $\FiniteSize(\FiniteSize-1)/2$ (\cite[arXiv v1, Theorem~A]{NakajimaYamauchiZava2025}).
Nakajima and Shioya constructed a compact metrizable space containing a dense topological copy of $\GHSpace$ (\cite[arXiv v1, Main Theorem~1.2]{NakajimaShioyaCompactification}).

A related classification concerns the Urysohn universal metric space.
Uspenskij proved that it is homeomorphic to $\SeparableHilbertSpace$ (\cite[arXiv v1, Theorem~2.1]{Uspenskij2004Urysohn}).

Our approximation theorem concerns a countable family of prescribed maps from compact metrizable domains.
We approximate all these maps so that every point of $\GHSpace$ has a neighborhood meeting at most one of their images.
The construction uses products with finite equilateral spaces, independently of Parts~\ref{part:measures}--\ref{part:topology}.
We combine this discrete approximation with the absolute retract property to prove the following theorem.

\begin{theorem}\label{thm:part-iv-main}
The space $\GHSpace$ is homeomorphic to the real Hilbert space $\SeparableHilbertSpace$.
\end{theorem}

This classification concerns the underlying topology.
The Gromov--Hausdorff metric carries further geometric information.
Ivanov and Tuzhilin proved that every surjective isometry of $\GHSpace$ onto itself is the identity (\cite[Main Theorem]{IvanovTuzhilin2019Isometries}).
Zava proved that the subspace of finite metric spaces admits no coarse embedding into any Hilbert space (\cite[arXiv v3, Theorem~B]{Zava2025Coarse}).
Consequently the homeomorphism in \Cref{thm:part-iv-main} cannot be bi-Lipschitz.

Let $\{\ApproximationDomain_\FamilyIndex\}_{\FamilyIndex\in\NonnegativeIntegers}$ be a sequence of compact metrizable spaces, let $\InputMap_\FamilyIndex\colon\ApproximationDomain_\FamilyIndex\to\GHSpace$ be continuous, and let $\OpenCover$ be an open cover of $\GHSpace$.
We construct continuous maps $\ApproximatingMap_\FamilyIndex\colon\ApproximationDomain_\FamilyIndex\to\GHSpace$ that are $\OpenCover$-close to $\InputMap_\FamilyIndex$ and whose images $\{\ApproximatingMap_\FamilyIndex(\ApproximationDomain_\FamilyIndex)\}_{\FamilyIndex\in\NonnegativeIntegers}$ form a discrete family in $\GHSpace$.
The first condition means that for every $\FamilyIndex$ and $\DomainPoint\in\ApproximationDomain_\FamilyIndex$ some member $U\in\OpenCover$ contains both $\InputMap_\FamilyIndex(\DomainPoint)$ and $\ApproximatingMap_\FamilyIndex(\DomainPoint)$.
The second means that every point of $\GHSpace$ has an open neighborhood $V$ for which $V\cap\ApproximatingMap_\FamilyIndex(\ApproximationDomain_\FamilyIndex)\ne\varnothing$ holds for at most one index $\FamilyIndex$.
This is \Cref{thm:discrete-approximation}.
Taking $\ApproximationDomain_\FamilyIndex=\HilbertCube
 =\UnitInterval^{\NonnegativeIntegers}$ for every $\FamilyIndex$ yields the approximation condition in Toru\'nczyk's characterization, stated in \Cref{thm:torunczyk} (Toru\'nczyk, \cite[p.~248, assertion~(i)]{Torunczyk1981} and \cite[Section~C]{Torunczyk1985}).
Together with completeness, separability, and the absolute retract property, this proves \Cref{thm:part-iv-main}.

The main theorem also implies that every point has a basis of open neighborhoods that strongly deformation retract onto that point.
We do not know whether these neighborhoods can be chosen to be balls for $\GHDistance$.

\medskip
\noindent \textbf{Dependence on the other parts.}
\Cref{thm:discrete-approximation} uses the common preliminaries and the metric constructions recalled in Section~\ref{sec:prelim-hilbert}.
The absolute retract result \Cref{thm:part-iii-main} enters only when we apply Toru\'nczyk's characterization in Subsection~\ref{subsec:approximation-hilbert-recognition}.
Thus the proof of the classification combines the approximation in this part with the result of Parts~\ref{part:measures}--\ref{part:topology}.

\medskip
\noindent \textbf{Organization.}
In Section~\ref{sec:prelim-hilbert}, we define discrete families and state the Hilbert space characterization.
In Subsection~\ref{subsec:approximation-products}, we recall the Gromov--Hausdorff estimate for metric products.
In Subsection~\ref{subsec:approximation-discrete}, we construct approximations of maps from countably many compact domains whose images form a discrete family, proving \Cref{thm:discrete-approximation}.
In Subsection~\ref{subsec:approximation-hilbert-recognition}, we apply Toru\'nczyk's characterization to prove \Cref{thm:part-iv-main} and obtain strong local contractions.
In \Cref{sec:questions}, we conclude with questions about related spaces and the geometry of Gromov--Hausdorff balls.

\medskip
\noindent \textbf{Conventions and notation.}
All compact metric spaces are nonempty, and Banach and Hilbert spaces are real.
We identify compact metric spaces with their isometry classes in $\GHSpace$.
An isometry is a surjective isometric embedding.
For a metric space $\MetricSpace{\BaseCarrier}{\MetricSymbol}$, we write $\CompactHyperspace(\BaseCarrier)$ for its nonempty compact subsets, equipped with the Hausdorff metric $\HausdorffDistance{\MetricSymbol}$.
For a compact space with an isometric embedding into an ambient metric space, we identify the compact space with its image and use the restricted ambient metric.
The notation $\IdentityMap$ denotes the identity map on its specified domain.
Sequences are indexed by $\NonnegativeIntegers=\{0,1,2,\ldots\}$.

\section{Preliminaries}\label{sec:prelim-hilbert}
We define discrete families and approximation with respect to an open cover.
We recall metric products and state Toru\'nczyk's characterization of Hilbert space.
\subsection{Discrete families and close maps}\label{subsec:hilbert-preliminaries-discrete}
To state the approximation condition in Toru\'nczyk's characterization (\Cref{thm:torunczyk}), we recall discreteness of a family of subsets and closeness of maps with respect to an open cover.

\begin{definition}\label{def:discrete-family}
A family $\{\ClosedDomain_\FamilyIndex\}_{\FamilyIndex\in\DomainIndexSet}$ of subsets of a space $\RecognitionTarget$ is \emph{discrete in $\RecognitionTarget$} if each $x\in\RecognitionTarget$ has an open neighborhood $U$ such that $U\cap\ClosedDomain_\FamilyIndex\ne\varnothing$ for at most one index $\FamilyIndex\in\DomainIndexSet$.
The family is \emph{locally finite} if each $x\in\RecognitionTarget$ has an open neighborhood meeting only finitely many members.
Both conditions concern all points of the space, including those outside $\bigcup_{\FamilyIndex\in\DomainIndexSet}\ClosedDomain_\FamilyIndex$.
For an open cover $\OpenCover$ of $\RecognitionTarget $, two maps $\InputMap ,\ApproximatingMap \colon \CommonMapDomain \to \RecognitionTarget $ are \emph{$\OpenCover$-close} if for each $\DomainPoint \in \CommonMapDomain $ some $U\in\OpenCover$ satisfies
\[
 \{\InputMap(\DomainPoint),\ApproximatingMap(\DomainPoint)\}\subset U.
\]
\end{definition}
\subsection{Metric products}\label{subsec:approximation-products}
We recall the Gromov--Hausdorff estimate for products used in the approximation.
For compact metric spaces $\MetricSpace{\FirstFactor }{\FirstFactorMetric }$ and $\MetricSpace{\SecondFactor }{\SecondFactorMetric }$, denote the \emph{maximum metric} on $\FirstFactor \times \SecondFactor $ by $\FirstFactorMetric \vee \SecondFactorMetric $.
For $\FirstFactorPoint,\OtherFirstFactorPoint\in\FirstFactor$ and $\SecondFactorPoint,\OtherSecondFactorPoint\in\SecondFactor$ put
\[
  (\FirstFactorMetric \vee \SecondFactorMetric )((\FirstFactorPoint ,\SecondFactorPoint ),(\OtherFirstFactorPoint ,\OtherSecondFactorPoint ))
  =\max\{\FirstFactorMetric (\FirstFactorPoint ,\OtherFirstFactorPoint ),\SecondFactorMetric (\SecondFactorPoint ,\OtherSecondFactorPoint )\}.
\]
\begin{lemma}\label{lem:metric-products}
Let $\MetricSpace{\FirstFactor }{\FirstFactorMetric }$, $\MetricSpace{\SecondFactor }{\SecondFactorMetric }$, $\MetricSpace{\ComparisonFirstFactor }{\ComparisonFirstFactorMetric }$, and $\MetricSpace{\ComparisonSecondFactor }{\ComparisonSecondFactorMetric }$ be nonempty compact metric spaces.
Then
\begin{equation}\label{eq:max-product}
 \begin{aligned}
  &\GHDistance(\MetricSpace{\FirstFactor \times \SecondFactor }{\FirstFactorMetric \vee \SecondFactorMetric },
        \MetricSpace{\ComparisonFirstFactor \times \ComparisonSecondFactor }{\ComparisonFirstFactorMetric \vee \ComparisonSecondFactorMetric })\\
  &\qquad\leq\max\{\GHDistance(\MetricSpace{\FirstFactor }{\FirstFactorMetric },\MetricSpace{\ComparisonFirstFactor }{\ComparisonFirstFactorMetric }),
                       \GHDistance(\MetricSpace{\SecondFactor }{\SecondFactorMetric },\MetricSpace{\ComparisonSecondFactor }{\ComparisonSecondFactorMetric })\}.
 \end{aligned}
\end{equation}
\end{lemma}

\begin{proof}
Let $\Correspondence_\FirstFactor\subset
 \FirstFactor\times\ComparisonFirstFactor$ and $\Correspondence_\SecondFactor\subset
 \SecondFactor\times\ComparisonSecondFactor$ be correspondences.
By rearranging coordinates, we regard their product $\Correspondence_\times
 =\Correspondence_\FirstFactor\times\Correspondence_\SecondFactor$ as a correspondence between $\FirstFactor\times\SecondFactor$ and $\ComparisonFirstFactor\times\ComparisonSecondFactor$.
For real numbers $a,b,c,d$, we have
\[
 |\max\{a,b\}-\max\{c,d\}|
 \leq\max\{|a-c|,|b-d|\}.
\]
Hence $\dis\Correspondence_\times
 \leq\max\{\dis\Correspondence_\FirstFactor,
            \dis\Correspondence_\SecondFactor\}$.
Taking the infimum over the two factor correspondences proves \eqref{eq:max-product} by \eqref{eq:correspondence}.
\end{proof}

\subsection{A characterization of Hilbert space}\label{subsec:hilbert-preliminaries-recognition}
We use the following characterization of Hilbert space, due to Toru\'nczyk (\cite[p.~248, condition~(i)]{Torunczyk1981}) with the correction in \cite[Section~C]{Torunczyk1985}.

\begin{theorem}[{Toru\'nczyk, \cite[p.~248, condition~(i)]{Torunczyk1981}, \cite[Section~C]{Torunczyk1985}}]\label{thm:torunczyk}
Let $\MetricSpace{\RecognitionTarget }{\TargetMetric }$ be a nonempty complete separable metric AR, and let $\HilbertCube =\UnitInterval^{\NonnegativeIntegers}$ be the Hilbert cube.
Equip $\NonnegativeIntegers$ with the discrete topology.
Then $\RecognitionTarget $ is homeomorphic to the real Hilbert space $\SeparableHilbertSpace$ if and only if, for every continuous map $\InputMap \colon \NonnegativeIntegers\times \HilbertCube \to \RecognitionTarget $ and every open cover $\OpenCover$ of $\RecognitionTarget $, there exists a continuous map $\ApproximatingMap \colon \NonnegativeIntegers\times \HilbertCube \to \RecognitionTarget $ that is $\OpenCover$-close to $\InputMap$ and whose image family $\{\ApproximatingMap (\{\FamilyIndex \}\times \HilbertCube )\}_{\FamilyIndex \in\NonnegativeIntegers}$ is discrete in $\RecognitionTarget $.
\end{theorem}

\section{Discrete approximation and Hilbert space topology}\label{sec:approximation}
To apply Toru\'nczyk's characterization of Hilbert space, we prove the discrete approximation property for maps from countably many Hilbert cubes.
\subsection{Discrete approximation for compact domains}\label{subsec:approximation-discrete}
We now construct the approximating maps and prove discreteness at every point of the ambient Gromov--Hausdorff space.

We first choose a positive Lipschitz scale adapted to an open cover.

\begin{lemma}\label{lem:cover-adapted-scale}
Let $\OpenCover$ be an open cover of $\GHSpace$.
Then there exists a $1/16$-Lipschitz function $\CoverScale\colon\GHSpace\to(0,1/16]$ with the following property.
For every $\MetricSpace{\BaseCarrier}{\MetricSymbol},
 \MetricSpace{\ComparisonCarrier}{\ComparisonMetric}\in\GHSpace$, if
\[
 \GHDistance(\MetricSpace{\BaseCarrier}{\MetricSymbol},
             \MetricSpace{\ComparisonCarrier}{\ComparisonMetric})
 <{\CoverScaleAt{\BaseCarrier}{\MetricSymbol}},
\]
then there exists $U\in\OpenCover$ with
\[
 \{\MetricSpace{\BaseCarrier}{\MetricSymbol},
   \MetricSpace{\ComparisonCarrier}{\ComparisonMetric}\}\subset U.
\]
\end{lemma}

\begin{proof}
Use the convention that distance to the empty set is $\infty$.
For every $\MetricSpace{\BaseCarrier}{\MetricSymbol}\in\GHSpace$, put
\begin{equation}\label{eq:cover-scale-2}
 \CoverScaleAt{\BaseCarrier}{\MetricSymbol}
 =\frac1{16}\sup_{\OpenSubset\in\OpenCover}
   \min\left\{1,\dist_{\GHDistance}
       \bigl(\MetricSpace{\BaseCarrier}{\MetricSymbol},\GHSpace\setminus\OpenSubset\bigr)\right\}.
\end{equation}
Each truncated distance function is $1$-Lipschitz, and their supremum is positive because $\OpenCover$ is an open cover.
Thus, for every $\MetricSpace{\BaseCarrier}{\MetricSymbol},
 \MetricSpace{\ComparisonCarrier}{\ComparisonMetric}\in\GHSpace$, we have
\begin{equation}\label{eq:scale-lipschitz-1}
 0<\CoverScaleAt{\BaseCarrier}{\MetricSymbol}\leq\frac1{16},
\end{equation}
and
\begin{equation}\label{eq:scale-lipschitz-2}
 \abs{\CoverScaleAt{\BaseCarrier}{\MetricSymbol}
      -\CoverScaleAt{\ComparisonCarrier}{\ComparisonMetric}}
 \leq\frac1{16}\GHDistance(\MetricSpace{\BaseCarrier}{\MetricSymbol},
                          \MetricSpace{\ComparisonCarrier}{\ComparisonMetric}).
\end{equation}
Now let $\MetricSpace{\BaseCarrier}{\MetricSymbol},
 \MetricSpace{\ComparisonCarrier}{\ComparisonMetric}\in\GHSpace$ satisfy $\GHDistance(\MetricSpace{\BaseCarrier}{\MetricSymbol},
             \MetricSpace{\ComparisonCarrier}{\ComparisonMetric})
 <\CoverScaleAt{\BaseCarrier}{\MetricSymbol}$.
By \eqref{eq:cover-scale-2}, there exists $\mathscr O\in\OpenCover$ such that
\[
 \GHDistance(\MetricSpace{\BaseCarrier}{\MetricSymbol},
             \MetricSpace{\ComparisonCarrier}{\ComparisonMetric})
 <\CoverScaleAt{\BaseCarrier}{\MetricSymbol}
 <\dist_{\GHDistance}
   \bigl(\MetricSpace{\BaseCarrier}{\MetricSymbol},\GHSpace\setminus\mathscr O\bigr).
\]
Both spaces therefore belong to $\mathscr O$.
\end{proof}

\begin{theorem}\label{thm:discrete-approximation}
Let $\{\ApproximationDomain _\FamilyIndex \}_{\FamilyIndex \in\NonnegativeIntegers}$ be a sequence of compact metrizable spaces.
For each $\FamilyIndex\in\NonnegativeIntegers$, let $\InputMap _\FamilyIndex \colon \ApproximationDomain _\FamilyIndex \to\GHSpace$ be continuous.
Let $\OpenCover$ be an open cover of $\GHSpace$.
Then there exist continuous maps $\ApproximatingMap _\FamilyIndex \colon \ApproximationDomain _\FamilyIndex \to\GHSpace$ such that each $\ApproximatingMap _\FamilyIndex $ is $\OpenCover$-close to $\InputMap _\FamilyIndex $ and the family $\{\ApproximatingMap _\FamilyIndex (\ApproximationDomain _\FamilyIndex )\}_{\FamilyIndex \in\NonnegativeIntegers}$ is discrete in $\GHSpace$.
\end{theorem}

\begin{proof}
Put $\DomainIndexSet=\{i\in\NonnegativeIntegers\mid\ApproximationDomain_i\ne\emptyset\}$.
For $i\notin\DomainIndexSet$ let $\ApproximatingMap_i$ be the empty map.
For each $i\in\DomainIndexSet$ we construct a continuous map $\ApproximatingMap_i$ by taking the product of each value of $\InputMap_i$ with a finite equilateral space.

\ProofStep{Product approximations.}\label{step:discrete-products}
By \Cref{lem:cover-adapted-scale}, choose $\CoverScale\colon\GHSpace\to(0,1/16]$ satisfying \eqref{eq:scale-lipschitz-2} and the cover condition in that lemma.

For an integer $\EquilateralSize\geq2$ put $\EquilateralCarrier_{\EquilateralSize}=\{0,\ldots,\EquilateralSize-1\}$.
For $\ProductScale>0$ denote by $\EquilateralMetric_{\EquilateralSize,\ProductScale}$ the metric on the \emph{equilateral space} $\EquilateralCarrier_{\EquilateralSize}$ whose off-diagonal distances are $\ProductScale$.
For each integer $\EquilateralSize\geq2$ define $\ProductApproximation_{\EquilateralSize}\colon\GHSpace\to\GHSpace$ as follows.
For $\MetricSpace{\BaseCarrier}{\MetricSymbol}\in\GHSpace$ set
\begin{equation}\label{eq:product-approximation-2}
 \ProductScale=\CoverScaleAt{\BaseCarrier}{\MetricSymbol}
\end{equation}
and
\begin{equation}\label{eq:product-approximation-1}
 \ProductApproximation_{\EquilateralSize}\MetricSpace{\BaseCarrier}{\MetricSymbol}
 =\MetricSpace{\BaseCarrier\times\EquilateralCarrier_{\EquilateralSize}}
              {\MetricSymbol\vee\EquilateralMetric_{\EquilateralSize,\ProductScale}}.
\end{equation}
The graph of the projection onto $\BaseCarrier$ is a correspondence of distortion at most $\ProductScale$.
Hence we have
\begin{equation}\label{eq:product-error}
 \GHDistance(\ProductApproximation_{\EquilateralSize}\MetricSpace{\BaseCarrier}{\MetricSymbol},
             \MetricSpace{\BaseCarrier}{\MetricSymbol})
 \leq\frac12\CoverScaleAt{\BaseCarrier}{\MetricSymbol}.
\end{equation}
For an integer $\EquilateralSize\geq2$ and spaces $\MetricSpace{\BaseCarrier}{\MetricSymbol},
 \MetricSpace{\ComparisonCarrier}{\ComparisonMetric}\in\GHSpace$, pair points with the same label in the equilateral factors.
As in the author's proof of \cite[Proposition~5.3]{Ishiki2022Branching}, the product correspondence and \eqref{eq:scale-lipschitz-2} yield
\begin{equation}\label{eq:equilateral-product-lipschitz}
 \begin{aligned}
 &\GHDistance(\ProductApproximation_{\EquilateralSize}\MetricSpace{\BaseCarrier}{\MetricSymbol},
              \ProductApproximation_{\EquilateralSize}\MetricSpace{\ComparisonCarrier}{\ComparisonMetric})\\
 &\quad\leq\max\left\{
     \GHDistance(\MetricSpace{\BaseCarrier}{\MetricSymbol},
                 \MetricSpace{\ComparisonCarrier}{\ComparisonMetric}),
     \frac12\left|\CoverScaleAt{\BaseCarrier}{\MetricSymbol}
                 -\CoverScaleAt{\ComparisonCarrier}{\ComparisonMetric}\right|
                  \right\}\\
 &\quad\leq\GHDistance(\MetricSpace{\BaseCarrier}{\MetricSymbol},
                      \MetricSpace{\ComparisonCarrier}{\ComparisonMetric}).
 \end{aligned}
\end{equation}
Thus every map $\ProductApproximation_{\EquilateralSize}$ is $1$-Lipschitz.

\ProofStep{Separation and local finiteness of the images.}\label{step:discrete-separation}

We use $\PackingSeparation{m}$ to separate each new image from all earlier images and to show that the family of images is locally finite at every point of $\GHSpace$.

Fix an integer $N\geq2$ and $\xi=\IsometryClass{X}{d}$.
Over each $x\in X$ the product $\ProductApproximation_N(\xi)$ contains an equilateral space of $N$ points and side length $\CoverScale(\xi)$.
Thus
\[
 \PackingSeparation{N}(\ProductApproximation_N(\xi))\geq\CoverScale(\xi).
\]
Take the indices $i\in\DomainIndexSet$ in increasing order and define $m_i$ and $\ApproximatingMap_i$ recursively.
Let the union of the previously chosen images be
\begin{equation}\label{eq:previous-approximation-images}
 \PreviousApproximationImages_{\FamilyIndex}
 =\bigcup_{\substack{k\in\DomainIndexSet\\k<\FamilyIndex}}
       \ApproximatingMap_k(\ApproximationDomain_k).
\end{equation}
This is a finite union of compact sets, empty for the least index in $\DomainIndexSet$.
Since $\ApproximationDomain_{\FamilyIndex}$ is nonempty and compact, the number
\begin{equation}\label{eq:minimum-product-scale}
 \ProductScale_{\FamilyIndex}
 =\min_{\DomainPoint\in\ApproximationDomain_{\FamilyIndex}}
       \CoverScale(\InputMap_{\FamilyIndex}(\DomainPoint))
\end{equation}
is positive.

By \Cref{cor:uniform-separated-sets}, the functions $\PackingSeparation{m}$ tend uniformly to zero on the compact set $\PreviousApproximationImages_i$.
We may therefore choose the number of points as follows.

Choose an integer $m_i$ larger than every previously chosen $m_k$ and satisfying
\begin{equation}\label{eq:equilateral-cardinality-choice}
 {m_i>i+1,\quad {}\PackingSeparation{m_i}(\zeta)<\ProductScale_i/2}
 \quad(\zeta\in\PreviousApproximationImages_i).
\end{equation}
This is possible also when $\PreviousApproximationImages_i=\emptyset$.
Define $\ApproximatingMap_i:\ApproximationDomain_i\to\GHSpace$ by
\begin{equation}\label{eq:discrete-maps}
 \ApproximatingMap_{\FamilyIndex}
 =\ProductApproximation_{\EquilateralSize_{\FamilyIndex}}\circ\InputMap_{\FamilyIndex}.
\end{equation}
The map $\ApproximatingMap_{\FamilyIndex}$ is continuous by \eqref{eq:equilateral-product-lipschitz}, and its image is compact.
For every $\DomainPoint\in\ApproximationDomain_{\FamilyIndex}$, equation \eqref{eq:product-error} implies
\begin{equation}\label{eq:approximation-error}
 \GHDistance(\ApproximatingMap_{\FamilyIndex}(\DomainPoint),
             \InputMap_{\FamilyIndex}(\DomainPoint))
 \leq\frac12\CoverScale(\InputMap_{\FamilyIndex}(\DomainPoint))
 <\CoverScale(\InputMap_{\FamilyIndex}(\DomainPoint)).
\end{equation}
Thus $\ApproximatingMap_{\FamilyIndex}$ and $\InputMap_{\FamilyIndex}$ are $\OpenCover$-close.

First estimate the distance from the new image to the earlier images.
For every $z\in\ApproximationDomain_i$ we have $\PackingSeparation{m_i}(\ApproximatingMap_i(z))\geq\CoverScale(\InputMap_i(z))\geq\ProductScale_i$.
For each $\zeta\in\PreviousApproximationImages_i$ the $2$-Lipschitz property of $\PackingSeparation{m_i}$ (\Cref{cor:uniform-separated-sets}) and \eqref{eq:equilateral-cardinality-choice} imply
\begin{equation}\label{eq:recursive-image-separation}
 2\GHDistance(\ApproximatingMap_i(z),\zeta)
 \geq\PackingSeparation{m_i}(\ApproximatingMap_i(z))-\PackingSeparation{m_i}(\zeta)
 >\ProductScale_i/2.
\end{equation}
Thus the new image is at distance at least $\ProductScale_i/4$ from every earlier image.

To prove local finiteness, suppose that every neighborhood of some $\zeta\in\GHSpace$ meets infinitely many images.
For each $j\in\NonnegativeIntegers$ choose an image meeting the open ball of radius $2^{-j}$ centered at $\zeta$, with strictly increasing indices $i_j$.
Choose $\zeta_j=\ApproximatingMap_{i_j}(z_j)$ in this intersection with $z_j\in\ApproximationDomain_{i_j}$.
Then $\zeta_j\to\zeta$.
Put $\xi_j=\InputMap_{i_j}(z_j)$ and $n_j=m_{i_j}$, so $\zeta_j=\ProductApproximation_{n_j}(\xi_j)$ and $n_j\to\infty$.
By the lower bound from the product fibers and the $2$-Lipschitz property of $\PackingSeparation{n_j}$,
\[
 0<\CoverScale(\xi_j)
 \leq\PackingSeparation{n_j}(\zeta_j)
 \leq\PackingSeparation{n_j}(\zeta)+2\GHDistance(\zeta_j,\zeta)
 \longrightarrow0.
\]
Here $\PackingSeparation{n_j}(\zeta)\to0$ by \Cref{cor:uniform-separated-sets}.
For the original spaces $\xi_j$, \eqref{eq:product-error} implies
\[
 \GHDistance(\xi_j,\zeta)
 \leq\tfrac12\CoverScale(\xi_j)+\GHDistance(\zeta_j,\zeta)\longrightarrow0.
\]
Continuity of $\CoverScale$ now implies $\CoverScale(\xi_j)\to\CoverScale(\zeta)>0$, a contradiction.

Each image is compact and hence closed.
At each point choose an open neighborhood meeting only finitely many images.
Removing those among them that omit the point yields an open neighborhood meeting at most one image.
Thus the family is discrete.

Figure~\ref{fig:discrete-approximation} illustrates the product approximation and the estimates for separation and local finiteness.
\begin{figure}[H]
\centering
\begin{tikzpicture}[>=Stealth,font=\small]
\node at (0,1.6) {$\BaseCarrier$};
\draw (0,0) ellipse (1.15 and 0.4);
\fill (-0.55,0.04) circle (1.1pt);
\fill (0.5,-0.12) circle (1.1pt);
\fill (0.05,0.08) circle (1.5pt);
\node[above,inner sep=3pt] at (0.05,0.08) {$\BasePoint$};
\draw[->] (1.5,0) -- node[above]
 {$\times\EquilateralCarrier_3$} (3.1,0);
\node at (4.65,1.6) {$\BaseCarrier\times\EquilateralCarrier_3$};
\fill[black!8] (4.5,-1.15) rectangle (4.9,1.3);
\foreach \level/\index in {0.95/0,0/1,-0.95/2}
{
 \draw (4.65,\level) ellipse (1.15 and 0.4);
 \fill (4.1,{\level+0.04}) circle (1.1pt);
 \fill (5.15,{\level-0.12}) circle (1.1pt);
 \fill (4.7,{\level+0.08}) circle (1.5pt);
 \node[right] at (5.95,\level)
  {$\BaseCarrier\times\{\index\}$};
}
\node[align=center] at (4.7,-1.95)
 {$\{\BasePoint\}\times\EquilateralCarrier_3$\\
  pairwise distance $\ProductScale$};
\end{tikzpicture}
\par\medskip
\begin{tikzpicture}[
 >=Stealth,
 font=\small,
 role/.style={draw,align=center,text width=3.4cm,
              minimum height=2.1cm,inner sep=5pt,anchor=north}
]
\node at (4,0.9)
 {$\MetricSpace{\ClosedDomain}{\SubspaceMetric}
   =\ProductApproximation_{\EquilateralSize_{\FamilyIndex}}
      \MetricSpace{\BaseCarrier}{\MetricSymbol},
   \qquad
   \ProductScale=\CoverScaleAt{\BaseCarrier}{\MetricSymbol}$};
\node (formula) at (4,0)
 {$\MetricSpace{\ClosedDomain}{\SubspaceMetric}
   =\MetricSpace{\BaseCarrier\times\EquilateralCarrier_{\EquilateralSize_{\FamilyIndex}}}
                {\MetricSymbol\vee\EquilateralMetric_{\EquilateralSize_{\FamilyIndex},\ProductScale}}$};
\node[role] (approximation) at (0,-1.7)
 {small scale $\ProductScale$\\
  $\begin{gathered}
    \GHDistance(\MetricSpace{\ClosedDomain}{\SubspaceMetric},
                 \MetricSpace{\BaseCarrier}{\MetricSymbol})\\
    \leq\ProductScale/2
   \end{gathered}$\\
  approximation};
\node[role] (separation) at (4,-1.7)
 {{\par \shortstack{$\PackingSeparation{m_i}<r_i/2$ (earlier)\\$\PackingSeparation{m_i}\geq r_i$ (new)\\distance $\ge r_i/4$}}};
\node[role] (packing) at (8,-1.7)
 {{{\par \shortstack{from distinct images\\$\zeta_j\to\zeta$\\$\CoverScale(\xi_j)\to0,\ \xi_j\to\zeta$\\contradicts $\CoverScale(\zeta)>0$}}}};
\draw[->,shorten <=2pt]
 ([xshift=-1.8cm]formula.south) .. controls +(0,-0.55) and +(0,0.55) .. (approximation.north);
\draw[->,shorten <=2pt]
 (formula.south) -- (separation.north);
\draw[->,shorten <=2pt]
 ([xshift=1.8cm]formula.south) .. controls +(0,-0.55) and +(0,0.55) .. (packing.north);
\end{tikzpicture}
\caption{The upper diagram shows three copies of $\BaseCarrier$ in a product with an equilateral space. For a general factor $\EquilateralCarrier_\EquilateralSize$, the indicated fiber consists of $\EquilateralSize$ points at pairwise distance $\ProductScale$. Below, the left panel bounds the approximation error. The middle panel compares $\PackingSeparation{m_i}$ on the earlier images and on the new image to give \eqref{eq:recursive-image-separation}. In the right panel, suppose that points $\zeta_j=\ApproximatingMap_{i_j}(z_j)$ from distinct images converge to $\zeta$. The original points $\xi_j=\InputMap_{i_j}(z_j)$ satisfy $\CoverScale(\xi_j)\leq\PackingSeparation{n_j}(\zeta_j)\to0$ and $\GHDistance(\xi_j,\zeta_j)\leq\CoverScale(\xi_j)/2$. Thus $\xi_j\to\zeta$, contradicting the continuity and positivity of $\CoverScale$.}
\label{fig:discrete-approximation}
\end{figure}
\end{proof}
\subsection{Proof of Main Theorem}\label{subsec:approximation-hilbert-recognition}
We apply Toru\'nczyk's characterization (\Cref{thm:torunczyk}) to $\GHSpace$ using the AR property from Part~\ref{part:topology} and the discrete approximation in \Cref{thm:discrete-approximation}.

\begin{proof}[{Proof of \Cref{thm:part-iv-main}}]
By Edwards's theorem (\Cref{lem:gh-basics}), $\GHSpace$ is complete and separable.
By \Cref{thm:part-iii-main} it is a nonempty AR.
Apply \Cref{thm:discrete-approximation} with $\ApproximationDomain_\FamilyIndex=\HilbertCube$ for every $\FamilyIndex\in\NonnegativeIntegers$ to obtain the approximation condition in Toru\'nczyk's characterization (\Cref{thm:torunczyk}).
It follows from \Cref{thm:torunczyk} that $\GHSpace$ is homeomorphic to $\SeparableHilbertSpace$.
\end{proof}

The proof of \Cref{thm:part-iv-main} also applies to the subspace of spaces of diameter one.

\begin{corollary}\label{cor:unit-diameter-hilbert-space}
The following subspace of $\GHSpace$
\[
 \UnitDiameterGHSpace
 =\{{\IsometryClass{\BaseCarrier}{\MetricSymbol}}\in\GHSpace
     \mid\diam\MetricSpace{\BaseCarrier}{\MetricSymbol}=1\}
\]
is a strong deformation retract of $\GHSpace$ and is homeomorphic to the real Hilbert space $\SeparableHilbertSpace$.
\end{corollary}

\begin{proof}
Fix a two-point space $B$ of diameter $1$.
For $\xi=\IsometryClass{X}{d}$ put $D=\diam(X,d)$, $\beta(D)=\max\{1-D,0\}$, and $b(D)=\max\{D,\beta(D)\}$.
For $c\geq0$ denote by $c\xi$ the isometry class of the metric quotient for the scaled distance $cd$.
In particular, $0\xi$ is the one-point space.
Equip products with the maximum metric.
Since $b(D)\geq1/2$, we can define
\[
 \DiameterNormalization\colon\GHSpace\longrightarrow\UnitDiameterGHSpace,\qquad
 \DiameterNormalization(\xi)=b(D)^{-1}\bigl(\xi\times\beta(D)B\bigr).
\]
Indeed, the product has diameter $b(D)$.
If $D=1$, then $\beta(D)=0$ and $\DiameterNormalization(\xi)=\xi$.
At the one-point space we have $\DiameterNormalization(\xi)=B$.

Define a homotopy from the identity to this map by
\[
 H(\xi,t)=\frac{1}{1-t+t\,b(D)}\,\bigl(\xi\times(t\beta(D))B\bigr)
 \qquad (0\leq t\leq1).
\]
The denominator is at least $1/2$.
Continuity of the diameter, of rescaling (\Cref{lem:scaling-contraction}), and of the product with the maximum metric (\Cref{lem:metric-products}) implies continuity of $H$.
Its endpoints satisfy $H(\xi,0)=\xi$ and $H(\xi,1)=\DiameterNormalization(\xi)$.
When $D=1$ we have $H(\xi,t)=\xi$ at every time.
Hence $\UnitDiameterGHSpace$ is a strong deformation retract of $\GHSpace$.
By \Cref{thm:absolute-extensor}, $\GHSpace$ is an AR, so its retract $\UnitDiameterGHSpace$ is also an AR.
Moreover, $\UnitDiameterGHSpace$ is a closed subspace of the complete separable metric space $\GHSpace$, so it is complete and separable.

Let $\OpenCover_{0}$ be a family of open subsets of $\GHSpace$ covering $\UnitDiameterGHSpace$, and take the open cover $\OpenCover=\OpenCover_{0}\cup\{\GHSpace\setminus\UnitDiameterGHSpace\}$ of $\GHSpace$.
Apply the construction in the proof of \Cref{thm:discrete-approximation} to maps into $\UnitDiameterGHSpace$ and the cover $\OpenCover$.
For each $\MetricSpace{\BaseCarrier}{\MetricSymbol}\in\UnitDiameterGHSpace$ and each integer $\EquilateralSize\geq2$, the approximating space constructed in Step~\ref{step:discrete-products} still has diameter $1$.
Indeed, since $\CoverScale\leq1/16$,
\[
 {\diam\bigl(\ProductApproximation_\EquilateralSize}
             \MetricSpace{\BaseCarrier}{\MetricSymbol}\bigr)
 =\max\{1,\CoverScaleAt{\BaseCarrier}{\MetricSymbol}\}=1.
\]
Thus each approximating map in the proof of \Cref{thm:discrete-approximation} takes values in $\UnitDiameterGHSpace$.
Its images form a discrete family in $\GHSpace$ and hence also in $\UnitDiameterGHSpace$, by restricting open neighborhoods.
Both the original and approximating maps take values in $\UnitDiameterGHSpace$, so $\OpenCover$-closeness implies $\OpenCover_0$-closeness.
As in the proof of \Cref{thm:part-iv-main}, Toru\'nczyk's characterization (\Cref{thm:torunczyk}) applies.
\end{proof}

\begin{remark}\label{rem:banach-homeomorphism}
Kadets proved that all separable infinite-dimensional Banach spaces are homeomorphic (\cite[Theorem]{Kadets1967}).
Thus \Cref{thm:part-iv-main} implies that $\GHSpace$ is homeomorphic to every real separable infinite-dimensional Banach space, including $\ContinuousFunctions([0,1],\RealNumbers)$ with the uniform norm and $\ell^p$ for every $1\leq p<\infty$.
\end{remark}

The main theorem also implies strong local contractions.

\begin{corollary}\label{cor:strong-local-contractions}
Every point $\MetricSpace{\BaseCarrier }{\MetricSymbol }\in \GHSpace$ has a basis of open neighborhoods which strongly deformation retract onto $\{\MetricSpace{\BaseCarrier }{\MetricSymbol }\}$.
\end{corollary}

\begin{remark}\label{rem:nonarchimedean-gh}
The author determined the topology of the non-Archimedean Gromov--Hausdorff space by identifying it isometrically with a function space (\cite[Theorems~1.1 and~1.3(3)--(4)]{Ishiki2024NonseparableUrysohn}).
Let $\UltrametricGHSpace$ denote the set of isometry classes of nonempty compact ultrametric spaces.
The non-Archimedean Gromov--Hausdorff distance $\NonArchimedeanGHDistance$ is the infimum of the Hausdorff distances between isometric images in a common ultrametric space.
Let $\UrysohnFunctionModel$ consist of all functions $f\colon[0,\infty)\to\NonnegativeIntegers$ satisfying
\[
 f(0)=0,
\]
and for every $\varepsilon>0$,
\[
 \Card\bigl(\{r\geq\varepsilon\mid f(r)\ne0\}\bigr)<\infty.
\]
Equip $\UrysohnFunctionModel$ with the ultrametric
\[
 \UrysohnFunctionDistance(f,g)
 =
 \begin{cases}
  \max\{r\in[0,\infty)\mid f(r)\ne g(r)\},&f\ne g,\\
  0,&f=g
 \end{cases}.
\]
By the author's theorem, $(\UltrametricGHSpace,\NonArchimedeanGHDistance)$ is isometric to $(\UrysohnFunctionModel,\UrysohnFunctionDistance)$.
\end{remark}

\section{Questions}\label{sec:questions}
We conclude with questions about related spaces and the geometry of GH balls.
We use the mm-isomorphism classes and box topology of \Cref{def:box-distance}.

\begin{question}\label{question:mm-hilbert-space}
Is the space of mm-isomorphism classes of metric measure spaces, equipped with the box topology, homeomorphic to the real Hilbert space $\SeparableHilbertSpace$?
\end{question}

We recall the concentration topology to distinguish it from the box topology.
For measurable real functions $f,g$ on $\BoxParameterInterval$, identified up to equality almost everywhere, define the \emph{Ky Fan distance} by
\[
 \KyFanDistance(f,g)=\inf\{\ErrorTolerance>0\mid
 \BoxLebesgueMeasure(\{s\in\BoxParameterInterval\mid
 |f(s)-g(s)|>\ErrorTolerance\})\leq\ErrorTolerance\}.
\]
Take mm-spaces $\mathbf X=(X,d_X,\mu_X)$ and $\mathbf Y=(Y,d_Y,\mu_Y)$.
Denote by $\LipschitzObservables(X,d_X)$ all real-valued $1$-Lipschitz functions on $(X,d_X)$.
For a parameter $\BoxFirstParameter$ of $\mathbf X$ put $\BoxFirstParameter^*\LipschitzObservables(X,d_X)
=\{f\circ\BoxFirstParameter\mid f\in\LipschitzObservables(X,d_X)\}$.
Let $\BoxFirstParameter$ and $\BoxSecondParameter$ range over all parameters of $\mathbf X$ and $\mathbf Y$, respectively.
The \emph{observable distance} is
\[
 \ConcentrationDistance(\mathbf X,\mathbf Y)
 =\inf_{\BoxFirstParameter,\BoxSecondParameter}
 \HausdorffDistance{\KyFanDistance}
 \bigl(\BoxFirstParameter^*\LipschitzObservables(X,d_X),
       \BoxSecondParameter^*\LipschitzObservables(Y,d_Y)\bigr).
\]
For nonempty sets $A,B$ of measurable real functions on $\BoxParameterInterval$, we use
\[
 \HausdorffDistance{\KyFanDistance}(A,B)
 =\max\left\{
   \sup_{f\in A}\inf_{g\in B}\KyFanDistance(f,g),
   \sup_{g\in B}\inf_{f\in A}\KyFanDistance(f,g)
 \right\}.
\]
This definition does not require compactness of $A$ or $B$.
It is finite because $\KyFanDistance\leq1$.
The metric $\ConcentrationDistance$ on $\MMClassSpace$ induces the \emph{concentration topology} (see Kazukawa, Nakajima, and Shioya~\cite[arXiv v1, Definitions~2.4, 2.13 and Theorem~2.14]{KazukawaNakajimaShioya2024}).

The analogous question for the concentration topology has a negative answer.
Kazukawa, Nakajima, and Shioya proved that $(\MMClassSpace,\ConcentrationDistance)$ is not a Baire space (\cite[arXiv v1, Theorem~1.3]{KazukawaNakajimaShioya2024}), whereas $\SeparableHilbertSpace$ is a Baire space.

Nakajima and Shioya constructed a compact metrizable space containing a dense topological copy of $\GHSpace$ (\cite[arXiv v1, Main Theorem~1.2]{NakajimaShioyaCompactification}).
This leads to the following question.

\begin{question}\label{question:gh-compactification-hilbert-cube}
Is the compactification of $\GHSpace$ constructed by Nakajima and Shioya homeomorphic to the Hilbert cube $\HilbertCube$?
\end{question}

The neighborhood basis in \Cref{cor:strong-local-contractions} need not consist of GH balls.
We use the notation $\GHball(\MetricSpace{\BaseCarrier}{\MetricSymbol},\Radius)$ for open Gromov--Hausdorff balls introduced in \Cref{sec:preliminaries}.

\begin{question}\label{question:gh-balls-contractible}
For every $\MetricSpace{\BaseCarrier}{\MetricSymbol}\in\GHSpace$, does there exist $\Radius_0>0$ such that for every $0<\Radius<\Radius_0$, the ball $\GHball(\MetricSpace{\BaseCarrier}{\MetricSymbol},\Radius)$ is contractible?
More strongly, can each such ball be strongly deformation retracted onto $\{\MetricSpace{\BaseCarrier}{\MetricSymbol}\}$?
\end{question}

One can also ask for a contraction whose image lies in a larger ball with a controlled radius.

\begin{question}\label{question:gh-linear-local-contraction}
For every $\MetricSpace{\BaseCarrier}{\MetricSymbol}\in\GHSpace$, do there exist $\ContractionFactor\geq1$ and $\Radius_0>0$ such that for every $0<\Radius<\Radius_0$, there is a continuous map $\ContractionHomotopy\colon
 \GHball(\MetricSpace{\BaseCarrier}{\MetricSymbol},\Radius)\times\UnitInterval
 \to
 \GHball(\MetricSpace{\BaseCarrier}{\MetricSymbol},\ContractionFactor\Radius)$ such that for every $\MetricSpace{\ComparisonCarrier}{\ComparisonMetric}\in\GHball(\MetricSpace{\BaseCarrier}{\MetricSymbol},\Radius)$, we have $\ContractionHomotopy(\MetricSpace{\ComparisonCarrier}{\ComparisonMetric},0)=\MetricSpace{\ComparisonCarrier}{\ComparisonMetric}$, and $\ContractionHomotopy(\MetricSpace{\ComparisonCarrier}{\ComparisonMetric},1)=\MetricSpace{\BaseCarrier}{\MetricSymbol}$, and for every $\HomotopyTime\in\UnitInterval$, $\ContractionHomotopy(\MetricSpace{\BaseCarrier}{\MetricSymbol},\HomotopyTime)=\MetricSpace{\BaseCarrier}{\MetricSymbol}$?
Can $\ContractionFactor$ be chosen independently of $\MetricSpace{\BaseCarrier}{\MetricSymbol}$, with $\Radius_0$ still allowed to depend on $\MetricSpace{\BaseCarrier}{\MetricSymbol}$?
\end{question}

To define the last compactification, use $\MMClassSpace$ and $\BoxDistance$ from \Cref{def:box-distance}.
For $\mathbf X,\mathbf Y\in\MMClassSpace$ choose representatives $(A,d_A,\mu_A)$ and $(B,d_B,\mu_B)$.
Write $\mathbf Y\prec \mathbf X$ if there exists a $1$-Lipschitz map $f\colon\supp\mu_A\to\supp\mu_B$ satisfying $f_*(\mu_A|_{\supp\mu_A})=\mu_B|_{\supp\mu_B}$.
This condition is independent of the representatives.
A \emph{pyramid} is a nonempty box-closed subset $\Pyramid\subset\MMClassSpace$ such that $\mathbf Y\in\Pyramid$ for every $\mathbf X\in\Pyramid$ and every $\mathbf Y\in\MMClassSpace$ with $\mathbf Y\prec \mathbf X$, and such that for every $\mathbf X,\mathbf Y\in\Pyramid$ some $\mathbf Z\in\Pyramid$ satisfies $\mathbf X\prec \mathbf Z$ and $\mathbf Y\prec \mathbf Z$ (Kazukawa, Nakajima, and Shioya, \cite[arXiv v1, Definitions~2.2 and~2.16]{KazukawaNakajimaShioya2024}).
Let $\PyramidSpace$ denote the set of pyramids.
Take a sequence $\{\Pyramid_n\}_{n\in\NonnegativeIntegers}$ in $\PyramidSpace$ and a point $\Pyramid\in\PyramidSpace$.
For $\mathbf X\in\MMClassSpace$ put $\BoxDistance(\mathbf X,\Pyramid_n)
 =\inf_{\mathbf Y\in\Pyramid_n}\BoxDistance(\mathbf X,\mathbf Y)$.
The weak topology on $\PyramidSpace$ is metrizable and is characterized by the following criterion for $\Pyramid_n\to\Pyramid$.
For every $\mathbf X\in\Pyramid$ one has $\BoxDistance(\mathbf X,\Pyramid_n)\longrightarrow0$, and for every $\mathbf X\in\MMClassSpace\setminus\Pyramid$ one has $\liminf_{n\to\infty}\BoxDistance(\mathbf X,\Pyramid_n)>0$.
Through the assignment $\mathbf X\mapsto\{\mathbf Y\in\MMClassSpace\mid \mathbf Y\prec \mathbf X\}$ the concentration space $(\MMClassSpace,\ConcentrationDistance)$ is identified with a dense subspace of the compact space $\PyramidSpace$ (Kazukawa, Nakajima, and Shioya, \cite[arXiv v1, Definition~2.17 and Theorem~2.18]{KazukawaNakajimaShioya2024}).
Kazukawa, Nakajima, and Shioya proved that $\PyramidSpace$ is contractible and locally path connected (\cite[arXiv v1, Theorems~1.4 and~1.6]{KazukawaNakajimaShioya2024}).
This compactification differs from the one in \Cref{question:gh-compactification-hilbert-cube}.

\begin{question}\label{question:pyramids-hilbert-cube}
Is $\PyramidSpace$ with the weak topology homeomorphic to the Hilbert cube $\HilbertCube$?
\end{question}

\Cref{cor:ghp-retract} embeds the Gromov--Hausdorff space $\GHSpace$ as a retract of the measured space $\MeasuredSpaces$.
We ask whether $\MeasuredSpaces$ has the same homeomorphism type as $\GHSpace$.

Equip $\MeasuredSpaces$ with the Gromov--Hausdorff--Prokhorov metric $\GHPDistance$ from \Cref{def:ghp-distance}.
A probability measure here can have proper support.
For example, concentrating all mass at one point of a two-point space yields the class of a singleton under the mm-isomorphism of \Cref{question:mm-hilbert-space}, since the support is a singleton.
In $\MeasuredSpaces$ the two-point space and the singleton represent distinct classes, since isomorphisms are surjective measure-preserving isometries.

\begin{question}\label{question:ghp-hilbert-space}
Is $(\MeasuredSpaces,\GHPDistance)$ homeomorphic to the real Hilbert space $\SeparableHilbertSpace$?
\end{question}

\printbibliography
\end{document}